\documentclass[11pt]{amsart}

\usepackage[T1]{fontenc}
\usepackage[utf8]{inputenc}
\usepackage[english]{babel}
\usepackage{amscd,amsmath,amsthm,amssymb}
\usepackage{mathtools}
\usepackage{mathrsfs}
\usepackage{comment}
\usepackage{xcolor}
\usepackage{array}
\usepackage{caption}
\usepackage{mathabx}
\usepackage{standalone}
\usepackage{scalerel}
\usepackage{csquotes}

\usepackage{enumitem}
\setlist[enumerate]{label=\textnormal{(\arabic*)}}
\newcommand{\defaultRoman}{{\textnormal{(\textit{\roman*})}}}

\definecolor{cadmiumgreen}{rgb}{0.0, 0.42, 0.24}
\usepackage[
  colorlinks,
  citecolor=cadmiumgreen,
  pdfstartview ={FitV},
]{hyperref}
\hypersetup{
  pdfauthor={Omid Amini, Matthieu Piquerez},
  pdftitle={Hodge theory for tropical varieties}}

\usepackage[
  alphabetic,
  msc-links,
  nobysame,
  lite,
]{amsrefs}

\usepackage{tikz, tikz-cd}
\usetikzlibrary{calc, backgrounds, quotes, arrows.meta}

\usepackage[left=3cm,top=3.5cm,right=3cm]{geometry}
\newtheorem{thm}{Theorem}[section]

\newtheorem{lemma}[thm]{Lemma}

\newtheorem{propdefi}[thm]{Proposition - Definition}
\newtheorem{prop}[thm]{Proposition}
\newtheorem{add}[thm]{Addendum}
\newtheorem{claim}[thm]{Claim}
\newtheorem{cor}[thm]{Corollary}

\theoremstyle{definition}
\newtheorem{defii}[thm]{Definition}
\newenvironment{defi}
  {\pushQED{\qed}\defii}
  {\popQED\enddefii}
\newtheorem{remm}[thm]{Remark}
\newenvironment{remark}
  {\pushQED{\qed}\remm}
  {\popQED\endremm}
\newtheorem{exx}[thm]{Example}
\newenvironment{example}
  {\pushQED{\qed}\exx}
  {\popQED\endexx}

\numberwithin{equation}{section}

\newcommand{\hodgediamond}{
\renewcommand{\d}{3} 
\renewcommand{\e}{1} 
\renewcommand{\u}{.7} 
\renewcommand{\v}{.7} 
\renewcommand{\a}{1.2} 
\renewcommand{\b}{.9} 
\renewcommand{\r}{.3} 
\begin{tikzpicture}[scale=.8]
\begin{scope}
\node (00) at (0,\d) {$\ST^{0,0}$};
\node (0d) at (\d,\d) {$\ST^{d,0}$};
\node (d0) at (-\d,-\d) {$\ST^{-d,2d}$};
\node (dd) at (0,-\d) {$\ST^{0,2d}$};
{ 
\node (p) at (-\b,\a-\b) {$\blacksquare$};
\node (p1) at (\b,\a+\b) {$\square$};
\node (p2) at (\b,-\a+\b) {$\square$};
\node (p3) at (-\b,-\a-\b) {$\square$};
}
\begin{scope}[on background layer]
  \begin{scope}
    \clip (00.center)--(0,0)--(d0.center)--cycle;
    \fill[gray!20] (p.center)--($(p)+(00)$)--($(p)+(d0)$)--cycle;
  \end{scope}
\end{scope}
\draw[loosely dotted] (00)--(0d)--(dd)--(d0)--(00);
\draw[gray] (00)--(dd) (0d)--(d0);
\node (center) at (0,0) {$\bullet$};
\draw[<->] (00) ++ (0,1) ++ (-\u,-\u) to node[near end, above] {$N$} node {\tikz{\draw[thick, -] (0,-.3)--(0,.3);}}  ++(2*\u, 2*\u);
\draw[<->] (-\d,0) ++ (-.5,0) ++ (0,\u) to node[very near end, left] {$\ell$} node {\tikz{\draw[thick, -] (-.2,-.2)--(.2,.2);}}  ++(0, -2*\u);
\draw[ultra thin, dashed] (p) -- (p1) -- (p2) -- (p3) -- (p);
\draw[-latex, thick] (150:\r) arc(150:-150:\r) --++ (130:.1);
\draw (45:\r) node[above right=-1mm] {$\pd$};
\draw (-\d, \d+1) node {\tikz[thick, -latex]{\draw edge["$a$"] ++ (\v, 0); \draw edge["$b$"'] ++ (0, -\v)}};
\end{scope}
\begin{scope}[shift={(2.8*\d,0)}]
\node (00) at (0,\d) {$\ST^{0,0}$};
\node (0d) at (\d,0) {$\ST^{d,0}$};
\node (d0) at (-\d,0) {$\ST^{-d,2d}$};
\node (dd) at (0,-\d) {$\ST^{0,2d}$};
{
\node (p)  at (-\b,\a) {$\blacksquare$};
\node (p1) at (\b,\a)  {$\square$};
\node (p2) at (\b,-\a) {$\square$};
\node (p3) at (-\b,-\a) {$\square$};
}
\coordinate (p) at (-\b,\a);
\begin{scope}[on background layer]
  \begin{scope}
    \clip (00.center)--(0,0)--(d0.center)--cycle;
    \fill[gray!20] (p.center)--($(p)+(00)$)--($(p)+(d0)$)--cycle;
  \end{scope}
\end{scope}
\draw[loosely dotted] (00)--(0d)--(dd)--(d0)--(00);
\draw[gray] (00)--(dd) (0d)--(d0);
\node (center) at (0,0) {$\bullet$};
\draw[<->] (00) ++ (0,1) ++ (-\u,0) to node[near end, above] {$N$}  ++(2*\u, 0);
\draw[<->] (\d,0) ++ (1,0) ++ (0,\u) to node[very near end, left] {$\ell$} ++(0, -2*\u);
\draw[ultra thin, dashed] (p) -- (p1) -- (p2) -- (p3) -- (p);
\draw[-latex, thick] (150:\r) arc(150:-150:\r) --++ (130:.1);
\draw (45:\r) node[above right=-1mm] {$\pd$};
\draw (\d, \d+1) node {\tikz[thick, -latex]{\draw edge["$a$"] ++ (\v, -\v); \draw edge["$b$"'] ++ (0, -\v)}};
\end{scope}
\draw[-Implies, double distance=4pt, thick] (1.9*\d, -.4*\d) -- (1.6*\d,-.6*\d);
\begin{scope}[shift={(1.2*\d,-\d*.9)}]
\node (00) at (0,\e) {$H^{0,0}$};
\node (0d) at (\e,0) {$H^{0,d}$};
\node (d0) at (-\e,0) {$H^{d,0}$};
\node (dd) at (0,-\e) {$H^{d,d}$};
\draw[dotted] (00)--(0d)--(dd)--(d0)--(00);
\draw[gray] (00)--(dd) (0d)--(d0);
\end{scope}
\end{tikzpicture}
}

\newcommand{\bigarray}{
\begin{array}{cccc|c|cccc|c|cccc|c|}
\cline{2-5} \cline{7-10} \cline{12-15}
\multicolumn{1}{c|}{\text{\Large\strut}}   & \multicolumn{4}{c|}{\trop}                                             & \multicolumn{1}{l|}{} & \multicolumn{4}{c|}{\ST}                                               & \multicolumn{1}{l|}{} & \multicolumn{4}{c|}{\D}                                               \\ \cline{2-5} \cline{7-10}  \cline{12-15}
\multicolumn{1}{c|}{\text{\Large\strut}}   & \multicolumn{1}{c|}{\d}   & \multicolumn{1}{c|}{\d'}  & \d^\i & \d^\pi & \multicolumn{1}{l|}{} & \multicolumn{1}{c|}{\d}   & \multicolumn{1}{c|}{\d'}  & \d^\i & \d^\pi & \multicolumn{1}{l|}{} & \multicolumn{1}{c|}{\d}   & \multicolumn{1}{c|}{\d'}  & \d^\i & \d^\pi \\ \cline{2-5} \cline{7-10}  \cline{12-15}
\multicolumn{1}{c|}{\d} & \multicolumn{1}{c|}{\taa} & \multicolumn{1}{c|}{\tab} & \tac  & \tad   & \multicolumn{1}{l|}{} & \multicolumn{1}{c|}{\saa} & \multicolumn{1}{c|}{\sab} & \sac  & \sad   & \multicolumn{1}{l|}{} & \multicolumn{1}{c|}{\aaa} & \multicolumn{1}{c|}{\aab} & \aac  & \aad   \\ \cline{2-5} \cline{7-10}  \cline{12-15}
\d'                     & \multicolumn{1}{c|}{ }    & \multicolumn{1}{c|}{\tbb} & \tbc  & \tbd   &                       & \multicolumn{1}{c|}{ }    & \multicolumn{1}{c|}{\sbb} & \sbc  & \sbd   &                       & \multicolumn{1}{c|}{ }    & \multicolumn{1}{c|}{\abb} & \abc  & \abd   \\ \cline{3-5} \cline{8-10}  \cline{13-15}
\d^\i                   &                           & \multicolumn{1}{c|}{ }    & \tcc  & \tcd   &                       &                           & \multicolumn{1}{c|}{ }    & \scc  & \scd   &                       &                           & \multicolumn{1}{c|}{ }    & \acc  & \acd   \\ \cline{4-5} \cline{9-10}  \cline{14-15}
\d^\pi                  &                           &                           &       & \tdd   &                       &                           &                           &       & \sdd   &                       &                           &                           &       & \add   \\ \cline{5-5} \cline{10-10} \cline{15-15}
\end{array}
}

\newcommand{\cf}[1]{cf.}
\newcommand{\ie}{i.e.}
\newcommand{\resp}{resp.\ }
\renewcommand{\~}{\widetilde}
\newcommand{\Q}{\mathbb{Q}}
\newcommand{\Z}{\mathbb{Z}}
\newcommand{\R}{\mathbb{R}}
\newcommand{\hooklongrightarrow}{\lhook\joinrel\longrightarrow}
\newcommand{\simto}{\xrightarrow{\raisebox{-3pt}[0pt][0pt]{\small$\hspace{-1pt}\sim$}}}
\newcommand{\longsimto}{\xrightarrow{\ \raisebox{-3pt}[0pt][0pt]{\small$\hspace{-1pt}\sim$\ }}}
\newcommand{\myand}{\text{ and }}
\newcommand{\bul}{\bullet} 

\renewcommand{\emptyset}{\varnothing}
\newcommand{\zint}[2]{\{#1,\dots,#2\}}
\newcommand{\rquot}[2]{#1\big/#2}
\newcommand{\rest}[1]{\raisebox{-1pt}{$\vert$}_{#1}}
\renewcommand{\pmod}[1]{\ (\mathrm{mod}\ #1)}
\newcommand{\ldot}{\,\cdot}
\newcommand{\ccdot}{\,\cdot\,}
\newcommand{\rdot}{\cdot\,}
\newcommand{\QED}{\hspace{\fill}{\qed}} 

\let\oldchi\chi
  \newcommand{\raisechi}[2]{\raisebox{.4ex}{$#1#2$}}
  \renewcommand{\chi}{{\mathpalette\raisechi\oldchi}}

\makeatletter
\let\oldsum\sum
\renewcommand{\sum}{\@ifnextchar_\@mysum\oldsum}
\def\@mysum_#1{\oldsum_{\substack{#1}}}
\let\oldbigoplus\bigoplus
\renewcommand{\bigoplus}{\@ifnextchar_\@mybigoplus\oldbigoplus}
\def\@mybigoplus_#1{\oldbigoplus_{\substack{#1}}}
\let\oldprod\prod
\renewcommand{\prod}{\@ifnextchar_\@myprod\oldprod}
\def\@myprod_#1{\oldprod_{\substack{#1}}}
\makeatother

\let\oldbigwedge\bigwedge
\renewcommand{\bigwedge}{{\textstyle\oldbigwedge\!}}

\newcommand{\cddots}{\raisebox{2pt}[8pt][4pt]{\makebox[20pt]{$\dots$}}}
\newcommand{\cdvdots}{\raisebox{0pt}[14pt][4pt]{\makebox[10pt]{$\vdots$}}}
\newcommand{\cdddots}{\raisebox{0pt}[14pt][4pt]{\makebox[20pt]{$\ddots$}}}

\let\Im\relax
  \DeclareMathOperator{\Im}{Im} 
\DeclareMathOperator{\cl}{cl} 
\DeclareMathOperator{\Hom}{Hom} 
\DeclareMathOperator{\sed}{sed} 
\DeclareMathOperator{\gys}{Gys} 
\DeclareMathOperator{\HR}{HR} 
\DeclareMathOperator{\HL}{HL} 
\DeclareMathOperator{\WPD}{WPD} 
\DeclareMathOperator{\PD}{PD} 
\DeclareMathOperator{\id}{id} 
\DeclareMathOperator{\Bl}{\mathcal B\ell} 
\DeclareMathOperator{\gr}{gr} 
\DeclareMathOperator{\Tot}{Tot} 
\DeclareMathOperator{\Vect}{Vect} 
\DeclareMathOperator{\dist}{dist} 
\DeclareMathOperator{\sign}{sign} 
\DeclareMathOperator{\pd}{pd} 
\DeclareMathOperator{\class}{cl} 

\newcommand{\G}{\mathbb G} 
\renewcommand{\u}{\mathbf u} 
\renewcommand{\v}{\mathbf v} 
\newcommand{\trop}{\mathrm{trop}} 
\newcommand{\Cl}{\mathcal C\!\ell} 
\newcommand{\hyp}{\mathbb H} 
\newcommand{\tm}{\mathscr{M}} 
\newcommand{\eR}{\mathbf R} 
\newcommand{\I}{\mathscr{I}} 
\newcommand{\C}{{C^\circ}} 
\newcommand{\e}{{\mathfrak e}} 
\newcommand{\f}{{\textnormal{\textsl{\textsf f}}}} 
\renewcommand{\O}{\mathcal O} 
\renewcommand{\P}{\mathbb P} 
\newcommand{\TP}{\mathbb{TP}} 
\newcommand{\corps}{\mathbb K} 
\newcommand{\prim}{\mathrm{prim}} 
\renewcommand{\div}{\mathrm{div}} 
\newcommand{\TT}{\mathrm{T}} 
\newcommand{\Tan}{\mathrm{T}_{\!\scriptscriptstyle\mathrm{aff}}} 
\newcommand{\dual}{\star} 
\newcommand{\CC}{\mathcal C} 
\newcommand{\cteun}{\mathbf 1} 
\let\i\relax
  \newcommand{\i}{{\mathop{}\mathrm{i}}} 
\renewcommand{\d}{{\rm d}} 
\newcommand{\dfrak}{\mathfrak d} 
\newcommand{\parr}{{\scaleto{/\!/}{6pt}}} 
\newcommand{\SF}{\textrm{\bf F}} 
\newcommand{\nvect}{\mathfrak n} 
\newcommand{\p}{{\rm p}} 
\newcommand{\cp}{\overline\p} 
\newcommand{\x}{\textsc{x}} 
\newcommand{\y}{\textsc{y}} 

\newcommand{\ooplus}{\raisebox{0pt}[0pt]{${}\overset\perp\oplus{}$}} 
\newcommand{\woplus}{\ \oplus\ } 
\newcommand{\wooplus}{\ \ooplus\ } 

\newcommand{\card}[1]{\lvert#1\rvert} 
\newcommand{\st}{\bigm|} 
\newcommand{\Bigst}{\Bigm|} 
\newcommand{\abs}[1]{\lvert #1\rvert} 
\newcommand{\comp}[1]{\overline{#1}} 
\newcommand{\modtrop}[2]{\Bl_{#1}(#2)} 
\newcommand{\pairing}[2]{\langle#1,#2\rangle} 
\newcommand{\ddelta}{{\underline\delta}} 
\newcommand{\eeta}{{\underline\eta}} 
\renewcommand{\ssum}{{\textstyle\sum}} 
\newcommand{\sminfty}{{\scaleobj{.7}{\infty}}} 

\newcommand{\X}{\mathfrak X}
\newcommand{\Y}{\mathscr Y}
\newcommand{\V}{\mathscr V}
\newcommand{\vZ}{\mathcal Z}
\newcommand{\U}{\mathscr U}
\newcommand{\W}{\mathscr W}

\newcommand{\dip}{\d'}
\newcommand{\dis}{\d''}
\newcommand{\Dol}{\mathrm{Dol}}
\newcommand{\A}{\mathcal A} 

\newcommand{\dimsaux}[2]{\raisebox{.2ex}{\scalebox{1}[.8]{$#1\lvert$}}#2\raisebox{.2ex}{\scalebox{1}[.8]{$#1\rvert$}}}
  \newcommand{\dims}[1]{\mathpalette\dimsaux{#1}}
\newcommand{\suppaux}[2]{\scalebox{1}[1.4]{$#1\lvert$}#2\scalebox{1}[1.4]{$#1\rvert$}}
  \newcommand{\supp}[1]{\mathpalette\suppaux{#1}}
\newcommand{\K}{\mathcal K} 
\newcommand{\cone}{\R_+} 
\newcommand{\pen}{\mathscr P} 
\newcommand{\conezero}{{\underline0}} 

\newcommand{\subface}{\prec}
\newcommand{\ssubface}{\mathbin{\mathchoice
  {\subface\!\!\!\cdot}%
  {\subface\!\!\!\cdot}%
  {\subface\!\cdot}%
  {\subface\!\cdot}%
}} 
\newcommand{\supface}{\succ}
\newcommand{\ssupface}{\mathbin{\mathchoice
  {\cdot\!\!\!\supface}%
  {\cdot\!\!\!\supface}%
  {\cdot\!\supface}%
  {\cdot\!\supface}%
}}

\DeclareMathOperator{\conv}{conv} 
\DeclareMathOperator{\adhop}{adh}
  \newcommand{\adh}[1]{\adhop(#1)} 
\DeclareMathOperator{\coneupop}{cone}
  \newcommand{\coneup}[1]{\coneupop(#1)} 
\DeclareMathOperator{\aff}{\mathcal L} 
\DeclareMathOperator{\lpm}{\mathcal L^{pm}} 
\DeclareMathOperator{\faceop}{face}
  \newcommand{\face}[1]{\faceop(#1)} 

\newcommand{\Ma}{\mathfrak M} 
\newcommand{\Fl}{\mathscr{F}}
\newcommand{\distel}{\ast}
\newcommand{\bases}{\mathfrak B}
\newcommand{\ind}{\mathfrak I}
\newcommand{\rkm}{\mathrm{rk}}
\newcommand{\contr}[1]{/#1} 
\newcommand{\del}{\backslash} 

\newcommand{\ws}{\scaleto{W}{4pt}}
\DeclareMathOperator{\maxsed}{maxsed}

\newcommand{\dd}{\d^\dual} 
\newcommand{\Lap}{\triangle} 
\newcommand{\w}{{\rm w}} 
\newcommand{\m}{\mathfrak m} 

\newcommand{\ST}{\textnormal{\textsf{S\!T}}} 
\newcommand{\HD}{\textnormal{\textsf{HD}}} 

\newcommand{\STpnop}{\textnormal{\textsf{St}}}
\newcommand{\STp}[1]{\prescript{}{#1}\STpnop}

\newcommand{\STpab}[2]{\displaystyle\bigoplus_{\delta\in X_\f \\ \dims{\delta}=#1} H^{#2}(\delta)}

\newcommand{\STi}{\overline\STpnop}
\newcommand{\STinf}[1]{\prescript{}{#1}\STi}
\newcommand{\STinfI}[1]{\prescript\downarrow{#1}\STi}

\newcommand{\STinfab}[4]{\displaystyle\bigoplus_{\dims{#1}=#2 \\ \dims{#1_\infty}\leq#3} H^{#4}(#1)}

\newcommand{\CCnop}{\textnormal{\textsf{C}}}
\newcommand{\CCp}[1]{\prescript{}{#1}{\CCnop}}
\newcommand{\CCab}[3]{\raisebox{0pt}[\height][2pt]{$\displaystyle \bigoplus_{\dims{\delta}=#1}\bigwedge^{#2}\TT^\dual\delta\otimes\SF^{#3}(\conezero^\delta)$}}

\newcommand{\Dnop}{\textnormal{\textsf{D}}}
\newcommand{\D}{\Dnop}
\newcommand{\Da}[1]{\prescript{}{#1}\Dnop}
\newcommand{\DI}[1]{\prescript\downarrow{#1}\Dnop}

\newcommand{\Dab}[4]{\displaystyle\bigoplus_{\dims\delta=#1} \bigwedge^{#2} \TT^\dual\delta \otimes \Bigl(\makebox[20pt]{$\displaystyle\bigoplus_{\eta \supface \delta \\ \dims\eta = #3 \\ \sed(\eta)=\sed(\delta)}$} H^{#4}(\eta)\Bigr)} 

\newcommand{\Enop}{\textnormal{\textsf{E}}}
\newcommand{\E}{\Enop}
\newcommand{\EI}{\prescript{\textup I}{}\Enop}

\renewcommand{\AA}{\textnormal{\textsf{A}}}
\newcommand{\AAa}[1]{\prescript{}{#1}\AA}
\newcommand{\AAp}[1]{\prescript{#1}{}\AA}

\newcommand{\BB}{\textnormal{\textsf{B}}}

\newcommand{\Cx}[2]{C^\bul_{#2}(#1)}
\newcommand{\Sxab}[2]{\displaystyle\bigoplus_{\delta\subface\eta \\ \dims\delta = #1}\bigwedge^{#2}\TT^\dual\delta} 
\newcommand{\Syab}[2]{\displaystyle\bigoplus_{\delta\subface\eta \\ \dims\delta = #1}\bigwedge^{#2}\TT\delta} 

\begin{document}
\title{Hodge theory for tropical varieties}

\author{Omid Amini}
\address{CNRS - CMLS, \'Ecole Polytechnique}
\email{\href{omid.amini@polytechnique.edu}{omid.amini@polytechnique.edu}}

\author{Matthieu Piquerez}
\address{CMLS, \'Ecole Polytechnique}
\email{\href{matthieu.piquerez@polytechnique.edu}{matthieu.piquerez@polytechnique.edu}}


\subjclass[2010]{
\href{https://mathscinet.ams.org/msc/msc2010.html?t=14T05}{14T05},
\href{https://mathscinet.ams.org/msc/msc2010.html?t=14C30}{14C30},
\href{https://mathscinet.ams.org/msc/msc2010.html?t=14D06}{14D06},
\href{https://mathscinet.ams.org/msc/msc2010.html?t=14D07}{14D07},
\href{https://mathscinet.ams.org/msc/msc2010.html?t=32G20}{32G20},
\href{https://mathscinet.ams.org/msc/msc2010.html?t=32S35}{32S35},
\href{https://mathscinet.ams.org/msc/msc2010.html?t=14C15}{14C15},
\href{https://mathscinet.ams.org/msc/msc2010.html?t=14B05}{14B05},
\href{https://mathscinet.ams.org/msc/msc2010.html?t=14F05}{14F05},
\href{https://mathscinet.ams.org/msc/msc2010.html?t=14F45}{14F45},
\href{https://mathscinet.ams.org/msc/msc2010.html?t=14C17}{14C17},
\href{https://mathscinet.ams.org/msc/msc2010.html?t=05B35}{05B35},
\href{https://mathscinet.ams.org/msc/msc2010.html?t=52B70}{52B70},
\href{https://mathscinet.ams.org/msc/msc2010.html?t=52B40}{52B40},
\href{https://mathscinet.ams.org/msc/msc2010.html?t=57Q15}{57Q15},
\href{https://mathscinet.ams.org/msc/msc2010.html?t=57P10}{57P10},
\href{https://mathscinet.ams.org/msc/msc2010.html?t=55T05}{55T05},
\href{https://mathscinet.ams.org/msc/msc2010.html?t=55N30}{55N30},
\href{https://mathscinet.ams.org/msc/msc2010.html?t=57R20}{57R20},
\href{https://mathscinet.ams.org/msc/msc2010.html?t=05E45}{05E45}
}


\date{\today}

\begin{abstract} In this paper we prove that the cohomology of smooth projective tropical varieties verify the tropical analogs of three fundamental theorems which govern the cohomology of complex projective varieties: Hard Lefschetz theorem, Hodge-Riemann relations and monodromy-weight conjecture.

\medskip

On the way to establish these results, we introduce and prove other results of independent interest. This includes a generalization of the results of Adiprasito-Huh-Katz, Hodge theory for combinatorial geometries, to any unimodular quasi-projective fan having the same support as the Bergman fan of a matroid, a tropical analog for Bergman fans of the pioneering work of Feichtner-Yuzvinsky on cohomology of wonderful compactifications (treated in a separate paper, recalled and used here), a combinatorial study of the tropical version of the Steenbrink spectral sequence, a treatment of K\"ahler forms in tropical geometry and their associated Hodge-Lefschetz structures, a tropical version of the projective bundle formula, and a result in polyhedral geometry on the existence of quasi-projective unimodular triangulations of polyhedral spaces.
\end{abstract}
\maketitle

\setcounter{tocdepth}{1}

\tableofcontents


\section{Introduction}\label{sec:intro}

The aim of this work and its forthcoming companions is to study Hodge theoretic aspects of tropical and non-archimedean geometry.

\medskip

Tropical geometry studies degenerations of algebraic varieties by enriching the theory of semistable models and their dual complexes by polyhedral geometry. This enrichment motivates the development of algebraic geometry for combinatorial and polyhedral spaces~\cites{Mik06, MS, GS11, MZ, IKMZ, Cart, JSS} and their hybrid counter-parts, which contain a mixture of algebraic and polyhedral components~\cites{Berkovich, AB15, BJ, HybridModuli, Mik06}.

\medskip

In dimension one, it is now well-understood that graphs and metric graphs are in many ways similar to Riemann surfaces and algebraic curves. For example, the analog of many classical theorems governing the geometry of Riemann surfaces, Riemann-Roch, Abel-Jacobi, Clifford and Torelli, are now established for graphs and metric graphs~\cites{BN07, GK08, MZ08, Coppens, CV10, AB15}.

This viewpoint and its subsequent extensions and refinements have been quite powerful in applications in the past decade, and have resulted in several advances ranging from the study of curves and their moduli spaces (e.g. Brill-Noether theory~\cite{CDPR}, maximum rank conjecture~\cite{JP}, the Kodaira dimension of the moduli spaces of curves~\cite{FJP}) to arithmetic geometry of curves (e.g. uniform bounds on the number of rational points on curves~\cite{KRZ} and equidistribution results~\cite{Ami14}). An overview of some of these results can be found in the survey paper~\cite{BJ16}.

\medskip

It is expected that higher dimensional combinatorial and polyhedral analogs of graphs and metric graphs should fit into the same picture and fundamental theorems of algebraic geometry should have their combinatorial and polyhedral analogs.

\medskip

Representative examples of this fresco in higher dimension have emerged in the pioneering work of McMullen and Karu~\cites{Mcmullen, Karu} on hard Lefschetz theorem for polytopes, and more recently, in connection with Hodge theory, in the work of Adiprasito-Huh-Katz~\cite{AHK} in the development of Hodge theory for matroids, in the work of Babaee-Huh~\cite{BH17} on Demailly's strong reformulation of the Hodge conjecture for currents, in the work of Elias and Williamson~\cite{EW} on Hodge theory of Soergel bimodules and positivity of Kazhdan-Lustig polynomials, and in the work of Adiprasito-Bj\"oner~\cite{AB} and Adiprasito~\cite{Adi18} on weak and combinatorial Lefschetz theorems, respectively.

\medskip

Our aim here is to show that cohomology of smooth projective tropical varieties verify the tropical analogs of three fundamental theorems which govern the cohomology of complex projective varieties: Hard Lefschetz theorem, Hodge-Riemann relations and monodromy-weight conjecture. As we will explain below, this provides a uniform treatment and a generalization of many of the above mentioned results within the framework of tropical geometry.

\medskip

In the rest of this introduction, we provide a brief discussion of the results presented in this article.

\subsection{Matroids} Central objects underlying our results in this paper are combinatorial structures called matroids. Discovered by Whitney~\cite{Whitney} in his study on classification of graphs with the same cycle structure, and independently by van der Waerden~\cite{vdw} and Nakasawa, matroids provide combinatorial axiomatization of central notions of linear algebra. In this way, they lead to a common generalization of other combinatorial structures such as graphs, simplicial complexes and hyperplane arrangements. They appear naturally in diverse fields of mathematics, e.g., in combinatorics and statistical physics~\cites{Tutte, MirPer, GGW, SS14, Piquerez}, optimization theory~\cite{schrijver}, topology~\cites{Mac91, GM92}, and algebraic geometry~\cites{Mnev, Laf, BB, BL,sturmfels}. In tropical geometry, they are combinatorial structures underlying the idea of maximal degenerations of algebraic varieties~\cites{deligne-md, KS}, and in this way give rise to a tropical notion of smoothness. We refer to~\cites{wilson, kung, Oxl11} for an introduction and a historical account of the theory of matroids.

\subsection{Bergman fans} To any matroid $\Ma$ over a ground set $E$ one associates a fan $\Sigma_\Ma$ called the \emph{Bergman fan of $\Ma$} which lives in the real vector space $\rquot{\R^E}{\R (1,\dots, 1)}$ and which is rational with respect to the lattice $\rquot{\Z^E}{\Z(1, \dots, 1)}$. In the case the matroid is given by an arrangement of hyperplanes, the Bergman fan can be identified with the tropicalization of the complement of the hyperplane arrangement, for the coordinates given by the linear functions which define the hyperplane arrangement~\cite{AK}.

The fan $\Sigma_\Ma$ has the following description. First, recall that a \emph{flat} of $\Ma$ is a set $F$ of elements of $E$ which is closed with respect to the linear dependency, \ie, an element $e\in E$ which is linearly dependent to $F$ already belongs to $F$. A non-empty flat $F \subsetneq E$  is called \emph{proper}. The Bergman fan $\Sigma_\Ma$ consists of cones $\sigma_\Fl$ in correspondence with \emph{flags of proper flats of $\Ma$}
\[\Fl\colon F_1 \subsetneq F_2 \subsetneq \dots \subsetneq \dots \subsetneq F_k, \qquad k \in \mathbb N \cup \{0\}. \]
For such a flag $\Fl$, the corresponding cone $\sigma_\Fl$ is the one generated by the characteristic vectors $\e_{F_1}, \dots, \e_{F_k}$ of $F_1, \dots, F_k$ in $\R^E / \R(1, \dots, 1)$.  In particular, the rays of $\Sigma_\Ma$ are in bijection with the proper flats $F$ of $\Ma$.

\smallskip
We call a \emph{Bergman support} or a \emph{tropical linear space} any subspace of a real vector space which is isomorphic via a linear map of real vector spaces to the support of the Bergman fan of a simple matroid. The terminology comes from~\cites{Bergman, Spe08}.

\smallskip
By a \emph{Bergman fan} we mean any fan which is supported on a Bergman support. In particular, Bergman fan $\Sigma_\Ma$ of a matroid $\Ma$ is an example of a Bergman fan.
Moreover, with this terminology, any complete fan in a real vector space is a Bergman fan. Indeed, for each integer $d$, the real vector space $\R^d$ is the support of the Bergman fan of the uniform matroid $U_{d+1}$ of rank $d+1$ on ground set $[d]:=\{0, 1, \dots, d\}$. (In $U_d$ every subset of $[d]$ is an independent set.) Thus, Bergman fans generalize both complete fans and Bergman fans of matroids. In addition, products of Bergman fans remain Bergman.

\smallskip
In this paper we consider \emph{unimodular Bergman fans} (although it should be possible to generalize with some extra effort our results to simplicial Bergman fans): these are Bergman fans $\Sigma$ in a real vector space $V$ which are moreover rational with respect to a full rank sublattice $N$ of $V$ and unimodular (smooth in the language of toric geometry). This means any cone $\sigma$ in $\Sigma$ is simplicial, that is generated by $\dims{\sigma}$, the dimension of $\sigma$,  rays. Moreover, denoting $\varrho_1, \dots, \varrho_{\dims{\sigma}}$ the rays of $\sigma$, each ray $\varrho_i$ is rational, and in addition, if $\e_{\varrho_i}$ denotes the primitive vector of the ray $\varrho_i$, the set of vectors $\e_{\varrho_1}, \dots, \e_{\varrho_{\dims\sigma}}$ is part of a basis of the lattice $N$. The product of two unimodular Bergman fans remains a unimodular Bergman fan.

\smallskip
Unimodular Bergman fans are building blocks for \emph{smooth tropical varieties} considered in this paper.

\subsection{Chow rings of matroids and more general Bergman fans} To a given simple matroid $\Ma$ over a ground set $E$, one associates its Chow ring $A^\bul(\Ma)$, which is a graded ring defined by generators and relations, as follows.

A \emph{flat} of $\Ma$ is a set $F$ of elements of $E$ which is closed with respect to the linear dependency: i.e., an element $e\in E$ which is linearly dependent to $F$ already belongs to $F$. To any flat $F$ of $\Ma$ one associates a variable $\x_F$. The Chow ring is then a quotient of the polynomial ring $\mathbb Z[\x_F \,\,: F \textrm{ flat of }\Ma]$ by the ideal $\I_\Ma$ generated by the polynomials
\begin{itemize}
\item $\x_F \x_H$ for any two \emph{incomparable} flats $F$ and $H$, i.e., with $F \not\subset H$ and $H \not\subset F$;
\item $\sum_{F \ni e} \x_F - \sum_{F \ni f} \x_F$ for any pair of elements $e,f\in E$.
\end{itemize}

In the case the matroid $\Ma$ is given by an arrangement of hyperplanes, the Chow ring gets identified with the Chow ring of the wonderful compactification of the complement of the hyperplane arrangement~\cites{CP95, FY}.

\medskip

Let now $\Sigma$ be any unimodular Bergman fan in a vector space $V$ and let $N$ be the corresponding full rank lattice of $V$, so $V = N_\R :=N\otimes\R $.

\smallskip
The Chow ring of $\Sigma$ denoted by $A^\bul(\Sigma)$ is defined as the quotient of the polynomial ring $\mathbb Z[\x_\varrho \mid \varrho\in \Sigma_1]$, with generators $\x_{\varrho}$ associated to the rays $\varrho$ of $\Sigma$, by the ideal $\I_1 +\I_2$ defined as follows:
\begin{itemize}
\item $\I_1$ is the ideal generated by all monomials of the form $\prod_{\varrho \in S} \x_\varrho$ for any subset $S \subseteq \Sigma_1$ which do not form a cone in $\Sigma$; and
\item $\I_2$ is the ideal generated by the elements of the form
\[\sum_{\varrho\in \Sigma_1} \langle m, \e_{\varrho}\rangle \x_\varrho\]
for any $m$ in the dual lattice $M = N^\vee$,
\end{itemize}
where $\e_\varrho$ denotes the primitive vector of $\varrho$, and $\langle \ccdot, \rdot\rangle$ is the duality pairing between $M$ and $N$.

\smallskip
This ring is known to be the Chow ring of the toric variety $\P_\Sigma$ defined by $\Sigma$, \cf.~\cites{Dan78, BDP90, Bri96, FS}.

\smallskip
In the case $\Sigma=\Sigma_\Ma$ is the Bergman fan of a simple matroid $\Ma$, by the description of the fan $\Sigma_\Ma$, it is easy to see that the ring $A^\bul(\Sigma_\Ma)$ coincides with the Chow ring $A^\bul(\Ma)$ of the matroid, introduced previously.

\subsection{Canonical compactifications} The \emph{canonical compactification} of a simplicial fan $\Sigma$ in a real vector space $V \simeq \R^n$ is defined as the result of completing any cone of $\Sigma$ to a \emph{hypercube}, by adding some faces at \emph{infinity}. This is defined more precisely as follows. The fan $\Sigma$ can be used to define a partial compactification of $V$ which we denote by $\TP_\Sigma$. The space $\TP_\Sigma$ is the tropical analog of the toric variety associated to $\Sigma$. In fact, in the case the fan $\Sigma$ is rational with respect to a full rank lattice $N$ in $V$, the tropical toric variety $\TP_\Sigma$ coincides with the tropicalization of the toric variety $\P_\Sigma$ associated to $\Sigma$, see e.g.~\cites{Payne, Thuillier}. These partial compactifications of $V$ coincide as well with the ones called \emph{manifolds with corners} considered in~\cites{AMRT, Oda}.\\
The canonical compactification of $\Sigma$ denoted by $\comp \Sigma$ is then defined as the closure of $\Sigma$ in the tropical toric variety $\TP_\Sigma$. The space $\comp\Sigma$ has the following description as the \emph{cubical completion} of $\Sigma$, see Section~\ref{sec:tropvar} for a more precise description. Any ray $\varrho$ of $\Sigma$ becomes a segment by adding a point $\infty_\varrho$ at infinity to $\varrho$. Any two dimensional cone $\sigma$ with two rays $\varrho_1$ and $\varrho_2$ becomes a parallelogram with vertex set the origin, $\infty_{\varrho_{1}}, \infty_{\varrho_2}$, and a new point $\infty_\sigma$ associated to $\sigma$. The higher dimensional cones are completed similarly to parallelepipeds.

\medskip

Canonical compactifications of fans serve as building blocks for the description and the study of more general tropical varieties, as we will explain below.

\subsection{Hodge isomorphism theorem and Poincar\'e duality} \label{sec:poincareduality}
Let $\Sigma$ be a unimodular rational Bergman fan and let $\comp \Sigma$ be the canonical compactification of $\Sigma$.
The compactification $\comp \Sigma$ is a \emph{smooth tropical variety} in the sense that it is locally modelled by unimodular rational Bergman fans, see Section~\ref{sec:smoothness_intro} below and Section~\ref{sec:tropvar} for a precise definition of smoothness. In this specific case, which concerns compactifications of Bergman fans, if the fan $\Sigma$ has support given by a simple matroid $\Ma$, which is not unique in general, the local fans associated to $\comp \Sigma$ have support a Bergman fan associated to \emph{restrictions} of \emph{quotients} of $\Ma$, or their products.

\medskip

Denote by $H^{p,q}_\trop(\comp \Sigma)$ the \emph{tropical cohomology group} of $\comp \Sigma$ of bidegree $(p,q)$ introduced by Itenberg-Katzarkov-Mikhalkin-Zharkov in~\cite{IKMZ}. We recall the definition of these groups in Section~\ref{sec:tropvar}. We prove in~\cite{AP} the following general theorem which allows in particular, when applied to the Bergman fan $\Sigma_\Ma$ of a matroid $\Ma$, to reinterpret the results of Adiprasito-Huh-Katz~\cite{AHK} as statements about the cohomology groups of a specific projective tropical variety, the canonical compactification of the Bergman fan of the matroid.

\begin{thm}[Hodge isomorphism for Bergman fans] \label{thm:HI} For any unimodular Bergman fan $\Sigma$ of dimension $d$, and for any non-negative integer $p\leq d$, there is an isomorphism
\[A^p(\Sigma) \longsimto H^{p,p}_\trop(\comp \Sigma, \Z).\]
These isomorphisms induce an isomorphism of $\Z$-algebras between the Chow ring of $\Sigma$ and the tropical cohomology ring $\bigoplus_{p=0}^d H_\trop^{p,p}(\comp \Sigma, \Z)$. Furthermore, the cohomology groups $H_\trop^{p,q}(\comp \Sigma, \Z)$ for $p\neq q$ are all trivial.
\end{thm}

Saying it differently, the theorem states that the \emph{cycle class map} from Chow groups to tropical cohomology is an isomorphism, and in particular, the tropical cohomology groups of $\comp \Sigma$ are generated by Hodge classes.

\smallskip
The theorem can be viewed as the tropical analog to a theorem of Feichtner and Yuzvinsky~\cite{FY}, which proves a similar theorem for wonderful compactifications of the complement of a hyperplane arrangement.

\smallskip
As a corollary of the above theorem, we obtain the following result.
\begin{thm}[Poincar\'e duality] \label{thm:pd} Let $\Sigma$ be a unimodular Bergman fan of dimension $d$. There is a canonical isomorphism
$ \deg\colon A^{d}(\Sigma)\to\Z$ called the degree map. With respect to this isomorphism, the Chow ring of $\Sigma$ satisfies the Poincar\'e duality, in the sense that the pairing
\[ \begin{tikzcd}[row sep=tiny, column sep=scriptsize]
A^k(\Sigma) \times A^{d-k}(\Sigma) \rar& A^{d}(\Sigma) \rar["\deg"', "\sim"]& \Z\\
(a,b) \arrow[mapsto]{rr} && \deg(ab)
\end{tikzcd} \]
is non-degenerate.
\end{thm}
This is in fact a consequence of Theorem~\ref{thm:HI}, the observation that the canonical compactification $\comp \Sigma$ is a smooth tropical variety, and the fact that smooth tropical varieties satisfy Poincar\'e duality, established by Jell-Shaw-Smacka~\cite{JSS} for rational coefficients and by Jell-Rau-Shaw~\cite{JRS19} for integral coefficients.

\medskip

We note that Theorem~\ref{thm:pd} is independently proved by Gross-Shokrieh~\cite{GS}.

\subsection{Hodge-Riemann and Hard Lefschetz for Bergman fans} \label{sec:local-intro} Theorem~\ref{thm:pd} generalizes Poincar\'e duality for matroids proved in~\cite{AHK} to any unimodular Bergman fan. In Section~\ref{sec:local} we will show that Chow rings of \emph{quasi-projective} unimodular Bergman fans verify Hard Lefschetz theorem and Hodge-Riemann relations, thus generalizing results of Adiprasito-Huh-Katz~\cite{AHK} to any quasi-projective unimodular Bergman fan. The results are as follows.

\smallskip
For a unimodular Bergman fan $\Sigma$ of dimension $d$ and an element $\ell \in A^1(\Sigma)$, we say that the pair $(\Sigma, \ell)$ verifies the \emph{Hard Lefschetz property $\HL(\Sigma, \ell)$} if the following holds:

\medskip

\noindent (Hard Lefschetz) for any non-negative integer $k \leq \frac d2$, the multiplication map by $\ell^{d-2k}$ induces an isomorphism
\[\begin{tikzcd}[column sep=large]
 A^k(\Sigma) \arrow{r}{\sim}& A^{d-k}(\Sigma)\,\,, \qquad a \arrow[mapsto]{r} & \ell^{d-2k}\cdot a.
\end{tikzcd}
\]

\medskip

We say the pair $(\Sigma, \ell)$ verifies \emph{Hodge-Riemann relations $\HR(\Sigma, \ell)$} if for any non-negative integer $k\leq \frac d2$, the following holds:

\medskip

\noindent (Hodge-Riemann relations) for any non-negative integer $k \leq \frac d2$, the symmetric bilinear form
\[\begin{tikzcd}[column sep=large]
 A^k(\Sigma) \times A^{k}(\Sigma) \arrow{r}& \Z \,\,, \qquad (a,b) \arrow[mapsto]{r} & (-1)^k\deg(\ell^{d-2k}\cdot a\cdot b)
\end{tikzcd}
\]
is positive definite on the \emph{primitive part} $P^k(\Sigma) \subseteq A^k(\Sigma)$, which by definition, is the kernel of the multiplication by $\ell^{d-2k+1}$ from $A^k(\Sigma)$ to $A^{d-k+1}(\Sigma)$.
\begin{thm}[K\"ahler package for Bergman fans]\label{thm-intro:HL-HR-local} Let $\Sigma$ be a quasi-projective unimodular Bergman fan of dimension $d$. For any ample element $\ell\in A^1(\Sigma)$, the pair $(\Sigma, \ell)$ verifies $\HL(\Sigma, \ell)$ and $\HR(\Sigma, \ell)$.
\end{thm}
Here $\ell$ is called \emph{ample} if $\ell$ can be represented in the form $\sum_{\varrho \in \Sigma_1} f(\e_\varrho) \x_\varrho$ for a strictly convex cone-wise linear function $f$ on the fan $\Sigma$, see Section~\ref{sec:tropvar} for the definition of strict convexity. A fan which admits such a function is called \emph{quasi-projective}. We note that the Bergman fan $\Sigma_\Ma$ of any matroid $\Ma$ is quasi-projective.

\smallskip
In addition to providing a generalization of~\cite{AHK} to any quasi-projective Bergman fan, our methods lead to arguably more transparent proofs of these results when restricted to Bergman fans of matroids treated in~\cite{AHK}. An alternate proof of the main result of~\cite{AHK} for the Chow ring of a matroid based on semi-simple decomposition of the Chow ring has been obtained independently by Braden-Huh-Matherne-Proudfoot-Wang~\cite{BHMPW}, as well as in the recent independent and simultaneous work of Ardila-Denham-Huh~\cite{ADH}.

\medskip

The main idea in our approach which allows to proceed by induction, basically by starting from trivial cases and deducing the theorem in full generality, are the \emph{ascent} and \emph{descent} properties proved in Section~\ref{sec:local}. These properties are also the final piece of arguments we will need in the global setting to conclude the proof of the main theorems of this paper. The possibility of having these nice properties is guaranteed by Keel's lemma for fans in the local setting, and by our projective bundle formula for tropical cohomology groups in the global setting, \cf. the discussion below and the corresponding sections for more details.

\subsection{Smoothness in tropical geometry} \label{sec:smoothness_intro} A tropical variety $\X$ is an \emph{extended polyhedral space} with a choice of an affine structure. A polyhedral complex structure on a tropical variety is a polyhedral complex structure which induces the affine structure. We refer to Section~\ref{sec:tropvar} which provides the precise definitions of these terminologies.

\smallskip
A tropical variety is called \emph{smooth} if it can be locally modeled by supports of Bergman fans \cites{IKMZ, MZ, JSS}. This means, locally, around any of its points, the variety looks like the support of a Bergman fan.

\smallskip
The smoothness is meant to reflect in polyhedral geometry the idea of \emph{maximal degeneracy} for varieties defined over non-Archimedean fields. Examples of smooth tropical varieties include smooth affine manifolds, canonical compactifications of Bergman fans, and tropicalizations of Mumford curves~\cite{Jell-MC}.

\smallskip
More generally, to get a good notion of smoothness in tropical geometry, one would need to fix, as in differential topology, a good class of fans and their supports which can serve as local charts. One such class which emerges from the results in our work and the work of Deligne~\cite{deligne-md} is that of fans which satisfy the Poincar\'e duality for their tropical cohomology groups and which have canonical compactifications with vanishing $H^{p,q}_\trop$ cohomology groups for all $p\neq q$. By Hodge isomorphism Theorem~\ref{thm:HI} and Poincar\'e duality~\cites{JSS, JRS19}, this class contains all Bergman fans, including therefore complete fans, and it would be certainly interesting to have a complete characterization of this class.

\subsection{Hodge-Riemann and Hard Lefschetz for smooth projective tropical varieties} In view of Hodge isomorphism Theorem~\ref{thm:HI} and K\"ahler package for Bergman fans, Theorem~\ref{thm-intro:HL-HR-local}, the results in Section~\ref{sec:local-intro} can be regarded as statements concerning the tropical cohomology ring of the canonical compactifications $\comp \Sigma$ of unimodular Bergman fans. Moreover, one can show that when $\Sigma$ is quasi-projective, the canonical compactification is projective~\cite{AP-geom}. One of our aim in this paper is to show that all the results stated above hold more generally for the tropical cohomology groups of any smooth projective tropical variety, and more generally, for K\"ahler tropical varieties, a concept we introduce in a moment.

\medskip

Let $V =N_\R$ be a vector space of dimension $n$ and $N$ a full rank sublattice in $V.$ For a fan $\Delta$ in $V$, we denote by $\TP_{\Delta}$ the corresponding \emph{tropical toric variety}, which is a partial compactification of $V$.

Let now $\Y$ be a polyhedral subspace of $V$, and let $Y$ be a unimodular polyhedral complex structure on $\Y$ with a quasi-projective unimodular \emph{recession fan} $Y_\infty$, \cf. Section~\ref{sec:tropvar} for the definitions of these terminologies. Let $\X$ be the closure of $\Y$ in $\TP_{Y_\infty}$, equipped with the induced polyhedral structure that we denote by $X$. Let $\ell$ be a strictly convex cone-wise linear function on $Y_\infty$. We will show that $\ell$ defines an element $\omega$ of $H_\trop^{1,1}(\TP_{Y_\infty}, \R)$, which lives in $H_\trop^{1,1}(\TP_{Y_\infty}, \Q)$ if $\ell$ has integral slopes. By restriction, this element gives an element in $H_\trop^{1,1}(\X, \Q)$ which, by an abuse of the notation, we still denote by $\omega$. Multiplication by $\omega$ defines a Lefschetz operator on $H_\trop^{\bul, \bul}(\X,\Q)$ still denoted by $\omega$.

\medskip

Here is the main theorem of this paper.

\begin{thm}[K\"ahler package for smooth projective tropical varieties]\label{thm:main2} Notations as above and working with $\Q$ coefficients, we have
\begin{itemize}
\item \emph{(weight-monodromy conjecture)} For any pair of non-negative integers $p\geq q$, we have an isomorphism
\[\phi\colon H_{\trop}^{p,q}(\X) \xrightarrow{\sim} H_{\trop}^{q,p}(\X).\]

\item \emph{(Hard Lefschetz)} For $p+q \leq d/2$, the Lefschetz operator $\ell$ induces an isomorphism
\[\ell^{d- p-q}\colon H_{\trop}^{p,q}(\X) \to H^{d-q, d-p}_{\trop}(\X).\]

\item \emph{(Hodge-Riemann)} The pairing $\bigl< \ccdot , \ell^{d-p-q} \phi\ccdot \bigr>$ induces a positive-definite symmetric bilinear form on the primitive part $P^{p,q}$ of $H_{\trop}^{p,q}$, where $\bigl< \ccdot,\rdot \bigr>$ denotes the pairing
\begin{equation*}
H^{p,q}_\trop(\X) \otimes H_\trop^{d-p,d-q}(\X) \longrightarrow H_\trop^{d,d}(\X) \xrightarrow[\raisebox{3pt}{$\sim$}]{\quad\deg\quad} \Q, \qquad
 \bigl< a,b \bigr> \longmapsto  (-1)^p\deg(a \cup b),
\end{equation*}
and $\phi$ is the isomorphism given in the first point.
\end{itemize}
\end{thm}

The first part of the theorem answers in positive a conjecture of Mikhalkin and Zharkov from~\cite{MZ} where they prove the isomorphism of $H^{p,q}_\trop(\X)$ and $H^{q,p}_\trop(\X)$ in the \emph{approximable case}. Recall that a projective tropical variety $\X$ is called \emph{approximable} if it is the \emph{tropical limit} of  a family of complex projective manifolds.

\medskip

Various versions of weak Lefschetz and vanishing theorems for projective tropical varieties were previously obtained by Adiprasito and Bj\"oner in~\cite{AB}. Our work answers several questions raised in their work, \cf.~\cite{AB}*{Section 11}, as well as~\cite{AB}*{Section 9}, and the results presented in the next sections.

\smallskip
We mention that an integral version of the weak Lefschetz theorem for hypersurfaces was obtained recently by Arnal-Renaudineau-Shaw in~\cite{ARS}.

\subsection{Hodge index theorem for tropical surfaces} Our main theorem in the case of smooth tropical surfaces implies the following tropical version of the Hodge index theorem. Let $\X$ be a smooth projective tropical surface. By Poincar\'e duality, we have $H^{2,2}_\trop(\X) \simeq \Q$, and the cup-product on cohomology leads to the pairing
\[\bigl< \ccdot , \rdot \bigr> \colon H^{1,1}_\trop(\X) \times H^{1,1}_\trop(\X) \to \Q,\]
which coincides with the intersection pairing on tropical homology group $H_{1,1}(\X)$.
\begin{thm}[Hodge index theorem for tropical surfaces]\label{thm:hodgeindex} The signature of the intersection pairing on $H^{1,1}_\trop(\X)$ is equal to $(1+b_2, h^{1,1}- 1- b_2)$ where $b_2$ is the second Betti number of $\X$ and $h^{1,1} =  \dim_\Q H^{1,1}_\trop(\X)$.
\end{thm}
This theorem was previously known  by explicit computations for an infinite family of surfaces in $\TP^3$, that of \emph{floor decomposed surfaces of varying degree $d$} in $\TP^3$, see~\cite{Sha13b}. In that case, the second Betti number counts the number of interior points in the standard simplex of width $d$, \ie, with vertices $(0,0,0), (d,0,0), (0,d,0), (0,0,d)$.

\medskip

\subsection{} What follows next gives an overview of the materials we need to establish in order to  conclude the proof of Theorem~\ref{thm:main2}. A sketch of our strategy appears in Section~\ref{sec:skecth_intro} and could be consulted at this point ahead of going through the details of the next coming sections.

\subsection{Tropical K\"ahler forms and classes} \label{sec:kahler-intro} The proof of the theorems above are all based on the notion of \emph{K\"ahler forms} and their corresponding \emph{K\"ahler classes} in tropical geometry that we introduce in this paper. Our definition is motivated by the study of the tropical Steenbrink sequence and its Hodge-Lefschetz structure, see below for a succinct presentation. A closely related notion of K\"ahler classes in non-Archimedean geometry has been introduced and studied in an unpublished work of Kontsevich-Tschinkel~\cite{KT} and in the work of Yu on non-Archimedean Gromov compactness~\cite{Yu} (see also~\cites{Bou, BFJ15, BGS, CLD, GK17, Zhang} for related work).

\medskip

Consider a smooth compact tropical variety $\X$ equipped with a unimodular polyhedral complex structure $X$. The polyhedral structure is naturally stratified according to \emph{sedentarity} of its faces, which is a measure of \emph{how far and in which direction at infinity} a point of $\X$ lies. We denote by $X_\f$ the \emph{finite part} of $X$ consisting of all the faces which do not touch the \emph{part at infinity}, \cf. Section~\ref{sec:tropvar} for the precise definition. A \emph{Kähler form} for $X$ is the data of ample classes $\ell^v$ in the Chow ring of the local fan around the vertex $v$, for each vertex $v\in X_\f$, such that for each edge $e=\{u,v\}$ of $X_\f$, the restriction of the two classes $\ell^v$ and $\ell^u$ in the Chow ring of the local fan around the edge $e$ are equal.

We show in Theorem~\ref{thm:KahlerClass} that any K\"ahler form defines a class in $H_\trop^{1,1}(\X)$. We call any class $\omega$ induced by a K\"ahler form coming from a unimodular triangulation a \emph{K\"ahler class}. A smooth compact tropical variety is then called \emph{K\"ahler} if it admits a Kähler class. In particular, our definition of K\"ahler requires the variety to admit a unimodular triangulation. We note however that it is possible to remedy this by working with tropical Dolbeault cohomology.

\smallskip
As in the classical setting, we prove the following theorem.
\begin{thm} A smooth projective tropical variety which admits a quasi-projective triangulation is K\"ahler.
\end{thm}
We believe that \emph{all smooth projective tropical varieties admit quasi-projective triangulations}, and plan to come back to this question in a separate publication. In the next section, we explain a weaker triangulation theorem we will prove in this paper which will be enough for our purpose.

\subsection{Quasi-projective and unimodular triangulations of polyhedral spaces} Our treatment of projective tropical varieties is based on the existence of what we call \emph{regular unimodular triangulations} of rational polyhedral spaces in $\R^n$. Regular triangulations are fundamental in the study of polytopes and their applications across different fields in mathematics and computer science, \cf. the book by De Lorea-Rambau-Santos~\cite{DRS} and the book by Gelfand-Kapranov-Zelvinsky~\cite{GKZ}, one of the pioneers to the field of tropical geometry, for the definition and its relevance in algebraic geometry.

\smallskip
Regularity for triangulations of a polytope is a notion of convexity, and as such, can be defined for any polyhedral subspace of $\R^n$. Generalizing the pioneering result of Kemp-Knudson-Mumford and Saint-Donat~\cite{KKMS} in the proof of the semistable reduction theorem, as well as previous variants proved by Itenberg-Kazarkov-Mikhalkin-Zharkov in~\cite{IKMZ} and W{\l}odarczyk~\cite{Wlo97}, we obtain the following theorem on the existence of unimodular  quasi-projective triangulations.

\begin{thm}[Triangulation theorem]\label{thm:regulartriangulations} Let $X$ be a rational polyhedral complex in $\R^n$. There exists a regular triangulation of $X$ which is quasi-projective and unimodular with respect to the lattice $\frac 1k \Z^n$, for some integer $k\in \mathbb N$.
\end{thm}

We expect stronger versions of the theorem, as well as a theory of \emph{secondary structures on regular triangulations}, and plan to come back to these questions in a future work.

\subsection{Canonical compactifications of polyhedral spaces}
Let $\Y$ be a polyhedral space in $V = N_\R$. The \emph{asymptotic cone} of $\Y$ which we denote by $\Y_\infty$ is defined as the pointwise limit of the rescaled subsets $\X/t$ when $t$ goes to $+\infty$. The term asymptotic cone is borrowed from geometric group theory, originating from the pioneering work of Gromov~\cites{Grom81, Grom83}, and exceptionally refers here to non-necessarily convex cones.

\medskip

Let $\Y$ be a smooth polyhedral space in $V = N_\R$. Let $\Delta$ be a unimodular fan structure on the asymptotic cone of $\Y_\infty$. Let $\X$ be the closure of $\Y$ in $\TP_\Delta$. The compactification $\X$ is then smooth. We call it the \emph{canonical compactification of $\Y$ with respect to $\Delta$}.

\smallskip
As in the theory of toric varieties, the tropical toric variety $\TP_\Delta$ has a natural stratification into  \emph{tropical toric orbits}. For each cone $\eta \in \Delta$, we have in particular the corresponding tropical toric subvariety denoted by $\TP_\Delta^\eta$ and defined as the closure of the torus orbit associated to $\eta$.

\medskip

For the canonical compactification $\X$ of $\Y$ as above, we define for any $\eta \in \Delta$ the closed stratum $D^\eta$ of $\X$ as the intersection $\X \cap \TP_\Delta^\eta$. With these definitions, we see that the compactification boundary $D := \X \setminus \Y$ is a \emph{simple normal crossing divisor} in $\X$, meaning that
\begin{itemize}
\item $D^{\conezero} = \X$, where $\conezero$ is the zero cone in $V$ consisting of the origin.
\item $D^\eta$ are all smooth of dimension $d-\dims{\eta}$.
\item We have \[D^\eta = \bigcap_{\substack{\varrho \subface \eta \\ \dims{\varrho}=1}} D^\varrho.\]
\end{itemize}

\subsection{Projective bundle formula} Notations as in the previous section, let $\delta$ be a cone in $\Delta$, and denote by $\Delta'$ the fan obtained by the barycentric star subdivision of $\Delta$. Denote by $\X'$ the closure of $\Y$ in $\TP_{\Delta'}$ and $\pi\colon \X' \to \X$ the projection map. The tropical variety $\X'$ is smooth again.

\smallskip
Let $\rho$ be the new ray in $\Delta'$ obtained after the star subdivision of $\delta \in \Delta$.  Consider the corresponding \emph{tropical divisor} $D'^\rho \subseteq \X'$, \ie, the closed stratum in $\X'$ associated to $\rho$.

\medskip

The divisor $D'^\rho$ defines by the tropical cycle class map an element $\class(D'^\rho)$ of $H^{1,1}_\trop(\X')$.

\smallskip
For any smooth compact tropical variety $\W$, and for each integer $k$, we define the $k$-th cohomology of $\W$ by
\[H^k(\W) := \sum_{p+q=k} H_\trop^{p,q}(\W).\]

\begin{thm}[Projective bundle formula] \label{thm:keelglobal-intro}
We have an isomorphism
\[H^k(\X' ) \simeq H^k(\X) \oplus T H^{k-2}(D^\delta) \oplus T^2 H^{k-4}(D^\delta) \oplus \dots \oplus T^{\dims{\delta}-1}H^{k-2\dims{\delta}+2}(D^\delta).\]
\end{thm}
Here the map from the right hand side to the left hand side restricts to $\pi^*$ on $H^k(\X)$ and sends $T$ to $-\class(D'^\rho)$. It is given on each factor $T^s H^{k-2s}(D^\delta)$ by
\[ \begin{array}{rcl}
T^s H^{k-2s}(D^\delta) & \longrightarrow & H^k(\X') \\[1em]
T^s \alpha & \longmapsto & (-1)^s \class(D'^\rho)^{s-1} \cup \pi^*\circ \gys(\alpha),
\end{array} \]
where $\pi^*$ and $\gys$ are the \emph{pull-back} and \emph{Gysin} maps for the tropical cohomology groups, with respect to the projection $\pi\colon \X' \to \X$ and the inclusion $D^\delta \hookrightarrow \X$, respectively.

\smallskip
The decomposition stated in the theorem provides for each pair of non-negative integers $(p,q)$,  a decomposition of the form
\[H^{p,q}(\X') \simeq H^{p,q}(\X) \oplus T H^{p-1, q-1}(D^\delta) \oplus T^2 H^{p-2, q-2}(D^\delta) \oplus \dots \oplus T^{\dims{\delta}-1}H^{p-\dims\delta+1,q - \dims\delta+1}(D^\delta).\]

\begin{remark}[Tropical Chern classes] The theorem describes the cohomology of the \emph{blow-up} $\X'$ of $\X$ along the subvariety $D^\delta \subset \X$. Our proof actually provides more generally a projective bundle theorem for the projective bundle associated to a vector bundle $E$ on a smooth tropical variety $\X$. As in the classical setting, it allows to define Chern classes of vector bundles in tropical geometry over smooth projective tropical varieties. \\
In the situation above, we get
\[T^{\dims{\delta}} + c_1(N_{D^\delta})T^{\dims\delta-1} +c_2(N_{D^\delta})T^{\dims\delta-2} +\dots + c_{\dims \delta}(N_{D^\delta}) =0 \]
where $c_i$ are the Chern classes of the \emph{normal bundle} $N_{D^\delta}$ in $\X$.\\
Chern classes of matroids were previously defined and studied by L{\'o}pez de Medrano, Rinc{\'o}n, and Shaw in~\cite{LRS}.
\end{remark}

\subsection{Tropical Steenbrink sequence} Let $\X$ be a smooth compact tropical variety and let $X$ be a unimodular polyhedral structure on $\X$. Denote by $X_\f$ the set of faces of $X$ whose closures do not intersect the boundary at infinity of $\X$, \ie, the set of compact faces of $X_\conezero$.

\medskip

Inspired by the shape of the first page of the Steenbrink spectral sequence~\cite{Ste76}, we define bigraded groups $\ST_1^{a,b}$ with a collection of maps between them as follows.

\medskip

For all pair of integers $a, b \in \mathbb Z$, we define
\[ \ST_1^{a,b} := \bigoplus_{s \geq \abs a \\ s \equiv a \pmod 2} \ST_1^{a,b,s} \]
where
\[ \ST_1^{a,b,s} = \bigoplus_{\delta \in X_\f \\ \dims\delta =s} H^{a+b-s}(\comp \Sigma^\delta). \]
Here $\Sigma^\delta$ is the \emph{star fan} of $X$ around $\delta$ (also called the \emph{transversal fan} of $\delta$ in the literature), and $\comp \Sigma^\delta = \comp{\Sigma^\delta}$ is the canonical compactification of $\Sigma^\delta$.

\medskip

The bigraded groups $\ST_1^{a,b}$ come with a collection of maps
\[\i^{a,b\,*} \colon \ST^{a,b}_1 \to \ST_1^{a+1, b} \qquad \textrm{and} \qquad \gys^{a,b}\colon \ST^{a,b}_1 \to \ST_1^{a+1, b}.\]
Both these maps are defined by our sign convention introduced in Section~\ref{sec:steenbrink} as an alternating sum of the corresponding maps on the level of cohomology groups appearing in the definition of $\ST_1^{a,b,s}$ above. In practice, we drop the indices and denote simply by $\i^*$ and $\gys$ the corresponding maps.

\medskip

Using these two maps, we define the \emph{differential} $\d\colon \ST_1^{a,b} \to \ST_1^{a+1,b}$ as the sum $\d = \i^*+ \gys$. For a unimodular triangulation $X$ of $\X$ and for any integer $b$, we will show that the differential $\d$ makes $\ST_1^{\bul,b}$ into a cochain complex.

\smallskip
For a cochain complex $(C^\bul, \d)$, denote by $H^a(C^\bul, \d)$ its $a$-th cohomology group, i.e.,
\[H^a(C^\bul, \d) = \frac{\ker\Bigl(\d\colon \, C^{a} \rightarrow C^{a+1}\Bigr)}{\Im\Bigl(\d\colon \, C^{a-1} \rightarrow C^{a}\Bigr)}. \]

We prove the following comparison theorem.
\begin{thm}[Steenbrink-Tropical comparison theorem] \label{thm:steenbrink-intro}
The cohomology of $(\ST_1^{\bul,b},\d)$ is described as follows. For $b$ odd, all the terms $\ST_1^{a,b}$ are zero, and the cohomology is vanishing. For $b$ even, let $b=2p$ for $p \in \mathbb Z$. Then for $q\in \mathbb Z$, we have a canonical isomorphism
\[H^{q-p}(\ST^{\bul,2p}_1, \d) \simeq H^{p,q}_{\trop}(\X).\]
\end{thm}

\medskip

In the approximable case, i.e., when $\X$ arises as the tropicalization of a family of complex projective varieties, this theorem was proved by Itenberg-Katzarkov-Mikhalkin-Zharkov in~\cite{IKMZ}; see also the work of Gross-Siebert~\cites{GS10, GS06} for a similar statement for the special case where $\X$ is an integral affine manifold with singularities.

\medskip

Our proof is inspired by the one given in the approximable case~\cites{IKMZ}, as well as by the use of the sheaf of logarithmic differentials on both Deligne's construction of a mixed Hodge structure on the cohomology of a smooth algebraic variety~\cite{Deligne-Hodge2} and Steenbrink's construction of the limit mixed Hodge structure on the special fiber of a semistable degeneration~\cite{Ste76}.

\medskip

In order to prove the theorem in this generality, we will introduce and develop a certain number of tools which we hope could be of independent interest and which will be used in our forthcoming work.

\smallskip
In particular, we will prove first an analogous result in the tropical setting of the \emph{Deligne resolution} which gives a resolution of the coefficient groups $\SF^p$ with cohomology groups of the canonically compactified fans $\comp \Sigma^\delta$. (The coefficient group $\SF^p$ is the discrete tropical analog of the sheaf of holomorphic forms of degree $p$, and is used to define the tropical cohomology groups $H^{p,q}(\X)$ for all non-negative integer $q$.) \\
In the approximable case, this resolution is a consequence of the Deligne spectral sequence in mixed Hodge theory as was observed in~\cite{IKMZ}. The proof we present for the general setting is based on the hypercohomology of a complex of sheaves which provides a resolution of the \emph{sheaf of tropical holomorphic forms} of a given order. It further uses Poincar\'e duality as well as our Hodge isomorphism Theorem~\ref{thm:HI}, which provides a description of the tropical cohomology groups of the canonically compactified fans.

\smallskip
Inspired by the weight filtration on the sheaf of logarithmic differentials, we define a natural weight filtration on the coefficient groups $\SF^p(\ccdot)$, also somehow explicitly present in~\cite{IKMZ},  and study the corresponding spectral sequence on the level of tropical cohomology groups. The resolution of the coefficient groups given by the Deligne exact sequence gives a double complex which allows to calculate the cohomology of the graded cochain complex associated to the weight filtration.

\smallskip
We then show that the spectral sequence associated to the double complex corresponding to the weight filtration which abuts to the tropical cohomology groups, abuts as well to the cohomology groups of the Steenbrink cochain complex.  The proof of this latter fact is based on an \emph{unfolding of the Steenbrink sequence into a double complex}, which allows to define new Steenbrink type spectral sequences on each degree, and a \emph{spectral resolution lemma} which allows to make a bridge between spectral sequences.

\medskip

While the treatment we give of these constructions might appear quite technical, we would like to note that this should be merely regarded as a manifestation of the rich combinatorial structure of the Steenbrink spectral sequence, and the geometric informations it contains. We hope the effort made to make a \emph{surgery} of this spectral sequence in this paper should paid off in further applications and developments, some of them will appear in our forthcoming work.

\subsection{Main theorem for triangulated tropical K\"ahler varieties}
In order to prove Theorem~\ref{thm:main2}, we first establish it for a unimodular triangulation of a smooth tropical variety which admits a K\"ahler form. The advantage of these triangulations is that we can use the tropical Steenbrink spectral sequence in order to study the tropical cohomology groups. We show that this spectral sequence can be endowed with a \emph{Hodge-Lefschetz structure}, where the monodromy operator is the analog of the one arising in the classical Steenbrink spectral sequence and the Lefschetz operator corresponds to multiplication with the K\"ahler form. This leads to the following theorem.

\begin{thm}\label{thm:main1-intro} Let $X$ be a unimodular triangulation of a smooth tropical variety which admits a K\"ahler form given by classes $\ell^v \in H^{1,1}(v)$ for $v$ vertices of the triangulation. Denote by $\ell$ the corresponding Lefschetz operator and by $N$ the corresponding monodromy operator. We have
\begin{itemize}
\item \emph{(Weight-monodromy conjecture)} For $q>p$ two non-negative integers, we get an isomorphism
\[N^{q-p} \colon H_{\trop}^{q,p}(X) \to H_{\trop}^{p,q}(X).\]

\item \emph{(Hard Lefschetz)} For $p+q \leq d/2$, the Lefschetz operator $\ell$ induces an isomorphism
\[\ell^{d- p-q}\colon H_{\trop}^{p,q}(X) \to H^{d-q, d-p}_{\trop}(X).\]

\item \emph{(Hodge-Riemann)} The pairing $(-1)^p \bigl< \ccdot, \ell^{d-p-q} N^{q-p}\ccdot\bigr>$ induces a positive-definite pairing on the primitive part $P^{q,p}$, where $\bigl< \ccdot, \rdot \bigr>$ is the natural pairing

\[\bigl< \ccdot , \rdot\bigr> \colon H^{q,p}(X) \otimes H^{d-q,d-p}(X) \to H^{d,d}(X) \simeq \mathbb Q.\]
\end{itemize}
\end{thm}

\subsection{Hodge-Lefschetz structures} In order to prove the theorem of the previous section, we show that $\ST_1^{\bul, \bul}$ and the tropical cohomology $H^{\bul,\bul}_\trop$ admit a Hodge-Lefschetz structure. To define this, we need to introduce the monodromy and Lefschetz operators on $\ST_1^{\bul, \bul}$.

The monodromy operator is of bidegree $(2,-2)$ and it is defined in such a way to mimic the one obtained from the Steenbrink spectral sequence if we had a semistable family of complex varieties.
So it is the data of maps $N^{a,b} \colon \ST_1^{a,b} \to \ST_1^{a+2, b-2}$ defined as
\[N^{a,b} = \bigoplus_{s \geq \abs a \\ s \equiv a \mod 2} N^{a,b,s}\]
with
\begin{equation}
N^{a,b, s} = \begin{cases} \mathrm{id} \colon\ST_1^{a,b,s} \to \ST_1^{a+2,b-2,s}   & \textrm{if } s \geq \abs{a+2}, \\
0 & \textrm{otherwise}.
\end{cases}
\end{equation}
Moreover, it coincides with the monodromy operator on the level of Dolbeault cohomology groups defined in~\cite{Liu19}, via our Steenbrink-Tropical comparison Theorem~\ref{thm:steenbrink-intro} and the comparison theorem between Dolbeault and tropical cohomology groups proved in~\cite{JSS}.

\medskip

The definition of the Lefschetz operator depends on the choice of the K\"ahler form $\ell$. Recall that a K\"ahler form is the data of ample classes $\ell^v$ at vertices $v$ which are moreover compatible along edges. This guaranties the compatibility of restrictions to higher dimensional faces giving an ample element $\ell^\delta \in H^2(\delta) = A^1(\delta)$ for all faces $\delta \in X_\f$.

\medskip

The Lefschetz operator $\ell^{a,b} \colon \ST_1^{a,b} \to \ST_1^{a, b+2} $ is then defined as the sum \[\ell^{a,b} = \bigoplus \ell^{a,b,s} \colon \bigoplus_{s \geq \abs{a} \\
s \equiv a \pmod 2} \ST_1^{a,b, s} \to \bigoplus_{s \geq \abs{a} \\
s \equiv a \pmod 2} \ST_1^{a,b+2, s}\]
where $\ell^{a,b,s}$ is the sum of the contributions of local ample elements
\[\ell^{a,b,s} = \bigoplus \ell^\delta \colon \bigoplus_{\delta \in X_\f \\ \dims{\delta} =s} H^{a+b-s}(\delta) \longrightarrow    \bigoplus_{\delta \in X_\f \\ \dims{\delta} =s} H^{a+b+2-s}(\delta).\]
In practice, the indices are dropped and this is denoted simply  by $\ell$.

\smallskip
The Monodromy and Lefschetz operators $N$ and $\ell$ verify the following properties
\begin{itemize}
\item $[\ell, N] = 0$, $[\ell, \i^*] = [\ell, \gys] =0$, and $[N,\i^*] = [N, \gys] =0$.

\smallskip
\item for a pair of integers $a,b$ with $a+b \geq d$, the map
\[\ell^{d-a-b} = \underbrace{\ell \circ \ell \circ \dots \circ \ell}_{(d-a-b) \,\textrm{times}} \colon \ST_1^{a,b} \longrightarrow \ST_1^{a, 2d - 2a-b}\]
is an isomorphism.

\smallskip
\item for a pair of integers $a,b$ with $a\leq 0$, the map
\[N^{-a} = \underbrace{N \circ N \circ \dots \circ N}_{(-a) \,\textrm{times}} \colon  \ST_1^{a,b} \longrightarrow \ST_1^{-a, b + 2a}\] is an isomorphism.
\end{itemize}

These properties show that the tropical Steenbrink spectral sequence with the two operators $\ell$ and $N$ form a Hodge-Lefschetz structure. Moreover, we show that local polarizations allow to define a polarization on $\ST_1^{\bullet, \bullet}$ leading to a polarized Hodge-Lefschetz structure.
From this, we deduce that the cohomology of the tropical Steenbrink spectral sequence admits a polarized Hodge-Lefschetz structure. Using this, and combining the corresponding \emph{primitive decomposition} with the comparison Theorem~\ref{thm:steenbrink-intro}, we can finish the proof of Theorem~\ref{thm:main1-intro}.

\medskip

Differential Hodge-Lefschetz structures are treated in Saito's work on Hodge modules~\cite{Saito}, in the paper by Guill\'en and Navarro Aznar on invariant cycle theorem~\cite{GNA90} and in the upcoming book by Sabbah and Schnell~\cite{SabSch} on mixed Hodge modules. We should however note that the tropical set-up is slightly different, in particular, the differential operator appearing in our setting is skew-symmetric with respect to the polarization. The proof we give of these results is elementary and does not make any recourse to representation theory, although it might be possible to recast in the language of representation theory the final combinatorial calculations we need to elaborate.

\medskip

We refer to Sections~\ref{sec:kahler} and \ref{sec:differential_HL_structure} for more details.

\subsection{Proof of Theorem~\ref{thm:main2}}\label{sec:skecth_intro}
Having explained all the needed ingredients, we now explain how the proof of the main theorem can be deduced. Starting with the projective tropical variety $\X$, and the quasi-projective unimodular fan structure $\Delta$ induced on $\X_\infty$ by the compactification we proceed as follows. The choice of the convex piecewise function $\ell$ on $\Delta$ gives an element $\ell \in H^{1,1}(\X)$. Using the ascent and descent properties, that we extend naturally from the local situation to the global setting using the projective bundle theorem, we show that it will be enough to treat the case where $\X$ admits a quasi-projective unimodular triangulation $X$ compatible with $\Delta$.  In this case, we show that the class $\ell \in H^{1,1}(\X)$ is K\"ahler, and so we can apply Theorem~\ref{thm:main1-intro} to conclude.

\subsection{Organization of the paper} The paper is organized as follows. In the next section, we provide the relevant background on tropical varieties, introduce canonical compactifications, and recall the definition of cohomology groups associated to them. This section introduces as well the corresponding terminology in polyhedral geometry we will be using all through the paper.

In Section~\ref{sec:local} we present the generalization of the results of Adiprasito-Huh-Katz to any unimodular Bergman fan, thus proving Theorem~\ref{thm-intro:HL-HR-local}. These local results are the basis for all the materials which will appear in the upcoming sections. Moreover, the ascent and descent properties used in the proof of local Hodge-Riemann relations will be again crucial in our proof of the global Hodge-Riemann relations.

Section~\ref{sec:triangulation}, which is somehow independent of the rest of the paper, is devoted to the proof of our triangulation theorem, cf. Theorem~\ref{thm:regulartriangulations}. The results are crucial in that they allow to introduce a tropical analog to the Steenbrink spectral sequence, which will be the basis for all the results which come thereby after.

In Section~\ref{sec:steenbrink} we study the tropical Steenbrink spectral sequence, introduce the weight filtration on tropical coefficient groups $\SF^p$, and prove the tropical Deligne resolution theorem as well as the comparison theorem between Steenbrink sequence and tropical cohomology groups, cf. Theorem~\ref{thm:steenbrink-intro}. For the ease of reading, few computational points of this section are treated separately in Appendix.

K\"ahler tropical varieties are introduced in Section~\ref{sec:kahler}. For a unimodular triangulation of a smooth tropical varieties which admits a K\"ahler form, the corresponding tropical Steenbrink spectral sequence can be endowed with a Hodge-Lefschetz structure, where the monodromy operator is the analog of the one arising in the classical Steenbrink spectral sequence and the Lefschetz operator corresponds to multiplication with the K\"ahler class. Using an adaptation of the theory of Hodge-Lefschetz structure to the tropical setting, Sections~\ref{sec:kahler} and~\ref{sec:differential_HL_structure} are devoted to the proof of the K\"ahler package for triangulated tropical varieties, cf. Theorem~\ref{thm:main1-intro}. The fact that the Hodge-Lefchetz structure will be inherited in the cohomology in the tropical setting is treated in Section~\ref{sec:differential_HL_structure}.

The last two sections are then devoted to the proof of Theorem~\ref{thm:main2}. In Section~\ref{sec:projective_bundle_theorem}, we  present the proof of the projective bundle theorem, cf. Theorem~\ref{thm:keelglobal-intro}. Section~\ref{sec:proofmaintheorem} then finishes the proof of the main theorem using the materials developed in the previous sections, ascent-descent property, projective bundle theorem and Theorem~\ref{thm:main1-intro}.

\medskip

The results presented in Appendix show that a certain triple complex constructed in Section~\ref{sec:steenbrink} provides a spectral resolution of the tropical spectral sequence, associated to the weight filtration on the tropical complex. Due to the calculatory nature of the content of this section, and necessity of introducing extra notations, we have decided to include it only as an appendix. It is also written somehow independently of the rest of the paper. This means we reproduce some of the materials from the paper there, and in few places, we adapt a terminology slightly different from the one used in the main body. We have found this more adapted both in terms of rigour and compactness to the purpose of this appendix which requires through case by case computations. Every time this change of terminology happens, we make a comment of comparison to the one used in the previous sections.


\section{Tropical varieties}
\label{sec:tropvar}

The aim of this section is to provide the necessary background on polyhedral spaces and tropical varieties. It contains a brief review of polyhedral geometry, matroids and their associated Bergman fans, tropical and Dolbeault cohomology, as well as a study of canonical compactifications of polyhedral spaces, which will be all used in the consequent sections. This will be also the occasion to fix the terminology and notations which will be used all through the text.

\subsection{Recollection in polyhedral geometry} \label{sec:recol}
A (closed convex) \emph{polyhedron} of a real vector space $V\simeq \R^n$ is a non-empty intersection of a finite number of affine half-spaces. All through this paper, a polyhedron means a \emph{strongly convex} one in the sense that we require it does not contain any affine line.

Let $P$ be a polyhedron. The \emph{affine tangent space} of $P$ denoted by $\Tan{P}$ is by definition the smallest affine subspace of $\R^n$ which contains $P$; it is precisely the set of all linear combinations of elements of $P$ whose coefficients sum up to one. The \emph{linear tangent space}, or simply \emph{tangent space}, of $P$ denoted by $\TT P$ is the linear subspace of $\R^n$ spanned by the differences $x-y$ for any pair of elements $x,y$ of $P$.

The \emph{dimension} of $P$ is denoted by $\dims{P}$, which is by definition that of its tangent space. A \emph{face} of $P$ is either $P$ or any nom-empty intersection $P\cap H$ for any affine hyperplane $H$ such that $P\subset H^+$ where $H^+$ is one of the two half-space delimited by $H$. Note that a face of a polyhedron $P$ is itself a polyhedron. We use the notation $\gamma\subface\delta$ for two polyhedra if $\gamma$ is a face of $\delta$. If moreover $\gamma$ is of codimension one in $\delta$, \ie, $\dims\gamma=\dims\delta-1$, then we write $\gamma\ssubface\delta$. A face of dimension zero is called a \emph{vertex} of $P$. A face of dimension one is called an \emph{edge}. The \emph{interior} of a polyhedron $P$ refers to the relative interior of $P$ in its affine tangent space, which is by definition the complement in $P$ of all the proper faces of $P$.

A \emph{cone} is a polyhedron with a unique vertex which is the origin of $\R^n$. (By our assumption, we only consider strongly convex cones, which means that the cone does not include any line.) A compact polyhedron is called a \emph{polytope}; equivalently, a polytope is a polyhedron which is the convex-hull of a finite number of points.

Suppose now that the vector space $V \simeq \R^n$ comes with a full rank lattice $N \simeq \Z^n $. A polyhedron $P$ in $V$ is called \emph{rational} if the half-spaces used to define $P$ can be defined in $N_\Q \simeq \Q^n$. We denote by $N_P$ the full-rank lattice $N\cap \TT(P)$ in $\TT(P)$.

A polyhedron is called \emph{integral with respect to $N$} if it is rational and all its vertices are in $N$. In the sequel, we simply omit to mention the underlying lattice $N$ if this is clear from the context.

The Minkowski sum of two polyhedra $P$ and $Q$ in $V$ is defined by
\[P+Q \colon = \bigl\{x+ y \st x\in P, y \in Q\Bigr\}.\]

By Minkowski-Weyl theorem, every polyhedron $P$ can be written as the Minkowski sum $Q+\sigma$ of a polytope $Q$ and a cone $\sigma$. Moreover, one can choose $Q$ to be the convex-hull of the vertices of $P$ and the cone $\sigma$ is uniquely characterized by the property that for any point $x\in P$, $\sigma$ is the maximal cone verifying $x+\sigma\subseteq P$. We denote this cone by $P_\infty$.

It follows that for any polyhedron $P$, one can choose points $v_0,\dots,v_k$ and vectors $u_1,\dots,u_l$ in $V$ such that
\[ P=\conv(v_0,\dots,v_k)+\sum_{i=1}^l\R_+u_i. \]
Here and in the whole article, $\R_+$ denotes the space of non-negative real numbers.

A polyhedron $P$ is called \emph{simplicial} if we can choose the points $v_i$ and the vectors $u_j$ in the above decomposition so that the collection of vectors $(v_1-v_0,\dots,v_k-v_0,u_1,\dots,u_l)$ be independent.
In this case, there is a unique choice for points $v_0,\dots, v_k$, which will be the vertices of $P$. (And the vectors $u_1, \dots, u_l$ are in bijection with the rays of $P_\infty$.) We define
\[P_\f:=\conv(v_0,\dots, v_k),\]
and we have $P = P_\f + P_\infty$. Furthermore, any point $x\in P$ can be decomposed in a unique way as the sum $x=x_\f+x_\infty$ of a point $x_\f\in P_\f$ and a point $x_\infty\in P_\infty$.

A simplicial polyhedron $P$ is called \emph{unimodular with respect to the lattice} $N$ in $V$ if it is integral with respect to $N$ and moreover, the points $v_i$ and the vectors $u_j$ can be chosen in $N$ so that $(v_1-v_0,\dots,v_k-v_0,u_1,\dots,u_l)$ is part of a basis of $N$ (or equivalently, if they form a basis of the lattice $N_P = N\cap\TT(P)$).

\medskip

A \emph{polyhedral complex} in a real vector space $V \simeq \R^n$ is a finite non-empty set $X$ of polyhedra in $\R^n$ called \emph{faces of $X$} such that for any pair of faces $\delta,\delta'\in X$
\begin{enumerate}[label=\defaultRoman]
\item \label{enum:subface} any face of $\delta$ belongs to $X$, and
\item the intersection $\delta\cap\delta'$ is either empty or is a common face of $\delta$ and $\delta'$.
\end{enumerate}
A finite collection of polyhedra $X$ which only verifies the first condition \ref{enum:subface} is called a \emph{polyhedral pseudo-complex}.

If $P$ is a polyhedron, the set of faces of $P$ form a polyhedral complex which we denote by $\face{P}$. A polyhedral complex which has a unique vertex is called a \enquote{\emph{fan}} if the unique vertex is the origin of $V$. In this case, all the faces of the complex are cones. The \emph{finite part of $X$} denoted by $X_\f$ is the set of all compact faces of $X$, which is itself again a polyhedral complex. The \emph{support of $X$} denoted by $\supp{X}$ is the union of the faces of $X$ in $V$. A polyhedral complex whose support is the entire vector space $V$ is said to be \emph{complete}.

A polyhedral complex is called \emph{simplicial}, \emph{rational}, \emph{integral}, or \emph{unimodular} (with respect to the lattice $N$), if all faces of $X$ are simplicial, rational, integral, or unimodular (with respect to $N$), respectively.

Let $X$ and $Y$ be two polyhedral complexes in a real vector space $V$. We say that $Y$ is a \emph{subdivision} of $X$ if $X$ and $Y$ have the same support, and each face of $Y$ is included in a face of $X$. We say that $Y$ is a \emph{triangulation} of $X$ if in addition $Y$ is simplicial. We say that $Y$ is a \emph{subcomplex} of $X$ or that $X$ \emph{contains} $Y$ if $Y \subseteq X$. In particular, the support of $X$ contains that of $Y$. We say that $X$ \emph{contains a subdivision of} $Y$ if there exists a subdivision $Z$ of $Y$ which is a subcomplex of $X$. The same terminology will be used for polyhedral pseudo-complexes.

\medskip

Let $S$ be a subset of a real vector space $V \simeq \R^n$. The set of \emph{affine linear functions on $S$} is defined as the restrictions  to $S$ of affine linear functions on $V$ and is denoted by $\aff(S)$. For a polyhedral complex, we simplify the notation and write $\aff(X)$ instead of $\aff(\supp X)$.

Let $f\colon \supp{X} \to \R$ be a continuous function. We say that $f$ is \emph{piecewise linear on $X$} if on each face $\delta$ of $X$, the restriction $f\rest \delta$ of $f$ to $\delta$ is affine linear. The set of piecewise linear functions on $X$ is denoted by $\lpm(X)$. ($\lpm$ reads \emph{lin\'eaire par morceaux}.)

A function $f\in\lpm(X)$ is called \emph{convex}, resp. \emph{strictly convex}, if for each face $\delta$ of $X$, there exists an affine linear function $\ell \in \aff(V)$ such that $f-\ell$ vanishes on $\delta$ and is non-negative, resp. strictly positive, on $U \setminus \delta$ for an open neighborhood $U$ of the relative interior of $\delta$ in $\supp{X}$.

We denote by $\K(X)$, \resp, $\K_+(X)$, the set of convex, resp. strictly convex, functions in $\lpm(X)$. A polyhedral complex $X$ is called \emph{quasi-projective} if the set $\K_+(X)$ is non-empty. We will treat this notion of convexity in more detail in Section~\ref{subsec:convexite}.

\subsection{Canonical compactifications of fans} \label{sec:can-compact-fans-2}

Let $\eR = \R \cup \{\infty\}$ be the extended real line with the topology induced by that of $\R$ and a basis of open neighborhoods of infinity given by intervals $(a, \infty]$ for $a\in \R$. Extending the addition of $\R$ to $\eR$ in a natural way, by setting $\infty + a = \infty$ for all $a \in \eR$, gives $\eR$ the structure of a topological monoid. We call $(\eR, +)$ the \emph{monoid of tropical numbers} and denote by $\eR_+ = \R_+ \cup\{\infty\}$ the submonoid of non-negative tropical numbers with the induced topology.

Both monoids admit a natural scalar multiplication by non-negative real numbers (setting $0\cdot\infty=0$). Moreover, the multiplication by any factor is continuous. As such, $\eR$ and $\eR_+$ can be seen as module over the semiring $\R_+$. Recall that modules of semirings are quite similar to classical modules over rings, except that such modules are commutative monoids instead of being abelian groups. Another important collection of  examples of topological modules over $\R_+$ are the cones.

We can consider the tensor product of two modules over $\R_+$. In this section, every tensor product will be over $\R_+$, thus we will sometimes omit to mention it.

\medskip

All through, $N$ will be a free $\Z$-module of finite rank and $M=N^{\vee}$ denotes the dual of $N$. We denote by $N_\R$ and $M_\R$ the corresponding real vectors spaces, so we have $M_\R =N_\R^\vee$. For a polyhedron $\delta$ in $N_\R$, we use both the notations $N_{\delta, \R}$ or $\TT\sigma$, depending on the context, to denote the \emph{linear tangent space} to $\delta$ which is the real vector subspace of $N_\R$ generated by differences $x-y$ for pairs of elements $x, y$ in $\delta$. Furthermore, we define the \emph{normal vector space} of $\delta$ denoted $N^{\delta}_{\R}$ by $N^{\delta}_{\R} := \rquot{N_\R}{N_{\delta, \R}}$. If the polyhedron $\delta$ is rational, then we naturally get lattices of full rank $N^\delta$ and $N_\delta$ in $N^\delta_\R$ and $N_{\delta, \R}$, respectively.

\medskip

By convention, all through this article, we most of the time use $\delta$ (or any other face) as a superscript to denote the quotient of some space by $\TT\delta$, or to denote elements related to this quotient. On the contrary, we use $\delta$ as a subscript for subspaces of $\TT\delta$ or elements associated to these subspaces.

\medskip

For a fan $\Sigma$ in $N_\mathbb R$, we denote by $\Sigma_k$ the set of $k$-dimensional cones of $\Sigma$; elements of $\Sigma_1$ are called \emph{rays}. In the case the fan is rational, for any ray $\varrho\in\Sigma_1$, we denote by $\e_\varrho$ the generator of $\varrho\cap N$.

\medskip

For any cone $\sigma$, denote by $\sigma^\vee \subseteq M_\R$ the dual cone defined by
\[\sigma^\vee := \Bigl\{m \in M_\mathbb R \st \langle m, a \rangle \geq 0 \ \textrm{ for all }\ a \in \sigma\Bigr\}. \]
We also define the orthogonal plane to $\sigma$ by
\[\sigma^\perp := \Bigl\{m \in M_\mathbb R \st \langle m, a \rangle = 0 \ \textrm{ for all }\ a \in \sigma\Bigr\}.  \]

\medskip

The \emph{canonical compactification} $\comp \sigma$ of a cone $\sigma$, called as well the \emph{extended cone} of $\sigma$, is defined by
\[\comp \sigma := \sigma \otimes_{\R_+} \eR_+. \]
The topology on $\comp\sigma$ is the finest one such that the endomorphisms
\[ z\mapsto z+z', \qquad z\mapsto x\otimes z \quad\myand\quad x\mapsto z\otimes a \]
are continuous for any $z'\in\comp\sigma$, $x\in\sigma$ and $a\in\eR_+$.
The space $\comp\sigma$ is a compact topological space, whose restriction to $\sigma$ coincides with the usual topology.

\begin{remark} \label{rem:definition_dual}
One can also define $\comp\sigma$ in the following way. Let $\tm$ be the category of topological monoids whose morphisms are continuous. Then
\[ \comp\sigma = \Hom_\tm(\sigma^\vee, \eR_+). \]
The advantage of this definition is that the link with compactification of toric varieties is more direct (cf. below). However, the description of the topology and of the general shape of $\comp\sigma$ is easier with the first definition.
\end{remark}

There is a distinguished point $\infty_\sigma$ in $\comp \sigma$ defined by $x\otimes\infty$ for any $x$ in the relative interior of $\sigma$. The definition does not depend on the chosen $x$. Note that for the cone $\conezero$, we have $\infty_{\conezero} = 0$.

For an inclusion of cones $\tau \subface \sigma$, we get a map $\comp \tau \subseteq \comp \sigma$. This inclusion identifies $\comp \tau$ as the topological closure of $\tau$ in $\comp \sigma$.

\medskip

Let $\Sigma$ be a rational fan in $N_\R$. The canonical compactification $\comp \Sigma$ is defined as the union of extended cones $\comp \sigma$ for any cone $\sigma\in\Sigma$ where we identify $\comp \tau$ with the corresponding subspace of $\comp \sigma$ for any inclusion of cones $\tau \subface \sigma$ in $\Sigma$. The topology of $\comp \Sigma$ is the induced quotient topology. Each extended cone $\comp \sigma$ naturally embeds as a subspace of $\comp \Sigma$.

\medskip

The canonical compactification $\comp \Sigma$ of $\Sigma$ naturally lives in a \emph{partial compactification} of $N_\R$ defined by $\Sigma$ that we denote by $\TP_\Sigma$. We define $\TP_\Sigma$ as follows. For any cone $\sigma$ in $\Sigma$, we consider the space $\~\sigma=(\sigma + \sigma^\perp) \otimes_{\R_+} \eR$, endowed with the finest topology such that sum and the wedge product are continuous. We have a natural inclusion of $\comp\sigma$ into $\~\sigma$. Notice that $(\sigma+\sigma^\perp)\otimes\R\simeq N_\R$. We set $N^\sigma_{\infty,\R}:=\infty_\sigma + N_\R \subset \~\sigma$. Clearly $N^\conezero_{\infty,\R}=N_\R$. Moreover, from the fact that for any $x\in\sigma$ and any $\lambda\in\R$, $\infty_\sigma + x\otimes \lambda=\infty_\sigma$, we can see the map
\[ N_\R \to N^\sigma_{\infty,\R}, \qquad z \mapsto z+\infty_\sigma \]
as a projection along $\sigma\otimes\R\simeq N_{\sigma,\R}$. Indeed, one can prove that $N^\sigma_{\infty,\R}\simeq N^\sigma_{\R}$.

Then, $\~\sigma$ is naturally stratified into a disjoint union of subspaces $N^\tau_{\infty,\R}$, each isomorphic to $N^\tau_\R$, for $\tau$ running over faces of $\sigma$. Because of this isomorphism, we might sometimes omit the $\infty$, when it is clear from the context that we are considering the stratum at infinity. Moreover, the inclusions $\~\tau\subseteq\~\sigma$ for pairs of elements $\tau \subface \sigma$ in $\Sigma$ allow to glue these spaces and define the space $\TP_\Sigma$.

The partial compactification $\TP_\Sigma$ of $N_\R$ is naturally stratified as the disjoint union of $N^\sigma_{\infty,\R} \simeq N^{\sigma}_\R$ for $\sigma \in \Sigma$.

\begin{remark}[Remark \ref{rem:definition_dual} continued] We can also define $\~\sigma$ as a set by
\[ \~\sigma = \Hom_\tm(\sigma^\vee, \eR_+). \]
This justifies the discussion of the next paragraph.
\end{remark}

For a rational fan $\Sigma$, we have $\TP_\Sigma = \mathrm{Trop}(\P_\Sigma)$ the extended tropicalization of the toric variety $\P_\Sigma$. In other words, $\TP_\Sigma$ can be viewed as the \emph{tropical toric variety associated to the fan $\Sigma$}. In particular, the tropical projective space $\TP^n$ coincides with $\TP_\Sigma$ for $\Sigma$ the standard fan of the projective space $\P^n$. We refer to~\cites{AP, BGJK, Kaj, Payne, Thuillier} for more details.

\medskip

The canonical compactification $\comp\Sigma$ of $\Sigma$ admits a natural stratification into cones and fans that we describe now. Moreover, this stratification can be enriched into an \emph{extended polyhedral structure}~\cite{AP-geom}, see Section~\ref{sec:extpc} for more details.

Consider a cone $\sigma \in \Sigma$ and a face $\tau$ of $\sigma$. Let $C^\tau_\sigma:=\infty_\tau+(\sigma\otimes 1) \subseteq \comp\sigma$. From this description, one sees that $C^\tau_\sigma$ is isomorphic to the projection of the cone $\sigma$ in the linear space $\rquot{N_{\sigma, \mathbb R}}{N_{\tau, \mathbb R}} \hookrightarrow N^{\tau}_{\infty,\R}$. We denote by $\C^\tau_\sigma$ the relative interior of $C^\tau_\sigma$.

One can show that the canonical compactification $\comp \Sigma$ is a disjoint union of (open) cones $\C^\tau_\sigma$ for pairs $\tau\subface\sigma$ of elements of $\Sigma$.

For a fixed cone $\tau \in \Sigma$, the collection of cones $C^\tau_\sigma$ for elements $\sigma \in \Sigma$ with $\tau \subface \sigma$ form a fan with origin $\infty_\tau$ that we denote by $\Sigma_\infty^\tau\subset{N^\tau_{\infty,\R}}$. Note that with this terminology, for the cone $\conezero$ of $\Sigma$ we have $\Sigma_\infty^\conezero = \Sigma$.

The fan $\Sigma_\infty^\tau$ is canonically isomorphic to the \emph{star fan} of $\tau$ in $\Sigma$ denoted by $\Sigma^\tau$ and defined by
\[\Sigma^\tau := \Bigl\{\, \pi_\tau (\sigma) \Bigst \sigma \supface \tau \textrm{ is a cone in $\Sigma$} \, \Bigr\}, \]
where $\pi_\tau\colon N_{\R} \to N^\tau_\R$. Note that our use of the star fan is consistent with the one used in~\cite{AHK} and differs from the usual terminology, e.g., in~\cites{Karu, BBFK} where this is called \emph{transversal fan}.

If there is no risk of confusion, we sometimes drop $\infty$ and write $\Sigma^\tau$ when referring to the fan $\Sigma_\infty^\tau$ based at $\infty_\tau$.

For any pair $(\tau, \sigma)$ with $\tau \subface \sigma$, the closure $\comp C^\tau_\sigma$ of $C^\tau_\sigma$ in $\comp \Sigma$ is the union of all the cones $\C^{\tau'}_{\sigma'}$ with $\tau \subface \tau' \subface \sigma' \subface \sigma$. The closure of $\Sigma_\infty^\tau$ is indeed canonically isomorphic to the canonical compactification of $\Sigma^\tau \subseteq N^{\tau}_\R$.

The stratification of $\comp \Sigma$ by cones $C^\tau_\sigma$ and their closures is called the \emph{conical stratification} of $\comp \Sigma$. Note that there is a second stratification of $\comp \Sigma$ into fans $\Sigma_\infty^\tau$ for $\tau \in \Sigma$.

\medskip

By definition, from the inclusion of spaces $\sigma\otimes\eR_+ \simeq (\sigma+\sigma^\perp)\otimes\eR_+ \subset (\sigma+\sigma^\perp)\otimes\eR$, we get an embedding $\comp \Sigma \subseteq \TP_\Sigma$. Since the strata $\Sigma^\tau_\infty$ of $\comp \Sigma$ lives in $N^\tau_{\infty, \R}$, this embedding is compatible with the stratification of both spaces. In particular, $\comp \Sigma$ is the compactification of the fan $\Sigma$ in the tropical toric variety $\TP_\Sigma$.

\subsection{Extended polyhedral structures}\label{sec:extpc} We now describe an enrichment of the category of polyhedral complexes and polyhedral spaces into \emph{extended} polyhedral complexes and \emph{extended} polyhedral spaces by replacing the ambient spaces $V\simeq \R^n$ with $\eR^n$, or more generally, with a partial compactification $\TP_\sigma$ for the fan $\face\sigma$ associated to a cone $\sigma$ in $V$, as described in the previous section. For more details we refer to~\cites{JSS, MZ, IKMZ}.

Let $n$ be a natural number and let $x = (x_1, \dots, x_n)$ be a point of $\eR^n$. The \emph{sedentarity} of $x$ denoted by $\sed(x)$ is the set $J$ of all integers $j \in [n]$ with $x_j = \infty.$ The space $\eR^n$ is stratified into subspaces $\R^n_J$ for $J \subseteq [n]$ where $\R^n_J$ consists of all points of $x$ of sedentarity $J$. Note that $\R^n_J \simeq \R^{n-\card{J}}$.

More generally, let $\sigma$ be a cone in $V$. By an abuse of the notation, we denote by $\sigma$ the fan in $V$ which consists of $\sigma$ and all its faces. By the construction described in the previous section, $\sigma$ provides a partial compactification $\TP_\sigma$ of $V$, where we get a bijection between strata $V^\tau_{\infty}$ of $\TP_\sigma$ and faces $\tau$ of the cone $\sigma$. In the context of the previous section, the stratum $V^\tau_\infty$ was denoted $N^\tau_{\infty,\R}$.

Generalizing the definition of the sedentarity, for any point $x\in \TP_\sigma$ which lies in the stratum $V^\tau_\infty$, with $\tau\subface\sigma$, the \emph{sedentarity} of $x$ is by definition $\sed(x):=\tau$. For $\sigma$ the positive quadrant in $\R^n$, these definitions coincide with the ones in the preceding paragraph, as we get $\TP_\sigma  =\eR^n$ and the subfaces of $\sigma$ can be identified with the subsets $[n]$.

\medskip

An \emph{extended polyhedron} $\delta$ in $\TP_\sigma$ is by definition the topological closure in $\TP_\sigma$ of any polyhedron included in a strata $V^\tau_\infty$ for some $\tau \subface \sigma$. A \emph{face} of $\delta$ is the topological closure of a face of $\delta \cap V^\zeta_\infty$ for some $\zeta\subface\sigma$. An \emph{extended polyhedral complex} in $\TP_\sigma$ is a finite collection $X$ of extended polyhedra in $\TP_\sigma$ such that the two following properties are verified.
\begin{itemize}
\item For any $\delta \in X$, any face $\gamma$ of $\delta$ is also an element of $X$.
\item For pair of elements $\delta$ and $\delta'$, the intersection $\delta \cap \delta'$ is either empty or a common face of $\delta$ and $\delta'$.
\end{itemize}
The support of the extended polyhedral complex $\supp X$ is by definition the union of $\delta \in X$. The space $\X = \supp X$ is then called an \emph{extended polyhedral subspace} of $\TP_\sigma$, and $X$ is called an \emph{extended polyhedral structure} on $\X$.

\medskip

We now use extended polyhedral subspaces of partial compactifications of vector spaces of the form $\TP_\sigma$ as local charts to define more general extended polyhedral spaces.

Let $n$ and $m$ be two natural numbers, and let $V_1\simeq \R^n$ and $V_2\simeq \R^m$ be two vector spaces with two cones $\sigma_1$ and $\sigma_2$ in $V_1$ and $V_2$, respectively. Let $\phi\colon V_1 \to V_2$ be an affine map between the two spaces and denote by $A$ be the linear part of $\phi$. Let $I$ be the set of rays $\varrho\subface\sigma_1$ such that $A\varrho$ lives inside $\sigma_2$. Let $\tau_I$ be the cone of $\sigma_1$ which is generated by the rays in $I$. The affine map $\phi$ can be then extended to a map $\bigcup_{\substack{\zeta \subface \sigma_1 \\ \zeta \subseteq \tau_I}} V^\zeta_{1,\infty} \to \TP_{\sigma_2}$ denoted by $\phi$ by an abuse of the notation. We call the extension an \emph{extended affine map}. More generally, for an open subset $U \subseteq \TP_{\sigma_1}$, a map $\phi \colon U \to \TP_{\sigma_2}$ is called an \emph{extended affine map} if it is the restriction to $U$ of an extended affine map $\phi\colon \TP_{\sigma_1} \to \TP_{\sigma_2}$ as above. (The definition thus requires that $U$ is a subset of $\bigcup_{\substack{\zeta \subface \sigma_1 \\ \zeta \subseteq \tau_I}} V^\zeta_{1,\infty}$.) In the case $V_1$ and $V_2$ come with sublattices of full ranks, the extended affine map is called \emph{integral} if the underlying map $\phi$ is integral, i.e., if the linear part $A$ is integral with respect to the two lattices.

\medskip

An \emph{extended polyhedral space} is a Hausdorff topological space $\X$ with a finite atlas of charts $\Bigl(\phi_i\colon W_i \to U_i \subseteq \X_i\Bigr)_{i\in I}$ for $I$ a finite set such that

\begin{itemize}
\item $\bigl\{\,W_i \st {i\in I}\, \bigr\}$ form an open covering of $\X$;
\item $\X_i$ is an extended polyhedral subspace of $\TP_{\sigma_i}$ for a finite dimensional real vector space $V_i\simeq\R^{n_i}$ and a cone $\sigma_i$ in $V_i$, and $U_i$ is an open subset of $\X_i$; and
\item $\phi_i$ are homeomorphisms between $W_i$ and $U_i$ such that for any pair of indices $i,j\in I$, the transition map
\[\phi_j\circ \phi_i^{-1} \colon \phi_i(W_i \cap W_j) \to \TP_{\sigma_j}\]
is an extended affine map.
\end{itemize}

The extended polyhedral space is called \emph{integral} if the transition maps are all integral with respect to the lattices $\Z^{n_i}$ in $\R^{n_i}$.

\medskip

Let $\X$ be an extended polyhedral space with an atlas of charts $\Bigl(\,\phi \colon  W_i \to U_i  \subseteq \X_i\Bigr)_{i\in I}$ as above. A \emph{face structure} on $\X$ is the choice of an extended polyhedral complex $X_i$ with $\supp{X_i} = \X_i$ for each $i$ and a finite number of closed sets $\theta_1, \dots, \theta_N$, for $N \in \mathbb N$, called \emph{facets} which cover $\X$ such that the following two properties are verified.

\begin{itemize}
\item Each facet $\theta_j$ is contained in some chart $W_i$ for $i\in I$ such that $\phi_i(\theta_j)$ is the intersection of the open set $U_i$ with a face $\eta_{j, i}$ of the polyhedral complex $X_i$.
\item For any subset $J \subseteq [N]$, any $j \in J$, and any chart $W_i$ which contains $\theta_j$, the image by $\phi_i$ of the intersection $\bigcap_{j\in J} \theta_j$ in $U_i$ is the intersection of $U_i$ with a face of $\eta_{j,i}$.
\end{itemize}
A \emph{face} of the face structure is the preimage by $\phi_i$ of a face of $\eta_{j,i}$ for any $j\in[N]$ and for any $i\in  I$ such that $\theta_j\subseteq W_i$. Each face is contained in a chart $W_i$  and the \emph{sedentarity of a face $\delta$ in the chart $W_i$} is defined as the sedentarity of any point in the relative interior of $\phi_i(\delta)$ seen in $X_i$.

\begin{prop} Let $\Sigma$ be a fan in $N_\R$. The canonical compactification $\comp \Sigma$ naturally admits the structure of an extended polyhedral space, and a face structure given by the closure $\comp C^\tau_\sigma$ in $\comp \Sigma$ of the cone $C^\tau_\sigma$ in the conical stratification of $\comp \Sigma$, for pairs of cones $\tau \subface \sigma$ in $\Sigma$. The extended polyhedral structure is integral if $\Sigma$ is a rational fan.
\end{prop}
\begin{proof} This is proved in~\cite{AP-geom}.
\end{proof}

Note that the sedentarity of the face $\comp C^\tau_\sigma$ of the conical stratification of $\comp \Sigma$ is the set of rays of $\tau$.

\subsection{Recession fan and canonical compactification of polyhedral complexes} \label{sec:can-compact-polycomp} Let $X$ be a polyhedral complex in a real vector space $N_\R \simeq \R^n$. The \emph{recession pseudo-fan} of $X$ denoted by $X_\infty$ is the set of cones $\{\delta_\infty\st\delta\in X\}$. In fact the collection of cones $\delta_\infty$ do not necessary form a fan, as depicted in the following example; we will however see later in Section~\ref{sec:triangulation} how to find a subdivision of $X$ whose recession pseudo-fan becomes a fan. In such a case, we will call $X_\infty$ the \emph{recession fan} of $X$. Recession fans are studied in~\cite{BS11}.

\begin{example} Let $V = \R^3$ and denote by $(e_1, e_2, e_3)$ the standard basis of $\R^3$. Define
\begin{gather*}
\sigma_1=\R_+e_1+\R_+e_2\subset\R^3, \\
\sigma_2=\R_+e_1+\R_+(e_1+e_2)\subset\R^3, \\
X=\face{\sigma_1}\cup\faceop\left(e_3+\sigma_2\right).
\end{gather*}
Note that $X_\infty$ contains the cones $\sigma_1$ and $\sigma_2$ whose intersection is not a face of $\sigma_1$.
\end{example}

Let $X$ be a polyhedral complex in $N_\R \simeq \R^n$ whose recession pseudo-fan $X_\infty$ is a fan. Let $\TP_{X_\sminfty}$ be the corresponding tropical toric variety. The \emph{canonical compactification of $X$} denoted by $\comp X$ is by definition the closure of $X$ in $\TP_{X_\sminfty}$.

\medskip

We now describe a natural stratification of $\comp X$. Let $\sigma \in X_\infty$, and consider the corresponding stratum $N^{\sigma}_{\infty, \R}$ of $\TP_\sigma$. Define $X^\sigma_{\infty}$ to be the intersection of $\comp X$ with $N^{\sigma}_{\infty, \R}$, that without risk of confusion, we simply denote $X^\sigma$ by dropping $\infty$. Note that $X^{\conezero} = X$, which is the \emph{open part} of the compactification $\comp X$.

\medskip

Let $D = \comp X \setminus X$ be the \emph{boundary at infinity} of the canonical compactification of $X$. For each non-zero cone $\sigma$ in $X_\infty$, we denote by $D^\sigma$ the closure of $X^\sigma$ in $\comp X$.

\begin{thm} [Tropical orbit-stratum correspondence]\label{thm:orbit-stratum-correspondence}
Notations as above, let $X$ be a polyhedral complex in $N_\R \simeq \R^n$. We have the following.
\begin{enumerate}
\item For each cone $\sigma \in X_\infty$, the corresponding strata $X^\sigma$ is a polyhedral complex in $N^{\sigma}_{\infty, \R}$.
\item The pseudo-recession fan $(X^{\sigma})_\infty$ of $X^\sigma$ is a fan which coincides with the fan $(X_\infty)_\infty^{\sigma}$ in $N_{\infty, \R}^\sigma$ in the stratification of the canonical compactification $\comp{X_\infty}$ of the fan $X_\infty$.
\item The closure $D^\sigma$ of $X^\sigma$ coincides with the canonical compactification of $X^{\sigma}$ in $N^{\sigma}_{\infty, \R}$, \ie, $D^\sigma = \comp{X^{\sigma}}$.
\item If $X$ has pure dimension $d$, then each stratum $X^{\sigma}$ as well as their closures $D^\sigma$ have pure dimension $d -\dims \sigma$.
\end{enumerate}
In particular, the canonical compactification $\comp X$ is an extended polyhedral structure with a face structure induced from that of $X$.
\end{thm}

\begin{proof} The theorem is proved in~\cite{AP-geom}.
\end{proof}

\subsection{Tropical homology and cohomology groups} We recall the definition of tropical cohomology groups and refer to~\cites{IKMZ, JSS, MZ, GS-sheaf} for more details.

Let $X$ be an extended polyhedral space with a face structure. We start by recalling the definition of the multi-tangent and multi-cotangent spaces $\SF_p(\delta)$ and $\SF^p(\delta)$, respectively, for the face $\delta$ of $X$. This leads to the definition of a chain, resp. cochain, complex which calculates the tropical homology, resp. cohomology, groups of $X$.

For a face $\delta$ of $X$ and for any non-negative integer $p$, the \emph{$p$-th multi-tangent} and the \emph{$p$-th multicotangent space} of $X$ at $\delta$ denoted by $\SF_p(\delta)$ and $\SF^p(\delta)$ are defined by
\[\SF_{p}(\delta)=  \!\!\sum_{\eta \supface \delta \\ \sed(\eta) = \sed(\delta) }\!\! \bigwedge^p \TT\eta, \qquad \textrm{and} \quad \SF^{p}(\delta) = \SF_p(\delta)^\dual. \]

For an inclusion of faces $\gamma \subface \delta$ in $X$, we get maps $\i_{\delta\supface\gamma} \colon \SF_p(\delta) \to \SF_p(\gamma)$ and $\i^*_{\gamma\subface\delta} \colon \SF^p(\gamma) \to \SF^p(\delta)$.

\medskip

For a pair of non-negative integers $p,q$, define
\[C_{p,q}(X) := \bigoplus_{\delta \in X \\ \dims{\delta} =q} \SF_p(\delta) \]
and consider the cellular chain complex, defined by using maps $\i_{\delta\supface\gamma}$ as above,
\[C_{p,\bul}\colon \quad \dots\longrightarrow C_{p, q+1}(X) \xrightarrow{\partial^\trop_{q+1}} C_{p,q}(X)  \xrightarrow{\ \partial^\trop_{q}\ } C_{p,q-1} (X)\longrightarrow\cdots\]
The tropical homology of $X$ is defined by
\[H_{p,q}^\trop(X) := H_q(C_{p,\bul}).\]
Similarly, we have a cochain complex
\[C^{p,\bul}\colon \quad \dots\longrightarrow C^{p, q-1}(X) \xrightarrow{\partial_\trop^{q-1}} C^{p,q}(X)  \xrightarrow{\ \partial_\trop^{q}\ } C^{p,q+1}(X) \longrightarrow\cdots\]
where
\[C^{p,q}(X) := \bigoplus_{\delta\in X \\
\dims{\delta}=q} \SF^p(\delta).\]
The tropical cohomology of $X$ is defined by
\[H^{p,q}_\trop(X) := H^q(C^{p,\bul}).\]

In the case $X$ comes with a rational structure, one can define tropical homology and cohomology groups with integer coefficients. In this case, for each face $\delta$, the tangent space $\TT\delta$ has a lattice $N_\delta$. One defines then
\[\SF_{p}(\delta,\Z)=  \!\!\sum_{\eta \supface \delta \\ \sed(\eta) = \sed(\delta) }\!\! \bigwedge^p N_\eta, \quad \textrm{and} \qquad \SF^{p}(\delta, \Z) = \SF_p(\delta,\Z)^{\vee}.\]
Similarly, we get complexes with $\Z$-coefficients $C_{p,\bul}^\Z$ and $C^{p, \bul}_\Z$, and define
\[H_{p,q}^\trop(X, \Z) := H_q(C_{p,\bul}^{\Z}) \qquad H^{p,q}_\trop(X, \Z) = H^q(C^{p,\bul}_\Z).\]

Similarly, working with the cochain $C^{p, \bullet}_c$ with \emph{compact support}, one can define tropical cohomology groups with compact support. We denote them by $H^{p,q}_{\trop,c}(X, A)$ with coefficients $A =\Z, \Q,$ or $\R$. If the coefficient field is not specified, it means we are working with real coefficients.

\subsection{Dolbeault cohomology and comparison theorem}
In this section we briefly recall the formalism of superforms and their corresponding Dolbeault cohomology on extended polyhedral complexes. The main references here are~\cites{JSS, JRS19}, see also~\cites{Lagerberg, CLD, GK17}.

Let $V \simeq \R^n$ be a real vector space and denote by $V^\dual$ its dual. We denote by $\A^{p}$ the sheaf of differential forms of degree $p$ on $V$. Note that the sheaf $\A^p$ is a module over the sheaf of rings of smooth functions $\CC^\infty$. The \emph{sheaf of (super)forms of bidegree} $(p,q)$ denoted by $\A^{p,q}$ is defined as the tensor product $\A^p\otimes_{\CC^{\infty}} \A^q$. Choosing a basis $\partial_{x_1}, \dots, \partial_{x_n}$ of $V$, we choose two copies of the dual vector space $V^\dual$ with the dual basis $\dip x_1, \dots, \dip x_n$ for the first copy and $\dis x_1, \dots, \dis x_n$ for the second. In the tensor product $\A^p\otimes \A^q$ the first basis is used to describe the forms in $\A^p$ and the second basis for the forms in $\A^q$. So for an open set $U \subseteq V$, a section $\alpha \in \A^{p,q}$ is uniquely described in the form
\[\alpha = \sum_{I, J \subseteq [n]\\ \card I = p , \card J = q} \alpha_{_{I, J}}\,  \dip x_I \wedge \dis x_J\]
for smooth functions $\alpha_{I, J} \in \CC^\infty(U)$. Here for a subset $I =\{i_1, \dots, i_\ell\} \subseteq [n]$ with $i_1< i_2<\dots<i_\ell$, we denote by $\dip x_I$ and $\dis x_I$ the elements of the two copies of $\bigwedge^{\ell} V^\dual$ defined by
\[ \dip x_I = \dip x_{i_1} \wedge \dots \wedge \dip x_{i_\ell} \qquad \textrm{and} \qquad \dis x_I = \dis x_{i_1} \wedge \dots \wedge \dis x_{i_\ell},\]
respectively. Moreover, for subsets $I, J \subseteq [n]$, $\dip x_I \wedge \dis x_J$ is the element of $\bigwedge^p V^\dual \otimes \bigwedge^q V^\dual$ given by $\dip x_I \otimes \dis x_J$. The usual rules of the wedge product apply. In particular, for any $i,j\in[n]$,
\[ \dis x_j \wedge \dip x_i=-\dip x_i \wedge \dis x_j. \]

\medskip

Observing that $\A^{p,q} \simeq \A^{p} \otimes_\R \bigwedge^q V^\dual$ and $\A^{p,q} \simeq \bigwedge^p V^\dual \otimes \A^q$, the collection of bigraded sheaves $\A^{p,q}$ come with differential operators $\dip\colon \A^{p,q} \to \A^{p+1,q}$ and $\dis \colon \A^{p,q} \to \A^{p,q+1}$ defined by $\dip = \d_p \otimes \id$ and $\dis = \id\otimes \d_q$ where for non-negative integer $\ell$, $\d_\ell \colon \A^\ell \to \A^{\ell+1}$ is the differentiation of smooth $\ell$-forms. More explicitly, the differential map $\dis$ is defined by its restriction to open sets $U \subseteq V$
\[\A^{p,q}(U) \to \A^{p,q+1}(U)\]
and sends a $(p,q)$-form $\alpha = \sum_{_{I, J}} \alpha_{_{I, J}}\,\dip x_I \wedge \dis x_J$ to
\[\dis \alpha = \sum_{I, J} \sum_{\ell=1}^n \partial_{x_\ell} \alpha_{_{I,J}} \,\,\dis x_\ell \wedge \dip x_I \wedge \dis x_J = (-1)^p \sum_{I, J} \sum_{\ell=1}^n \partial_{x_\ell} \alpha_{_{I,J}} \,\,  \dip x_I \wedge (\dis x_\ell\wedge \dis x_J).  \]
The operator $\dip$ has a similar expression.

\medskip

Consider now a fan $\Sigma$ in $V$, and let $\TP_\Sigma$ be the corresponding partial compactification of $V$. We define the  \emph{sheaf of bigraded $(p,q)$-forms} $\A^{p,q}$ on $\TP_\Sigma$ as follows. For each pair of cones $\tau\subface\sigma$ of $\Sigma$, we have a natural projection $\pi_{\tau\subface\sigma}\colon V^\tau_\infty \to V^\sigma_\infty$. This projection induces a map $\pi_{\sigma\supface\tau}^*\colon \A^{p,q}(V^\sigma_\infty)\to\A^{p,q}(V^\tau_\infty)$. Let $U \subseteq \TP_\Sigma$ be an open set. A section $\alpha \in \A^{p,q}(U)$ is the data of the collection of $(p,q)$-forms $\alpha^\sigma \in \A^{p,q}(U \cap V^{\sigma}_\infty)$, where $U \cap V^\sigma_\infty$ is viewed as an open subset of the real vector space $V^\sigma_\infty \simeq \rquot{V}{\TT\sigma}$, so that the following compatibility at infinity is verified among these sections. For each pair of cones $\tau\subface\sigma$ in $\Sigma$, there exists a neighborhood $U'$ of $U\cap V^\sigma_\infty$ in $U$ such that the restriction of $\alpha^\tau$ to $U'$ is given by the pullback ${\pi^*_{\sigma\supface\tau}} \bigl(\alpha^{\sigma}\bigr)\rest{U'}$.

The differential operators $\dip$ and $\dis$ naturally extend to $\TP_\Sigma$ and define operators $\dip\colon \A^{p, q}\to \A^{p+1, q}$ and $\dis\colon \A^{p, q} \to \A^{p, q+1}$ for any pair of non-negative integers $(p, q)$.

\medskip

Let now $X$ be a closed polyhedral complex in $\TP_\Sigma$, and denote $\iota\colon X  \hookrightarrow \TP_\Sigma$ the inclusion. The sheaf of $(p,q)$-forms $\A^{p,q}_X$ is defined as a quotient of the restriction $\iota^* \A^{p,q}_{\TP_\Sigma}$ of the sheaf of $(p,q)$-forms on $\TP_\Sigma$ to $X$ modulo those sections which vanish on the tangent space of $X$. Thus for an open set $U$ in $X$, a section $\alpha \in \A^{p,q}_X(U)$ is locally around each point $x$ given by a section $\alpha_x \in \A^{p,q}(U_x)$ for an open subset $U_x \subseteq \TP_\Sigma$ with $U_x \cap X  \subseteq U$, and two sections $\alpha$ and $\beta$ in $\A_X^{p,q}(X)$ are equivalent if for any point $x$ living in a face $\sigma$ of $X$, the difference $\alpha(x) - \beta(x)$ vanishes on any $(p,q)$-multivector in $\bigwedge^p \TT \sigma \otimes \bigwedge^q \TT \sigma$.

The definition above can be extended to any extended polyhedral complex.

The two differential operators $\dip$ and $\dis$ naturally extend to $\TP_{X_\sminfty}$ and we get $\dip\colon \A_X^{p, q}\to \A_X^{p+1, q}$ and $\dis\colon \A_X^{p, q} \to \A_X^{p, q+1}$ for any pair of non-negative integers $(p, q)$.

\medskip

The above definition only depends on the support of $X$ and not on the chosen polyhedral structure on this support. As a consequence, if $\X$ is any polyhedral space, the sheaf $\A^{p,q}_\X$ is well-defined.

\medskip

The \emph{Dolbeault cohomology group $H^{p,q}_{\Dol}(X)$ of bidegree $(p,q)$} is by definition the $q$-th cohomology group of the cochain complex
\[\bigl(\A_X^{p,\bullet}(X), \dis\bigr)\colon \qquad \A_X^{p,0}(X) \to \A^{p,1}_X(X) \to \dots \to \A^{p,q-1}_X(X) \to \A^{p,q}_X(X) \to \A^{p,q+1}_X(X) \to \cdots \]

Similarly, we define Dolbeault cohomology groups with compact support $H^{p,q}_{\Dol,c}(X)$ as the $q$-th cohomology group of the cochain complex
\[\bigl(\A_{X,c}^{p,\bullet}(X), \dis\bigr)\colon \qquad \A_{X,c}^{p,0}(X) \to \A^{p,1}_{X,c}(X) \to \cdots \]
where $\A^{p,q}_{X, c}(X)$ is the space of $(p,q)$-forms whose support forms a compact subset of $X$.

Denote by $\Omega_X^p$ the kernel of the map $\dis\colon \A_X^{p,0} \to \A^{p,1}_X$. We call this the sheaf of \emph{holomorphic $p$-forms} on $X$.
By Poincar\'e lemma for superforms, proved in~\cites{Jell, JSS}, the complex $(\A^{p,\bullet}, \dis)$ provides an acyclic resolution of the sheaf $\Omega_X^p$. We thus get the isomorphism
\[H^{p,q}_\Dol(X) \simeq H^q(X, \Omega_X^p),\]
and similarly,
\[H^{p,q}_{\Dol, c}(X) \simeq H^q_c(X, \Omega_X^p).\]
Moreover, we have the following comparison theorem from~\cite{JSS}.
\begin{thm}[Comparison theorem]\label{thm:comparison}
For any extended polyhedral complex $X$, we have
\[H^{p,q}_\trop(X) \simeq H^{p,q}_\Dol(X) , \qquad \textrm{and} \qquad H^{p,q}_{\trop,c}(X) \simeq H^{p,q}_{\Dol,c}(X).\]
\end{thm}

\bigskip

\subsection{Matroids} \label{sec:matroid} We first recall basic definitions involving matroids, and refer to relevant part of \cite{Oxl11} for more details.

\medskip

For a set $E$, a subset $A \subseteq E$ and an element $a\in E$, if $a\in A$, we denote by $A-a$ the set $A \setminus \{a\}$; if $a\notin A$, we denote by $A+a$ the set $A \cup\{a\}$.

\medskip

A matroid can be defined in different \emph{equivalent} manners, for example by specifying its collection of \emph{independent sets}, or its collection of \emph{bases}, or its collection of \emph{flats}, or still by giving its \emph{rank function}. In each case, the data of one of these determines all the others. The definition of a matroid with respect to independent sets is as follows.

\begin{defi}[Matroid: definition with respect to independent sets] \label{defi:matroid}
A \emph{matroid} $\Ma$ is a pair consisting of a finite set $E$ called the \emph{ground set of $\Ma$} and a family $\ind$ of subsets of $E$ called the \emph{family of independent sets of $\Ma$} which verifies the following three axiomatic properties:
\begin{enumerate}
\item $\emptyset\in\ind$.
\item (hereditary property) $\ind$ is stable under inclusion, in the sense that if $J\subset I$ and $I\in\ind$, then $J\in\ind$.
\item (augmentation property) if $I, J\in\ind$ and if $\card J<\card I$, then there exists an element $i\in I\setminus J$ such that $J+i\in\ind$. \label{defi:matroid:augmentation_property}  \qedhere
\end{enumerate}
\end{defi}
If needed, we write $\Ma=(E, \ind)$ in order to emphasize that $\Ma$ is defined on the ground set $E$ and has independent sets $\ind \subseteq 2^E$.

\begin{example}\label{ex:representable} An example of a matroid is given by a collection of vectors $v_1, \dots, v_m$ in a vector space $H \simeq k^n$, over a field $k$. Such a collection defines a matroid $\Ma$ on the ground set $[m]$ whose independent $\ind$ are all the subsets $I \subseteq [m]$ such that the corresponding vectors $v_i$ for $i\in I$ are linearly independent. A matroid $\Ma$ which is defined by a collection of vectors in a vector space is said to be \emph{representable}. By an interesting recent result of Peter Nelson~\cite{Nelson}, we know that \emph{almost any matroid is non representable over any field}. More precisely, if one denotes by $\epsilon_m$ the number of representable matroids on the ground set $[m]$ divided by the total number of matroids on the same ground set, then the numbers $\epsilon_m$ tend to zero when $m$ goes to infinity; Nelson gives estimates for all these values.
\end{example}

For a matroid $\Ma=(E,\ind)$, the \emph{rank function} $\rkm\colon 2^E \to \mathbb N \cup\{0\}$ is defined as follows. For a subset $A \subseteq E$, the \emph{rank of $A$} is defined by
\[ \rkm(A):=\max_{\substack{I \in \ind \textrm{ with }I \subseteq A}} \, \card I. \]
We define the \emph{rank of $\Ma$} to be the integer number $\rkm(E)$. The rank function satisfies what is called the \emph{submodularity property}, meaning
\[\forall A, A' \subseteq E, \qquad \rkm(A) + \rkm(A') \leq \rkm(A\cap A') + \rkm(A \cup A').\]
The data of $E$ and rank function $\rkm$ determines $\ind$ as the set of all subsets $I \subseteq E$ which verify $\rkm(I) = \card{I}$.

A maximal independent set in $\Ma$ is called a \emph{basis} of $\Ma$. The collection of bases of $\Ma$ is denoted by $\bases(\Ma)$.

The \emph{closure} of a subset $A \subseteq E$ is the set
\[ \cl(A):=\{a \in E \textrm{ such that }  \rkm(A + a)=\rkm(A)\}. \]
A subset $F \subseteq E$ is called a \emph{flat} of $\Ma$ if $\cl(F) = F$. Flats are also equally called \emph{closed set} and the collection of flats of the matroid $\Ma$ is denoted by $\Cl(\Ma)$.

\begin{example}[Example~\ref{ex:representable} continued] For a matroid $\Ma$ representable by a collection of vectors $v_1, \dots, v_m$ in a vector space $ H \simeq k^n$, the above notions coincide with the usual linear algebra notions. For example the rank $\rkm(A)$ is the dimension of the vector subspace generated by the set of vectors $v_i$ for $i\in A$. A flat $F$ of $\Ma$ is a subset of $[m]$ such that $\{v_i\}_{i\in F}$ is the set of the vectors in the collection which live in a vector subspace of $H$.
\end{example}

An element $e\in E$ is called a \emph{loop} if $\{e\}$ is not an independent set, i.e., if $\rkm(\{e\})=0$. Two elements $e$ and $e'$ of $E$ are called \emph{parallel} if $\rkm(\{e, e'\}) = \rkm(\{e\}) = \rkm(\{e'\}) = 1$. A matroid is called \emph{simple} if it neither contains any loop nor any parallel elements.

An element $e\in E$ is called a \emph{coloop} if $\rkm(E-e)=\rkm(E)-1$, equivalently, if $E-e$ is a flat. A \emph{proper element} of $\Ma$ is an element which is neither a loop nor a coloop.

\subsection{Bergman fans} For a simple matroid $\Ma$ on a ground set $E$, the Bergman fan of $\Ma$ denoted by $\Sigma_\Ma$ is defined as follows.

First, we fix some notations. We denote by $\{\e_i\}_{i\in E}$ the standard basis of $V$. For a subset $A \subseteq E$, we denote by $\e_A$ the sum $\sum_{i\in A} \e_i$ in $V$.

Consider the real vector space $V := \rquot{\R^E}{\R \e_E}$. It has a full rank lattice $N = \rquot{\Z^E}{\Z \e_E}$. By an abuse of the notation, we denote by $\e_A$ the projection of $\e_A=\sum_{i\in A}\e_i$ in $V$. In particular, we have $\e_E =0$.

The \emph{Bergman fan of $\Ma$} denoted by $\Sigma_\Ma$ is a rational fan in $N_\R =V$ of dimension $\rkm(\Ma) -1$ defined as follows. A \emph{flag of proper flats $\Fl$ of $\Ma$ of length $\ell$} is a collection
\[\Fl\colon \quad \emptyset \neq F_1 \subsetneq F_2 \subsetneq \dots \subsetneq F_{\ell} \neq E\]
consisting of flats $F_1, \dots, F_\ell$ of $\Ma$.

For each such flag $\Fl$, we define the rational cone $\sigma_\Fl \subseteq V$ as the cone generated by the vectors $\e_{F_1}, \e_{F_2}, \dots, \e_{F_\ell}$, so
\[\sigma_\Fl := \Bigl\{ \lambda_1 \e_{F_1} + \dots + \lambda_\ell \e_{F_\ell} \Bigst \lambda_1, \dots, \lambda_\ell \geq0\Bigr\}.\]
Note that the dimension of $\sigma_\Fl$ is equal to the length of $\Fl$.

The Bergman fan of $\Ma$ is the fan consisting of all the cones $\sigma_\Fl$ for $\Fl$ a flag of proper lattices of $\Ma$:
\[\Sigma_\Ma := \Bigl\{\, \sigma_\Fl \Bigst \Fl \textrm{ flag of proper flats of }\Ma\,\Bigr\}.\]
From the definition, it follows easily that $\Sigma_\Ma$ has pure dimension $\rkm(\Ma)-1$.

\medskip

A subset $S$ in a real vector space $W$ is called a \emph{Bergman support} if there exists a matroid $\Ma$ on a ground set $E$ and a linear map $\phi\colon W \to \rquot{\R^E}{\R \e_E}$ which induces a bijection $\phi\rest S\colon S \to \supp{\Sigma_\Ma}$. A fan $\Sigma$
in a real vector space is called a \emph{Bergman fan} if the support $\supp{\Sigma}$ of $\Sigma$ is a Bergman support.

A \emph{rational Bergman fan} is a fan in a vector space $W$ which is rational with respect to a full rank lattice $N_W$ in $W$ such that there exists a $\Z$-linear map $\phi\colon N_W \to \rquot{\Z^E}{\Z\e_E}$ which, once extended to a linear map $W \to \rquot{\R^E}{\R \e_E}$, induces an isomorphism
$\supp{\Sigma} \to \supp{\Sigma_\Ma}$.

\begin{example} Any complete fan $\Sigma$ in a real vector space $\R^n$ is a Bergman fan. In fact, denoting by $U_{n+1}$ the \emph{uniform matroid on the ground set $E=[n+1]$} whose collection of independent sets is $\ind = 2^{[n+1]}$, we see that $\R^n \simeq \rquot{\R^{n+1}}{\R \e_{E}}$ is a Bergman support. Therefore, $\Sigma$ is a Bergman fan. Moreover, $\Sigma$ is a rational Bergman fan if it is rational (in the usual sense) with respect to the lattice $\Z^n$.
\end{example}

\begin{example}[Example~\ref{ex:representable} continued] Let $\Ma$ be a matroid on the ground set $[m]$ representable by a collection of non-zero vectors $v_1, \dots, v_m$ in a vector space $V$ over $k$. Viewing $v_1, \dots, v_m$ as elements of the dual vector space of $V^\dual$, we get a collection of hyperplanes $H_1, \dots, H_m \subset V^\dual$ with $H_j = \ker\Bigl( v_j \colon  V^\dual \to k\Bigr)$ for each $j$, which leads to a hyperplane arrangements denoted, by abusing of the notation, by $H_1, \dots, H_m \subseteq \P(V^\dual)$ in the projective space $\P(V^\dual)$.

Denote by $X$ the complement in $\P(V^\dual)$ of the union $H_1 \cup \dots \cup H_m$. The collection of linear forms $v_j$ gives an embedding of $X\hookrightarrow \rquot{\G_m^{m}}{\G_m}$, where $\G_m$ is the algebraic torus over $k$ whose set of $K$-points for each field extension $K$ of $k$ is $K^\times$, and the quotient $\rquot{\G_m^{m}}{\G_m}$ means the quotient of $\G_m^{m}$ by the diagonal action of $\G_m$. The tropicalization of $X$ in the torus $\rquot{\G_m^{m}}{\G_m}$ coincides with the Bergman fan $\Sigma_\Ma$ for the representable matroid $\Ma$ which is associated to the collection of vectors $v_1, \dots, v_m$. See~\cites{AK, MS} for more details.
\end{example}

\medskip

The category of Bergman fans is closed under product as the following proposition shows.

\begin{prop}\label{prop:productBergman} The product of any two Bergman supports is a Bergman support. In particular, the product of two Bergman fans is again Bergman.

More precisely, for two matroids $\Ma$ and $\Ma'$, we have
\[\supp{\Sigma_{\Ma}} \times \supp{\Sigma_{\Ma'}} = \supp{\Sigma_{\Ma\vee\Ma'}}\]
where $\Ma\vee \Ma'$ is a parallel connection of $\Ma$ and $\Ma'$.
\end{prop}

We recall the definition of \emph{parallel connection} between two matroids, \cf.~\cite{Oxl11}*{Chapter 7}, and refer to~\cite{AP} for the proof of this proposition.

A \emph{pointed matroid} is a matroid $\Ma$ on a ground set $E$ with a choice of a distinguished element $\distel = \distel_\Ma$ in $E$.

For two pointed matroids $\Ma$ and $\Ma'$ on ground sets $E$ and $E'$ with distinguished elements $\distel_\Ma \in E$ and $\distel_{\Ma'} \in E'$, respectively, the \emph{parallel connection} of $\Ma$ and $\Ma'$, that we can also call the \emph{wedge sum} of $\Ma$ and $\Ma'$, denoted by $\Ma \vee \Ma'$, is by definition the pointed matroid on the ground set $E \vee E' = \rquot{E \sqcup E'}{(\distel_\Ma=\distel_{\Ma'})}$, the wedge sum of the two pointed sets $E$ and $E'$, with distinguished element $\distel = \distel_\Ma = \distel_{\Ma'}$, and with the following collection of bases:

\begin{align*}
\bases( \Ma \vee \Ma') =& \Bigl\{ \,B \cup B' \st B \in \bases(\Ma), B' \in \bases(\Ma') \,\textrm{ with }\,\distel_\Ma\in B\, \textrm{ and }\,  \distel_{\Ma'}\in B'\Bigr \} \\
& \cup  \Bigl\{ \,B \cup B' \setminus\{\distel_\Ma\} \st B \in \bases (\Ma), B' \in \bases(\Ma') \, \textrm{ with }\, \distel_\Ma\in B \, \textrm{ and }\, \distel_{\Ma'}\notin B'\Bigr \}\\
& \cup  \Bigl\{ \,B \cup B' \setminus\{\distel_{\Ma'}\}\st B \in \bases(\Ma), B' \in \bases(\Ma')\, \textrm{ with }\, \distel_\Ma \notin B\, \textrm{ and }\,  \distel_{\Ma'}\in B'\Bigr \}.
\end{align*}

A parallel connection of two matroids $\Ma$ and $\Ma'$ is the wedge sum of the form $(\Ma, e) \vee (\Ma', e')$ for the choice of two elements $e\in \Ma$ and $e'\in \Ma'$, which make pointed matroids out of them.

\subsection{Tropical cycles and divisors} Let $V =N_\R \simeq \R^n$ be a real vector space of dimension $n$ with a full rank lattice $N$. Let $Y$ be a rational polyhedral complex of pure dimension $d$ in $V$. For each facet $\sigma$ of $Y$, suppose we are given a \emph{weight} which is an integer $w(\sigma)$. Let $C := (Y, w)$ be the weighted polyhedral complex $Y$ with the weight function on its facets $w$. The weighted polyhedral complex $C$ is called a \emph{tropical cycle} if the following \emph{balancing condition} is verified:
\[\forall \tau\in Y\text{ of dimension $d-1$},\qquad  \sum_{\sigma \ssupface \tau}  w(\sigma) \e_{\sigma/\tau}  = 0 \in N^{\tau}.\]
Here $\e_{\sigma/\tau}$ is the primitive vector of the ray $\rho_{\sigma/\tau}$ corresponding to $\sigma$ in $N^\tau_\R$. The following is a direct consequence of the the definition of the Bergman fan given in the previous section.

\begin{prop} The Bergman fan $\Sigma_\Ma$ of any matroid with weights equal to one on facets is a tropical cycle. More generally any Bergman fan is a tropical cycle.
\end{prop}

A \emph{tropical regular function} $f$ on a tropical cycle $C$ is a piecewise affine function with integral slopes which is the restriction to $C$ of a concave function on $\R^n$. In other words, it is a function of the form
\[f = \min_{i \in I} \,\bigl\{ \langle \alpha_i, \ccdot \rangle + \beta_i\bigr\}\]
for a finite index set $I$ and $\alpha_i \in \Z^n$, $\beta_i\in \Q$.

Any tropical regular function induces a subdivision of the tropical cycle $C$. On each face $\delta$ of $C$, the new polyhedral structure consists of set of the form $\{x\in\delta \st f(x) = g(x) \}$ where $g$ is an affine linear function such that $g(x)\geq f(x)$ for every $x$. Keeping the same weight on facets, this leads to a new tropical cycle that we denote by $C_1$ and note that $C$ and $C_1$ have the same support and same weight function.

Given a tropical regular function $f$ on $C$, we can consider the graph $\Gamma_f$ of $f$ in $V \times \R \simeq \R^{n+1},$ with the polyhedral structure obtained by that of the subdivided $C_1$. The graph $\Gamma_f$ is not necessary a tropical cycle in $V\times \R$ anymore, as it might have unbalanced codimension one faces. In this case, we can add new faces to $\Gamma_f$ to obtain a tropical cycle. For this, consider a codimension one face $\delta$ of $C_1$ which turns to be unbalanced in the graph $\Gamma_f$. Consider the  polyhedron $E_\delta$ defined by
\[E_\delta := \Bigl\{ (x, t) \st x \in \delta \myand t \geq f(x)\Bigr\} \subset V \times \R.\]
Define $\widetilde C$ as the union of $\Gamma_f$ with the polyhedra $E_\delta$ and there faces, one polyhedron for each unbalanced codimension one face $\delta$ of $\Gamma_f$.  Note that for each such face $\delta$, the primitive vector $\e_{E_\delta/\delta}$ is the one given by the last vector $\e_{n+1}$ in the standard basis of $V\times \R \simeq \R^{n+1}$.

One verifies there is a unique integer weight $w(E_\delta)$ which we can associate to new faces $E_\delta$ in order to make $\widetilde C$ a tropical cycle in $V\times \R$.  The new tropical cycle in $V\times \R$ is called tropical modification of $C$ along $\div_C(f)$, where $\div_C(f)$ is the union of the $\delta$ in $C_1$ of codimension one which happens to be unbalanced when seen in $\Gamma_f$.

\subsection{Tropical modifications of Bergman fans} \label{sec:modification}

For a matroid $\Ma$ on a ground set $E$ and for an element $a \in E$, we denote by $\Ma \del a$ the matroid obtained by the \emph{deletion of $a$}. This is the matroid on ground set $E - a$ whose independent sets are those of $\Ma$ which do not contain $a$. The rank function of $\Ma\del a$ is thus that of $\Ma$ restricted to the set $E-a$. The flats of $\Ma\del a$ are all the sets $A = F\setminus\{a\}$ for a flat $F$ of $\Ma$, so
\[\Cl(\Ma \del a) := \Bigl\{F \subseteq E - a \st \textrm{either $F$ or $F+a$ is a flat of $\Ma$} \Bigr\}.\]

For an element $a \in E$, the matroid obtained by \emph{contraction of $a$ in $\Ma$} is the matroid denoted by $\Ma\contr a$ on the ground set $E-a$ whose independent sets are precisely those subsets $I \subseteq E -a$ such that $I+a$ is in independent set in $\Ma$.
It follows that the flats of $\Ma\contr a$ are obtained by removing $a$ from the flats of $\Ma$ which contain $a$, i.e.,
\[\Cl(\Ma \contr a) := \Bigl\{F \subseteq E - a \st \textrm{$F+a$ is a flat of $\Ma$} \Bigr\}.\]

Let $N=\rquot{\Z^E}{\Z \e_E}$, $M = N^\vee$ the dual lattice, and let $N^a := \rquot{\Z^{E-a}}{\Z \e_{_{E-a}}} = \rquot{N}{\Z \e_a}$ and $M^a = N^{a\vee}$.
The two Bergman fans $\Sigma_{\Ma \del a}$ and $\Sigma_{\Ma \contr a}$ both live in the real vector space $N^a_\R.$ Moreover, from the description of the flats of these matroids, it is clear that $\supp{\Sigma_{\Ma \contr a}} \subseteq \supp{\Sigma_{\Ma \del a}}$. In addition, if $a$ is proper, then we have $\rkm(\Ma \del a) = \rkm(\Ma \contr a)-1$. In this case, $\supp{\Sigma_{\Ma \contr a}}$ is of codimension one in $\supp{\Sigma_{\Ma \del a}}$.

\begin{prop} Assume $\Ma$ is a simple matroid and $a\in\Ma$ is a proper element. Let $\Delta'$ be a Bergman fan with support $ \supp{\Delta'} = \supp{\Sigma_{\Ma \contr a}} $ and let $\Sigma'$ be a Bergman fan with support $\supp{\Sigma'} = \supp{\Sigma_{\Ma \del a}}$ such that $\Delta' \subseteq \Sigma'$. Denote by $\Sigma$ the tropical modification of $\Sigma'$ along the cycle $\Delta'$. Then $\Sigma$ is a Bergman fan with support $\supp{\Sigma} = \supp{\Sigma_\Ma}$. Moreover, if $\Delta'$ and $\Sigma'$ are unimodular, then so is $\Sigma'$.
\end{prop}
\begin{proof}
This follows from the observation that $\Sigma_{\Ma}$ is a subdivision the tropical modification of $\Sigma_{\Ma\del a}$ along $\Sigma_{\Ma \contr a}$. More detail is given in \cite{AP-geom}.
\end{proof}

\medskip

\subsection{Tropical smoothness of the support} An extended polyhedral complex $X$ is said to have a \emph{tropically smooth,} or simply \emph{smooth, support} if its support $\supp X$ is locally isomorphic to a product $\eR^a \times \supp \Sigma$ for a non-negative integer $a$ and a Bergman fan $\Sigma$.

\medskip

We prove in~\cite{AP-geom} that the canonical compactification of any unimodular Bergman fan has smooth support. In fact, we have the following more general result.

\begin{thm}\label{thm:smoothness-compact} Let $X$ be a rational polyhedral complex in $N_\R$ with smooth support. Suppose that the recession fan $X_\infty$ is unimodular. Then the canonical compactification $\comp X$ of $X$ has smooth support. Moreover, if the polyhedral structure on $X$ is unimodular, then the induced extended polyhedral structure on $\comp X$, given in Theorem~\ref{thm:orbit-stratum-correspondence}, is unimodular as well.
\end{thm}

This theorem is proved in~\cite{AP-geom}.


\section{Local Hodge-Riemann and Hard Lefschetz} \label{sec:local}
The aim of this section is to study the local situation and to prove Theorem~\ref{thm-intro:HL-HR-local}.

The results in this section are valid for integral, rational, or real Chow rings of simplicial Bergman fans when they are definable. Since our primary objective in this paper is to establish global versions of these results, and we do not have tried to control the torsion part of the tropical cohomology for more general tropical varieties, to simplify the presentation, we will only consider the case of rational fans and we suppose the coefficient ring $\corps$ be either $\Q$ or $\R$.

\subsection{Basic definitions and results}

In this section, we review some basic results concerning Hodge theory of finite dimensional algebras.

\medskip

Let $A^\bul = \bigoplus_{i\in \mathbb Z} A^i$ be a graded vector space of finite dimension over $\corps$ such that the graded pieces of negative degree are all trivial.

\medskip

\begin{defi}[Poincar\'e duality]  We say that $A^\bul$ verifies
\begin{itemize}[label=-]
\item  \emph{weak Poincar\'e duality for an integer $d$} denoted by $\WPD^d(A^\bul)$ if for any integer $k$, we have the equality $\dim_{\corps} A^k = \dim_{\corps} A^{d-k}$.
\end{itemize}
Let $\Phi \colon A^\bul\times A^\bul \to \corps$ be a symmetric bilinear form on $A^\bul$. We say that the pair $(A^\bul, \Phi)$ verifies
\begin{itemize}[label=-]
\item  \emph{Poincar\'e duality for an integer $d$} that we denote by $\PD^d(A^\bul, \Phi)$ if $\Phi$ is non-degenerate and for any integer $k$, $A^k$ is orthogonal to $A^l$ for all $l \neq d-k$. \qedhere
\end{itemize}
\end{defi}

\medskip

Let $A^\bul$ be a graded $\corps$-vector space of finite dimension as above, let $d$ be a positive integer, and let $Q$ be a symmetric bilinear form on the truncation $\bigoplus_{k\leq \frac d2} A^k$ such that $A^i$ and $A^j$ are orthogonal for distinct pair of non-negative integers $i,j \leq \frac d2$.

\begin{defi}[Hard Lefschetz] Notations as above, we say that the pair $(A^\bul,Q)$ verifies
\begin{itemize}[label=-]
\item  \emph{Hard Lefschetz property} denoted $\HL^d(A^\bul, Q)$ if $A^\bul$ verifies $\WPD^d(A^\bul)$ and $Q$ is non-degenerate.
\end{itemize}
\end{defi}

In the sequel, for any $k\leq d/2$, we denote by $Q_k$ the restriction of $Q$ to $A^k$.

\begin{defi}[Hodge-Riemann] We say the pair $(A^\bul, Q)$ verifies
\begin{itemize}[label=-]
\item \emph{Hodge-Riemann relations} $\HR^d(A^\bul, Q)$ if $A^\bul$ verifies $\WPD^d(A^\bul)$ and moreover, for each non-negative integer $k\leq \frac d2$, the signature of $Q_k$ is given by the sum
\[ \sum_{i=0}^k (-1)^i\bigl(\dim(A^i)-\dim(A^{i-1})\bigr). \qedhere \]
\end{itemize}
\end{defi}

\begin{remark} As we will show in Proposition \ref{prop:HR} below, the above definitions are in fact equivalent to the ones given in the introduction.
\end{remark}

\subsubsection{Bilinear form associated to a degree map} Assume that $A^\bul$ is a unitary graded $\corps$-algebra with unit $1_{A^\bul}$ and suppose $A^0 = \corps \cdot 1_{A^\bul}$. Moreover, suppose that we are given a linear form $\deg\colon A^d \to \corps$ called the \emph{degree map} which we extend to the entire algebra $A^\bul$ by declaring it to be zero on each graded piece $A^k$ for $k\neq d$. In this way, we get a bilinear form on $A^\bul$ denoted by $\Phi_{\deg}$ and defined by
\[ \begin{array}{rccl}
  \Phi_{\deg} \colon & A^\bul \times A^\bul & \to & \corps, \\
  & (a, b) & \mapsto & \deg(ab).
\end{array} \]

\subsubsection{Lefschetz operator and its associated bilinear form}\label{sec:lefschetz} Notations as in the preceding paragraph, let now $L\colon A^\bul\to A^{\bul+1}$ be a morphism of degree one which sends any element $a$ of $A^\bul$ to the element $La\in A^{\bul+1}$. Let $\Phi$ be a symmetric bilinear form on $A^\bul$. For example, if we have a degree map on $A^\bul$, we can take $\Phi = \Phi_{\deg}$.

\medskip

We define a bilinear form $Q_{L,\Phi}$ on $A^\bul$ by setting
\[\forall\,\, a\in A^k,\, b\in A^j, \qquad Q_{L,\Phi}(a, b):=\Phi(a,L^{d-2k}b)\]
for all pairs of integers $k,j\leq d/2$. We furthermore require $A^j$ and $A^k$ to be orthogonal for $j\neq k$. In the case $\Phi = \Phi_{\deg}$, we will have $Q_{L,\Phi}(a, b) = \deg(a\cdot L^{d-2k}b)$ and this latter requirement will be automatic.

\medskip

Assuming that $L$ is auto-adjoint with respect to $Q_{L, \Phi}$, i.e., $\Phi(a,L b) = \Phi(La, b)$ for all $a,b \in A^\bul$, we get a symmetric bilinear form $Q = Q_{L, \Phi}$ on $A^\bul$.

\medskip

If $L$ is given by multiplication by an element $\ell\in A^1$ and $\Phi=\Phi_{\deg}$, we write $\ell$ instead of $L$ for the corresponding map, and write $Q_{\ell}$ for the corresponding bilinear form, which is then automatically symmetric: in fact, we have
\[Q_{\ell}(a,b)=\deg(\ell^{d-2k}ab)\]
for any pair of elements $a,b \in A^k$. In this case, we write $\HL^d(A^\bul, \ell)$ and $\HR^d(A^\bul, \ell)$ instead of $\HL^d(A^\bul, Q_{\ell,\Phi_{\deg}})$ and $\HR^d(A^\bul, Q_{\ell,\Phi_{\deg}})$, respectively. If there is no risk of confusion, we furthermore omit the mention of the top dimension $d$ in the above notations. Moreover, we sometimes use $\deg$ instead of $\Phi_{\deg}$ in order to simplify the notation.

\subsubsection{Primitive parts induced by a Lefschetz operator} For a Lefschetz operator $L\colon A^\bul \to A^{\bul+1}$, as in the preceding section, and for a non-negative integer $k\leq \frac d2$,
the \emph{primitive part of $A^k$ with respect to $L$} denoted by $A^k_{\prim, L}$ is by definition the kernel of the map
\[L^{r-2k+1} \colon A^{k} \longrightarrow A^{k-r+1}.\]

If $L$ is given by an element $\ell \in A^1$, we denoted by $A^k_{\prim, \ell}$ the primitive part of $A^k$ with respect to $\ell$.

If there is no risk of confusion, and $L$ or $\ell$ are understood from the context, we simply drop them and write $A^k_{\prim}$ for the primitive part of $A^k$.

\subsubsection{Basic implications and consequences} For future use, we state here basic results around these notions and some of their consequences.

We start by the following proposition which gathers several useful implications and relations between the properties introduced above.

We use the notations of the preceding paragraph. We suppose the Lefschetz operator $L$ is auto-adjoint with respect to $\Phi$ so that the corresponding bilinear form becomes symmetric.

\begin{prop}  \label{prop:HR} Notations as above, the following properties are verified.

\let\olditem\item
\renewcommand{\item}{\stepcounter{equation}\olditem}
\begin{enumerate}[label=\textnormal{(\theequation)}, ref={\theequation}, labelindent={0pt}, leftmargin=*]
\item \label{prop:HR:implication} \ \vspace*{-1.3em} \[ \HR^d(A^\bul, L, \Phi)\Longrightarrow\HL^d(A^\bul, L, \Phi)\Longrightarrow\PD^d(A^\bul, \Phi)\Longrightarrow\WPD^d(A^\bul). \]

\item \label{prop:HL:implication} {\setlength{\abovedisplayskip}{-1em}
  \[ \HL^d(A^\bul, L, \Phi)\iff \begin{cases}
    \PD^d(A^\bul, \Phi)\myand\\
    \forall\, k\leq d/2, \,\, L^{d-2k}\colon A^k\to A^{d-k}\text{ is an isomorphism.}
  \end{cases} \]}

\item \label{prop:HR-sign} {\setlength{\abovedisplayskip}{-.5em}
  \[ \HR^d(A^\bul, L, \Phi)\iff
  \begin{cases}
    \WPD^d(A^\bul) \\
    \forall k\leq d/2, \text{ the restriction $(-1)^kQ_{L,\Phi}\rest{A^k_{\prim}}$ is positive definite.}
  \end{cases} \]}

\item \label{prop:HR:automorphism}
  For $S\colon A^\bul \to A^\bul$ an automorphism of degree 0,
  \begin{align*}
  \PD(A^\bul,\Phi) &\iff \PD(A^\bul,\Phi\circ S), \\
  \HL(A^\bul,Q)    &\iff \HL(A^\bul,Q\circ S),    \\
  \HR(A^\bul,Q)    &\iff \HR(A^\bul,Q\circ S).
  \end{align*}

\item \label{prop:HR:oplus}
  For two graded $\corps$-algebras $A^\bul$ and $B^\bul$, we have
  \begin{align*}
  \WPD(A^\bul)       \myand \WPD(B^\bul)       &\iff \WPD(A^\bul)       \myand \WPD(A^\bul\oplus B^\bul),                    \\
  \PD(A^\bul,\Phi_A) \myand \PD(B^\bul,\Phi_B) &\iff \PD(A^\bul,\Phi_A) \myand \PD(A^\bul\oplus B^\bul,\Phi_A\oplus \Phi_B), \\
  \HL(A^\bul,Q_A)    \myand \HL(B^\bul,Q_B)    &\iff \HL(A^\bul,Q_A)    \myand \HL(A^\bul\oplus B^\bul,Q_A\oplus Q_B),       \\
  \HR(A^\bul,Q_A)    \myand \HR(B^\bul,Q_B)    &\iff \HR(A^\bul,Q_A)    \myand \HR(A^\bul\oplus B^\bul,Q_A\oplus Q_B).
  \end{align*}

\item \label{prop:HR:otimes}
  For two graded $\corps$-algebras $A^\bul$ and $B^\bul$, we have
  \begin{align*}
  \WPD(A^\bul)\myand\WPD(B^\bul)             &\Longrightarrow \WPD(A^\bul\otimes B^\bul), \\
  \PD(A^\bul,\Phi_A)\myand\PD(B^\bul,\Phi_B) &\Longrightarrow \PD(A^\bul\otimes B^\bul,\Phi_A\otimes \Phi_B), \\
  \HL(A^\bul,Q_A)\myand\HL(B^\bul,Q_B)       &\Longrightarrow \HL(A^\bul\otimes B^\bul,Q_A\otimes Q_B), \\
  \HR(A^\bul,Q_A)\myand\HR(B^\bul,Q_B)       &\Longrightarrow \HR(A^\bul\otimes B^\bul,Q_A\otimes Q_B).
  \end{align*}

\item \label{prop:HR:poincare_plus}
  For two graded $\corps$-algebras $A^\bul$ and $B^\bul$
  \begin{align*}
  \Phi_{\deg_A}\otimes\Phi_{\deg_B}   &= \Phi_{\deg_A\otimes\deg_B}, \text{ and } \\
  Q_{L_A,\Phi_A}\oplus Q_{L_B,\Phi_B} &= Q_{L_A\oplus L_B,\Phi_A\oplus\Phi_B}.
  \end{align*}

\item \label{prop:HR:otimes_lefschetz}
  For two graded $\corps$-algebras $A^\bul$ and $B^\bul$, we have
  \[ \HR(A^\bul,L_A,\Phi_{A})\myand\HR(B^\bul,L_B,\Phi_B)\Longrightarrow\HR(A^\bul\otimes B^\bul, L_A\otimes\id+\id\otimes L_B,\Phi_A\otimes \Phi_B). \]
\end{enumerate}
\end{prop}

Most of these statements are routine. We only give the proof of the first three and the very last statements.
\begin{proof}[Proof of~\eqref{prop:HR:implication}]
The last implication is trivial. If $\HL^d(A^\bul, L, \Phi)$ holds, then the bilinear form $Q=Q_{L, \Phi}$ is non-degenerate. This implies that for each non-negative integer $k \leq d/2$, the map $L^{d-2k}$ is injective on $A^k$. By $\WPD^d(A^\bul)$, this implies that $L^{d-2k}$ induces an isomorphism between $A^k$ and $A^{d-k}$, and so, as $A^j$ and $A^k$ are orthogonal for distinct pair of integers $j$ and $k$, by our requirement, this implies that $\Phi$ is non-degenerate and so we have $\PD^d(A^\bul, \Phi)$.

It remains to prove the first implication. We prove by induction on $k$ that $Q_k$ is non-degenerate. The case $k=0$ is immediate by the assumption on the signature and the fact that $A^0 =\corps\cdot 1_{A^\bul}$. Assume this holds for $k-1$ and $k \leq d/2$. Since $L$ is auto-adjoint, we have for any $a,b \in A^{k-1}$,
\[Q_{k-1}(a,b) = Q_{k}(La,Lb)\]
which implies that $L\colon A^{k-1} \to A^k$ is injective, and $Q_k$ restricted to $L(A^{k-1})$ is non-degenerate and has signature that of $Q_{k-1}$. Since the difference between the signature of $Q_k$ and that of $Q_{k-1}$ is $(-1)^k \bigl(\dim(A^k) - \dim(A^{k-1})\bigr) = (-1)^k \bigl(\dim(A^k) - \dim(LA^{k-1})\bigr)$, it follows that $Q_k$ is non-degenerate.
\end{proof}
\begin{proof}[Proof of~\eqref{prop:HL:implication}] The implication $\Rightarrow$ follows from the proof of~\eqref{prop:HR:implication}. The reverse implication is a consequence of our requirement that $A^j$ and $A^k$ are orthogonal for $\Phi$ if $j$ and $k$ are distinct, which combined with the non-degeneracy of $\Phi$ and the isomorphism between $A^k$ and $A^{d-k}$ induced by $L^{d-2k}$, implies that $Q_k$ is non-degenerate.
\end{proof}
\begin{proof}[Proof of~\eqref{prop:HR-sign}] Let $Q = Q_{L,\Phi}$. By~\eqref{prop:HR:implication}, $\HR^d(A^\bul, L, \Phi)$ implies $\HL^d(A^\bul, L, \Phi)$. By Proposition~\ref{prop:lefschetzdecomposition} below this leads to an orthogonal decomposition
\[ A^{k}=\bigoplus_{i=0}^k L^{k-i}A^{i}_{\prim} = A^k_\prim \oplus LA^{k-1}\]
where the restriction of $Q_k$ to $LA^{k-1}$ coincides with that of $Q_{k-1}$ on $A^{k-1}$. Proceeding by induction and using $\HR^d(A^\bul, L, \Phi)$ which gives the signature of $Q_k$, we get the implication $\Rightarrow$.

\medskip

It remains to prove the reverse implication $\Leftarrow$. We proceed by induction and show that the bilinear form $Q_k$ is non-degenerate, that $L^{d-2k}\colon A^{k}\to A^{d-k}$ is an isomorphism, and that we have an orthogonal decomposition $A^k = A^k_\prim \bigoplus LA^{k-1}$. This shows that $L$ is an isometry with respect to $Q_{k-1}$ and $Q_k$ on the source and target, respectively, and the result on signature follows by the positivity of $(-1)^k Q_k$ on the primitive part combined with $\dim A^k_\prim = \dim(A^{k}) -\dim(A^{k-1})$.

The case $k=0$ comes from the assumption on the positivity of $Q_0$, since $A^0=A^0_\prim$, which implies $L^d \colon A^0 \to A^d$ is injective, and the property $\WPD^d(A^\bul)$. Suppose the property holds for $k-1$ and that $k\leq d/2$. By the assumption of our induction, $L\colon A^{k-1} \to A^k$ is injective, and in addition, it preserves the bilinear forms $Q_{k-1}$ and $Q_k$. Moreover, $A^k_\prim$ is orthogonal to $LA^{k-1}$. We now explain why an element $a \in A^k $ which is orthogonal to $LA^{k-1}$ is in the primitive part $A^k_\prim$: for such an element, and for an element $b\in A^{k-1}$, we have
\[Q_k(a, Lb) = \Phi(a, L^{d-2k+1}b) = \Phi(L^{d-2k+1}a, b) =0.\]
Since $Q_{k-1}$ is non-degenerate, using $\WPD^{d}$, we infer the bilinear form $\Phi$ restricted to $A^{d-k+1} \times A^{k-1}$ is non-degenerate. This implies that $L^{d-2k+1}a$ is zero and so $a\in A^k_\prim$.

\medskip

We thus get the orthogonal decomposition $A^k = A^k_\prim \oplus LA^{k-1}$. We infer that $Q_k$ is non-degenerate, and using $\WPD^d(A^\bul)$, we get the isomorphism $L^{d-2k} \colon A^k \to A^{d-k}$.
\end{proof}

\begin{proof}[Proof of \eqref{prop:HR:otimes_lefschetz}]  By~\eqref{prop:HR:implication} and Proposition~\ref{prop:lefschetzdecomposition}, we can decompose each graded piece $A^k$ and $B^l$ for $k\leq d_A/2$ and $l \leq d_B/2$ into the orthogonal sum
\[ A^{k}=\bigoplus_{i=0}^k L_A^{k-i}A^{i}_{\prim} \qquad B^{l}=\bigoplus_{i=0}^l L_B^{l-i}B^{i}_{\prim}\]
such that $Q_{A, k}$ and $Q_{B,j}$ are definite on $A^k_\prim$ and $B^j_\prim$, respectively. Taking an orthonormal basis of the primitive parts $A^k_\prim$ and $B^j_\prim$, and making a shift in the degree of graded pieces, by Proposition \ref{prop:HR}, we can reduce to verifying the statement in the case where $A^\bul =\rquot{\corps[\x]}{\x^{r+1}}$ and $B^\bul = \rquot{\corps[\y]}{\y^{s+1}}$, for two non-negative integers $r$ and $s$, and $L_A$ and $L_B$ are multiplications by $\x$ and $\y$, respectively. The proof in this case can be obtained either by using Hodge-Riemann property for the projective  complex variety $\mathbb C\P^{r} \times \mathbb C\P^s$, or by the direct argument given in~\cite{BBFK}*{Proposition 5.7}, or still by the combinatorial argument given in~\cite{McD}*{Lemma 2.2} and~\cite{AHK}*{Lemma 7.8} based on the use of Gessel-Viennot-Lindstr\"om lemma~\cites{GV, Lin}.
\end{proof}

\begin{prop}[Lefschetz decomposition]\label{prop:lefschetzdecomposition}
Assume that $\HL^d(A^\bul, L, \Phi)$ holds. Then for $k\leq\frac{d}2$, we have the orthogonal decomposition \[ A^{k}=\bigoplus_{i=0}^k L^{k-i}A^{i}_{\prim}, \]
where, for each $i\in\zint0k$, the map $L^{k-i}\colon A^{i}_{\prim} \simto L^{k-i}A^{i}_{\prim}$ is an isomorphism which preserves respective bilinear forms \emph{($Q_i$ on $A^i$ and $Q_k$ on $A^{k}$)}.
\end{prop}
\begin{proof} First, we observe that for $i\leq k$, the map $L^{k-i}$ induces an isomorphism between $A^{i}$ and its image $L^{k-i}A^{i}$ in
$A^k$. This follows from the fact previously stated that for $a,b \in A^{j}$,
\[Q_{j}(a,b) = Q_{j+1}(La,Lb)\]
which implies that $L\colon A^{j} \to A^{j+1}$ is injective provided that $j+1 \leq d/2$. Note that $A_\prim$ is orthogonal to $LA^{k-1}$. Indeed for $a\in A_\prim$ and $b\in A^{k-1}$, we have
\[\Phi(a , L^{d-2k} Lb) = \Phi(L^{d-2k+1}a , b)=0. \]
Moreover, an element $a \in A^{k}$ orthogonal to $LA^{k-1}$ is necessarily in the primitive part, \cf. the proof of~\eqref{prop:HR-sign}.
We thus get the orthogonal decomposition
\[A^k = A^k_{\prim} \oplus L A^{k-1}.\]
Proceeding by induction on $k$, this leads to the orthogonal decomposition
\[A^k = A_\prim^k \oplus L A^{k-1} = \bigoplus_{i=0}^k L^{k-i}A^{i}_{\prim}. \qedhere \]
\end{proof}

The following proposition is straightforward.

\begin{prop}\label{prop:HRbis} The properties $\HR$ and $\HL$ are open conditions in $Q$ (and so in $\Phi$ and $L$). The property $\PD$ is open in $\Phi$. Moreover, $\HR$ is a closed condition in $Q$ when restricted to the space $\Bigl\{Q\st\HL^d(A^\bul, Q)\Bigr\}$ consisting of those $Q$ which verify $\HL$.
\end{prop}

\subsection{Chow rings of fans} For a simplicial fan $\Sigma$ of dimension $d$ in a real vector space $V$, we denote by $A^\bul(\Sigma)$ its Chow ring with $\corps$ coefficients. It is defined as follows.

For each ray $\varrho \in \Sigma_1$, pick a vector $\e_\varrho$ which generates $\varrho$.
If $\Sigma$ is a rational fan with respect to a lattice, a natural choice for the generating vector is the primitive vector in the lattice which generates the ray.

Consider the polynomial ring $\corps[\x_\varrho]_{\varrho\in \Sigma_1}$ with indeterminate variables $\x_\varrho$ associated to rays $\varrho$ in $\Sigma_1$. The Chow ring $A^\bul(\Sigma)$ is by definition the quotient ring
\[A^\bul(\Sigma) := \rquot{\corps[\x_\varrho]_{\varrho\in \Sigma_1}}{\bigl(\I_1 + \I_2\bigr)}\]
where
\begin{itemize}
\item $\I_1$ is the ideal generated by the products $\x_{\varrho_1}\x_{\varrho_2} \dots \x_{\varrho_k}$ for any positive integer $k$ and rays $\varrho_1, \dots, \varrho_k$ which do not form a cone in $\Sigma$;
\item $\I_2$ is the ideal generated by the homogeneous polynomials of degree one of the form $\sum_{\varrho\in \Sigma_1} f(\e_\varrho) \x_\varrho$ where $f \in V^\dual$ is a linear form on the vector space $V$, and $f(\varrho)$ is the evaluation of $f$ at the generator $\e_\varrho$ of $\varrho$.
\end{itemize}

Note that since the ideal $\I_1+\I_2$ is homogeneous, the Chow ring inherits a graded ring structure. Moreover, the definition is well-posed even in the case of non-rational fans, in which case $\corps=\R$, as any two different choices of generating vectors for the rays lead to an isomorphism of the quotient rings. (This isomorphism sends the generator $\x_\varrho$ in the Chow ring defined with respect to the first generator set to $c_\varrho \x_\varrho$ in the second, for a constant $c_\varrho \in \R_{>0}$.)

\medskip

For each ray $\varrho$ of $\Sigma$, we denote by $x_\varrho$ the image of $\x_\varrho$ in $A^1(\Sigma)$.\\
Any cone-wise linear function $f$ on $\Sigma$ gives an element of $A^1(\Sigma)$ defined as
\[\sum_{\varrho \in \Sigma_1} f(\e_\varrho) x_\varrho.\]
In fact, all elements of $A^1(\Sigma)$ are of this form, and $A^1(\Sigma)$ can be identified with the space of cone-wise linear functions on $\Sigma$ modulo linear functions.

\medskip

In the case the fan $\Sigma$ is rational and unimodular with respect to a lattice $N$ in $V$, one can define Chow rings with integral coefficients by choosing $\e_\varrho$ to be the primitive vector of $\varrho$, for each ray $\varrho$, and by requiring $f$ in the definition of $\I_2$ to be integral.
In this case, we have the following characterization of the Chow ring, cf.~\cites{Dan78, BDP90, Bri96, FS}.
\begin{thm}
Let $\Sigma$ be a unimodular fan, and denote by $\P_\Sigma$ the corresponding toric variety. The Chow ring $A^\bullet(\Sigma)$ is isomorphic to the Chow ring of $\P_\Sigma$.
\end{thm}

In the sequel, unless otherwise stated, we will be only considering Bergman fans, which are thus fans on the support of the Bergman fan of a matroid. Moreover, we suppose that all the matroids are simple which means they are loopless and do not have parallel elements. We call as before \emph{Bergman support} the support of the Bergman fan of any matroid.

\subsection{Bergman supports are smooth} Let $\Sigma$ be a Bergman fan. Recall that for any $\sigma \in \Sigma$, we denote by $\Sigma^\sigma$ the star fan of $\sigma$ which lives in the vector space $N^{\sigma}_\R = \rquot{N_\R}{N_{\sigma, \R}}$.

\begin{prop} \label{prop:loc-smooth}For any Bergman fan $\Sigma$, all the star fans $\Sigma^\sigma$ are Bergman. It follows that Bergman supports are smooth. Moreover, $\Sigma$ is connected in codimension one.
\end{prop}
\begin{proof} See e.g.~\cites{AP, Sha13a}.
\end{proof}

\subsection{The canonical element} The following proposition allows to define a degree map for the Chow ring of a Bergman fan.

\begin{prop} \label{prop:canonical}Let $\Sigma$ be a unimodular Bergman fan of dimension $d$. For each cone $\sigma\in\Sigma_d$, the element
\[\prod_{\varrho\subface \sigma \\ \dims \varrho =1}x_\varrho\in A^d(\Sigma)\]
does not depend on the choice of $\sigma$.
\end{prop}
We call this element the \emph{canonical element} of $A^d(\Sigma)$ and denote it by $\omega_\Sigma$.

\begin{proof} By Proposition~\ref{prop:loc-smooth}, $\Sigma$ is smooth. Consider a face $\tau$ of codimension one in a $d$-dimensional cone $\sigma$ of $\Sigma$. The star fan $\Sigma^\tau$ is Bergman. It has dimension one, and it follows by unimodularity that the primitive vectors $\e_{\eta/\tau}$ of the rays corresponding to $\eta$ in $\Sigma^\tau$, for $\eta \ssupface \tau$, form a circuit in $N^\tau_\R$: this means we have
\[\sum_{\eta \ssupface \tau} \e_{\eta/\tau} =0\]
and, moreover, this equation and its scalar multiples are the unique linear relations between these vectors. Combined with unimodularity, this implies that for any facet $\zeta \neq \sigma$ in $\Sigma$ which contains $\tau$, we can find a linear function $f \in {N^{\tau}}^\vee$ such that, viewing $f$ as an element of $N^\vee$ which vanishes on $N_\tau$, we have
\begin{itemize}
\item $f$ vanishes on $\nvect_{\eta/\tau}$ for any facet $\eta \ssupface \tau$ distinct from $\sigma$ and $\zeta$, and
\item $f$ takes values $f(\nvect_{\sigma/\tau}) = -f(\nvect_{\zeta/\tau}) =1$.
\end{itemize}
If we denote by $\rho_{\eta/\tau}$ the unique ray of $\eta\in\Sigma$ which is not in $\tau$, we get
\[\Bigr(\sum_{\eta \ssupface \tau} f(\nvect_{\eta/\tau}) x_{\rho_{\eta/\tau}}\Bigl) \prod_{\varrho \subface \tau \\ \dims \varrho=1}x_\varrho = \Bigr(\sum_{\varrho\in\Sigma_1} f(\varrho) x_{\varrho}\Bigl) \prod_{\varrho \subface \tau \\ \dims \varrho=1}x_\varrho  = 0.\]
We infer that
\[\prod_{\varrho \subface \sigma \\ \dims\varrho=1}x_\varrho  = \prod_{\varrho \subface \zeta \\ \dims \varrho=1}x_\varrho. \]

This shows that for any two facets $\sigma$ and $\zeta$ which share a face of dimension $d-1$, the two associated elements
$\prod_{\varrho\subface \sigma \\ \dims \varrho =1}x_\varrho$ and $\prod_{\varrho\subface \zeta \\ \dims \varrho =1}x_\varrho$ coincide in $A^d(\Sigma)$. Since $\Sigma$ is connected in codimension one, this proves the proposition.
\end{proof}

\begin{prop}
Let $\Sigma$ be a unimodular Bergman fan. Then, for each $k\in\{0,\dots,d\}$, $A^k(\Sigma)$ is generated by square-free monomials.
\end{prop}
\begin{proof}
To prove this statement, one can proceed by a lexicographic induction and use the relations defining the Chow ring in order to replace monomials which have exponents larger than one by a linear combination of square-free monomials, as in~\cite{AHK} and \cite{Ami}.
\end{proof}

\begin{cor}\label{prop:non-vanishing} Let $\Sigma$ be a unimodular Bergman fan. Then we have $\omega_\Sigma \neq 0$ and $A^d(\Sigma)$ is generated by $\omega_\Sigma$.
\end{cor}
\begin{proof} This follows from the fact that $A^d(\Sigma)$ is generated by square-free monomials associated to top-dimensional cones, which by Poincar\'e duality $A^0(\Sigma) \simeq A^{d}(\Sigma)^\dual$ and Proposition \ref{prop:canonical} ensures the non-vanishing of $\omega_\Sigma$.
\end{proof}

\medskip

Let $\deg\colon A^d(\Sigma)\to\corps$ be the linear map which takes value 1 on the canonical element. This degree map and the corresponding bilinear map $\Phi_{\deg}$ is systematically used below, so we use our convention and omit the mention of $\Phi$ and $\deg$ in $\HL$ and $\HR$.

\subsection{Chow rings of products} Let $\Sigma$ and $\Sigma'$ be two Bergman fans. By Proposition~\ref{prop:productBergman} the product $\Sigma \times \Sigma'$ is Bergman. We have the following proposition for the Chow ring of the product fan whose proof is straightforward and is thus omitted.

\begin{prop} \label{prop:cartesian_product}
Let $\Sigma$ and $\Sigma'$ be two fans. There exists a natural isomorphism of rings
\[ A^\bul(\Sigma\times\Sigma')\simeq A^\bul(\Sigma)\otimes A^\bul(\Sigma'). \]
Moreover, under this isomorphism, one has
\begin{gather*}
\omega_{\Sigma\times\Sigma'}=\omega_\Sigma\otimes\omega_{\Sigma'}, \\
\deg_{\Sigma\times\Sigma'}=\deg_\Sigma\otimes\deg_{\Sigma'}.
\end{gather*}
\end{prop}

\subsection{Restriction and Gysin maps} For a pair of cones $\tau \subface \sigma$ in $\Sigma$, we define the restriction and Gysin maps $\i^*_{\tau \subface \sigma}$ and $\gys_{\sigma \supface \tau}$ between the Chow rings of $\Sigma^\tau$ and $\Sigma^\sigma$.

\medskip

The restriction map
\[\i^*_{\tau \subface \sigma}\colon A^\bul(\Sigma^\tau) \to A^\bul(\Sigma^\sigma)\]
is a graded $\corps$-algebra homomorphism given on generating sets
\[\i^*_{\tau \subface \sigma}\colon A^1(\Sigma^\tau) \to A^1(\Sigma^\sigma)\]
by
\[
\forall \, \rho \textrm{ ray in } \Sigma^\tau, \qquad \i^*_{\tau \subface \sigma} (x_\rho) =
\begin{cases}
  x_\rho & \qquad  \textrm{if $\sigma+\rho$ is a cone $\zeta \supface \sigma$}, \\
  -\sum_{\varrho \in \Sigma^\sigma_1} f(\e_\varrho) x_\varrho & \qquad \textrm{if $\rho \in \sigma$} \\
  0 & \qquad \textrm{otherwise},
\end{cases}
\]
where in the second equality, $f$ is any linear function on $N^\tau$, viewed naturally in $N^\vee$, which takes value $1$ on $\e_\rho$ and value zero on all the other rays of $\sigma$, and the sum is over all the rays $\varrho$ of $\Sigma^\sigma$. Note that any two such choices of $f$ and $f'$ differ by a linear function on $N^\sigma$, and so the element $\sum_{\varrho \in \Sigma^\sigma} f(\e_\varrho) x_\varrho$ in $A^1(\Sigma^\sigma)$ does not depend on the choice of $f$.

\begin{remark} For a pair of cones $\tau \subface \sigma$ as above, we get an inclusion of canonical compactifications $\i \colon \comp\Sigma^\sigma \hookrightarrow \comp\Sigma^{\tau}$. The map $\i^*_{\tau \subface \sigma}\colon A^k(\Sigma^\tau) \to A^k(\Sigma^\sigma)$ coincides with the restriction map
\[\i^* \colon H^{k,k}_\trop(\comp\Sigma^{\tau}) \longrightarrow H^{k,k}_\trop(\comp\Sigma^{\sigma}) \]
via the isomorphism between the Chow and tropical cohomology groups given in Theorem~\ref{thm:HI}.
\end{remark}

The Gysin map is the $\corps$-linear morphism
\[ \gys_{\sigma \supface \tau} \colon A^{\bul}(\Sigma^\sigma) \longrightarrow A^{\bul+ \dims \sigma -\dims \tau}(\Sigma^\tau) \]
defined as follows. Let $r = \dims \sigma -\dims \tau$, and denote by $\rho_1, \dots, \rho_r$ the rays of $\sigma$ which are not in $\tau$.

Consider the $\corps$-linear map
\[\corps[\x_\varrho]_{\substack{\varrho \in \Sigma_1^\sigma}} \longrightarrow \corps[\x_\varrho]_{\substack{\varrho \in \Sigma^\tau_1}}\]
defined by multiplication by $\x_{\rho_1}\x_{\rho_2} \dots \x_{\rho_r}$. Obviously, it sends an element of the ideal $\I_2$ in the source to an element of the ideal $\I_2$ in the target. Moreover, any linear function on $N^\sigma$ defines a linear function on $N^\tau$ via the projection
\[N^\tau = \rquot{N}{N_\tau} \longrightarrow N^\sigma = \rquot{N}{N_{\sigma}}.\]
This shows that the elements of $\I_1$ in the source are sent to elements of $\I_1$ in the target as well. Passing to the quotient, we get a $\corps$-linear map
\[ \gys_{\sigma \supface \tau} \colon A^{k}(\Sigma^\sigma) \longrightarrow A^{k+ \dims \sigma -\dims \tau}(\Sigma^\tau). \]
\begin{remark} For a pair of cones $\tau \subface \sigma$ as above, from the inclusion of canonical compactifications $\i \colon \comp\Sigma^\sigma \hookrightarrow \comp\Sigma^{\tau}$, we get a restriction map
\[\i^*\colon H^{d-\dims \sigma -k,d-\dims \sigma - k}_\trop(\comp\Sigma^{\tau}) \longrightarrow H^{d-\dims \sigma-k,d-\dims \sigma -k}_\trop(\comp\Sigma^{\sigma}) \]
which by Poincar\'e duality on both sides, gives the Gysin map
\[\gys_{\sigma>\tau} \colon H^{k,k}_\trop(\comp\Sigma^{\sigma})  \longrightarrow H^{k+\dims \sigma -\dims\tau,k +\dims \sigma -\dims \tau}_\trop(\comp\Sigma^{\tau}).\]
This map coincides with the Gysin map between Chow groups defined in the preceding paragraph, via the isomorphism of Chow and tropical cohomology given in Theorem~\ref{thm:HI}.
\end{remark}

The following proposition gathers some basic properties of the restriction and Gysin maps.
\begin{prop} \label{lem:i_gys_basic_properties-local}
Let $\tau\ssubface\sigma$ be a pair of faces, and let $x\in A^\bul(\Sigma^\tau)$ and $y\in A^\bul(\Sigma^\sigma)$. Denote by $\rho_{\sigma/\tau}$ the unique ray associated to $\sigma$ in $\Sigma^\tau$, and by $x_{\sigma/\tau}$ the associated element of $A^1(\Sigma^\tau)$. The following properties hold.
\begin{gather}
\text{$\i^*_{\tau\ssubface\sigma}$ is a surjective ring homomorphism.} \label{eqn:i_surjective_homeo-local} \\
\gys_{\sigma\ssupface\tau}\circ\i^*_{\tau\ssubface \sigma}(x)=x_{\sigma/\tau}\cdot x. \label{eqn:gys_circ_i-local} \\
\gys_{\sigma\ssupface\tau}(\i^*_{\tau\ssubface\sigma}(x)\cdot y)=x\cdot\gys_{\sigma \ssupface\tau}(y).  \label{eqn:gys_i_simplification-local}
\end{gather}
Moreover, if $\deg_\tau\colon A^{d-\dims\tau}(\Sigma^\tau)\to\corps$ and $\deg_\sigma\colon A^{d-\dims\sigma}(\Sigma^\sigma)\to\corps$ are the corresponding degrees maps, then
\begin{equation} \label{eqn:deg_circ_gys-local}
\deg_\sigma=\deg_\tau\circ\gys_{\sigma\ssupface\tau}.
\end{equation}
Finally, $\gys_{\sigma\ssupface\tau}$ and $\i^*_{\tau\ssubface\sigma}$ are dual in the sense that
\begin{equation} \label{eqn:i_gys_dual-local}
\deg_\tau(x\cdot\gys_{\sigma\ssupface\tau}(y))=\deg_\sigma(\i^*_{\tau\ssubface\sigma}(x)\cdot y).
\end{equation}
\end{prop}

\begin{proof} In order to simplify the presentation, we drop the indices of $\gys$ and $\i^*$. Properties \eqref{eqn:i_surjective_homeo-local} and \eqref{eqn:gys_circ_i-local} follow directly from the definitions. From Equation \eqref{eqn:gys_circ_i-local}, we can deduce Equation \eqref{eqn:gys_i_simplification-local} by the following calculation. Let $\~y$ be a preimage of $y$ by $\i^*$. Then,
\[ \gys(\i^*(x)\cdot y)=\gys(\i^*(x\cdot \~y))=x_{\delta/\gamma}\cdot x\cdot\~y=x\cdot\gys\circ\i^*(\~y)=x\cdot\gys(y). \]
For Equation \eqref{eqn:deg_circ_gys-local}, let $\eta$ be a maximal cone of $\Sigma^\sigma$. Let $\~\eta$ be the corresponding cone containing $\rho_{\sigma/\tau}$ in $\Sigma^\tau$. The cone $\~\eta$ is maximal in $\Sigma^\tau$. We have the respective corresponding generators $x_{\~\eta}\in A^{d-\dims\tau}(\Sigma^\tau)$ and $x_\eta\in A^{d-\dims\sigma}(\Sigma^\sigma)$. By definition of the  degree maps, $\deg_\tau(x_{\~\eta})=\deg_\sigma(x_\eta)$. Using the definition of $\gys$, we get that $x_{\~\eta}=\gys(x_\eta)$. Thus, we can conclude that
\[ \deg_\sigma = \deg_\tau\circ\gys. \]
Finally, we get Equation \eqref{eqn:i_gys_dual-local}:
\[ \deg_\tau(x\cdot\gys(y))=\deg_\tau(\gys(\i^*(x)\cdot y))=\deg_\sigma(\i^*(x)\cdot y). \qedhere \]
\end{proof}

\subsection{Primitive parts of the Chow ring} Let $\Sigma$ be a Bergman fan of dimension $d$.
For an element $\ell \in A^1(\Sigma)$ and any non-negative integer number $k \leq \frac d2$, the primitive part $A^k_{\prim, \ell}(\Sigma)$ of $A^k(\Sigma)$ is the kernel of the Lefschetz operator given by multiplication with $\ell^{d-2k+1}$
\[\ell^{d-2k+1}\cdot - \colon A^k(\Sigma) \longrightarrow A^{d-k+1}(\Sigma).\]
Let $Q$ be the bilinear form defined on $A^k$, for $k \leq \frac{d}2$, by
\[ \forall \, a,b\in A^{k}(\Sigma), \qquad Q (a,b) = \deg(\ell^{d-2k}ab).\]
Note that by Propositions~\ref{prop:HR} and~\ref{prop:lefschetzdecomposition}, we have the following properties.
\begin{itemize}
\item $\HL(\Sigma, \ell)$ is equivalent to the following: For each $k \leq \frac{d}2$, the map
\[\ell^{d-2k} \cdot - \colon A^{k}(\Sigma) \to A^{d-k}(\Sigma)\]
is an isomorphism.

\item $\HR(\Sigma, \ell)$ is equivalent to the following: The bilinear form $(-1)^k Q(\ccdot,\rdot)$ restricted to the primitive part $A^{k}_{\prim}(\Sigma)$ is positive definite.

\item Assume that $\HL(\Sigma, \ell)$ holds. Then for $k\leq\frac{d}2$, we have the Lefschetz decomposition \[ A^{k}(\Sigma)=\bigoplus_{i=0}^k \ell^{k-i}A^{i}_{\prim}(\Sigma). \]
\end{itemize}

\subsection{Hodge-Riemann for star fans of rays implies Hard Lefschetz} Let $\ell \in A^1(\Sigma)$. We assume that $\ell$ has a representative in $A^1(\Sigma)$ with strictly positive coefficients, i.e.,
\[\ell =\sum_{\varrho \in \Sigma_1} c_\varrho x_\varrho\]
for scalars $c_\varrho >0$ in $\corps$. For each $\varrho \in \Sigma_1$, define $\ell^\varrho = \i^*_{\conezero \subface \varrho}(\ell) \in A^1(\Sigma^\varrho)$.

\medskip

The following is well-known, see e.g.~\cite{CM05} or~\cite{AHK}*{Proposition 7.15}.
\begin{prop} \label{prop:local_HR}
If $\HR(\Sigma^\varrho, \ell^\varrho)$ holds for all rays $\varrho\in\Sigma_1$, then we have $\HL(\Sigma, \ell)$.
\end{prop}

\begin{proof} Let $k \leq \frac d2$ be a non-negative integer where $d$ denotes the dimension of $\Sigma$. We need to prove that the
bilinear form $Q_\ell$ defined by
\[\forall \, a, b \in A^{k}(\Sigma), \qquad Q_\ell (a, b)=  \deg(\ell^{d-2k}ab)\]
is non-degenerate. By Poincar\'e duality for $A^\bul(\Sigma)$, Theorem \ref{thm:pd}, this is equivalent to showing the multiplication map
\[\ell^{d-2k}\cdot - \colon A^{k}(\Sigma) \longrightarrow A^{d-k}(\Sigma)\]
is injective. Let $a \in A^k(\Sigma)$ such that $\ell^{d-2k}\cdot a =0$. We have to show that $a=0$.

\medskip

There is nothing to prove if $k=d/2$, so assume $2k < d$. For each $\varrho \in \Sigma_1$, define $a_\varrho :=\i^*_{\conezero \subface \varrho}(a)$.
It follows that
\[(\ell^\varrho)^{d-2k} \cdot a_\varrho = \i^*_{\conezero \subface \varrho}(\ell^{d-2k}a) =0,\]
and so $a_\varrho \in A^{k}_{\prim, \ell^\varrho}(\Sigma^\varrho)$ for each ray $\varrho \in \Sigma_1$.

\medskip

Using now Proposition~\ref{lem:i_gys_basic_properties-local}, for each $\varrho\in \Sigma_1$, we get
\[\deg_\varrho\bigl((\ell^\varrho)^{d-2k-1}\cdot a_\varrho \cdot a_{\varrho}\bigr) = \deg_\varrho\bigl(\i^*_{\conezero \subface \varrho}(\ell^{d-2k-1}\cdot a) \cdot a_{\varrho}\bigr) = \deg\bigl(\ell^{d-2k-1}\cdot a \cdot x_\varrho \cdot a\bigr).\]

We infer that
\[\sum_{\varrho \in \Sigma_1} c_\varrho \deg_\varrho\bigl((\ell^\varrho)^{d-2k-1}\cdot a_\varrho \cdot a_{\varrho}\bigr) = \deg\bigl(\ell^{d-2k-1}\cdot a \cdot \bigl(\,\sum_{\varrho }c_\varrho x_\varrho\bigr) \cdot a\bigr) = \deg(\ell^{d-2k}\cdot a\cdot a)=0.\]

By $\HR(\Sigma^\varrho,\ell^\varrho)$, since $a_\varrho\in A^{k}_{\prim, \ell^\varrho}(\Sigma^\varrho)$, we have $(-1)^k \deg_\varrho\bigl((\ell^\varrho)^{d-2k-1}\cdot a_\varrho \cdot a_{\varrho}\bigr) \geq 0$ with equality if and only if $a_\varrho =0$. Since $c_\varrho >0$ for all $\varrho$, we conclude that $a_\varrho =0$ for all $\varrho \in \Sigma_1$.

\medskip

Applying Proposition~\ref{lem:i_gys_basic_properties-local} once more, we infer that
\[x_\varrho a =  \gys_{\varrho \supface \conezero}\circ \i^*_{\conezero \subface \varrho} (a) = \gys_{\varrho \supface \conezero}(a_\varrho) =0.\]

By Poincar\'e duality for $A^\bul(\Sigma)$, and since the elements $x_\varrho$ generate the Chow ring, this implies that $a=0$, and the proposition follows.
\end{proof}

\subsection{Isomorphism of Chow groups under tropical modification} In this section, we prove that tropical modifications as defined in Section~\ref{sec:modification} preserve the Chow groups. This will be served as the basis of our induction process to prove $\HR$ and $\HL$ properties for Bergman fans.

\begin{prop} \label{prop:modification} Let $\Ma$ be a simple matroid of rank $d+1$ on a ground set $E$, and let $a$ be a proper element of $E$. Let $\Delta'$ be a unimodular Bergman fan with support $\supp{\Sigma_{\Ma\contr a}}$, and let $\Sigma'$ be a unimodular Bergman fan with support $\supp{\Sigma_{\Ma\del a}}$ which contains $\Delta'$ as a subfan. Let $\Sigma$ be the unimodular fan obtained by the tropical modification of $\Sigma'$ along $\Delta'$. Then we have
\[A^\bul(\Sigma) \simeq A^\bul(\Sigma').\]
\end{prop}

\begin{proof} We follow the notation of Section~\ref{sec:modification}. So $\Sigma$ is embedded in $N_\R$ for $N = \rquot{\Z^E}{\Z \e_E}$, and $(\e_b)$ for $b\in E$ represents the standard basis of $\Z^E$. We denote by $M = N^\vee$ the dual lattice. We also denote by $N^a$ the lattice $\rquot{\Z^{E - a}}{\Z \e_{_{E-a}}} = \rquot{N}{\Z \e_a}$ and $M^a$ the corresponding dual lattice. Thus, we have $\Sigma'\subset N^a_\R$. The tropical modification induces an injective map of fans $\Gamma\colon \Sigma'\to\Sigma$ which preserves the fan structure. We denote by $\rho_a$ the ray in $\Sigma$ associated to $a$; note that with our convention, we have $\e_{\rho_a} =\e_a$. The ray $\rho_a$ is the only one which is not in $\Gamma(\Sigma')$.

\medskip

Let $h$ be a linear form in $M_\corps$ which takes value one on $\e_a$. Then $h$ induces a piecewise linear function $f := \Gamma^{-1}\circ h\rest{\Gamma(\Sigma')}$ on $\Sigma'$ and we have $\div(f) = -\Delta'$, in the sense that the sum of slopes of $f$ around all codimension one faces $\sigma\in\Sigma'_{d-1}\setminus\Delta'_{d-1}$ is zero while the sum of slopes is $-1$ for all $\sigma\in\Delta'_{d-1}$.

Consider the surjective $\corps$-algebra homomorphism
\[\phi\colon \corps[\x_\varrho]_{\varrho\in \Sigma_1} \longrightarrow \corps[\x_\varrho]_{\varrho\in \Sigma'_1}\]
which on the level of generators is defined by sending $\x_{\rho_a}$ to $- \sum_{\varrho \in\Sigma_1 - \rho_a} f(\e_\varrho)\x_\varrho$, and which maps $\x_\varrho$ to $\x_{\Gamma^{-1}(\varrho)}$ for rays $\varrho$ of $\Sigma_1$ different from $\rho_a$. We will show that $\phi$ induces an isomorphism between $A^\bul(\Sigma)$ and $A^\bul(\Sigma')$. For this, we first show that we have $\phi(\I_2) \subseteq \I_2$ and $\phi(\I_1) \subseteq \I_1+\I_2$, which allows to pass to the quotient and to get a morphism of $\corps$-algebra between the two Chow rings.

\medskip

To show that $\phi(\I_2) \subseteq \I_2$, consider a linear form $l$ on $N_\corps$. Let $l' := l - l(\e_a)h$. We have $l'(\e_a)=0$ and so $l'$ gives a linear form on $N^a_\corps$. We have
\[ \phi(\sum_{\varrho\in\Sigma_1} l(\e_\varrho)\x_\varrho) = \sum_{\varrho\in \Sigma_1-\rho_a} l(\e_\varrho)\x_\varrho+l(\e_a)\phi(\x_{\rho_a})=\sum_{\varrho\in \Sigma'_1} l'(\e_\varrho)x_\varrho. \]
This shows that
\[\phi(\I_2) \subseteq \I_2.\]

We now consider the image of $\I_1$. Consider a collection of distinct rays $\varrho_0, \dots, \varrho_k$, $k\in \mathbb N$, and suppose they are not comparable, which translates to $\x_{\varrho_0} \dots \x_{\varrho_k} \in \I_1$.

\medskip

If $\varrho_0, \dots, \varrho_k \neq \rho_a$, then they all belong to $\Sigma'$ and we have $\phi(\x_{\varrho_0} \dots \x_{\varrho_k}) = \x_{\varrho_0} \dots \x_{\varrho_k} \in \I_1$.

\medskip

Otherwise, one of the rays, say $\varrho_0$, is equal to $\rho_a$. Two cases can happen:

\begin{enumerate}
\item \label{prop:topmodif_case1} Either, $\varrho_1, \dots, \varrho_k$ belong to $\Delta'$,
\item \label{prop:tropmodif_case2} Or, one of the rays, say $\varrho_k$, does not belong to $\Delta'$.
\end{enumerate}

In the first case, since $\rho_a, \varrho_1, \dots, \varrho_k$ are not comparable, we deduce that the rays $\varrho_1, \dots, \varrho_k$ are not comparable in $\Delta'$, and again, we get
\[\phi(\x_{\varrho_0}\x_{\varrho_1} \dots \x_{\varrho_k}) = \phi(\x_{\varrho_0}) \x_{\varrho_1} \dots \x_{\varrho_k} \in \I_1.
\]
Suppose we are in case~\ref{prop:tropmodif_case2}, so $\varrho_k \not \in \Delta_1'$. It follows from Proposition~\ref{prop:voisinage} below applied to $\varrho =\varrho_k$ that there exists a neighborhood of $\e_{\varrho_k}$ in $\Sigma$ which is entirely included in a hyperplane $H$ in $N_\corps$ which does not contain $\rho_a$. From this we infer that any ray comparable with $\varrho_k$ should belong to $H$. Moreover, the projection map $N \to N^a$ induces an isomorphism $H\to N^a_\corps$. It follows that $h\rest H$ induces a linear function $l$ on $N^a_\corps$, which coincides with $h$ on all rays $\xi$ which are comparable with $\varrho_k$.\\
Therefore,
\begin{align*}
\phi(\x_{\rho_a}\x_{\varrho_1} \dots \x_{\varrho_k})
&= - \sum_{\xi\in\Sigma_1-\rho_a} h(\e_\xi)\x_\xi \x_{\varrho_1} \dots \x_{\varrho_k} \\
&=- \sum_{\xi\in\Sigma'_1}l(\e_\xi) \x_\xi \x_{\varrho_1} \dots \x_{\varrho_k}\in \I_2.
\end{align*}

From the discussion above, we obtain a surjective $\corps$-algebra morphism
\[\phi\colon A^\bul(\Sigma) \to A^\bul(\Sigma').\]
This map is alos injective. Indeed, a reasoning similar to the above shows that the injection \[\corps[\x_\varrho]_{\varrho\in \Sigma'_1} \hookrightarrow \corps[\x_\varrho]_{\varrho\in \Sigma_1}\]
passes to the quotient and provides a morphism
\[A^\bul(\Sigma') \to A^\bul(\Sigma)\]
which is clearly inverse to $\phi$.
\end{proof}

\begin{prop} \label{prop:voisinage} Notations as in the previous proposition, let $\rho_a$ be the unique ray in $\Sigma\setminus \Gamma(\Sigma')$ and let $N = \rquot{\Z^{E}}{\Z {\e_{E}}}$. Let $\varrho$ be a ray in $\Gamma(\Sigma') \setminus \Gamma(\Delta')$. There exists a hyperplane $H$ in $N_\corps$ which does not contain $\rho_a$ such that there exists a neighborhood $U$ of $\e_\varrho$ in $\Sigma$ which is entirely included in $H_\R$.
\end{prop}

\begin{proof} Let as before $N^a$ be the lattice $\rquot{\Z^{E - a}}{\Z \e_{_{E-a}}} = \rquot{N}{\Z \e_a}$ so that both the fans $\Sigma'$ and $\Delta'$ live in $N^a_\R$.

Let $T$ be the tangent space to $\Sigma$ at $\e_\varrho$. Viewing $\e_\varrho$ in the Bergman fan $\Sigma_{\Ma}$ of the matroid $\Ma$, there exists a flag of flats $\Fl$ of $\Ma$ consisting of proper flats $\emptyset \neq F_1 \subsetneq F_2 \subsetneq \dots \subsetneq F_i \neq E$ such that $\e_\varrho$ lies in the relative interior of $\sigma_\Fl$ in $\Sigma_{\Ma}$. Note that since $\varrho\in\Gamma(\Sigma')$, the flats $F_j$ are in fact elements of $\Cl(\mu\del a)$; in particular, they do not contain $a$. The linear subspace $T \subseteq N_\R$ is then generated by the vectors $\e_F$ for $F$ a flat of $\Ma$ which is compatible with $\Fl$, in the sense that $\Fl \cup \{F\}$ forms again a flag of flats of $\Ma$ ($F$ may belong to $\Fl$).

Since $\varrho$ is not a ray of $\Gamma(\Delta')$, the cone $\sigma_\Fl$ is not in $\supp{\Gamma(\Delta')} =\supp{\Gamma(\Sigma_{\Ma\contr a})}$. This means there is an element $F_j \in \Fl$ which is not a flat of $\Ma \contr a$. By characterization of flats of $\Ma \del a$ and $\Ma \contr a$ given in Section~\ref{sec:modification}, we infer that $F_j \in \Cl(\Ma)$ but $F_j+a \notin \Cl(\Ma)$. In particular, the closure $\cl(F_j+a)$ in $\Ma$ of $F_j+a$ contains an element $b \in E \setminus (F_j+a)$.

\medskip

Now, for any flat $F$ in $\Ma$ which is compatible with $\Fl$, we have
\begin{itemize}
\item Either, $F$ contains both $a$ and $b$.
\item Or, $F$ contains none of $a$ and $b$.
\end{itemize}
Indeed, if $F \subseteq F_j$, then $F$ contains none of $a,b$. Otherwise, $F_j \subseteq F$, and in this case, if $F$ contains one of $a$ or $b$, then it must contain the other by the property that $\cl(F_j+a) =\cl(F_j+b)$.

\medskip

We infer that for any such $F$, we have $\e_F \in H_\Z:=\ker(\e_a^\dual -\e_b^\dual)$, where $\e_a^\dual -\e_b^\dual \in M = N^\vee$. This implies $T \subseteq H$ for $H = H_\Z \otimes_\Z \corps$. Obviously, $\e_a \notin H$ and the result follows.
\end{proof}

\subsection{Keel's lemma} \label{sec:Keel} Let $\Sigma$ be a unimodular fan and let $\sigma$ be a cone in $\Sigma$. Let $\Sigma'$ be the star subdivision of $\Sigma$ obtained by star subdividing $\sigma$. Denote by $\rho$ the new ray in $\Sigma'$. In this paper we only consider barycentric star subdivisions of rational fans. This means that $\rho=\R_+(\e_1+\dots+\e_k)$ where $\e_1, \dots, \e_k$ are the primitive vectors of the rays of $\sigma$.

\begin{thm}[Keel's lemma] \label{thm:keel}
Let $J$ be the kernel of the surjective map $\i^*_{\conezero\subface\sigma}\colon A^\bul(\Sigma)\to A^\bul(\Sigma^\sigma)$ and let
\[ P(T):=\prod_{\varrho\subface \sigma \\ \dims\varrho=1}(x_\varrho+T). \]
We have
\[ A^\bul(\Sigma')\simeq \rquot{A^\bul(\Sigma)[T]}{(JT+P(T))}. \]
The isomorphism is given by the map
\[\chi\colon \rquot{A^\bul(\Sigma)[T]}{(JT+P(T))} \simto A^\bul(\Sigma')\]
which sends $T$ to $-x_\rho$ and which verifies
\[ \forall \varrho \in \Sigma_1, \qquad
\chi(x_\varrho) = \begin{cases}
x_\varrho+x_\rho & \textrm{ if $\varrho \subface \sigma$}\\
x_\varrho & \textrm{otherwise}.
\end{cases}
\]
In particular this gives a vector space decomposition of $A^\bul(\Sigma')$ as
\begin{align}\label{eq:keel}
A^\bul(\Sigma')\simeq A^\bul(\Sigma)\oplus A^{\bul-1}(\Sigma^\sigma)T \oplus \dots \oplus A^{\bul-\dims{\sigma}+1}(\Sigma^\sigma)T^{\dims{\sigma}-1}.
\end{align}
\end{thm}

\begin{proof} This follows from~\cite{Keel}*{Theorem 1 in the appendix} for the map of toric varieties $\P_{\Sigma'} \to \P_{\Sigma}$. Here $P(T)$ is the polynomial in $A^\bul(\P_\Sigma)$ whose restriction in $A^\bul(\Sigma^\sigma)$ is the Chern polynomial of the normal bundle for the inclusion of toric varieties $\P_{\Sigma^\sigma} \hookrightarrow \P_{\Sigma}$.
\end{proof}

\subsection{Ascent and Descent} Situation as in Section~\ref{sec:Keel}, let $\Sigma$ be a unimodular fan and $\sigma$ a cone in $\Sigma$. Let $\Sigma'$ be the star subdivision of $\Sigma$ obtained by star subdividing $\sigma$. Denote by $\rho$ the new ray in $\Sigma'$. Let $\ell$ be an element of $A^1(\Sigma)$, and denote by $\ell^\sigma$ the restriction of $\ell$ to $\Sigma^\sigma$, i.e., the image of $\ell$ under the restriction map $A^\bul(\Sigma) \to A^\bul(\Sigma^\sigma)$.

\begin{thm} \label{thm:barycentric_subdivision} We have the following properties.
\begin{itemize}
\item \emph{(Ascent)} Assume the property $\HR(\Sigma^\sigma,\ell^\sigma)$ holds. Then $\HR(\Sigma, \ell)$ implies $\HR(\Sigma',\ell'-\epsilon x_\rho)$ for any small enough $\epsilon>0$, where
\[ \ell' = \sum_{\varrho\in\Sigma}\ell(\e_\varrho)x_\varrho + \bigl(\sum_{\varrho \subface \sigma \\ \dims\varrho=1}\ell(\e_\varrho)\bigr)x_\rho. \]

\item \emph{(Descent)} We have the following partial inverse: if both the properties $\HR(\Sigma^\sigma,\ell^\sigma)$ and $\HL(\Sigma,\ell)$ hold, and if we have the property $\HR(\Sigma',\ell'+\epsilon T)$ for any small enough $\epsilon>0$, then we have $\HR(\Sigma, \ell)$.
\end{itemize}
\end{thm}

Notice that the definition of $\ell'$ does not depend on the chosen representative of the class of $\ell$ as a cone-wise linear function on $\Sigma$.

\begin{remark} \label{rem:keel}
One can check in the proof of Theorem \ref{thm:barycentric_subdivision} below that we only use Keel's lemma, Poincaré duality and Proposition \ref{prop:HR}. For instance, we do not use that we are dealing with fans. Therefore, Theorem \ref{thm:barycentric_subdivision} might be useful in other contexts where an analog of Keel's lemma holds. In particular, this remark will be useful in Section \ref{sec:proofmaintheorem} where we need a similar ascent-descent theorem in the global case.
\end{remark}

In preparation for the proof, we introduce some preliminaries. Let $d$ be the dimension of $\Sigma$. By Keel's lemma, we have $A^d(\Sigma) \simeq A^d(\Sigma')$, and under this isomorphism, the canonical element $\omega_\Sigma$ corresponds to the canonical element $\omega_{\Sigma'}.$

Let $\alpha$ be a lifting of $\omega_{\Sigma^\sigma}$ in $A^{d-\dims \sigma}(\Sigma)$ for the restriction map $A^\bul(\Sigma) \to A^\bul(\Sigma^{\sigma})$. Such a lifting can be obtained for example by fixing a maximum dimensional cone $\tau$ in $\Sigma^{\sigma}$; in this case we have
\[\omega_{\Sigma^\sigma} = \prod_{\varrho \subface \tau \\ \dims\varrho=1} x_\varrho\]
and a lifting $\alpha$ of $\omega_{\Sigma^\sigma}$ is given by taking the product $\prod_{\varrho \subface \~\tau}x_\varrho$ in $A^{d-\dims\sigma}(\Sigma)$, where $\~\tau$ is the face of $\Sigma$ corresponding to $\tau$.

Using this lifting, the canonical element $\omega_{\Sigma'}$ can be identified with $-T^{\dims \sigma} \alpha$ in $A^\bul(\Sigma') \simeq \rquot{A^\bul(\Sigma)[T]}{(JT+P(T))}$. Indeed, since $P(T) = 0$, we get
\[-T^{\dims \sigma} = \sum_{j=0}^{\dims \sigma -1} S_{j} T^{\dims \sigma -j},\]
with $S_j$ the $j$-th symmetric function in the variables $\x_\varrho$ for $\varrho$ a ray in $\sigma$, seen as an element of $A^\bul(\Sigma')$. Therefore we get
\begin{align*}
-T^{\dims \sigma} \alpha &= \sum_{j=1}^{\dims \sigma} S_{j} \alpha T^{\dims \sigma -j} \\
&= \prod_{\varrho \subface \sigma\\ \dims\varrho =1} x_\varrho \cdot \alpha + \sum_{j=1}^{\dims \sigma -1} S_{j} \alpha T^{\dims \sigma -j}.
\end{align*}
Since $\omega_{\Sigma^{\sigma}}=\i^*_\sigma(\alpha)$ lives in the top-degree part of $A^\bul(\Sigma^\sigma)$, the products $S_j \omega_{\Sigma^\sigma}$ belongs to $J$. Thus the terms of the sum are all vanishing for $j=1,\dots, \dims \sigma -1$. Therefore, we get
\[-T^{\dims \sigma} \alpha = \prod_{\varrho \subface \sigma\\ \dims\varrho =1} x_\varrho\cdot \alpha = \prod_{\varrho \subface \sigma\\ \dims\varrho =1} x_\varrho \cdot \prod_{\varrho \subface \~\tau\\ \dims\varrho =1} x_\varrho = \omega_{\Sigma} = \omega_{\Sigma'}.\]

We infer that the following commutative diagram holds.

\[ \begin{tikzcd}
A^{d-\dims\sigma}(\Sigma^\sigma) \arrow[r, "\cdot \left(-T^{\dims{\sigma}}\right)", "\sim"'] \arrow[rd, "\sim"{sloped, above}, "\deg"'] & A^d(\Sigma') \arrow[r, "\sim"'] \arrow[d, "\deg", "\sim"{sloped, below}] & A^d(\Sigma) \arrow[ld, "\deg", "\sim"{sloped, above}] \\
& \corps
\end{tikzcd} \]

We are now ready to give the proof of the main theorem of this section.

\begin{proof}[Proof of Theorem~\ref{thm:barycentric_subdivision}]

(Ascent)\quad The idea for the proof of the ascent property is well-known, e.g., it is used as a way to derive Grothendieck's standard conjecture of Lefschetz type for a blow-up from the result on the base. It is also used in~\cite{AHK}. We give the proof here. Let $\epsilon>0$ and define the linear automorphism
\[S_\epsilon\colon A^\bul(\Sigma')\to A^\bul(\Sigma')\]
of degree $0$ which is the identity map $\id$ on $A^\bul(\Sigma)$ and multiplication by $\epsilon^{-\dims{\sigma}/2+k}\id$ on each $A^\bul(\Sigma^\sigma)T^k$, in the direct sum decomposition in Theorem~\ref{thm:keel}.

We admit for the moment that $Q_{\ell'+\epsilon T}\circ S_\epsilon$ admits a limit $Q$ such that the pair $(A^\bul(\Sigma'),Q)$, consisting of a graded vector space and a bilinear form on this vector space, is naturally isomorphic to the pair
\[ \Xi:=\Bigl(A^\bul(\Sigma)\oplus \big(T\cdot(\rquot{\R[T]}{T^{\dims\sigma-1}})\otimes A^\bul(\Sigma^\sigma)\big), \ell \oplus(T\otimes\id+\id\otimes \ell^\sigma)\Bigr). \]
It follows from the statement~\eqref{prop:HR:otimes_lefschetz} in Proposition~\ref{prop:HR}, the hypothesis of the theorem, and $\HR(\rquot{\R[T]}{T^{\dims\sigma-1}}, T)$, that the pair $\Xi$ verifies $\HR(\Xi)$.
Using Keel's lemma, it follows that we have $\HR(A^\bul(\Sigma'),Q)$. By Proposition~\ref{prop:HRbis}, the property $\HR(A^\bul(\Sigma'),Q_{\ell'+\epsilon T}\circ S_\epsilon)$ holds for any small enough value of $\epsilon>0$. Finally, applying Proposition~\ref{prop:HR}, we get $\HR(A^\bul(\Sigma'),Q_{\ell'+\epsilon T})$, which finishes the proof of the ascent part of our theorem.

\medskip

Before turning to the proof of the existence of the limit, let us explain the proof of the descent part.

\bigskip

\noindent(Descent)

\smallskip
We again use Proposition~\ref{prop:HR}. By our assumption, we have the properties $\HL(A^\bul(\Sigma), \ell)$ and $\HR(A^\bul(\Sigma^\sigma,\ell^\sigma))$, from which we get the property
\[ \HR\Bigl((T\cdot(\rquot{\R[T]}{T^{\dims\sigma-1}})\otimes A^\bul(\Sigma^\sigma)\big), (T\otimes\id+\id\otimes \ell^\sigma)\Bigr). \]

Defining $\Xi$ as in the first part of the proof, we deduce $\HL(\Xi)$ by applying Proposition~\ref{prop:HR}. By the statement~\eqref{prop:HR:oplus} in that proposition, it will be enough to prove
$\HR(\Xi)$ in order to get $\HR(A^\bul(\Sigma), \ell)$.

Note that by the hypothesis in the theorem, we have $\HR(A^\bul(\Sigma'),Q_{\ell'+\epsilon T})$ for small enough values of $\epsilon>0$. Using Proposition~\ref{prop:HR}, we deduce that $\HR(A^\bul(\Sigma'),Q_{\ell'+\epsilon T}\circ S_\epsilon)$ holds for $\epsilon>0$ small enough. It follows that we have as well $\HL(A^\bul(\Sigma'),Q_{\ell'+\epsilon T}\circ S_\epsilon)$.

Applying now Proposition~\ref{prop:HRbis}, since $\HR$ is a closed condition restricted to the space of bilinear forms which verify $\HL$, we deduce that $\HR(\Xi)$ holds: indeed, $Q$
is the limit of $Q_{\ell'+\epsilon T}\circ S_\epsilon$ and $\HR(A^\bul(\Sigma'),Q_{\ell'+\epsilon T}\circ S_\epsilon)$ and $\HL(A^\bul(\Sigma'),Q_{\ell'+\epsilon T}\circ S_\epsilon)$ hold for $\epsilon >0 $ small enough.

\vspace{.7cm}

We are thus left to prove the last point, namely the existence of the limit $\lim_{\epsilon \to 0}\, Q_{\ell'+\epsilon T}\circ S_\epsilon$.

Consider the decomposition~\eqref{eq:keel} given in Keel's lemma. The decomposition for degree $k$ part $A^{k}(\Sigma')$ of the Chow ring has pieces $A^{k}(\Sigma)$ and $A^{k-i}(\Sigma^\sigma) T^{i}$ for $1 \leq i \leq \min\{k, \dims \sigma -1\}$. Consider the bilinear form $\Phi_{A^\bul(\Sigma')}$ given by the degree map of $A^\bul(\Sigma')$.

For this degree bilinear form, each piece $A^{k-i}(\Sigma^\sigma)T^i$ is orthogonal to the piece $A^{d-k}(\Sigma)$ as well as to the piece $A^{d-k-j}(\Sigma^\sigma)T^j$ for $j < \dims \sigma -i$, in the decomposition of $A^{d-k}(\Sigma')$ given in Keel's lemma.

We now work out the form of the matrix of the bilinear forms $Q_{\ell'+\epsilon T}\circ S_\epsilon$ in degree $k$ with respect to the decomposition given in \eqref{eq:keel}.

First, note that for $a,b\in A^k(\Sigma)$, we have
\[ Q_{\ell'+\epsilon T}(a,b)=\deg(a\ell^{\prime d-2k}b)+\O(\epsilon^{\dims \sigma}). \]

Second, for an element $a\in A^k(\Sigma)$ and an element $b\in A^{k-i}(\Sigma^\sigma)T^i$, we get the existence of an element $c\in A^{d-2k-\dims\sigma+i}(\Sigma')$ such that
\[ Q_{\ell'+\epsilon T}(a,b)=\deg(ac(\epsilon T)^{\dims\sigma-i}\epsilon^{-\dims\sigma/2+i}b)=\O(\epsilon^{\dims \sigma/2}). \]

Third, for an element $a\in A^{k-i}(\Sigma^\sigma)T^i$ and an element $b\in A^{k-j}(\Sigma^\sigma)T^j$, with $i+j\leq\dims\sigma$, we get
\begin{align*}
Q_{\ell'+\epsilon T}(a,b) &= \deg\Big(\epsilon^{-\dims\sigma/2+i}a\big(C(\epsilon T)^{\dims\sigma-i-j}\ell'^{d-2k-\dims\sigma+i+j}+\O(\epsilon^{\dims\sigma-i-j+1})\big)\epsilon^{-\dims\sigma/2+j}b\Big) \\
  &=C\deg_{A^\bul(\Sigma^\sigma)}(a(\ell^\sigma)^{d-\dims\sigma-2k+i+j}b)+\O(\epsilon),
\end{align*}
with $C=\binom{d-2k}{\dims\sigma-i-j}$.

Finally, if $a\in A^{k-i}(\Sigma^\sigma)T^i$ and $b\in A^{k-j}(\Sigma^\sigma)T^j$, with $i+j>\dims\sigma$, we get
\[ Q_{\ell'+\epsilon T}(a,b)=\deg(\epsilon^{-\dims\sigma/2+i}a\O(1)\epsilon^{-\dims\sigma/2+j}b)=\O(\epsilon^{i+j-\dims\sigma}). \]

\medskip

Doing the same calculation for the bilinear form $Q_{\~L}$ associated to the linear map of degree one $\~L\colon A^\bul(\Sigma')\to A^{\bul+1}(\Sigma')$ given by the matrix
\[ \begin{pmatrix}
\ell \\
0 &\ell^\sigma & & \makebox(0,0){\huge0} \\
 & T & \ell^\sigma \\[1ex]
 & \makebox(0,0){\huge0} & \ddots & \ddots \\[1ex]
 & & & T & \ell^\sigma
\end{pmatrix} \]
relative to the decomposition~\eqref{eq:keel} given in Keel's lemma, we see that
\[\lim_{\epsilon \to 0}\,Q_{\ell'+\epsilon T} = Q_{\~L}. \]
To conclude, we observe that the pair $(A^\bul(\Sigma'), Q_{\~L})$ is isomorphic to the previously introduced pair,
\[ \Xi=\Bigl(A^\bul(\Sigma)\oplus \big(T\cdot(\rquot{\R[T]}{T^{\dims\sigma-1}})\otimes A^\bul(\Sigma^\sigma)\big), \ell\oplus(T\otimes\id+\id\otimes \ell^\sigma)\Bigr), \]
as claimed.
\end{proof}

\subsection{Weak factorization theorem} Two unimodular fans with the same support are called \emph{elementary equivalent} if one can be obtained from the other by a barycentric star subdivision. The \emph{weak equivalence} between unimodular fans with the same support is then defined as the transitive closure of the elementary equivalence relation.

\begin{thm}[Weak factorization theorem] \label{thm:equivalent_fan}
Two unimodular fans with the same support are always weakly equivalent.
\end{thm}
\begin{proof}
This is Theorem A of~\cite{Wlo97} proved independently by Morelli~\cite{Mor96} and expanded by Abramovich-Matsuki-Rashid, see~\cite{AMR}.
\end{proof}

We will be interested in the class of quasi-projective unimodular fans on the same support. In this class, we define the \emph{weak equivalence} between quasi-projective unimodular fans as the transitive closure of the elementary equivalence relation between quasi-projective unimodular fans with the same support. This means all the intermediate fans in the star subdivision and star assembly are required to remain quasi-projective.

\begin{thm}[Weak factorization theorem for quasi-projective fans] \label{thm:equivalent_fan2}
Two quasi-projective unimodular fans with the same support are weakly equivalent within the class of quasi-projective fans.
\end{thm}
\begin{proof}
This can be obtained from relevant parts of~\cites{Wlo97, Mor96, AKMW} as  discussed and generalized by Abramovich and Temkin in~\cite{AT}*{Section 3}.
\end{proof}

We thank Kenji Matsuki and Kalle Karu for helpful correspondence and clarification on this subject.

\subsection{Ample classes}

Let $\Sigma\subseteq N_\R$ be a unimodular Bergman fan. Let $f$ be a strictly convex cone-wise linear function on $\Sigma$. Recall that if $\varrho$ is a ray in $\Sigma_1$, we denote by $\e_\varrho$ the primitive element of $\varrho$ and by $x_\varrho$ the corresponding element of $A^1(\Sigma)$. We denote by $\ell(f)$ the element of $A^1(\Sigma)$ defined by
\[ \ell(f):=\sum_{\varrho\in\Sigma^1}f(\e_\varrho)x_\varrho. \]

\begin{prop} \label{prop:ell_f_independent}
Notations as above, let $\sigma$ be a cone of $\Sigma$ and let $\phi$ be a linear form on $N_\R$ which coincides with $f$ on $\sigma$. The function $f-\phi$ induces a cone-wise linear function on $\Sigma^\sigma$ which we denote by $f^\sigma$. Then
\[ \ell(f^\sigma)=\i^*_{\conezero\subface\sigma}(\ell(f)). \]
In particular, $\ell(f^\sigma)$ does not depend on the chosen $\phi$.

Moreover, if $f$ is strictly convex on $\Sigma$, then $f^\sigma$ is strictly convex on $\Sigma^\sigma$.
\end{prop}

\begin{proof}
We write $\varrho\sim\sigma$ if $\varrho$ is a ray in the link of $\sigma$, \ie, $\varrho\not\in\sigma$ and $\sigma+\varrho$ is a face of $\Sigma$. Such rays are in one-to-one correspondence with rays of $\Sigma^\sigma$. Recall the following facts:
\[ \ell(\phi)=0, \qquad \i^*_{\conezero\subface\sigma}(x_\varrho)=0 \text{ if $\varrho\not\sim\sigma$ and $\varrho\not\in\sigma$}, \qquad f(\e_\varrho)=\phi(\e_\varrho) \text{ if $\varrho\subface\sigma$}. \]
Then we have
\begin{align*}
\i^*_{\conezero\subface\sigma}(\ell(f))
  &= \i^*_{\conezero\subface\sigma}(\ell(f-\phi)) \\
  &= \sum_{\varrho\in\Sigma_1}(f-\phi)(\e_\varrho)\i^*_{\conezero\subface\sigma}(x_\varrho) = \sum_{\varrho\in\Sigma_1 \\ \varrho\sim\sigma}(f-\phi)(\e_\varrho)\i^*_{\conezero\subface\sigma}(x_\varrho) \\
  &= \sum_{\varrho\in\Sigma^\sigma\\ \dims\varrho=1}f^\sigma(\e_{\varrho})x_{\varrho} = \ell(f^\sigma). \qedhere
\end{align*}

For the last point, let $f$ be a strictly convex cone-wise linear function on $\Sigma$. We have to prove that $f^\sigma$ is strictly convex. Let $\zeta$ be a cone of $\Sigma^\sigma$. The cone $\zeta$ corresponds to a unique cone $\~\zeta$ in $\Sigma$ that contains $\sigma$. Since $f$ is strictly convex, there exists a linear form $\psi$ on $N_\R$ such that $f-\psi$ is strictly positive around $\~\zeta$ (in the sense of Section \ref{subsec:convexite}). Then $\psi-\phi$ is zero on $\sigma$, thus induces a linear form $\psi^\sigma$ on $N^\sigma_\R$. We get that $f^\sigma-\psi^\sigma$ is strictly positive around $\zeta$, thus $f^\sigma$ is strictly convex around $\zeta$. We infer that $f^\sigma$ is strictly convex on $\Sigma^\sigma$.
\end{proof}

\begin{defi}
Let $\Sigma$ be a Bergman fan. An element $\ell\in A^1(\Sigma)$ is called an \emph{ample class} if there exists a strictly convex cone-wise linear function $f$ on $\Sigma$ such that $\ell=\ell(f)$.
\end{defi}

\subsection{Main Theorem}

Here is the main theorem of this section.

\begin{thm}\label{thm:mainlocal} Any quasi-projective unimodular Bergman fan $\Sigma$ verifies $\HR(\Sigma, \ell)$ for any ample class $\ell \in A^1(\Sigma)$.
\end{thm}

The following proposition is the base of our proof. It allows to start with a fan on any Bergman support which verifies the Hodge-Riemann relations, and then use the ascent and descent properties to propagate the property to any other fan on the same support.

\begin{prop}\label{prop:baseHR}
For any Bergman support, there exists a quasi-projective unimodular fan $\Sigma$ on this support and a strictly convex element $\ell \in A^1(\Sigma)$ which verifies $\HR(\Sigma, \ell)$.
\end{prop}

We prove both statements, in the theorem and proposition, by induction on the dimension of the ambient space, i.e., on the size of the ground set of the matroid underlying the support of the Bergman fan.
So suppose by induction that the theorem and the proposition hold for all Bergman support of $\R^k$ with $k<n$, and for any quasi-projective unimodular fan with a Bergman support lying in $\R^k$ for $k<n.$

\begin{proof}[Proof of Proposition~\ref{prop:baseHR}]
Let $\Ma$ be a simple matroid on the ground set $\{0,\dots,n\}$.

If $\Ma$ is the free matroid, then we have $\supp{\Sigma_\Ma}=\R^n$. In this case, we can take any projective unimodular fan. For example, the fan of the projective space $\P^n$, equivalently of $\TP^{n}$, $\Sigma_{\TP^n}$ has Chow ring isomorphic to $\rquot{\R[\x]}{\x^{n+1}}$. In this case, we can take $\ell = \x$, and we have $\HR(\rquot{\R[\x]}{\x^{n+1}},\x)$.

\medskip

Otherwise, $\Ma$ is not the free matroid. Let $i$ be a proper element of $\Ma$. By the induction hypothesis, there exists a quasi-projective unimodular fan $\Sigma'$ of support $\supp{\Sigma_{\Ma\del{i}}}$. Doing some blow-ups if necessary, we can assume that $\Sigma'$ contains a subdivision $\Delta'$ of $\Sigma_{\Ma\contr{i}}$ (cf. Section \ref{sec:triangulation} for more details on this argument). Let $\ell'$ be a strictly convex function on $\Sigma'$. Modifying $\ell'$ by adding a linear function if necessary, we can suppose without loss of generality that $\ell'$ is strictly positive on $\supp{\Sigma'}\setminus\{0\}$. Let $\Sigma = \modtrop{\Delta'}{\Sigma'}$, the unimodular fan with support $\supp{\Sigma_\Ma}$ obtained as the tropical modification of $\Delta'$ in $\Sigma'$. Using Proposition~\ref{prop:modification}, we deduce $\HR(\Sigma, \ell)$ from $\HR(\Sigma_{\Ma\del i},\ell')$ where $\ell$ is the cone-wise linear function on $\Sigma$ which takes value $0$ on the new ray $\varrho_i$ of $\Sigma$ and which coincides on other rays with $\Gamma_*(\ell')$ where $\Gamma\colon \Sigma'\hookrightarrow\Sigma$ is the natural injective map. One easily verifies that $\ell$ is strictly convex. So the proposition follows.
\end{proof}

\begin{prop} \label{prop:HR-oneall}
For any unimodular fan $\Sigma$ of dimension $d$, we have the equivalence of the following statements.
\begin{itemize}
\item $\HR(\Sigma, \ell)$ is true for any strictly convex element $\ell \in A^1(\Sigma)$.
\item $\HR(\Sigma, \ell)$ is true for a strictly convex element $\ell \in A^1(\Sigma)$.
\end{itemize}
\end{prop}
\begin{proof} Let $\ell$ be a strictly convex piecewise linear function on $\Sigma$. For all $\sigma\in \Sigma$, we get a strictly convex function $\ell^\sigma$ on $\Sigma^\sigma$. By the hypothesis of our induction, we know that $\HR(\Sigma^\sigma, \ell^\sigma)$ holds. By Proposition~\ref{prop:local_HR} we thus get $\HL(\Sigma, \ell)$.

By Proposition \ref{prop:HRbis} we know that $\HR(\Sigma,\ell)$ is an open and closed condition on the set of all $\ell$ which satisfy $\HL(\Sigma,\ell)$. In particular, if there exists $\ell_0$ in the cone of strictly convex elements which verifies $\HR(\Sigma,\ell_0)$, any $\ell$ in this cone should verify $\HR(\Sigma,\ell)$.
\end{proof}

Let $\Sigma$ be a unimodular fan of dimension $d$, $\ell$ a strictly convex element on $\Sigma$, and let $\Sigma'$ be the fan obtained from $\Sigma$ by star subdividing a cone $\sigma\in\Sigma$. Denote by $\rho$ the new ray in $\Sigma'$, and let
\[\ell':=\ell+\bigl(\sum_{\varrho \subface \sigma \\ \dims\varrho=1}\ell(\e_\varrho)\bigr)x_\rho.\]

The following is straightforward.

\begin{prop} For any small enough $\epsilon>0$, the element $\ell' - \epsilon x_\rho$ of $A^1(\Sigma')$ is strictly convex.
\end{prop}

\begin{prop} \label{prop:HR-trans}
Notations as above, the following statements are equivalent.
\begin{enumerate}[label=\defaultRoman]
\item \label{hr:trans1} We have $\HR(\Sigma, \ell)$.
\item \label{hr:trans2} The property $\HR(\Sigma',\ell'-\epsilon x_\rho)$ holds for any small enough $\epsilon>0$.
\end{enumerate}
\end{prop}

\begin{proof}
It will be enough to apply Theorem~\ref{thm:barycentric_subdivision}. A direct application of the theorem, and the fact that $T$ is identified with $-x_\rho$ via Keel's lemma, leads to the implication $\ref{hr:trans1} \Rightarrow \ref{hr:trans2}$.

We now explain the implication $\ref{hr:trans2} \Rightarrow \ref{hr:trans1}$. By the hypothesis of our induction, we have $\HR(\Sigma^\sigma, \ell^\sigma)$. Moreover, Proposition \ref{prop:local_HR} and our induction hypothesis imply that we have $\HL(\Sigma, \ell)$.
Applying now the descent part of Theorem~\ref{thm:barycentric_subdivision} gives the result.
\end{proof}

\begin{proof}[Proof of Theorem~\ref{thm:mainlocal}] Let $\Ma$ be a matroid on the ground set $\{0,\dots, n\}$. By Proposition~\ref{prop:baseHR} there exists a quasi-projective unimodular fan $\Sigma_0$ with support $\supp{\Sigma_\Ma}$ and a strictly convex element $\ell_0 \in A^1(\Sigma_0)$ such that $\HR(\Sigma_0, \ell_0)$ holds. By Propositions~\ref{prop:HR-trans} and~\ref{prop:HR-oneall}, for any quasi-projective fan $\Sigma$ in the quasi-projective weak equivalence class of $\Sigma_0$, so with support $\supp{\Sigma_\Ma}$, and any strictly convex element $\ell \in A^1(\Sigma)$, we get $\HR(\Sigma, \ell)$. The theorem now follows in dimension $n$ by applying Theorem~\ref{thm:equivalent_fan2}
\end{proof}

\subsection{Examples} We finish this section with some concrete examples of fans verifying the Hodge-Riemann relations with respect to ample and non-ample classes.\label{ex:convex_U33}

\medskip

\subsubsection{} 
\begin{figure}[ht]
\caption{The Bergman fan of $U_{3,3}$} \label{fig:Bergman_fan_U33}
\begin{tikzpicture}
\clip (0, 0) circle (2);
\fill[color=gray!10] (-2,-2) rectangle (2,2);
\draw[color=gray!25, scale=.6] (-4,-4) grid (4,4);
\draw  (0 , 0)
  edge["$\rho_1$", near end]     (2 , 0)
  edge["$\rho_{12}$"{above left=-4pt}]  (2 , 2)
  edge["$\rho_{2}$", near end]   (0 , 2)
  edge["$\rho_{02}$", near end]  (-2, 0)
  edge["$\rho_{0}$"{below right=-4pt}]   (-2,-2)
  edge["$\rho_{01}$", near end]  (0 ,-2);
\end{tikzpicture}
\end{figure}
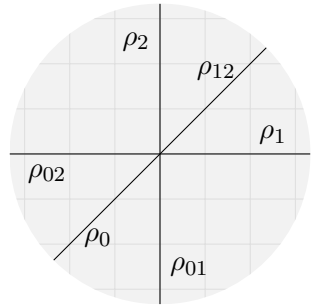

Let $\Sigma$ be the complete Bergman fan associated to the uniform matroid $U_{3,3}$ on the ground set $\{0,1,2\}$ as illustrated in Figure \ref{fig:Bergman_fan_U33}.
The flats of $U_{3,3}$ are $\Cl(U_{3,3}):=\{1,12,2,02,0,01\}$ (where $ij$ denote the set $\{i,j\}$). The rays of $\Sigma$ correspond to the flats of $U_{3,3}$ and are denoted $\rho_F$ with $F\in\Cl(U_{3,3})$. These rays corresponds to elements of $A^1(\Sigma)$ denoted $x_F$, for $F\in\Cl(U_{3,3})$. We have
\[ \dim(A^1(\Sigma))=\card{\Cl(U_{3,3})}-\dim\bigl(\R^{2^\dual}\bigr)=4. \]
We have a natural degree map $\deg\colon A^2(\Sigma)\to\R$ which verifies the following.
\[ \deg(x_Fx_G)=\begin{cases}
  -1 & \text{if $F=G$,} \\
  1  & \text{if $F\neq G$ and $F$ and $G$ are comparable,} \\
  0  & \text{if $F$ and $G$ are not comparable.}
\end{cases} \]

\medskip

Let $\ell:=x_{1}+x_{12}+x_{2}+x_{02}+x_{0}+x_{01}\in A^1(\Sigma)$. This element corresponds to a strictly convex cone-wise linear function on $\Sigma$, thus, it must verify the Hodge-Riemann relations $\HR(\Sigma,\ell)$. Let us check that this is the case. We have to check that $\deg(\ell^2)>0$ and that
\[ \begin{array}{rccc}
Q_1\colon & A^1(\Sigma)  \times  A^1(\Sigma) & \to     & \R,       \\
          &( x            ,       y )        & \mapsto & \deg(xy),
\end{array} \]
has signature $-2=2\dim(A^0)-\dim(A^1)$. We have
\[ \deg(\ell x_1)=\deg(x_1^2+x_1x_{12}+x_1x_{01})=1, \]
and, by symmetry, $\deg(\ell x_F)=1$ for every $F\in\Cl(U_{3,3})$. Thus, $\deg(\ell^2)=6$.

For $Q_1$ we have an orthogonal basis $(x_0, x_1, x_2, x_0+x_{01}+x_1)$. Moreover $\deg(x_0^2)=\deg(x_1^2)=\deg(x_2^2)=-1$ and
\[ \deg((x_0+x_{01}+x_1)^2)=\deg(x_0^2+x_{01}^2+x_1^2+2x_0x_{01}+2x_{01}x_1)=1. \]
Thus, the signature of $Q_1$ is $-2$ and the Hodge-Riemann relations are verified.

\medskip

\subsubsection{} Let us give another example in the same Chow ring treated in the previous paragraph. Let $\ell'=x_1+x_{12}+x_2$. The cone-wise linear function associated to $\ell'$ is not convex. It is even zero around $\rho_0$, which implies that $\i^*_{\conezero\ssubface\rho_0}(\ell')=0$. Nevertheless, $\ell'$ verifies the Hodge-Riemann relations $\HR(\Sigma, \ell')$. We only have to check that $\deg(\ell'^2)>0$, since $Q_1$ does not depend on $\ell'$ and we have already checked that $Q_1$ has signature $-2$. We have
\[ \deg((x_1+x_{12}+x_2)^2)=\deg(x_1^2+x_{12}^2+x_2^2+2x_1x_{12}+2x_{12}x_2)=1>0. \]
This gives us an example of a non-ample element verifying the Hodge-Riemann relations.

\medskip

\subsubsection{} As a last example, let $\Sigma$ be the Bergman fan of the uniform matroid $U_{3,4}$ on the ground set $\{0,1,2,3\}$, and let $\ell=(2+2\epsilon)x_0+x_2+x_3+\epsilon x_{01}+(3+\epsilon)(x_{02}+x_{03})-x_{23}\in A^1(\Sigma)$ where $\epsilon$ is a small positive number. Then $\ell$ is not ample, but we have $\HR(\Sigma, \ell)$. Indeed, we have something stronger: $\ell^\sigma$ is ample for any $\sigma\in\Sigma\setminus\{\conezero\}$ where $\ell^\sigma = \i^*_{\conezero \subface \sigma}(\ell)$.


\section{Quasi-projective unimodular triangulations of polyhedral complexes}\label{sec:triangulation}

In this section we show that any rational polyhedral complex in $\R^n$ admits a simplicial subdivision which is both \emph{quasi-projective} and \emph{unimodular} with respect to the lattice $\frac 1k \Z^n$, for some natural number $k$, in a sense which we will precise in a moment. Moreover, we will show this property holds for any \emph{tropical compactification} of $X$. The novelty here is in ensuring the \emph{convexity property} underlying the definition of quasi-projective triangulations which will be the crucial property used in what will follow later, as well as in ensuring the unimodularity of unbounded polyhedra which form the polyhedral subdivision of $X$. The theorem thus extends the pioneering result proved in~\cite{KKMS}, in Mumford's proof of the semistable reduction theorem, as well as previous variants proved in~\cite{IKMZ} and~\cite{Wlo97}.

\subsection{Statement of the triangulation theorem.}

This theorem corresponds to Theorem \ref{thm:regulartriangulations} stated in the introduction. However, we add a technical result about the existence of strictly convex functions to the statement which  will be needed later.

\begin{thm}[Triangulation theorem] \label{thm:triangulation_unimodulaire_convexe}
Let $X$ be a rational polyhedral complex in $\R^n$. There exists an integer $k$ and a triangulation $X'$ of $X$ which is quasi-projective and unimodular with respect to the lattice $\frac 1k \Z^n$.

Moreover we ensure that $X'_\infty$ is a fan, and that there exists a strictly convex function $f$ on $X'$ such that $f_\infty$ is well-defined and strictly convex.
\end{thm}
The terminology will be introduced in the next section.

\medskip

Furthermore, we will prove the following theorem which will be later used in Section \ref{sec:projective_bundle_theorem}.

\begin{thm}[Unimodular triangulation preserving the recession fan] \label{thm:unimodular_preserving_recession}
Let $X$ be a polyhedral complex in $\R^n$ such that its pseudo-recession fan $X_\infty$ is a unimodular fan. Then there exists a positive integer $k$ and a subdivision $Y$ of $X$ which is unimodular with respect to $\frac1k\Z^n$ and which verifies $Y_\infty=X_\infty$.
\end{thm}

\subsection{Preliminary definitions and constructions}

This section presents certain constructions which will be needed later. While definitions are given here for the case of polyhedral complexes and fans, we note that the majority of the constructions naturally extend to the case of polyhedral pseudo-complexes and pseudo-fans.

\subsubsection{Restriction, intersection, and hyperplane cut}

Let $X$ be a polyhedral complex of $\R^n$ and let $S$ be a subset of $\R^n$. The \emph{restriction of $X$ to $S$} denoted by $X\rest S$ is by definition the polyhedral complex consisting of all the faces of $X$ which are included in $S$. If $S$ is an affine hyperplane, an affine half-space, or more generally a polyhedron of $\R^n$, the \emph{intersection of $X$ with $S$} denoted by $X\cap S$ is the polyhedral complex whose faces are $\delta \cap S$ for any face $\delta$ of $X$, i.e.,
\[ X\cap S:=\bigl\{\delta\cap S\st \delta\in X\bigr\}. \]
For an affine hyperplane $H$ in $\R^n$, the \emph{hyperplane cut of $X$ by $H$} denoted by $X\cdot H$ is the subdivision of $X$ defined by
\[ X\cdot H:=(X\cap H)\cup(X\cap H_+)\cup(X\cap H_-), \]
where $H_+$ et $H_-$ are the two half-spaces of $\R^n$ defined by $H$.

\subsubsection{Slicing with respect to a pencil of hyperplanes} Let $\Sigma$ be a fan of $\R^n$. For any cone $\sigma \in \Sigma$, choose a collection of half-spaces $H^+_{\sigma,1}, \dots, H^+_{\sigma,k_\sigma}$ for $k_\sigma \in \Z_{>0}$ each defined by hyperplanes $H_{\sigma,1}, \dots, H_{\sigma,k_\sigma}$, respectively, such that $\sigma$ is the intersection of these half-spaces.
We call the collection of hyperplanes
\[\pen_\Sigma:= \Bigl\{H_{\sigma, j} \, \st \, \sigma \in \Sigma \quad \textrm{and} \quad j =1, \dots, k_\sigma\Bigr\}\]
\emph{a pencil of hyperplanes defined by $\Sigma$}.

\medskip

Let now $\Delta$ be a complete fan in $\R^n$ and let $\Sigma$ be a second (arbitrary) fan in $\R^n$. The \emph{slicing of $\Delta$ with a pencil of hyperplanes $\pen_\Sigma$} defined by $\Sigma$ is the fan denoted by $\Delta \cdot \pen_\Sigma$ and defined by
\[ \Delta\cdot \pen_{\Sigma} := \Delta \cdot H_{\sigma_1,1} \cdot H_{\sigma_1,2}\cdots H_{\sigma_1,k_{\sigma_1}}\cdots H_{\sigma_m,k_{\sigma_m}}, \]
\ie, defined by starting from $\Delta$ and successively taking hyperplane cuts by elements $H_{\sigma,j}$ of the pencil $ \pen_\Sigma$.

\begin{prop} \label{prop:decoupe_subdivision}
The slicing of a complete fan with a pencil of hyperplanes $\pen_\Sigma$ defined by a fan $\Sigma$ necessarily contains a subdivision of $\Sigma$.
\end{prop}

\begin{proof} First note that if $\Delta$ is a fan in $\R^n$ and $H$ is a hyperplane, all the cones $\delta \in \Delta \cdot H$ live in one of the two half-spaces $H^+$ and $H^-$ defined by $H$. This property holds as well for any subdivision of $\Delta\cdot H$.

Let now $\Delta$ be complete. Let $\pen_\Sigma$ be a pencil of hyperplanes defined by the fan $\Sigma$ consisting of hyperplanes $H_{\sigma, j}$ for $\sigma \in \Sigma$ and $j=1, \dots, k_\sigma$ as above. Let $\Delta' = \Delta \cdot \pen_\Sigma$. In order to prove the proposition, we will need to show the following two properties:
\begin{enumerate}
\item \label{enum:pencil:1} the support of the restriction $\Delta'\rest{\supp{\Sigma}}$ is the entire support $\supp{\Sigma}$; and
\item \label{enum:pencil:2} any face of $\Delta'\rest{\supp{\Sigma}}$ is included in a face of $\Sigma$.
\end{enumerate}

In order to prove \ref{enum:pencil:1}, let $x$ be a point in $\supp{\Sigma}$. It will be enough to show the existence of a face of $\Delta'$ which contains $x$ and which is entirely in $\supp{\Sigma}$.

Let $\delta$ be the smallest face of $\Delta'$ which contains $x$; $\delta$ is the face which contains $x$ in its interior. Let $\sigma \in \Sigma$ be a face which contains $x$. We will show that $\delta \subseteq \sigma$.

By what preceded, for each $i \in \{1, \dots, k_\sigma\}$, there exists $\epsilon_i \in \{+, -\} $ such that $\delta$ is included in $H^{\epsilon_i}_{\sigma,i}$. Recall that by definition, we have $\sigma=\cap_{i=1}^{k_\sigma}H^+_{\sigma,i}$. We will show that we can take $\epsilon_i = +$ for each $i$, which gives the result. Suppose $\epsilon_i=-$ for an $i\in \{1, \dots, k_\sigma\}$. In particular, $x$ belongs to both half-spaces $H^-_{\sigma,i}$ and $H^+_{\sigma,i}$, and so $x\in H_{\sigma,i}$. Since $\delta\cap H_{\sigma,i}$ is a non-empty face of $\delta$ which contains $x$, by minimality of $\delta$, we infer that $\delta\subset H_{\sigma,i}$ which shows that we can choose $\epsilon_i = +$. This proves \ref{enum:pencil:1}.

In order to prove the second assertion, let $\delta$ be a face of the restriction $\Delta'\rest{\supp{\Sigma}}$. Let $x$ be a point in the interior of $\delta$, so $\delta$ is the smallest face of $\Delta'$ which contains $x$. Let $\sigma$ be a face of $\Sigma$ which contains $x$. As we showed in the proof of \ref{enum:pencil:1}, we have the inclusion $\delta \subseteq \sigma$, which is the desired property in \ref{enum:pencil:2}.

It follows that $\Delta'\rest{\supp{\Sigma}}$ is a subdivision of $\Sigma$ and the proposition follows.
\end{proof}

\subsubsection{Blow-ups and existence of triangulations} Let $\sigma$ be a cone in $\R^n$ and let $x$ be a vector of $\R^n \setminus \{0\}$. The \emph{blow-up of $\sigma$ at $x$} is the fan denoted by $\sigma_{(x)}$ and defined by
\[ \sigma_{(x)}:=\begin{cases}
\displaystyle \bigcup_{\substack{\tau\subface\sigma \\ x\not\in\tau}} \bigl\{\,\tau,\tau+\R_+x\, \bigr\}& \text{if $x \in\sigma$,} \\[2em]
\face{\sigma} & \text{otherwise.}
\end{cases} \]
Note that $\sigma_{(x)}$ is a polyhedral subdivision of $\face{\sigma}$. Moreover, the definition depends only on the half line $\R_+x$.

More generally, if $\Sigma$ is a fan, we define the \emph{blow-up of $\Sigma$ at $x$} to be the fan denoted by $\Sigma_{(x)}$ and defined as the union of blow-ups $\sigma_{(x)}$ for $\sigma \in \Sigma$, \ie,
\[ \Sigma_{(x)}=\bigcup_{\sigma\in\Sigma}\sigma_{(x)}. \]
Then $\Sigma_{(x)}$ is a subdivision of $\Sigma$.

\medskip

For a pair of points $x_1, x_2$, we denote by $\Sigma_{(x_1)(x_2)}$ the fan obtained by first blowing-up $\Sigma$ at $x_1$, and then blowing-up the resulting fan at $x_2$. The definition extends to any ordered sequence of points $x_1, \dots, x_k$ in $\R^n$.

\begin{prop} \label{prop:triangulation}
Let $\Sigma$ be a fan. Let $\sigma_1, \dots, \sigma_k$ be the minimal faces of $\Sigma$ which are not simplicial. For each such face $\sigma_i$, $i\in\zint1k$, pick a vector $x_i$ in its interior. Then
\[ \Sigma_{(x_1)(x_2)\cdots(x_k)} \]
is a triangulation of $\Sigma$.
\end{prop}

\begin{proof}
Faces of $\Sigma_{(x_1)}$ which were not in $\Sigma$ are of the form $\tau+\R_+x_1$ where $\tau\in\Sigma$ and $x_1\not\in\TT\tau$. Let $\tau'=\tau+\R_+x_1$ be such a new face. Clearly $\tau'$ is simplicial if and only if $\tau$ is simplicial. By contrapositive, if $\tau'$ is not simplicial, then $\tau\subface\tau'$ is not simplicial. In any case, $\tau'$ is not a non-simplicial minimal cone of $\Sigma_{(x_1)}$. Therefore, the non-simplicial minimal cones of $\Sigma_{(x_1)}$ are included in (in fact, equal to) $\{\sigma_2,\dots,\sigma_k\}$. We concludes the proof by repeating this argument $k$-th time.
\end{proof}

\subsubsection{External cone over a polyhedral complex} We now define external cones over polyhedral complexes.

First, let $\phi\colon\R^n\to\R^m$ be an affine linear application and let $X$ be a polyhedral complex in $\R^n$. The \emph{image of $X$ by $\phi$} denoted by $\phi(X)$ is the polyhedral pseudo-complex $\{\phi(\delta)\st\delta\in X\}$ of $\R^m$. If the restriction $\phi\rest{\supp X}$ is injective, $\phi(X)$ is a polyhedral complex.

\medskip

For a subset $S$ of $\R^n$, we denote by $\adh{S}$ the topological closure of $S$ in $\R^n$.

\begin{prop}[Coning over a polyhedron] \label{prop:cone_polyedre}
Let $\delta$ be an arbitrary polyhedron which does not contain the origin in $\R^n$. Let $H$ be any affine hyperplane containing $\Tan\delta$ but not $0$, and denote by $H_0$ the linear hyperplane parallel to $H$. Let $H_0^+$ be the half-space defined by $H_0$ which contains $\delta$ in its interior. Then
\begin{enumerate}
\item\label{enum:coning:1} $\sigma:=\adh{\R_+\delta}$ is a cone which verifies $\delta=\sigma\cap H$ and $\delta_\infty=\sigma\cap H_0$.
\item\label{enum:coning:2} Moreover, $\sigma$ is the unique cone included in $H_0^+$ which verifies $\delta=\sigma\cap H$.
\item\label{enum:coning:3} In addition, if $\delta$ is rational, resp. simplicial, then so is $\sigma$.
\end{enumerate}
\end{prop}
Note that by the terminology adapted in this paper, part \ref{enum:coning:1} means $\sigma$ is strongly convex, i.e., it does not contain any line.
\begin{proof} For part \ref{enum:coning:1}, we will show more precisely that if
\[\delta=\conv(v_0,\dots,v_l) +\R_+u_1+\dots+\R_+u_m,\]
for points $v_0, \dots, v_l$ in $\R^n$ and vectors $u_1, \dots, u_m$ in $H_0$, then $\sigma=\sigma'$ where
\[\sigma':=\R_+v_0+\dots+\R_+v_l+\R_+u_1+\dots+\R_+u_m.\]

Note that obviously we have $\delta\subset\sigma'$, which implies the inclusion $\sigma\subseteq\sigma'$. To prove the inclusion $\sigma'\subseteq\sigma$, since $v_0,\dots,v_l \in \sigma$, we need to show $u_1, \dots, u_m \in \sigma$. Let $x$ be an element of $\delta$. For all $\lambda\geq0$, all the points $x+\lambda u_j$ are in $\delta$. Dividing by $\lambda$, we get $x/\lambda+u_j\in\sigma$. Making $\lambda$ tend to infinity, we get $u_j\in\sigma$ for all $j=1, \dots, m$. This proves $\sigma = \sigma'$.

We now show that $\sigma$ is a cone, i.e., it does not contain any line. To show this, let $\lambda_0,\dots,\lambda_l,\mu_1,\dots,\mu_m\in\R_+$ such that $\lambda_0v_0+\dots+\lambda_lv_l+\mu_1u_1+\dots+\mu_mu_m=0$. It will be enough to show all the coefficients $\lambda_j$ and $\mu_i$ are zero.

First, if one of the coefficients $\lambda_j $ was non-zero, dividing by the linear combination above with the sum $\sum_{j=1}^l \lambda_j$, we would get $ 0 \in \delta$, which would be a contradiction. This shows all the coefficients $\lambda_j$ are zero. Moreover, $\delta_\infty$ is a cone, i.e., it is strongly convex, which proves that all the coefficients $\mu_i$ are zero as well. (Note that in our definition of polyhedra, we assume that a polyhedron is strongly convex, that is it does not contain any affine line.)

\medskip

To finish the proof of part \ref{enum:coning:1}, it remains to show $\delta = \sigma \cap H$ and that $\delta_\infty = \sigma \cap H_0$.

Obviously, we have the inclusion $\delta\subseteq\sigma\cap\Tan\delta\subseteq\sigma\cap H$. Proceeding by absurd, suppose there exists $y\in(\sigma\cap H)\setminus\delta$. We can write $y=\lambda_0v_0+\dots+\lambda_lv_l+\mu_1u_1+\dots+\mu_mu_m$ with non-negative coefficients $\lambda_j, \mu_i$.

We first show that all the coefficients $\lambda_j$ are zero. Otherwise, let $\Lambda=\sum_i\lambda_i$, which is non-zero. Since we are assuming $y\not\in\delta$, we get $\Lambda\neq 1$. The point $x=y/\Lambda$ belongs to $\delta$, and so both points $x,y\in \delta \subset H$. It follows that any linear combination $(1-\lambda) x +\lambda y$ for $\lambda\in \R$ belongs to $ H$. In particular,
\[  0 = \frac\Lambda{1-\Lambda}x-\frac1{1-\Lambda}y\in H, \]
which is a contradiction. This shows $\lambda_0=\dots =\lambda_l =0$.

We thus infer that $y = \mu_1 u_1 + \dots+\mu_m u_m$ belongs to $H_0$. But this is a contradiction because $H \cap H_0 =\emptyset$. We conclude with the equality $\delta = \sigma \cap H$.

We are left to show $\delta_\infty=\sigma\cap H_0$. Take a point $y\in\sigma\cap H_0$, and write $y=\lambda_0v_0+\dots+\lambda_nv_n+\mu_1u_1+\dots+\mu_mu_m$ with non-negative coefficients $\lambda_j, \mu_i$. Then by a similar reasoning as above, all the coefficients $\lambda_i$ should be zero, which implies $y\in\delta_\infty$. This implies the inclusion $\sigma \cap H_0 \subseteq \delta_\infty$. The inclusion $\delta_\infty \subset H_0$ comes from the inclusion $\delta_\infty \subset \TT\delta \subseteq H_0$ and the inclusion $\delta_\infty \subset \sigma$.

\medskip

We now prove part \ref{enum:coning:2}. Let $\sigma'$ be a cone included in $H_0^+$ such that $\sigma'\cap H=\delta$. Obviously, $\sigma=\adh{\R_+\delta}\subseteq\sigma'$. To prove the inverse inclusion, let $y$ be a point $\sigma'$. If $y\not\in H_0$, then there exists a scalar $\lambda>0$ such that $\lambda y\in H$. It follows that $\lambda y\in\delta$ and so $y\in\sigma$. If now $y\in H_0$, then for all points $x\in\delta$ and any real $\lambda>0$, the point $x+\lambda y$ belongs to $\sigma\cap H$, which is equal to $\delta$ by the first part of the proposition. We infer that $y\in\delta_\infty\subset\sigma$. This implies that $\sigma'=\sigma$ and the unicity follows.

\medskip

Part \ref{enum:coning:3} follows from the description of the cone $\sigma$ given in the beginning of the proof.
\end{proof}

Let now $X$ be a polyhedral complex in $\R^n$. We define the \emph{cone over $X$} denoted by $\cone{X}$ by
\[ \cone{X}:=\Bigl\{\,\adh{\R_+\delta}\st\delta\in X\,\Bigr\}\cup\Bigl\{\conezero\Bigr\}. \]

In particular, we note that this is not always a pseudo-fan. However, if one assumes that $\supp{X}$ does not contain $0$, then $\cone{X}$ is indeed a pseudo-fan which contains $X_\infty$.

\medskip

We define the \emph{external cone over $X$} denoted by $\coneup{X}$ as the pseudo-fan of $\R^{n+1}$ defined as follows. This is also called sometimes the \emph{homogenization of $X$},

Let $(e_0,\dots,e_n)$ be the standard basis of $\R^{n+1}$, and let $x_0,\dots,x_n$ be the corresponding coordinates.
Consider the projection $\pi\colon\R^{n+1}\to\R^n$ defined by
\[ \pi(x_0,\dots,x_n)=(x_1,\dots,x_n). \]

Consider the affine hyperplane $H_1$ of $\R^{n+1}$ given by the equation $x_0=1$. The restriction of $\pi$ to $H_1$ gives an isomorphism $\pi_1\colon H_1\simto\R^n$, and we define
\[ \coneup X:=\cone{\pi_1^{-1}(X)}. \]

In particular, for any polyhedron $\delta$, we get the external cone
\[ \coneup \delta:=\cone{\pi_1^{-1}(\delta)}. \]
By Proposition~\ref{prop:cone_polyedre}, intersecting with $H_1$ gives
\[ \pi(\coneup{X}\cap H_1)=X. \]

The pseudo-fan $\coneup{X}$ is included in $\{x_0\geq 0\}$. Moreover, denoting by $H_0$ the linear hyperplane $\{x_0 = 0\}$, we get
\[\coneup{X}\rest{H_0}=\coneup{X}\cap H_0.\]
Applying again Proposition~\ref{prop:cone_polyedre}, we deduce that $\pi(\coneup{X}\cap H_0)$ is the recession pseudo-fan $X_\infty$ of $X$.

\subsubsection{From $X_\infty$ to $X$ and vice-versa}
\label{sec:X_infty_X}

Let $X$ be any polyhedral complex such that $X_\infty$ is a fan. Let $f$ be a piecewise linear function on $X$. The function $f$ naturally induces a function on $\cone(X)\cap H_1$. This function can be extended by linearity to $\cone{X}\setminus H_0$. We denote this extension by $\~f$. We say that \emph{$f_\infty$ is well-defined} if $\~f$ can be extended by continuity to $H_0$. In this case, we define $f_\infty$ to be the restriction of $\~f$ to $X_\infty=\coneup{X}\cap H_0$.

\begin{remark}
Notice that if $f$ is strictly convex and if $f_\infty$ is well-defined, then $f_\infty$ is strictly convex.
\end{remark}

\medskip

Assume $X$ is simplicial. In this case, recall that for any polyhedron $\delta\in X$ and any point $x$ in $\delta$, there is a unique decomposition $x=x_\f+x_\infty$ with $x_\f\in\delta_\f$ and $x_\infty\in\delta_\infty$. We define the \emph{projection onto the finite part} $\pi_\f\colon\supp{X}\to\supp{X_\f}$ as the map $x\mapsto x_\f$. Similarly, we define the \emph{projection onto the asymptotic part} $\pi_\infty\colon \supp X\to \supp{X_\infty}$ by $x\mapsto x_\infty$. If $Y_\f$ is subdivision of $X_\f$, we define
\[ \pi_\f^{-1}(Y_\f):=\bigcup_{\delta\in X}\Bigl\{\,\gamma+\sigma\,\st\,\sigma\subface\delta_\infty, \gamma\in Y_\f, \gamma\subset\delta_\f\, \Bigr\}. \]

\begin{remark} \label{rem:extension_compatible}
Setting $Z:=\pi_\f^{-1}(Y_\f)$, one can show that $Z$ is a polyhedral complex, which is moreover simplicial if $Y_\f$ is simplicial. In addition, we have $Z_\f=Y_\f$ et $Z_\infty=X_\infty$.
\end{remark}

\subsection{Regular subdivisions and convexity}
\label{subsec:convexite}

Let now $X$ be a polyhedral complex in $\R^n$, and let $\delta$ be a face of $X$. Recall that a function $f$ is called strictly positive around $\delta$ if there exists a a neighborhood $V$ of the relative interior of $\delta$ such that $f$ vanishes on $\delta$ and it is strictly positive on $V \setminus \delta$. Recall as well from Section~\ref{sec:tropvar} that a function $f\in\lpm(X)$ is called \emph{strictly convex around $\delta$} if there exists a function $\ell\in\aff(\R^n)$ such that $f-\ell$ is strictly positive around the relative interior of $\delta$ in $\supp{X}$. We also note that $\K_+(X)$ denotes the set of strictly convex functions on $X$, i.e., the set of those functions in $\lpm(X)$ which are strictly convex around any face $\delta$ of $X$.

\begin{remark} For a polyhedral complex $X$, the set $\K_+(X)$ is an open convex cone in the real vector space $\lpm(X)$, \ie, it is an open set in the natural topology induced on $\lpm(X)$, and moreover, if $f,g\in\K_+(X)$, then we have $f+g \in\K_+(X)$ and $\lambda f \in\K_+(X)$ for any positive real number $\lambda>0$. In addition, for $\ell\in\aff(X)$, we have $f+\ell\in\K_+(X)$.
\end{remark}

\begin{example} If $P$ is a polyhedron in $\R^n$, then we have $0\in\K_+(\face P)$. Indeed, for each face $\tau$ of $P$, by definition of faces $\delta$ in a polyhedron, there exists a linear form $\ell$ on $\R^n$ which is zero on $\delta$ and which is strictly negative on $P\setminus \delta$, so $0-\ell = -\ell$ is strictly positive around $\delta$.
\end{example}

We call a polyhedral complex $X$ \emph{quasi-projective} if the set $\K_+(X)$ is non-empty. Quasi-projective polyhedral complexes are sometimes called \emph{convex} in the literature.

\begin{remark}
In the case of fans, the notion coincides with the notion of quasi-projectivity (convexity) used in the previous section.
\end{remark}

We are now going to make a list of operations on a polyhedral complex which preserves the quasi-projectivity property.
To start, we have the following which is straightforward.

\begin{prop} \label{prop:sous-complexe_projectif}
A subcomplex of a quasi-projective polyhedral complex is quasi-projective.
\end{prop}

Let now $Y$ be a subdivision of $X$. For each face $\delta$ of $X$, recall that the restriction $Y\rest \delta$ consists of all the faces of $Y$ which are included in $\delta$.

\medskip

Let $X$ be a polyhedral complex. A \emph{regular subdivision of $X$} is a subdivision $Y$ such that there exists a function $f\in\lpm(Y)$ whose restriction to $Y\rest\delta$ for any face $\delta\in X$ is strictly convex. In this case, we say $Y$ is \emph{regular relative to $X$}, and that the function $f$ is \emph{strictly convex relative to $X$}. For example, a strictly convex function $f$ on $Y$ is necessary strictly convex relative to $X$ (but the inverse is not necessarily true).

\begin{remark} The terminology is borrowed from the theory of triangulations of polytopes and secondary polytopes~\cites{GKZ, DRS}. In fact, a polyhedral subdivision of a polytope is regular in our sense if and only it is regular in the sense of~\cite{GKZ}.
\end{remark}

\begin{remark} To justify the terminology, we should emphasize that regularity for a subdivision is a relative notion while quasi-projectivity for a polyhedral complex is absolute. In particular, the underlying polyhedral complex of a regular subdivision does not need to be quasi-projective.
\end{remark}

\begin{remark}
A function $f \in \lpm(Y)$ which is strictly convex relative to $X$ allows to recover the subdivision $Y$ as the polyhedral complex whose collection of faces is
\[ \bigcup_{\delta\in X}\bigcup_{\substack{\ell\in(\R^n)^\dual \\ \ell\rest\delta\leq f\rest\delta}}\Bigl\{x\in\delta\st\ell(x)=f(x)\Bigr\}\ \setminus\Bigl\{\emptyset\Bigr\}. \qedhere \]
\end{remark}

\begin{prop} \label{prop:projectivite_convexite}
A regular subdivision of a quasi-projective polyhedral complex is itself quasi-projective.
\end{prop}

\begin{proof}
Let $X$ be a quasi-projective polyhedral complex in $\R^n$ and let $f\in\K_+(X)$. Let $Y$ be a regular subdivision of $X$ and let $g\in\lpm (Y)$ be strictly convex relative to $X$. We show that for any small enough $\epsilon>0$, the function $f+ \epsilon g \in \lpm(Y)$ is strictly convex on $Y$, which proves the proposition.

Let $\gamma$ be a face of $Y$. It will be enough to show that $f+ \epsilon g$ is strictly convex around $\gamma$ for any sufficiently small value of $\epsilon>0$.

Let $\delta$ be the smallest face of $X$ which contains $\gamma$. Let $\ell\in\aff(X)$ and $\ell' \in \aff(Y\rest\delta)$ be two affine linear functions such that $f-\ell$ is strictly positive around $\delta$ and $g\rest{\delta} -\ell'$ is strictly positive around $\gamma$. Without loss of generality, we can assume that $\ell' \in \aff(X)$.

Let $V$ be a small enough neighborhood of the interior of $\gamma$ in $\supp{X}$. Equipping $\R^n$ with its Euclidean norm, we can in addition choose $V$ in such a way that for each point $x \in V$, the closest point of $x$ in $\delta$ belongs to $V$.

Since $f$ and $g$ are piecewise linear, we can find $c, c'>0$ such that for any $x\in V$,
\begin{gather*}
f(x)-\ell(x)\geq c\,\dist(x,\delta)\quad\text{and} \\
g(x)-\ell'(x)\geq \min_{y\in V\cap\delta}\big\{g(y) -\ell'(y)\bigr\}-c'\,\dist(x,\delta)=-c'\,\dist(x,\delta),
\end{gather*}
where $\dist(x,\delta)$ denotes the Euclidean distance between the point $x$ and the face $\delta$.

We infer that the function $f+\epsilon g-(\ell+\epsilon \ell')$ is strictly positive around $\gamma$ for any value of $\epsilon<c/c'$. It follows that $f+\epsilon g$ is strictly convex around $\gamma$ for small enough values of $\epsilon>0$, as required. This shows that $Y$ is quasi-projective.
\end{proof}

As an immediate consequence we get the following.

\begin{cor} \label{cor:transitivite_projectivite}
The relation of \enquote{being a regular subdivision} is transitive.
\end{cor}

\begin{proof}
Let $Y$ be a regular subdivision of a polyhedral complex $X$ with $f \in \lpm(Y)$ strictly convex relative to $X$, and let $Z$ be a regular subdivision of $Y$ with $g \in \lpm(Z)$ strictly convex relative to $Y$.
It follows that for a face $\delta \in X$, the restriction $f\rest{\delta}$ is strictly convex on $Y\rest{\delta}$, and the restriction $g\rest{\delta}$, seen as a function in $\lpm (Z\rest\delta)$, is strictly convex relative to $Y \rest \delta$.

By the previous proposition, more precisely by its proof, for sufficiently small values of $\epsilon>0$, the function $f\rest{\delta}+\epsilon g\rest{\delta}$ is strictly convex on $Z\rest{\delta}$. Taking $\epsilon>0$ small enough so that this holds for all faces $\delta$ of $X$, we get the result.
\end{proof}

\begin{prop} \label{prop:eclate_projectif}
Let $\Sigma$ be a fan in $\R^n$ and take a point $x\in \R^n\setminus\{0\}$. The blow-up $\Sigma_{(x)}$ of $\Sigma$ at $x$ is regular relative to $\Sigma$.
\end{prop}

\begin{proof} There exists a unique function $ f \in \lpm (\Sigma_{(x)}) $ which vanishes on any face of $\Sigma$ that does not contain $x$ and which takes value $ -1 $ at $ x $. We show that this function is strictly convex relative to $\Sigma$, which proves the proposition.

Let $\sigma \in \Sigma$. If $\sigma$ does not contain $x$, then $\sigma$ belongs to $\Sigma_{(x)}$, and $f$ is zero, therefore convex, on $\Sigma_{(x)}\rest\sigma=\face{\sigma}$. Consider now the case where $ x \in \sigma $. Let $ \tau $ be a face of $ \sigma$ not containing $x$. By the definition of the blow-up, we have to show that $ f \rest \sigma$ is convex around $ \tau $ and around $ \tau + \R_+ x $. Let $ \ell \in \aff (\R^n)$ be an affine linear function which vanishes on $\tau $ and which is strictly negative on $ \sigma \setminus \tau $. Moreover, we can assume without losing generality that $ \ell(x) = -1 $. The function $f-\ell$ is therefore zero on $\tau + \R_+ x $. In addition, if $\nu$ is a face of $\sigma$ not containing $x$, $f$ vanishes on $\nu$, so $f-\ell$ is strictly positive on $\nu \setminus \tau$. If now $y$ is a point lying in $\sigma \setminus (\tau + \R_+ x)$, there is a face $\nu$ of $\sigma$ such that $ y$ belongs to $(\nu \setminus \tau) + \R_+ x$, so we can write $y = \lambda x + z $ for some real number $\lambda \geq 0$ and an element $z \in \nu \setminus \tau$. Given the definition of $f$ and the choice of $\ell$, we see that
\[ (f- \ell) (y) = \lambda \cdot0 + (f- \ell) (z)> 0. \]
Thus, $f-\ell$ is strictly positive around $\tau + \R_+ x$ in $\sigma$. In the same way, we see that $f-2\ell$ is strictly positive around $\tau$ in $\sigma$. So $f\rest \sigma$ is strictly convex around $\tau + \R_+ x$ and around $\tau$ in $\Sigma_{(x)}\rest \sigma$, and the proposition follows.
\end{proof}

\begin{prop} \label{prop:polyedre_convexe}
Let $P$ be a polytope in $\R^n$ which contains $0$. The fan $\Sigma:=\cone{\face{P}\rest{\R^n\setminus\{0\}}}$ consisting of cones of the form $\R_+\sigma$ for $\sigma$ a face of $P$ is quasi-projective.
\end{prop}

\begin{proof} We omit the proof that $\Sigma$ is a fan. In order to prove the quasi-projectivity, let $\tau$ be a face of $P$ which does not contain $0$. There exists an affine linear $\ell^\tau$ function on $\R^n$ which takes value $0$ on $\tau$ and which takes negative values on $P \setminus \tau$. Moreover, we can assume that $\ell^\tau(0)=-1$.
The restrictions $\ell^\tau\rest{\cone{\tau}}$ glue together to define a piecewise affine linear function $f$ on $\Sigma$. One can show that $f$ is strictly convex on $\Sigma$.
\end{proof}

\begin{prop} \label{prop:coupe_projective}
Let $X$ be a polyhedral complex and let $H$ be an affine hyperplane in $\R^n$. The hyperplane cut $X\cdot H$ is regular relative to $X$.
\end{prop}

\begin{proof} Let $\dist(\ccdot,\rdot)$ be the Euclidean distance in $\R^n$. It is easy to see that the distance function $x\mapsto \dist(x, H)$ is strictly convex on $X\cdot H$ relative to $X$.
\end{proof}

\begin{cor} \label{prop:intersection_projective}
Suppose that $X$ is quasi-projective. Then the intersection $X\cap H$ is quasi-projective as well.
\end{cor}

\begin{proof} If $X$ is quasi-projective, by combining the preceding proposition with Proposition~\ref{prop:projectivite_convexite}, we infer that the hyperplane cut $X\cdot H$ is quasi-projective as well. Since $X\cap H=(X\cdot H)\rest H$ is a polyhedral subcomplex of $X\cdot H$, it will be itself quasi-projective by Proposition~\ref{prop:sous-complexe_projectif}.
\end{proof}

\begin{prop} \label{prop:extension_projectivite_partie_finie}
Let $X$ be a polyhedral complex in $\R^n$. If $X$ is simplicial and $Y_\f$ is a regular subdivision of $X_\f$, the subdivision $Y:=\pi_\f^{-1}(Y_\f)$ of $X$ is regular relative to $X$. Moreover, there exists a function $f\in\lpm(Y)$ which is strictly convex relative to $X$ and such that $f_\infty=0$.
\end{prop}

\begin{proof}
\renewcommand{\int}{\textrm{int}}

By Remark~\ref{rem:extension_compatible}, the polyhedral complex $Y$ is a subdivision of $X$ whose finite part coincides with $Y_\f$.

\medskip

Let $f_\f$ be a piecewise affine linear function on $Y_\f$ which is strictly convex relative to $X_\f$.
Let $f:=\pi_\f^*(f_\f) = f_\f \circ \pi_\f$ be the extension of $f$ to $\supp Y=\supp X$. It is clear that $f_\infty=0$. For a face $\gamma\in Y$, since both the projection $\pi_\f\rest\gamma\colon\gamma\to\gamma_\f$ and the restriction $f_\f\rest{\gamma_\f}$ are affine linear, $f\rest\gamma$ is also affine linear. This shows that $f \in \lpm(Y)$. We now show that $f$ is strictly convex relative to $X$, which proves the proposition.

\medskip

Let $\delta$ be a face of $X$ and let $\gamma$ be a face in $Y$ which is included in $\delta$. We need to show that $f\rest\delta$, as a piecewise linear function on $Y\rest\delta$, is strictly convex around $\gamma$.
By the choice of $f_\f$, we know that $f_\f$ is strictly convex around $\gamma_\f$ in $Y_\f\rest{\delta_\f}$. Let $\ell_\f$ be an element of $\aff(\delta_\f)$ such that $f_\f\rest{\delta_\f}-\ell_\f$ is strictly positive around $\int(\gamma_\f)$ in $\delta_\f$, where $\int(\gamma_\f)$ denotes the relative interior of $\gamma_\f$.

The function $\ell:=\pi_\f^*(\ell_\f)$ is an element of $\aff(\delta)$. Moreover, the difference $f-\ell$ is strictly positive around $\int(\gamma_\f)+\delta_\infty$ in $\delta$. Notice that the relative interior of $\gamma$ is included in $\int(\gamma_\f)+\delta_\infty$. Let $\ell'\in\aff(\delta)$ be an affine linear function which is zero on $\gamma$ and which takes negative values on $\bigl(\gamma_\f+\delta_\infty\bigr)\setminus\gamma$. Then, for $\epsilon>0$ small enough, $f-\ell-\epsilon\ell'$ is strictly positive around $\int(\gamma)$. Thus, $f$ is strictly convex around $\gamma$, and the proposition follows.
\end{proof}

\subsection{Unimodularity} In this section, we discuss the existence of unimodular subdivisions of a rational polyhedral complex in $\R^n$.

\begin{prop} \label{prop:subdivision_unimodulaire_projective}
Let $\Sigma$ be a rational simplicial fan in $\R^n$. There exists a regular subdivision $\~\Sigma$ of $\Sigma$ which is in addition unimodular with respect to $\Z^n$. Moreover, we can assume that unimodular cones of $\Sigma$ are still in $\~\Sigma$.
\end{prop}

\begin{proof}
This follows from the results of~\cite{Wlo97}*{Section 9}, where it is shown how to get a unimodular subdivision $\~\Sigma$ of $\Sigma$ by a sequence of blow-ups on some non-unimodular cones. By Proposition~\ref{prop:eclate_projectif} and Corollary~\ref{cor:transitivite_projectivite}, the resulting fan $\~\Sigma$ is a projective subdivision of $\Sigma$, from which the result follows.
\end{proof}

For a rational vector subspace $W$ of $\R^n$, we denote by $\p_W\colon\R^n\to\rquot{\R^n}{W}$ the resulting projection map. We endow the quotient $\rquot{\R^n}W$ by the quotient lattice $\rquot{\Z^n}{\Z^n\cap W}$.

\medskip

Let now $P$ be a simplicial polyhedron in $\R^n$ and $W:=\TT(P_\f)$. We say that \emph{$P_\infty$ is unimodular relative to $P_\f$} if the projection $\p_W(P_\infty)$ is unimodular in $\rquot{\R^n}{W}$.

\begin{remark} \label{rem:modulaire_subdivision}
If $P_\infty$ is unimodular relative to $P_\f$, then for any polytope $Q$ included in $\TT(P_\f)$, for the Minkowski sum $R:=Q+P_\infty$ we see that $R_\infty$ is unimodular relative $R_\f$
\end{remark}

\begin{prop} \label{prop:decomposition_unimodularite}
Let $P$ be an integral simplicial polyhedron in $\R^n$. Then $P$ is unimodular if and only if both the following conditions are verified:
\begin{enumerate}
\item $P_\f$ is unimodular; and
\item $P_\infty$ is unimodular relative to $P_\f$.
\end{enumerate}
\end{prop}

For a collection of vectors $v_1, \dots, v_m$ in a real vector space, we denote by $\Vect(v_1, \dots, v_m) := \R v_1+ \dots +\R v_m$ the vector subspace generated by the vectors $v_i$.
\begin{proof}
We can write $P=\conv(v_0,\dots,v_l)+\R_+v_{l+1}+\dots+\R_+v_m$, as in Section~\ref{sec:recol}. Moreover, we can suppose without loss of generality that $v_0=0$, that $v_i$ is primitive in $\R_+v_i$ for $i\in\zint{l+1,m}$, and that $\TT(P)=\R^n$. The claimed equivalence can be rewritten in the form: \emph{$(v_1,\dots,v_m)$ is a basis of $\Z^n$ if and only if $(v_1,\dots,v_l)$ is a basis of $\Z^n\cap\Vect(v_1,\dots,v_m)$ and $(v_{l+1},\dots,v_m)$ is a basis of $\rquot{\Z^n}{\Z^n\cap\Vect(v_1,\dots,v_m)}$.} This is clear.
\end{proof}

\begin{cor} \label{cor:coneup_modulaire}
If a polyhedron $P$ is \emph{rational} and its coning $\coneup{P}$ is unimodular in $\R^{n+1}$, with respect to the lattice $\Z^{n+1}$, then $P_\infty$ is unimodular relative to $P_\f$.
\end{cor}

\begin{proof}
We use the notations used in the definition of $\coneup P$. So let $\pi\colon\R^{n+1}\to\R^n$ be the projection, and let $H_0$ be the hyperplane $\{x_0=0\}$ and $\pi_0:=\pi\rest{H_0}$. Let $u_1, \dots, u_l, u_{l+1} \dots, u_m\in\R^{n+1}$ be so that $\coneup P=\R_+u_1+\dots+\R_+u_m$ and $u_1,\dots,u_l\not\in H_0$ while $u_{l+1}, \dots, u_m\in H_0$, and $u_1,\dots,u_m$ is part of a basis of $\Z^n$.

We apply the preceding proposition to $P':=\conv(0,u_1,\dots, u_l)+\R_+u_{l+1}+\dots+\R_+u_m$ which is clearly unimodular. Setting $W'=\TT(P'_\f)$, we infer that $\p_{W'}(P'_\infty)$ is a unimodular cone in $\rquot{\R^{n+1}}{W'}$. Note that $\rquot{\R^{n+1}}{W'}=\rquot{H_0}{W'\!\cap\!H_0}$. Moreover, the linear application $\pi_0\colon H_0\to\R^n$ is an isomorphism which preserves the corresponding lattices, and which sends $P'_\infty=\coneup P\cap H_0$ onto $P_\infty$, and $W'\cap H_0$ onto $W := \TT(P_\f)$. Therefore, $\pi_0$ induces an isomorphism from $\rquot{H_0}{W'\!\cap\!H_0}$ to $\rquot{\R^n}{W}$ which sends $\p_{W'}(P'_\infty)$ onto $\p_W(P_\infty)$. This finally shows that $\p_W(P_\infty)$ is unimodular.
\end{proof}

\begin{thm} \label{thm:triangulation_unimodulaire_projective}
Let $X$ be a \emph{compact} rational polyhedral complex in $\R^n$ (so $X_\infty=\{\conezero\}$). There exists a positive integer $k$ and a regular subdivision $Y$ of $X$ such that $Y$ is unimodular with respect to the lattice $\frac 1k \Z^n$.
\end{thm}

\begin{proof}
This is a consequence of Theorem 4.1 of \cite{KKMS}*{Chapter III}.
\end{proof}

\begin{remark} \label{rem:stabilite_eventail_unimodulaire}
Note that if a fan $\Sigma$ in $\R^n$ is unimodular with respect to $\Z^n$, then it is as well unimodular with respect to any lattice of the form $\frac1k\Z^n$.
\end{remark}

\subsection{Proof of Theorem~\ref{thm:triangulation_unimodulaire_convexe}}
We are now in position to prove the main theorem of this section. Let $\Sigma_X:=\coneup{X}$, which is a pseudo-fan, and define $H_0$, $H_1$ and $\pi\colon\R^{n+1}\to\R^n$ as in the definition of the external cone. We thus have $X=\pi(\Sigma_X\cap H_1)$ and $X_\infty=\pi(\Sigma_X\cap H_0)$.

\medskip

We first prove the existence of a regular subdivision of $\Sigma_X$.

\medskip

Let $P$ be a polytope containing $0$ in its interior. Define the fan $\Sigma_1:=\cone{\face{P}\rest{\R^{n+1}\setminus\{0\}}}$. This is a complete fan which moreover, by Proposition \ref{prop:polyedre_convexe}, is projective. Define $\Sigma_2$ as the fan obtained by slicing $\Sigma_1$ with respect to a pencil of hyperplanes associated to $\Sigma_X$. (We defined the slicing with respect to a fan but as mentioned previously, the definition extends to pseudo-fans.)

By Proposition \ref{prop:decoupe_subdivision}, $\Sigma_2$ contains a subdivision of $\Sigma_X$.
Moreover, since $\Sigma_1$ is quasi-projective, and since by Proposition \ref{prop:coupe_projective}, $\Sigma_2$ is regular relative to $\Sigma_1$ (being obtained by a sequence of hyperplane cuts), it follows that $\Sigma_2$ is quasi-projective.

\medskip

Let $\Sigma_3:=\Sigma_2\cdot H_0$, and note that this is a regular subdivision of $\Sigma_2$ by the same proposition. In particular, $\Sigma_3$ itself is quasi-projective. The fan $\Sigma_3$ is special in that if $\Sigma'$ if a subdivision of $\Sigma_3$, then any cone in $\Sigma'$ is included in one of the two half-spaces $H_0^+$ or $H_0^-$ defined by $H_0$.

\medskip

Let $\Sigma_4$ be a rational triangulation of $\Sigma_3$ obtained from $\Sigma_3$ by applying Proposition \ref{prop:triangulation}. Note that $\Sigma_4$ is obtained from $\Sigma_3$ by blow-ups, thus, by Proposition \ref{prop:eclate_projectif}, it is regular relative to $\Sigma_3$. In particular, it is quasi-projective.

\medskip

Let $\Sigma_5$ be a unimodular subdivision of $\Sigma_4$, which exists by Proposition~\ref{prop:subdivision_unimodulaire_projective}. The fan $\Sigma_5$ is therefore unimodular and quasi-projective, and it contains moreover a subdivision of $\Sigma_X$.

\medskip

We now use $\Sigma_5$ to get a subdivision of $X$. Let $X_1:=\pi(\Sigma_5\cap H_1)$.

\medskip

Since $\Sigma_5$ is quasi-projective, the intersection $\Sigma_5\cap H_1$ is quasi-projective as well. Since $\Sigma_5$ contains a subdivision of $\Sigma_X$, it follows that $\Sigma_5\cap H_1$ contains a subdivision of $\Sigma_X\cap H_1=X$.

Passing now to the projection, we infer that $X_1$ is quasi-projective and contains a subdivision of $X$. More precisely, there exists a strictly convex function $f$ on $X_1$ such that $f_\infty$ is well-defined.

\medskip

Pick $\delta$ in $X_1$, and let $\sigma$ be a cone in $\Sigma_5$ such that $\pi(\sigma\cap H_1)=\delta$. Applying Proposition~\ref{prop:cone_polyedre} to $\delta':=\sigma\cap H_1$ and to the hyperplane $H_1$, given that $\sigma$ is included in the half-space defined by $H_0$ which contains $\delta'$, we infer that $\sigma=\adh{\R_+\delta}=\coneup\delta$. Since $\sigma$ is unimodular, by Corollary~\ref{cor:coneup_modulaire} we infer that $\delta_\infty$ is unimodular relative $\delta_\f$.

\medskip

Let $X_2:=X_1\rest{\supp X}$. It follows that $X_2$ is a subdivision of $X$ which is rational, simplicial, and quasi-projective, and moreover, for any $\delta\in X_2$, the cone $\delta_\infty$ is unimodular relative to $\delta_\f$.

\medskip

At this moment, we are left to show the existence of a unimodular regular subdivision $Y$ of $X_2$. By Theorem~\ref{thm:triangulation_unimodulaire_projective}, there exists an integer $k\geq1$ and a regular subdivision $Y_\f$ of $X_{2,\f}$ such that $Y_\f$ is unimodular with respect to the lattice $\frac1k\Z^n$. Define $Y:=\pi_\f^{-1}(Y_\f)$. By Proposition~\ref{prop:extension_projectivite_partie_finie}, $Y$ is regular relative to $X_2$, and so by quasi-projectivity of $X_2$, it is quasi-projective itself. More precisely, there exists a function $f$ strictly convex on $Y$ such that $Y_\infty$ is well-defined.

\medskip

Using Remark~\ref{rem:modulaire_subdivision}, we infer that for each face $\delta \in Y$, the cone $\delta_\infty$ is unimodular relative to $\delta_\f$, with respect to the lattice $\Z^n$ of $\R^n$. By Remark~\ref{rem:stabilite_eventail_unimodulaire}, this remains true for the lattice $\frac 1k\Z^n$ as well.

\medskip

Applying Proposition~\ref{prop:decomposition_unimodularite}, we finally conclude that $Y$ is unimodular with respect to the lattice $\frac 1k\Z^n$. At present, we have obtained a regular subdivision $Y$ of $X$ which is moreover quasi-projective and unimodular with respect to the lattice $\frac 1k \Z^n$. Moreover, it is clear that $Y_\infty$ is a fan. This finishes the proof of Theorem~\ref{thm:triangulation_unimodulaire_convexe}.

\medskip

\subsection{Unimodular triangulations preserving the recession fan}
We finish this section by presenting the proof of Theorem~\ref{thm:unimodular_preserving_recession}. We use the same notations as in the previous section for external cones.

\begin{proof}[Proof of Theorem~\ref{thm:unimodular_preserving_recession}]
Let $\Sigma_1=\coneup X$. We know that $\Sigma_1\rest{H_0}=X_\infty$ is a unimodular fan. Thus $\Sigma_1$ is a fan.

\medskip

By Proposition~\ref{prop:triangulation}, we can construct a simplicial subdivision $\Sigma_2$ of $\Sigma_1$. Moreover, since faces of $\Sigma_1\rest{H_0}$ are already simplicial, the construction described in Proposition~\ref{prop:triangulation} does not change this part of the fan, \ie, $\Sigma_2\rest{H_0}=\Sigma_1\rest {H_0}$.

\medskip

In the same way, by Proposition~\ref{prop:subdivision_unimodulaire_projective}, there exists a unimodular subdivision $\Sigma_3$ of $\Sigma_2$. Moreover, since $\Sigma_2\rest{H_0}$ is unimodular, the same proposition ensures we can assume $\Sigma_3\rest{H_0}=\Sigma_2\rest{H_0}$.

\medskip

Set $X_1=\pi(\Sigma_3\cap H_1)$. Then $X_1$ is a triangulation of $X$ with $X_{1,\infty}=X_\infty$.
\medskip

By Theorem~\ref{thm:triangulation_unimodulaire_projective}, there exists an integer $k\geq1$ and a subdivision $Y_\f$ of $X_{1,\f}$ such that $Y_\f$ is unimodular with respect to the lattice $\frac1k\Z^n$. Define $Y:=\pi_\f^{-1}(Y_\f)$. Clearly, $Y_\infty=X_{1,\infty}=X_\infty$.

\medskip

Using the same argument as in the final part of the proof of Theorem~\ref{thm:triangulation_unimodulaire_convexe} in the previous section, we infer that $Y$ is unimodular, which concludes the proof.
\end{proof}


\section{Tropical Steenbrink spectral sequence}\label{sec:steenbrink}

Let $V:=\R^n$ be a vector space. Denote by $N$ the lattice $\Z^n\subset V$, thus $V=N_\R$. Let $Y$ be a unimodular rational simplicial complex of dimension $d$ in $N_\R$ whose support $\supp{Y}$ is a smooth tropical variety. Assume that the recession pseudo-fan of $Y$ is a fan $Y_\infty$. This fan induces a partial compactification $\TP_{Y_\sminfty}$ of $\R^n$. Denote by $X$ the extended polyhedral complex induced by $Y$ on the closure of $Y$ in $\TP_{Y_\sminfty}$. Then, $\supp{X}$ is a smooth compact tropical variety. To simplify the notations, we will denote by $X_\infty$ the recession fan $Y_\infty$ of $Y$, and by $\TP_{X_{\sminfty}}$ the partial compactification $\TP_{Y_{\sminfty}}$ of $\R^n$.

\medskip

Recall from Section~\ref{sec:tropvar} that the space $\TP_{X_\sminfty}$ is naturally stratified by open strata $N^\sigma_\R$ for $\sigma\in X_\infty$. If $\sigma\in X_\infty$, we denote by $X^\sigma$ the intersection of $X$ with $N^\sigma_\R$. Notice that $X^\sigma$ is a unimodular simplicial complex of dimension $d-\dims{\sigma}$. In particular, $X^\conezero=Y$. If $\delta\in X$ is a face, we denote by $\sed(\delta)\in X_\infty$ its sedentarity: this is the unique cone $\sigma\in X_\infty$ such that the relative interior of $\delta$ is in $N^\sigma_\R$.

\medskip

Recall that for a simplicial polyhedron $\delta$ in $\R^n$, we defined $\delta_\infty$ and $\delta_\f$ in Section~\ref{sec:tropvar}. We extend the definition to faces $\delta$ of $X$ of higher sedentarity as follows. By $\delta_\infty$ (\resp $\delta_\f$) we mean $(\delta\cap N^{\sed(\delta)}_\R)_\infty$ (\resp $(\delta\cap N^{\sed(\delta)}_\R)_\f$) which belongs to the fun $(X^{\sed(\delta)})_\infty$ (\resp $(X^{\sed(\delta)})_\f$).

\medskip

The \emph{max-sedentarity} of a face $\delta \in X$ that we denote by $\maxsed(\delta)$ is defined as
\[\maxsed(\delta) :=\max\Bigl\{\sed(\gamma) \st \gamma \subface \delta\Bigr\}.\]
One verifies directly that we have $\maxsed(\delta) = \sed(\delta)+\delta_\infty$ and that this is the maximum cone $\sigma\in X_\infty$ such that $\delta$ intersects $N^\sigma_\R$.

\medskip

For each $\delta \in X$, let $\Sigma^\delta$ be the fan in $\rquot{N^{\sed(\delta)}_\R\!}{\TT\delta}$ induced by the star of $\delta$ in $X$. Unimodularity of $X$ implies that the fan $\Sigma^\delta$ is unimodular. We denote as before by $\comp\Sigma^\delta$ the canonical compactification of $\Sigma^\delta$.

\medskip

We define the \emph{$k$-th cohomology} $H^{k}(\delta)$ of $\comp\Sigma^\delta$ by
\[H^k(\delta):= \bigoplus_{r+s = k} H_{\trop}^{r,s}(\comp \Sigma^\delta) = \begin{cases} H_{\trop}^{k/2,k/2} (\comp \Sigma^\delta) \simeq A^{k/2}(\Sigma^\delta) & \textrm{ if $k$ is even} \\
0 & \textrm{ otherwise}
\end{cases}\]
where the last equalities follows from Theorem~\ref{thm:HI}, proved in~\cite{AP}.

\subsection{Basic maps} \label{sec:basic_maps}

We define some basic maps between the cohomology groups associated to faces of $X$.

\subsubsection{The sign function} We will need to make the choice of a \emph{sign function} for $X$. This is a function
\[\sign\colon \Bigl\{(\gamma,\delta)\in X\times X   \,\st\, \gamma\ssubface\delta \Bigr\}\longrightarrow \{-1,+1\}\]
which takes a value $\sign(\gamma, \delta)\in\{-1,+1\}$ for a pair of faces $\gamma,\delta\in X$ with $\gamma\ssubface\delta$, and which verifies the following property that we call the \emph{diamond property}.

\begin{defi}[Diamond property for a sign function]\rm Let $\delta$ be a face of dimension at least 2 in $X$ and $\nu$ be a face of codimension 2 in $\delta$. The diamond property of polyhedra ensures that there exist exactly two different faces $\gamma, \gamma'\in X$ such that $\nu\ssubface\gamma,\!\gamma'\ssubface\delta$. Then, we require the sign function to verify the following compatibility
\[ \sign(\nu, \gamma)\sign(\gamma, \delta)=-\sign(\nu, \gamma')\sign(\gamma', \delta). \]
Moreover, if $e$ is an edge in $X$ and if $u$ and $v$ are its extremities, then we require the sign function to verify
\[ \sign(u,e)=-\sign(v,e). \]
\end{defi}
For the existence, note that such a sign function can be easily defined by choosing an orientation on each face.

From now on, we assume that a sign function for $X$ has been fixed once for all.

\subsubsection{The restriction and Gysin maps}
If $\gamma$ and $\delta$ are two faces with $\gamma\ssubface\delta$, we have the \emph{Gysin map}
\[ \gys_{\delta\ssupface\gamma}\colon H^k(\delta)\to H^{k+2}(\gamma), \]
and the \emph{restriction map}
\[ \i^*_{\gamma\ssubface\delta}\colon H^k(\gamma)\to H^k(\delta), \]
where we define both maps to be zero if $\gamma$ and $\delta$ do not have the same sedentarity.
If $\gamma$ and $\delta$ have the same sedentarity, these maps correspond to the corresponding maps in the level of Chow groups defined in Section~\ref{sec:local}.

\begin{remark} We recall that if $\gamma$ and $\delta$ have the same sedentarity, using the isomorphism of the Chow groups with the tropical cohomology groups, the two maps have the following equivalent formulation. First, we get an inclusion $\comp\Sigma^\delta \hookrightarrow \comp \Sigma^\gamma$. To see this, note that $\delta$ gives a ray $\rho_{\delta/\gamma}$ in $\Sigma^\gamma$. The above embedding identifies $\comp\Sigma^\delta$ with the closure of the part of sedentarity $\rho_{\delta/\gamma}$ in $\comp \Sigma^\gamma$.
The restriction map $\i^*_{\gamma\ssubface\delta}$ is then the restriction map $H^k(\comp \Sigma^\gamma) \to H^k(\comp\Sigma^\delta)$ induced by this inclusion. The Gysin map is the dual of the restriction map $\i^*_{\gamma\ssubface\delta} \colon H^{2d-2\dims{\delta}-k}(\comp \Sigma^\gamma) \to H^{2d - 2\dims{\delta}-k}(\comp\Sigma^\delta)$, which by Poincar\'e duality for $\comp \Sigma^\delta$ and $\comp \Sigma^\gamma$, is a map from $H^k(\comp \Sigma^\delta)$ to $H^{k+2}(\comp \Sigma^\gamma)$.
\end{remark}

In addition to the above two maps, if two faces $\gamma$ and $\delta$ do not have the same sedentarity, the projection map $\pi_{\delta\ssupface\gamma} \colon N^{\sed(\delta)}_\R \to N^{\sed(\gamma)}_\R$
induces an isomorphism $\Sigma^\delta \simeq \Sigma^\gamma$, which leads to an isomorphism of the canonical compactifications, and therefor of the corresponding cohomology groups that we denote by
\[ \pi^*_{\gamma\ssubface\delta}\colon H^k(\gamma)\simto H^k(\delta). \]

\begin{remark} All the three above maps can be defined actually even if $\gamma$ is not of codimension 1 in $\delta$. Note that in this case, the Gysin map is a map $\gys_{\delta \supface \gamma} \colon H^k(\delta)\simto H^{k+ 2\dims{\delta} - 2\dims{\gamma}}(\gamma)$.
\end{remark}

Let now $S$ and $T$ be two collections of faces of $X$. We can naturally define the map
\[ \begin{array}{rccc}
i_{S,T}^*\colon & \bigoplus_{\gamma\in S}H^k(\gamma) & \longrightarrow & \bigoplus_{\delta\in T}H^k(\delta), \\
\end{array} \]
by sending, for any $\gamma\in S$, any element $x$ of $H^k(\gamma)$ to the element
\[ \i_{S,T}^* (x):=\sum_{\delta\in T \\ \delta\ssupface\gamma}\sign(\gamma, \delta)\i^*_{\gamma\ssubface\delta}(x), \]
and extend it by linearity to the direct sum $\bigoplus_{\gamma\in S} H^k(\gamma)$.

In the same way, we define
\[ \begin{array}{rccc}
\gys_{S,T}\colon & \bigoplus_{\delta\in S}H^k(\delta) &  \longrightarrow & \bigoplus_{\gamma \in T}H^{k+2}(\gamma), \\
\end{array} \]
by setting for $x \in H^k(\delta)$, for $\delta \in S$,
\[ \gys_{S,T}(x) := \sum_{\gamma \in T \\ \delta \ssupface \gamma} \sign(\gamma,\delta)\gys_{\delta \ssupface \gamma}(x). \]
Similarly, we get
\[ \begin{array}{rccc}
\pi^*_{S,T}\colon & \bigoplus_{\gamma\in S}H^k(\gamma) &  \longrightarrow & \bigoplus_{\delta\in T}H^k(\delta). \\
\end{array} \]

In what follows, we will often extend maps in this way. Moreover, we will drop the indices if they are clear from the context.

\subsection{Combinatorial Steenbrink and the main theorem}

Denote by $X_\f$ the set of faces of $X$ whose closures do not intersect the boundary of $\TP_{X_\sminfty}$, \ie, the set of compact faces of $X_\conezero$. Note that this was previously denoted by $X_{\conezero,\f}$ that we now simplify to $X_\f$ for the ease of presentation.

In what follows, inspired by the shape of the first page of the Steenbrink spectral sequence~\cite{Ste76}, we define bigraded groups $\ST_1^{a,b}$ with a collection of maps between them. In absence of a geometric framework reminiscent to the framework of degenerating families of smooth complex projective varieties, the heart of this section is devoted to the study of this combinatorial shadow in order to obtain interesting information about the geometry of $X$.

\medskip

For all pair of integers $a, b \in \mathbb Z$, we define
\[ \ST_1^{a,b} := \bigoplus_{s \geq \dims{a} \\ s \equiv a \pmod 2} \ST_1^{a,b,s} \]
where
\[ \ST_1^{a,b,s} = \bigoplus_{\delta \in X_\f \\ \dims\delta =s} H^{a+b-s}(\delta). \]

The bigraded groups $\ST_1^{a,b}$ come with a collection of maps
\[\i^{a,b\,*} \colon \ST^{a,b}_1 \to \ST_1^{a+1, b} \qquad \textrm{and} \qquad \gys^{a,b}\colon \ST^{a,b}_1 \to \ST_1^{a+1, b},\]
where both these maps are defined by our sign convention introduced in the previous section.

In practice, we drop the indices and denote simply by $\i^*$ and $\gys$ the corresponding maps.
\begin{remark} More precisely, we could view the collection of maps $\i^*$ and $\gys$ bi-indexed by $a,b$ as maps of bidegree $(1,0)$
\[\i^* = \bigoplus_{a,b}\i^{a,b\,*} \colon \bigoplus_{a,b} \ST_1^{a,b} \longrightarrow \bigoplus_{a,b} \ST_1^{a,b},\]
and
\[\gys = \bigoplus_{a,b}\gys^{a,b} \colon \bigoplus_{a,b} \ST_1^{a,b} \longrightarrow \bigoplus_{a,b} \ST_1^{a,b}. \qedhere \]
\end{remark}

\begin{prop}\label{prop:app1}
The two collections of maps $\i^*$ and $\gys$ have the following properties.
\begin{align*}
\i^* \circ \i^* =0, \qquad \gys \circ \gys =0, \qquad \i^* \circ \gys + \gys \circ \i^* =0.
\end{align*}
\end{prop}
\begin{proof}
Section~\ref{sec:proofapp1} is devoted to the proof of this proposition.
\end{proof}

\medskip

We now define the \emph{differential} $\d\colon \ST_1^{a,b} \to \ST_1^{a+1,b}$ as the sum $\d = \i^*+ \gys$. It follows from the properties given in the above proposition that we have the following.

\begin{prop}
For a unimodular triangulation of $\X$ and for any integer $b$, the differential $\d$ makes $\ST_1^{\bul,b}$ into a cochain complex.
\end{prop}

\begin{remark} This proposition suggests $\ST_1^{\bul, \bul}$ might be the first page of a spectral sequence converging to the cohomology of $X$. We could think that such a spectral sequence could be defined if we could \emph{deform infinitesimally the tropical variety}. In the sequel we will only use the property stated in the proposition, that the lines of $\ST_1^{\bul, \bul}$ with the differential $d$ form a cochain complex, as well as certain combinatorial properties of these cochain complexes.
\end{remark}

For a cochain complex $(C^\bul, \d)$, denote by $H^a(C^\bul, \d)$ its $a$-th cohomology group, i.e.,
\[H^a(C^\bul, \d) = \frac{\ker\Bigl(\d\colon \, C^{a} \rightarrow C^{a+1}\Bigr)}{\Im \Bigl(\d\colon \, C^{a-1} \rightarrow C^{a}\Bigr)}. \]

In this section we prove the following theorem, generalizing the main theorem of~\cite{IKMZ} from approximable setting (when $X$ arises as the tropicalization of a family of complex projective varieties) to any general smooth tropical variety.
\begin{thm}[Steenbrink-Tropical comparison theorem] \label{thm:steenbrink}
The cohomology of $(\ST_1^{\bul,b},\d)$ is described as follows. For $b$ odd, all the terms $\ST_1^{a,b}$ are zero, and the cohomology is vanishing. For $b$ even, let $b=2p$ for $p \in \mathbb Z$, then for $q\in \mathbb Z$, we have
\[H^{q-p}(\ST^{\bul,2p}_1, \d) = H^{p,q}_{\trop}(\X).\]
\end{thm}

The rest of this section is devoted to the proof of Theorem~\ref{thm:steenbrink}.

In absence of a geometric framework reminiscent to the framework of degenerating families of smooth complex projective varieties, the heart of our proof is devoted to introducing and developing certain tools, constructions, and results about the combinatorial structure of the Steenbrink spectral sequence.

\medskip

Broadly speaking, the proof which follows is inspired by the one given in the approximable case~\cite{IKMZ}, as well as by the ingredients in Deligne's construction of the mixed Hodge structure on the cohomology of a smooth algebraic variety~\cite{Deligne-Hodge2} and by Steenbrink's construction of the limit mixed Hodge structure on the special fiber of a semistable degeneration~\cite{Ste76}, both based on the use of the sheaf of logarithmic differentials and their corresponding spectral sequences.

In particular, we are going to prove first an analogous result in the tropical setting of the Deligne exact sequence which gives a resolution of the coefficient groups $\SF^p$ with cohomology groups of the canonically compactified fans $\comp \Sigma^\delta$. The proof here is based on the use of Poincar\'e duality and our Theorem~\ref{thm:HI}, the tropical analog of Feichtner-Yuzvinsky theorem~\cite{FY}, which provides a description of the tropical cohomology groups of the canonically compactified fans. We then define a natural filtration that we call the \emph{tropical weight filtration} on the coefficient groups $\SF^p(\cdot)$, inspired by the weight filtration on the sheaf of logarithmic differentials, and study the corresponding spectral sequence. The resolution of the coefficient groups given by the Deligne spectral resolution gives a double complex which allows to calculate the cohomology of the graded cochain complex associated to the weight filtration. We then show that the spectral sequence associated to the double complex corresponding to the weight filtration which abuts to the tropical cohomology groups, abuts as well to the cohomology groups of the Steenbrink cochain complex, which concludes the proof. The proof of this last result is based on the use of our Spectral Resolution Lemma~\ref{lem:spectral_resolution} which allows to make a bridge between different spectral sequences.

\subsection{Tropical Deligne resolution}

Let $\Sigma$ be a unimodular fan structure on the support of a Bergman fan. We follow the notations of the previous section. In particular, for each cone $\sigma$, we denote by $\Sigma^\sigma$ the induced fan on the star of $\sigma$ in $\Sigma$, and by $\comp \Sigma^\sigma$ its canonical compactification.

\begin{thm}[Tropical Deligne resolution]\label{thm:deligne}
We have the following exact sequence of $\mathbb Q$-vector spaces:
\[0 \rightarrow \SF^p(\conezero) \rightarrow \bigoplus_{\sigma \in \Sigma \\
\dims{\sigma} =p} H^0(\sigma) \rightarrow \bigoplus_{\sigma \in \Sigma \\
\dims{\sigma} =p-1} H^2(\sigma) \rightarrow \dots \rightarrow \bigoplus_{\sigma \in \Sigma \\
\dims{\sigma} =1} H^{2p-2}(\sigma) \rightarrow H^{2p} (\conezero) \to 0, \]
where the maps between cohomology groups are given by $\gys$.
\end{thm}

\begin{remark} We should justify the name given to the theorem. In the realizable case when $\Sigma$ has the same support as the Bergman fan $\Sigma_\Ma$ of a matroid $\Ma$ realizable over the field of complex numbers, this can be obtained from the Deligne spectral sequence which describes the mixed Hodge structure on the cohomology of a complement of complex hyperplane arrangement. This is done in~\cite{IKMZ}. Our theorem above states this is more general and holds for any Bergman fan.
\end{remark}
\begin{remark} The theorem suggests this should be regarded as a cohomological version of the inclusion-exclusion principle. The cohomology groups are described in terms of the coefficient sheaf and the coefficient sheaf can be recovered from the cohomology groups. It would be certainly interesting to see whether a cohomological version of the M\"obius inversion formula could exist.
\end{remark}

The rest of this section is devoted to the proof of this theorem. First using the Poincar\'e duality for canonical compactifications $\comp \Sigma^\sigma$, it will be enough to prove the exactness of the following complex for each $k$ (here $k=d-p$):

\begin{equation}
0 \rightarrow H^{2k}(\conezero) \to \bigoplus_{\sigma \in \Sigma \\
\dims{\sigma} =1} H^{2k}(\sigma) \to \bigoplus_{\sigma \in \Sigma \\
\dims{\sigma} =2} H^{2k}(\sigma) \to \dots \to \bigoplus_{\sigma \in \Sigma \\
\dims{\sigma} =d-k} H^{2k}(\sigma) \to \SF_{d-k}(\conezero) \to 0.
\end{equation}

Recall from Section~\ref{sec:tropvar} that for any tropical variety $Z$, we denote by $\Omega^k_Z$ the sheaf of tropical holomorphic $k$-forms on $Z$, defined as the kernel of the second differential operator from Dolbeault $(k,0)$-forms to Dolbeault $(k,1)$-forms on $Z$. Since we are going to use some results from~\cite{JSS}, we note that this sheaf is denoted $\mathcal L^k_Z$ in \emph{loc. cit.}

\medskip

We will need the following alternative characterization of this sheaf, which also shows that it can be defined over $\R$ as the sheafification of the combinatorial sheaf $\SF^k$. Let $Z$ be a tropical variety with an extended polyhedral structure as in Section~\ref{sec:tropvar}. If $U$ is an open set of $Z$, we say that $U$ is \emph{nice} if either $U$ is empty, or there exists a face $\gamma$ of $Z$ intersecting $U$ such that, for each face $\delta$ of $Z$, every connected components of $U\cap\delta$ contains $U\cap\gamma$. (Compare with the \emph{basic open sets} from~\cite{JSS}.) This condition implies that $U$ is connected, and for each face $\delta\in Z$ which intersects $U$, $\gamma$ is a face of $\delta$ and $\delta\cap U$ is connected. We call $\gamma$ the \emph{minimum face of $U$}. Nice open sets form a basis of open sets on $Z$.

The sheaf $\Omega^k_Z$ is then the unique sheaf on $Z$ such that, for each nice open set $U$ of $Z$ with minimum face $\gamma$, we have
\[ \Omega^k_Z(U)=\SF^k(\gamma). \]

\medskip

Suppose now $Z \subseteq \comp Z$ is a compactification of $Z$. We denote by $\Omega^k_{Z,c}$ the \emph{sheaf of holomorphic $k$-forms on $\comp Z$ with compact support in $Z$} defined on connected open sets by
\[ \Omega^k_{Z,c}(U):=\begin{cases}
\Omega^k_Z(U) & \text{if $U\subseteq Z$} \\
0             & \text{otherwise.}
\end{cases} \]

\begin{remark} \label{rem:direct_image} To justify the name of the sheaf $\Omega^k_{Z,c}$, denote by $\iota\colon Z\hookrightarrow\comp Z$ the inclusion. Then we have
\[ \Omega^k_{Z,c}=\iota_!\Omega^k_Z. \]
In other words, this is the \emph{direct image with compact support} of the sheaf $\Omega^k_Z$. In particular, in the case we study here, the cohomology with compact support of $\Omega^k_Z$ is computed by the usual sheaf cohomology of $\Omega^k_{Z,c}$, \ie, we have
\[ H_c^\bul(\comp Z, \Omega^k_Z)=H^\bul(\comp Z, \Omega^k_{Z,c}). \qedhere \]
\end{remark}

We now specify the above set-up for $Z = \Sigma$ and $\comp \Sigma$ the canonical compactification of $\Sigma$.

\medskip

For a cone $\sigma \in \Sigma$, the canonical compactification $\comp\Sigma^\sigma$ naturally leaves at infinity of $\comp \Sigma$. Using the notations of Section~\ref{sec:tropvar}, the fan $\Sigma_\infty^\sigma$, which is based at the point $\infty_\sigma$ of $\comp \Sigma$, is naturally isomorphic to $\Sigma^\sigma$. Via this isomorphism, the closure $\comp{\Sigma_\infty^\sigma}$ of ${\Sigma_\infty^\sigma}$ in $\comp \Sigma$ coincides with the canonical compactification $\comp \Sigma^\sigma$ of $\Sigma^\sigma$. In the remaining of this section, to simplify the notations, we identify $\comp \Sigma^\sigma$ as $\comp{\Sigma_\infty^\sigma}$ living in $\comp \Sigma$.

\medskip

To each $\sigma \in \Sigma$ corresponds the sheaf $\Omega^k_{\comp \Sigma^\sigma}$ of holomorphic $k$-forms on $\comp \Sigma^\sigma$ which by extension by zero leads to a sheaf on $\comp \Sigma$. We denote this sheaf by $\Omega^k_\sigma$. The following proposition describes the cohomology of these sheaves.

\begin{prop} \label{prop:cohomology_Omega_fan} Notations as above, for each pair of non-negative integers $m,k$,
we have
\[ H^m(\comp \Sigma, \Omega^k_\sigma) = \begin{cases} H^{k,k}(\sigma) = H^{2k}(\sigma) & \textrm{if $m = k$} \\
0  & \textrm{otherwise}.
\end{cases}\]
\end{prop}
\begin{proof}
We have
\[H^m(\comp \Sigma, \Omega^k_\sigma) \simeq H^m(\comp \Sigma^\sigma, \Omega^k_{\comp \Sigma^\sigma}) \simeq H^{k,m}_{\textrm{Dolb}}(\comp\Sigma^\sigma) \simeq H^{k,m}_{\trop}(\comp\Sigma^\sigma),\]
by the comparison Theorem~\ref{thm:comparison}, proved in~\cite{JSS}, and the result follows from Theorem~\ref{thm:HI}, proved in~\cite{AP}.
\end{proof}

For a pair of faces $\tau\subface\sigma$ in $\Sigma$, we get natural inclusion maps $\comp \Sigma^\sigma \hookrightarrow \comp\Sigma^\tau \hookrightarrow \comp \Sigma$, leading to natural restriction maps of sheaves $\i^*_{\tau\subface\sigma}\colon\Omega^k_\tau \to\Omega^k_\sigma$ on $\comp \Sigma$. Here the map $\i_{\tau\subface\sigma}$ ($= \i_{\sigma \supface\tau}$, using our convention from Section~\ref{sec:intro}) denotes the inclusion $\comp\Sigma^\sigma\hookrightarrow\comp\Sigma^\tau$.

\medskip

We consider now the following complex of sheaves on $\comp \Sigma$:
\begin{equation}
\Omega^k_\bul\colon \qquad \Omega_\conezero^k \to \bigoplus_{\sigma \in \Sigma \\ \dims{\sigma} =1} \Omega^k_{\sigma} \to \dots \to \bigoplus_{\sigma \in \Sigma \\ \dims{\sigma} =d-k} \Omega^k_{\sigma}
\end{equation}
concentrated in degrees $0, 1, \dots, d-k$, given by the dimension of the cones $\sigma$ in $\Sigma$, whose boundary maps are given by
\[ \alpha\in\Omega_\sigma^k \mapsto \d\alpha = \sum_{\zeta\ssupface\sigma}\sign(\sigma,\zeta)\i^*_{\sigma\ssubface\zeta}(\alpha). \]
We will derive Theorem~\ref{thm:deligne} by looking at the hypercohomology groups $\hyp^\bul(\comp\Sigma, \Omega^k_\bul)$ of this complex and by using the following proposition.

\begin{prop} \label{prop:exactness_Omega}
The following sequence of sheaves is exact
\[ 0 \to \Omega_{\Sigma, c}^k \to \Omega_\conezero^k \to \bigoplus_{\varrho \in \Sigma \\ \dims{\varrho}=1}\Omega^k_\varrho \to \bigoplus_{\sigma \in \Sigma \\ \dims{\sigma}=2}\Omega^k_\sigma \to \dots \to \bigoplus_{\sigma \in \Sigma \\ \dims{\sigma}=d-k}\Omega^k_\sigma \to 0. \]
\end{prop}

\begin{proof} It will be enough to prove that taking sections over nice open sets $U$ give exact sequences of $\R$-vector spaces.

If $U$ is included in $\Sigma$, clearly by definition we have $\Omega^k_{\Sigma,c}(U)\simeq\Omega^k_{\conezero}(U)$, and the other sheaves of the sequence have no nontrivial section over $U$. Thus, the sequence is exact over $U$.

It remains to prove that, for every other nice open set $U$ intersecting $\comp\Sigma\setminus\Sigma$, the sequence
\[ 0 \to \Omega_\conezero^k(U) \to \bigoplus_{\varrho \in \Sigma \\ \dims{\varrho}=1}\Omega^k_\varrho(U) \to \bigoplus_{\sigma \in \Sigma \\ \dims{\sigma}=2}\Omega^k_\sigma(U) \to \dots \to \bigoplus_{\sigma \in \Sigma \\ \dims{\sigma}=d-k}\Omega^k_\sigma(U) \to 0 \]
is exact. Let $\gamma\in\comp\Sigma$ be the minimum face of $U$. Let $\sigma\in\Sigma$ be the sedentarity of $\gamma$. The closed strata of $\comp \Sigma$ which intersect $U$ are exactly those of the form $\comp{\Sigma_\infty^\tau} \simeq \comp\Sigma^{\tau}$ with $\tau\subface\sigma$. Moreover, if $\tau$ is a face of $\sigma$,
\[ \Omega^k_\tau(U)=\SF^k(\gamma). \]
Thus, the previous sequence can be rewritten in the form
\[ 0 \to \SF^k(\gamma) \to \bigoplus_{\sigma\subface\sigma \\ \dims{\sigma}=1}\SF^k(\gamma) \to \bigoplus_{\tau\subface\sigma \\ \dims{\tau}=2}\SF^k(\gamma) \to \dots \to \bigoplus_{\tau\subface\sigma \\ \dims{\tau}=\dims{\sigma}}\SF^k(\gamma) \to 0. \]
This is just the cochain complex of the simplicial cohomology (for the natural simplicial structure induced by the faces) of the cone $\sigma$ with coefficients in the group $\SF^k(\gamma)=\SF_k(\gamma)^\dual$. This itself corresponds to the reduced simplicial cohomology of a simplex shifted by $1$. This last cohomology is trivial, thus the sequence is exact. That concludes the proof of the proposition.
\end{proof}

\begin{proof}[Proof of Theorem~\ref{thm:deligne}]
By~\cite{JSS}, we have
\[ H^m(\comp\Sigma,\Omega^k_{\Sigma,c}) = H^m_c(\Sigma,\Omega^k_\Sigma) = H^{k,m}_{\trop,c}(\Sigma)=
\begin{cases} \SF^{d-k}(\conezero)^\dual = \SF_{d-k}(\conezero) & \text{if $m=d$}\\
0 & \text{otherwise}.
\end{cases} \]
By proposition \ref{prop:exactness_Omega}, the cohomology of $\Omega^k_{\Sigma,c}$ becomes isomorphic to the hypercohomology $\hyp(\Sigma, \Omega^k_\bul)$. Thus, we get
\[ \hyp(\Sigma, \Omega^k_\bul)\simeq \SF_{d-k}(\conezero)[d], \]
meaning
\[\hyp^m(\comp\Sigma, \Omega^k_\bul) = \begin{cases}   \SF_{d-k}(\conezero) & \textrm{ for $m =d$}\\
0 & \textrm{ otherwise.}
\end{cases}\]

On the other hand, using the hypercohomology spectral sequence, combined with Proposition \ref{prop:cohomology_Omega_fan}, we infer that the hypercohomology of $\Omega^k_\bul$ is given by the cohomology of the following complex:
\[ 0 \rightarrow H^{2k}(\conezero)[k] \to \bigoplus_{\sigma \in \Sigma \\
\dims{\sigma} =1} H^{2k}(\sigma)[k+1] \to \bigoplus_{\sigma \in \Sigma \\
\dims{\sigma} =2} H^{2k}(\sigma)[k+2] \to \dots \to \bigoplus_{\sigma \in \Sigma \\
\dims{\sigma} =d-k} H^{2k}(\sigma)[d] \to 0. \]
We thus conclude the exactness of the sequence
\begin{equation}
0 \rightarrow H^{2k}(\conezero) \to \bigoplus_{\sigma \in \Sigma \\
\dims{\sigma} =1} H^{2k}(\sigma) \to \bigoplus_{\sigma \in \Sigma \\
\dims{\sigma} =2} H^{2k}(\sigma) \to \dots \to \bigoplus_{\sigma \in \Sigma \\
\dims{\sigma} =d-k} H^{2k}(\sigma) \to \SF_{d-k}(\conezero) \to 0,
\end{equation}
and the theorem follows.
\end{proof}

\subsection{Weight filtration on $\SF^p$}

In this section, we introduce a natural filtration on the structure sheaf $\SF^p$. This filtration could be regarded as the analog in polyhedral geometry of the weight filtration on the sheaf $\Omega^p(\log \mathfrak X_0)$ of logarithmic differentials around the special fiber $\mathfrak X_0$ of a semistable degeneration $\mathfrak X$ over the unit disk.

\medskip

Let $\delta$ be a face of the triangulation $X$ of the tropical variety $\X$. First we define a filtration denoted by $\~ W_\bul$ on $\SF^p(\delta)$, and then define the weight filtration $W_\bul$ by shifting the filtration induced by $\~ W_\bul$.

\medskip

The tangent space $\TT\delta$ of $\delta$ is naturally included in $\SF_1(\delta)$. More generally, for any integer $s$, we have an inclusion $\bigwedge^s\TT\delta \subseteq \SF_s(\delta)$, which in turn gives an inclusion
\[\bigwedge^s\TT\delta  \wedge \SF_{p-s}(\delta) \subseteq \SF_p(\delta)
\]
for any non-negative integer $s$.

\begin{defi} For any integer $s\geq 0$, define $\~ W_s \SF^p(\delta)$ as follows
\begin{align*}
\~W_s \SF^p(\delta) &:= \\
  & \hspace{-2em}\Bigl\{ \, \alpha \in \SF^p(\delta) \st \textrm{restriction of $\alpha$ on subspace } \bigwedge^{s+1}\TT\delta  \wedge \SF_{p-s-1}(\delta) \textrm{ of } \SF_p(\delta) \textrm{ is trivial}\, \Bigr\}. \qedhere
\end{align*}
\end{defi}

Recall that for a pair of faces $\gamma \ssubface \delta$, we denote by $\i^*_{\gamma\ssubface\delta}$ the restriction map $\SF^p(\gamma) \to \SF^p(\delta)$.

\begin{prop} Let $\gamma$ be a face of codimension one of $\delta$. We have the following two cases.
\begin{itemize}
\item If $\gamma$ and $\delta$ have the same sedentarity, then the natural map $\i^*_{\gamma\ssubface\delta}\colon\SF^p(\gamma) \to \SF^p(\delta)$ sends $\~ W_{s-1}\SF^p(\gamma)$ into $\~ W_s\SF^p(\delta)$.
\item If $\gamma$ and $\delta$ do not have the same sedentarity, then the natural map $\pi^*_{\gamma\ssubface\delta}\colon\SF^p(\gamma) \to \SF^p(\delta)$ sends $\~ W_s\SF^p(\gamma)$ into $\~ W_s\SF^p(\delta)$.
\end{itemize}
\end{prop}
\begin{proof}
For the first point, let $\nvect_{\delta/\gamma}$ be a primitive normal vector to $\TT\gamma$ in $\TT\delta$. The map $\SF_p(\delta) \to \SF_p(\gamma)$ restricted to the subspace $\bigwedge^{s+1} \TT\delta  \wedge \SF_{p-s-1}(\delta)$ can be decomposed as the composition of the following maps
\begin{align*}
\bigwedge^{s+1} \TT\delta & \wedge \SF_{p-s-1}(\delta) = \bigwedge^{s} \TT\gamma \wedge \nvect_{\delta/\gamma} \wedge \SF_{p-s-1}(\delta) + \bigwedge^{s+1} \TT\gamma \wedge \SF_{p-s-1}(\delta) \\
&\longrightarrow \bigwedge^{s} \TT\gamma \wedge \nvect_{\delta/\gamma} \wedge \SF_{p-s-1}(\gamma) + \bigwedge^{s+1} \TT\gamma \wedge \SF_{p-s-1}(\gamma)\hooklongrightarrow \bigwedge^{s} \TT\gamma  \wedge \SF_{p-s}(\gamma).
\end{align*}

An element of $\~W_{s-1}\SF^p(\gamma)$ is thus sent to an element of $\SF^p(\delta)$ which restricts to zero on the subspace $\bigwedge^{s+1} \TT\delta \wedge \SF_{p-s-1}(\delta)$, in other words, to an element of $\~ W_s\SF^p(\delta)$.

\medskip

For the second point, the projection $\pi_{\delta\ssupface\gamma}$ maps $\TT\delta$ onto $\TT\gamma$. Thus, it also maps $\bigwedge^{s+1}\TT\delta$ onto $\bigwedge^{s+1}\TT\gamma$. The image of an element of $\~W_s\SF^p(\gamma)$ thus restricts to zero on $\bigwedge^s\TT\delta$, \ie, it is an element of $\~W_s\SF^p(\delta)$.
\end{proof}

We now consider the graded pieces $\gr_{s}^{\~ \ws} \SF^p(\delta) := \rquot {\~ W_s \SF^p(\delta)} {\~ W_{s-1}\SF^p(\delta)}$.

\medskip

For $\delta \in X$, we denote as before by $\conezero^\delta$ the zero cone of the fan $\Sigma^\delta$ induced by the triangulation $X$ around $\delta$. The following proposition gives the description of the graded pieces of the filtration.

\begin{prop}\label{prop:grading1}
For each face $\delta$, we have
\[ \gr_s^{\~ \ws}\SF^p(\delta) \simeq \bigwedge^s \TT^\dual\delta \otimes \SF^{p-s}(\conezero^\delta). \]
Here $\TT^\dual\delta$ is the cotangent space of the face $\delta$.
\end{prop}

In order to prove this proposition, we introduce some useful maps. Denote by $\pi_\delta$ the natural projection $N^{\sed(\delta)}_\R \to N^\delta_\R$. We choose a projection $\p_\delta\colon N^{\sed(\delta)}_\R \to \TT\delta$. Both projections naturally extend to exterior algebras, and for any integer $p$, we get maps
\begin{align*}
\p_\delta\colon\,\,& \bigwedge^p N^{\sed(\delta)}_\R \to \bigwedge^p \TT\delta, \\
\pi_\delta\colon\,\,& \bigwedge^p N^{\sed(\delta)}_\R \to \bigwedge^p N^\delta_\R.
\end{align*}
Furthermore, we get a map
\[ \pi_\delta\colon \SF_p(\delta) \to \SF_p(\conezero^\delta). \]

\medskip

We have the corresponding pullback $\p_\delta^*$ and $\pi_\delta^*$ between the corresponding dual spaces $\SF^p(\delta)$ and $\SF(\conezero^\delta)$. Moreover, if $\alpha\in\SF^p(\delta)$ is zero on $\TT\delta\wedge\SF_{p-1}(\delta)$, since $\TT\delta=\ker(\pi_\delta)$, we can define a natural pushforward $\pi_{\delta\,*}(\alpha)$ of $\alpha$ in $\SF_p(\conezero^\delta)$. This leads to a map (in fact an isomorphism)
\[ \pi_{\delta\,*}\colon \~W_0(\SF^p(\delta)) \to \SF^p(\conezero^\delta). \]

\begin{proof}[Proof of Proposition~\ref{prop:grading1}]
Let $\alpha$ be an element of $\~W_s\SF^p(\delta)$. Let $\u\in \bigwedge^s\TT\delta$. Recall that the contraction of $\alpha$ by $\u$ is the multiform $\beta:=\alpha(\u\wedge {}\cdot{})\in \SF^{p-s}(\delta)$.

Since $\alpha$ is in $\~W_s\SF^p(\delta)$, and $ \u \wedge \TT\delta\wedge \SF_{p-s-1}(\delta) \subseteq \bigwedge^{s+1} \TT\delta\wedge \SF_{p-s-1}(\delta)$, by definition of the filtration $\~W_\bul$, the contracted multiform $\beta$ is zero on $\TT\delta\wedge \SF_{p-s-1}(\delta)$. Thus, $\beta\in\~W_0(\delta)$. Hence, we get a morphism
\[ \begin{array}{rrcl}
\Psi\colon \~W_s\SF^p(\delta) & \longrightarrow & \Hom\Bigl(\,\bigwedge^s\TT\delta\,,\, \SF^{p-s}(\conezero^\delta)\,\Bigr) &\simeq \,\,\bigwedge^s \TT^\dual\delta \otimes \SF^{p-s}(\conezero^\delta), \\[.4em]
\alpha  & \longmapsto &\quad\bigl(\,  \u \longmapsto \pi_{\delta\,*}(\alpha(\u \wedge {}\cdot{}))  \, \bigr).
\end{array} \]
Notice that the kernel of $\Psi$ is $\~W_{s-1}(\SF^p(\delta))$, \ie, its cokernel is $\gr^{\~W}_s(\SF^p(\delta))$.

Set
\[ \begin{array}{rccl}
\Phi\colon& \bigwedge^s \TT^\dual\delta \otimes \SF^{p-s}(\conezero^\delta), & \longrightarrow & \SF^p(\delta) \\[.4em]
& x \otimes y  & \longmapsto & \p_\delta^*(x)\wedge\pi_\delta^*(y),
\end{array} \]
and extended by linearity.

\medskip

One verifies directly that $\Im(\Psi)\subseteq\~W_s\SF^p(\delta)$, and $\Psi\circ\Phi$ is identity on $\bigwedge^s \TT^\dual\delta \otimes \SF^{p-s}(\conezero^\delta)$. In particular, $\Psi$ is surjective. We thus infer that $\Psi$ induces an isomorphism between its cokernel and $\bigwedge^s \TT^\dual\delta \otimes \SF^{p-s}(\conezero^\delta)$. This concludes the proof as we have seen that its cokernel of $\Psi$ is $\gr^{\~W}_s(\SF^p(\delta))$.
\end{proof}

\medskip

\begin{defi}[Weight filtration]
For each $\delta \in X$, we define the \emph{weight filtration} $W_\bul$ on $\SF^p(\delta)$ and its \emph{opposite filtration} $W^\bul$ by
\[ W_s \SF^p(\delta) = \~W_{s+\dims\delta}\SF^p(\delta), \textrm{ and }\]
\[W^s \SF^p(\delta) = \~ W_{-s+\dims{\delta}}\SF^p(\delta),\]
for each integer $s$.
\end{defi}

Proposition~\ref{prop:grading1} directly translates into the following facts about the filtration $W^\bul$.

\begin{cor} \label{cor:grading} The following holds.
\begin{enumerate}
\item For each face $\delta$, we have
\[ \gr^s_{\ws} \SF^p(\delta) \simeq \bigwedge^{\dims{\delta}-s} \TT^\dual\delta \otimes \SF^{p+s-\dims{\delta}}(\conezero^\delta).\]
\item For inclusion of faces $\gamma\subface\delta$, the map $\SF^p(\gamma) \to \SF^p(\delta)$ respects the filtration $W^\bul$. In particular, we get an application at each graded piece $\gr^s_{\ws}\SF^p(\gamma) \to \gr^s_{\ws} \SF^p(\delta)$.
\end{enumerate}
\end{cor}

\subsubsection{Description of the restriction maps on graded pieces of the weight filtration}
We would now like to explicitly describe the map induced by $\i^*$ on the level of graded pieces. In what follows, we identify $\gr^s_{\ws} \SF^p(\delta)$ with the decomposition $\bigwedge^{\dims{\delta}-s} \TT^\dual\delta \otimes \SF^{p+s-\dims{\delta}}(\conezero^\delta)$ and also, by duality, with $\Hom\Bigl(\,\bigwedge^{\dims\delta-s}\TT\delta\,,\, \SF^{p+s-\dims\delta}(\conezero^\delta)\,\Bigr)$. We call $\bigwedge^{\dims{\delta}-s} \TT^\dual\delta$ and $\SF^{p+s-\dims{\delta}}(\conezero^\delta)$ the \emph{parallel part} and the \emph{transversal part} of the graded piece of the filtration, respectively. We use the notation $\alpha_\parr, \beta_\parr$, etc. when referring to the elements of the parallel part, and use $\alpha_\perp, \beta_\perp$, etc. when referring to the elements of the transversal part. In particular, each element of the graded piece is a sum of elements of the form $\alpha_\parr\otimes \alpha_\perp$.

\medskip

We denote by $\Psi_{s,\delta}\colon W_s\SF^p(\delta)\to\gr^s_{\ws}\SF^p(\delta)$ the projection. As we have seen in the proof of Proposition \ref{prop:grading1}, this map is explicitly described by
\[ \begin{array}{rrl}
\Psi_{s,\delta}\colon W^s\SF^p(\delta) & \longrightarrow & \gr_{\ws}^s\SF^p(\delta) \\[.4em]
\alpha  & \longmapsto & \bigl(\,  \u \longmapsto \pi_{\delta\,*}(\alpha(\u \wedge {}\cdot{}))  \, \bigr).
\end{array} \]
We also have a section $\Phi_{s,\delta}$, which this time depends on the chosen projection $\p_\delta$, and which is given by
\[ \begin{array}{rrcl}
\Phi_{s,\delta}\colon& \gr_{\ws}^s\SF^p(\delta) & \longrightarrow & W^s\SF^p(\delta) \\[.4em]
& \alpha_\parr \otimes \alpha_\perp  & \longmapsto & \p_\delta^*(\alpha_\parr)\wedge\pi_\delta^*(\alpha_\perp).
\end{array} \]

\medskip

Consider now two faces $\gamma\ssubface\delta$ in $X$. We define two maps $\i^*_\parr$ and $\i^*_\perp$ between parallel and transversal parts of the graded pieces as follows. For each non-negative integer $t$, the map
\[\i^*_\parr \colon \bigwedge^{t} \TT^\dual\gamma \longrightarrow   \bigwedge^{t+1} \TT^\dual\delta \]
sends an element $\alpha_\parr \in  \bigwedge^{t} \TT^\dual\gamma$ to the element $\beta_\parr =  \i^*_\parr(\alpha_\parr)$ of $\bigwedge^{t+1} \TT^\dual\delta$ defined as follows. The multiform $\beta_\parr$ is the unique element of $\bigwedge^{t+1} \TT^\dual\delta$ which restricts to zero on $\bigwedge^{t+1}\TT\gamma$ under the inclusion map $\bigwedge^{t+1} \TT \gamma \hookrightarrow \bigwedge^{t+1} \TT\delta$, and which verifies
\[ \beta_\parr(\u\wedge \nvect_{\delta/\gamma})=\alpha_\parr(\u) \quad\text{for any $\u\in\bigwedge^{t}\TT\gamma\subset\bigwedge^{t}\TT\delta$}. \]
Here $\nvect_{\delta/\gamma}$ is any primitive normal vector to $\TT\gamma$ in $\TT\delta$. In other words,
\[ \beta_\parr=\p_\gamma^*(\alpha_\parr)\rest{\TT\delta} \wedge \nvect_{\delta/\gamma}^\dual, \]
where $\nvect_{\delta/\gamma}^\dual$ is the form on $\TT\delta$ which vanishes on $\TT\gamma$ and which takes value 1 on $\nvect_{\delta/\gamma}$.

For any positive integer $t$, the other map
\[\i^*_\perp\colon \SF^{t}(\conezero^\gamma) \longrightarrow \SF^{t-1}(\conezero^\delta)\]
between transversal parts is similarly defined as follows. This is the unique linear map which sends an element $\alpha_\perp$ of $\SF^{t}(\conezero^\gamma)$ to the element $\beta_\perp = \i^*_\perp(\alpha_\perp)$ which verifies
\[ \alpha_\perp(\e_{\delta/\gamma}\wedge\v)=\beta_\perp(\pi_{\gamma\subface\delta}(\v)) \quad\text{for any $\v\in \SF_{t-1}(\conezero^\gamma)$}. \]
Here, as in the previous sections, $\e_{\delta/\gamma}$ is the primitive vector of the ray $\rho_{\delta/\gamma}$ corresponding to $\delta$ in $\Sigma^\gamma$, and $\pi_{\gamma\subface\delta}$ is the projection from $\comp\Sigma^\gamma$ to $\comp\Sigma^\delta$ (along $\R \e_{\delta/\gamma}$), which naturally induces a surjective map from $\SF_{t-1}(\conezero^\gamma)\to \SF_{t-1}(\conezero^\delta)$. In other words, $\beta_\perp$ is the pushforward by $\pi_{\gamma\subface\delta}$ of the contraction of $\alpha_\perp$ by $\e_{\delta/\gamma}$ (which is well-defined because $\ker{\pi_{\gamma\subface\delta}}=\R \e_{\delta/\gamma}$).

\begin{prop}\label{prop:grading} Let $\gamma\ssubface\delta$ be two faces in $X$. Notations as above, the induced map on graded pieces $\gr^s_{\ws} \SF^p(\gamma) \to \gr^s_{\ws} \SF^p(\delta)$ is zero provided that $\sed(\gamma)\neq\sed(\delta)$. Otherwise, it coincides with the map
\[ \begin{array}{rrclcrcl}
i^*_\parr \otimes \i^*_\perp\colon & \bigwedge^{\dims{\gamma}-s} \TT^\dual\gamma &\otimes& \SF^{p+s-\dims{\gamma}}(\conezero^\gamma) & \longrightarrow & \bigwedge^{\dims{\delta}-s} \TT^\dual\delta &\otimes& \SF^{p+s - \dims{\delta}}(\conezero^\delta), \\
& \alpha_\parr &\otimes& \alpha_\perp & \longmapsto & \sign(\gamma,\delta) \, \beta_\parr &\otimes& \beta_\perp,
\end{array} \]
where $\beta_\parr = \i^*_\parr(\alpha_\parr)$ and $\beta_\perp =\i^*_\perp(\alpha_\perp)$.
\end{prop}

\begin{proof}
If $\sed(\delta)\neq\sed(\gamma)$, then the natural map $\SF^p(\gamma)\to\SF^p(\delta)$ is $\pi_{\gamma\ssubface\delta}^*$, and it sends $W^s\SF^p(\gamma)$ onto $W^{s+1}\SF^p(\delta)$. Thus, it induces the zero map on the graded pieces.

\medskip

Assume $\sed(\delta)=\sed(\gamma)$. Let $\alpha=\alpha_\parr\otimes\alpha_\perp\in \gr_{\ws}^s\SF^p(\gamma)$. By definition, the map induced by $\i^*$ on the graded pieces maps $\alpha$ onto $\Psi_{s+1,\delta}\circ\i^*_{\gamma\ssubface\delta}\circ\Phi_{s,\gamma}(\alpha)$. To understand this image, take an element $\u_\parr\otimes\u_\perp$ in $\bigwedge^{\dims\delta-s}\TT\delta\otimes \SF_{p+s-\dims\delta}(\conezero^\delta)$. By definition of $\Psi_{s+1,\delta}$ we get
\[ \Psi_{s+1,\delta}\circ \i^*_{\gamma\ssubface\delta}\circ \Phi_{s,\gamma}(\alpha)\,\bigl(\u_\parr\otimes\u_\perp\bigr) = \i^*_{\gamma\ssubface\delta}\circ \Phi_{s,\gamma}(\alpha)\,\bigl(\u_\parr\wedge\~\u_\perp\bigr) , \]
where $\~\u_\perp$ is any preimage of $\u_\perp$ by $\Phi_{s,\gamma}$. Then,
\begin{align*}
\i^*_{\gamma\ssubface\delta}\circ \Phi_{s,\gamma}(\alpha)\,\bigl(\u_\parr\wedge\~\u_\perp\bigr)
  &= \Phi_{s,\gamma}(\alpha)\,\bigl(\u_\parr\wedge\~\u_\perp\bigr) \\
  &= \bigl(\p^*_\gamma(\alpha_\parr)\wedge\pi^*_\gamma(\alpha_\perp)\bigr)\bigl(\u_\parr\wedge\~\u_\perp\bigr).
\end{align*}
Here, $\u_\parr \in \bigwedge^{\dims\delta-s} \TT\delta$, and we have the decomposition
\[ \bigwedge^{\dims\delta-s} \TT\delta = \bigwedge^{\dims\gamma-s+1}\TT \gamma \oplus \bigwedge^{\dims\gamma-s} \TT\gamma \wedge\nvect_{\delta, \gamma}, \]
where $\nvect_{\delta/\gamma}$ is a primitive normal vector to $\TT\gamma$ in $\TT\delta$.
Thus, it suffices to study the following two cases.
\begin{itemize}
\item Assume that $\u_\parr=\u'_\parr\wedge\nvect_{\delta/\gamma}$ with $\u'_\parr\in\bigwedge^{\dims\gamma-s}\TT\gamma$. Since $\TT\gamma=\ker(\pi_\gamma)$, we get
\begin{align*}
\bigl(\p^*_\gamma(\alpha_\parr)\wedge\pi^*_\gamma(\alpha_\perp)\bigr)\bigl(\u'_\parr\wedge\nvect_{\delta/\gamma}\wedge\~\u_\perp\bigr)
  &= \alpha_\parr(\p_\gamma(\u_\parr)) \cdot \alpha_\perp(\pi_\gamma(\nvect_{\delta/\gamma} \wedge \~\u_\perp)) \\
  &= \alpha_\parr(\u'_\parr) \cdot \alpha_\perp(\e_{\delta/\gamma} \wedge \pi_\gamma(\~\u_\perp)) \\
  &= \beta_\parr(\u_\parr)\cdot\beta_\perp(\u_\perp),
\end{align*}
since $\u_\parr=\u'_\parr\wedge\nvect_{\delta/\gamma}$ and $\u_\perp=\pi_{\gamma\subface\delta}(\pi_\gamma(\~\u_\perp))$.
\item If $\u_\parr\in\bigwedge^{\dims\gamma-s+1} \TT\gamma$, since $\p_\gamma^*(\alpha_\parr)\in\bigwedge^{\dims\gamma-s} N^{\sed(\gamma)}_\R$, a part of $\u_\parr$ must be evaluated by $\pi_\gamma^*(\alpha_\perp)$. But this evaluation will be zero since $\TT_\gamma=\ker(\pi_\gamma)$. Thus,
\begin{align*}
0
  &= \bigl(\p^*_\gamma(\alpha_\parr)\wedge\pi^*_\gamma(\alpha_\perp)\bigr)(\u'_\parr\wedge\nvect_{\delta/\gamma}\wedge\~\u_\perp) \\
  &= \beta_\parr(\u_\parr)\cdot\beta_\perp(\u_\perp).
\end{align*}
\end{itemize}
In any case, the statement of the proposition holds.
\end{proof}

Using the above proposition, we identify in the sequel the graded pieces of the opposite weight filtration and the maps between them by $\bigwedge^{\bul} \TT^\dual\gamma \otimes \SF^{\bul}(\conezero^\gamma)$ and maps between them.

\subsection{Yoga of spectral sequences and proof of the main theorem} \label{sec:steenbrinkdoublecomplex}

In this section, we study the combinatorics behind the Steenbrink spectral sequence, and use this to prove Theorem~\ref{thm:steenbrink}.
We introduce two spectral sequences. One of these spectral sequences computes the tropical cohomology $H_\trop^{p,\bul}$, and the another one computes the cohomology of the $2p$-th row of the tropical Steenbrink spectral sequence. We show that the two spectral sequences are isomorphic in page one. Then using a \emph{spectral resolution lemma} we generalize the isomorphism in page one to all pages, from which we will deduce that these two cohomologies coincide as stated by Theorem \ref{thm:steenbrink}.

\medskip

Consider first the cochain complex $C^\bul(X, \SF^p)$ which calculates the tropical cohomology groups $H_\trop^{p,q}(X)$. By Corollary~\ref{cor:grading}, the decreasing filtration $W^\bul$ on $\SF^p$ induces a decreasing filtration on the cochain complex $C^\bul(X, \SF^p)$. By an abuse of the notation, we denote this filtration by $W^\bul$. This leads to a spectral sequence
\begin{equation} \label{eqn:abutment_trop}
\CCp{p}_0^{\bul,\bul} \Longrightarrow H^\bul(X, \SF^p) = H^{p,\bul}_{\trop}(X)
\end{equation}
where
\[ \CCp{p}_0^{a,b} = \gr_{\ws}^a C^{a+b}(X, \SF^p) = \bigoplus_{\delta\in X \\ \dims{\delta}=a+b}\bigwedge^b\TT^\dual\delta\otimes\SF^{p-b}(\conezero^\delta) \]
and the differentials in page zero, which are of bidegree $(0,1)$, are given by Proposition~\ref{prop:grading}. We call this the \emph{tropical spectral sequence}. The zero-th page of this spectral sequence is given in Figure \ref{fig:tropical_spectral_sequence}. The dashed arrows correspond to the maps of the first page. The explicit form of all these maps appear later in this section.

\begin{figure}
\caption{The zero-th page of the tropical spectral sequence $\CCp{p}_0^{\bul,\bul}$ over $\SF^p$.} \label{fig:tropical_spectral_sequence}
\renewcommand\S\CCab
\newcommand{\scddots}{\hspace{-2ex}\cddots\hspace{-2ex} }
\small
\[ \begin{tikzcd}
\S00p      \dar["\i^*_\parr\otimes\i^*_\perp"]\rar[dashed, "\pi^*+\d^\i"]&  \S10p      \dar&  \scddots  &  \S{d-p}0p        \dar\\
\S11{p-1}                                                 \dar\rar[dashed]&  \S21{p-1}  \dar&  \scddots  &  \S{d-p+1}1{p-1}  \dar\\
\cdvdots                                                                      \dar&  \cdvdots   \dar&  \scddots  &  \cdvdots         \dar\\
\S{p}p0                                                       \rar[dashed]&  \S{p+1}p0      &  \scddots  &  \S{d}p0
\end{tikzcd}
\]
\end{figure}

Before introducing the second spectral sequence, let us consider the $2p$-th row $\ST^{\bul,2p}_1$ of the Steenbrink spectral sequence. This row can be decomposed into a double complex as follows. We define the double complex $\STp{p}^{\bul, \bul}$ by
\[ \STp{p}^{a,b}:=\begin{cases}
  \bigoplus_{\delta\in X_\f \\ \dims\delta=p+a-b} H^{2b}(\delta) & \text{if $a\geq0$ and $b\leq p$}, \\
  0 & \text{otherwise,}
\end{cases} \]
whose differential of bidegree $(1,0)$ is $\i^*$ and whose differential of bidegree $(0,1)$ is $\gys$. From the identification
\[ \STp{p}^{a,b}=\ST^{a+b-p,2p,p+a-b}, \]
and the following equivalence
\[ p+a-b\geq\abs{a+b-p} \Longleftrightarrow a\geq0 \myand b\leq p, \]
we deduce that
\[ \ST^{\bul,2p}_1=\Tot^\bul(\STp{p}^{\bul,\bul})[-p]. \]
The double complex $\STp{p}^{\bul,\bul}$ is represented in Figure \ref{fig:STp}.

\begin{figure}
\caption{The $p$-th unfolded Steenbrink double complex $\STp{p}^{\bul,\bul}$.} \label{fig:STp}
\small
{ \renewcommand\S\STpab
\[ \begin{tikzcd}
\S{p}0     \dar["\gys"']\rar["\i^*"]&  \S{p+1}0   \dar\rar&  \cddots   \rar&  \S{d}0       \dar\\
\S{p-1}2                    \dar\rar&  \S{p}2     \dar\rar&  \cddots   \rar&  \S{d-1}2     \dar\\
\cdvdots                        \dar&  \cdvdots       \dar&  \cdddots      &  \cdvdots     \dar\\
\S{0}{2p}                       \rar&  \S{1}{2p}      \rar&  \cddots   \rar&  \S{d-p}{2p}
\end{tikzcd} \] }
\end{figure}

\medskip

Finally, we introduce another double complex $\STinf{p}^{\bul,\bul}$, which will play the role of a bridge between $\STp{p}^{\bul,\bul}_0$ and $\CCp{p}^{\bul,\bul}_0$ via Proposition~\ref{prop:STinf_to_ST} and Theorem~\ref{thm:spectral_isomorphism} below. Recall that for $\delta\in X$, the max-sedentarity of $\delta$ is defined by
\[ \maxsed(\delta) := \max\{\sed(x) \st x\in\delta\} \in X_\infty. \]
The double complex $\STinf{p}^{\bul,\bul}$ is defined by
\[ \STinf{p}^{\bul,\bul}=\bigoplus_{\sigma\in X_\infty}\STinf{p}^{\sigma,\bul,\bul} \]
where, for $\sigma\in X_\infty$, we set
\[ \STinf{p}^{\sigma,a,b}:=
\begin{cases}
  \bigoplus_{\delta\in X \\ \maxsed(\delta)=\sigma \\ \dims\delta=p+a-b \\ \dims{\delta_\infty}\leq a} H^{2b}(\delta) & \text{if $b\leq p$}, \\
  0 & \text{otherwise.}
\end{cases} \]
On $\STinf{p}^{\sigma,\bul,\bul}$, the differentials of bidegree $(0,1)$ are given by $\gys$, and those of bidegree $(1,0)$ are given by $\i^*+\pi^*$. Note that, on the whole complex $\STinf{p}^{\bul,\bul}$, no differentials goes from $\STinf{p}^{\sigma,\bul,\bul}$ to $\STinf{p}^{\sigma',\bul,\bul}$ if $\sigma\neq\sigma'$. In other words, the differentials $\gys$ and $\i^*$ are restricted to those pairs of faces $\gamma\ssubface\delta$ such that $\maxsed(\gamma)=\maxsed(\delta)$ and $\sed(\gamma)=\sed(\delta)$, or equivalently such that $\dims{\gamma_\infty}=\dims{\delta_\infty}$ and $\sed(\gamma)=\sed(\delta)$.

The total double complex $\STinf{p}^{\bul,\bul}$ is represented in Figure \ref{fig:STinf}.

\begin{figure}
\caption{The $p$-th extended Steenbrink double complex $\STinf{p}^{\bul,\bul}$. The filtration by columns gives the spectral sequence $\STinfI{p}^{\bul,\bul}_\bul$. Note that here, $\gys$ has to be restricted to pairs of faces with the same max-sedentarity.} \label{fig:STinf}
\small
{ \renewcommand{\S}[3]{\STinfab{\delta}{#1}{#3}{#2}}
\[ \begin{tikzcd}
\S{p}00     \dar["\gys"']\rar["\i^*+\pi^*"]&  \S{p+1}01   \dar\rar&  \cddots   \rar&  \S{d}0{d-p}       \dar\\
\S{p-1}20                          \dar\rar&  \S{p}21     \dar\rar&  \cddots   \rar&  \S{d-1}2{d-p}     \dar\\
\cdvdots                               \dar&  \cdvdots        \dar&  \cdddots      &  \cdvdots          \dar\\
\S{0}{2p}0                             \rar&  \S{1}{2p}1      \rar&  \cddots   \rar&  \S{d-p}{2p}{d-p}
\end{tikzcd} \] }
\end{figure}

\medskip

For the cone $\conezero$ of $X_\infty$, the double complex $\STinf{p}^{\,\conezero,\bul,\bul}$ is described as follows. Faces $\delta$ such that $\maxsed(\delta)=\conezero$ are exactly faces of $X_\f$. Moreover, this implies that $\dims{\delta_\infty}=0$. Thus, $\dims{\delta_\infty}\leq a$ is equivalent to $a\geq0$. Hence,
\[ \STinf{p}^{\,\conezero,\bul,\bul} = \STp{p}^{\bul,\bul}. \]

For any non-zero cone $\sigma\in X_\infty$, we have the following result whose proof will be given at the end of this section.
\begin{prop} \label{prop:STinf_to_ST}
Let $\sigma\in X_\infty\setminus\{0\}$. Then the cohomology of $\Tot^\bul(\STinf{p}^{\sigma,\bul,\bul})$ is zero.
\end{prop}

\begin{proof}
Section \ref{sec:STinf_to_ST} is devoted to the proof of this proposition.
\end{proof}

Filtration by columns of the double complex $\STinf{p}^{\bul,\bul}$ gives a spectral sequence $\STinfI{p}_0^{\bul,\bul}$, which abuts to the cohomology of the total complex of $\STinf{p}^{\bul,\bul}$. By the previous discussion and Proposition~\ref{prop:STinf_to_ST}, this cohomology is just the cohomology of $\Tot^\bul(\STp{p}^{\bul,\bul})$, \ie,
\begin{equation} \label{eqn:abutment_STinf}
\STinfI{p}_0^{\bul,\bul} \Longrightarrow H^\bul(\Tot^\bul(\STp{p}^{\bul,\bul})) = H^{\bul}(\ST^{\bul,2p}_1)[p].
\end{equation}

\medskip

The link between the two spectral sequences $\CCp{p}_\bul^{\bul,\bul}$ and $\STinfI{p}_\bul^{\bul,\bul}$ is summarized by the following theorem.
\begin{thm} \label{thm:spectral_isomorphism}
The Tropical Deligne exact sequence induces compatible canonical isomorphisms
\[ \CCp{p}_k^{\bul,\bul}\simeq\STinfI{p}_k^{\bul,\bul}, \]
between the $k$-th pages for any $k\geq1$.
\end{thm}
\begin{proof}
Section \ref{sec:CC_to_STinf} is devoted to the proof of the isomorphism between first pages, and Section \ref{sec:every_pages} extends the isomorphism to further pages.
\end{proof}

\smallskip
We can now present the proof of the main theorem of this section.
\begin{proof}[Proof of Theorem \ref{thm:steenbrink}]
We recall that we have to prove
\[ H^\bul(\ST^{\bul,2p}_1, \d)[p] = H^{p,\bul}_{\trop}(\X). \]

By Theorem \ref{thm:spectral_isomorphism}, we have
\[ \CCp{p}_k^{\bul,\bul} \simeq \STinfI{p}_k^{\bul,\bul}, \]
for any $k\geq0$.
In particular, both spectral sequences have the same abutment. This implies y \eqref{eqn:abutment_trop} and \eqref{eqn:abutment_STinf} that we have a (non-canonical) isomorphism
\[ H^\bul(\ST^{\bul,2p}_1, \d)[p] \simeq H^{p,\bul}_{\trop}(\X). \]
\end{proof}

\begin{remark}
We note that though we cannot expect to have a canonical isomorphism $H^\bul(\ST^{\bul,2p}_1, \d)[p] \simeq H^{p,\bul}_{\trop}(\X)$, the isomorphism
\[ \gr_F^s H^\bul(\ST^{\bul,2p}_1, \d)[p] \simeq \gr_{\ws}^s H^{p,\bul}_{\trop}(\X) \]
between the graduations is canonical for any $s$. Here the filtration $F^\bul$ is the filtration induced by columns of $\STp{p}^{\bul,\bul}$. The isomorphism between the two cohomology groups above should depend on the data of a smooth deformation of $\X$.
\end{remark}

\medskip

In the remaining of this section, we present the proofs of Theorem~\ref{thm:spectral_isomorphism}, Proposition \ref{prop:STinf_to_ST} and Proposition~\ref{prop:app1}.

\subsection{From $\CCp{p}_1^{\bul,\bul}$ to $\STinfI{p}_1^{\bul,\bul}$} \label{sec:CC_to_STinf}

In this section and the next one, we prove Theorem~\ref{thm:spectral_isomorphism}, which states that $\CCp{p}^{\bul,\bul}_k$ is canonically isomorphic to $\STinf{p}^{\bul,\bul}_k$ for $k\geq1$. As a warm-up, in this section, we study the case of the first page $k=1$. Some of the results and constructions in this section will be used in the next section to achieve the isomorphism of other pages. In order to keep the reading flow, a few technical points of this section will be treated in Appendix~\ref{sec:technicalities}.

\medskip

We have to prove that, for every integer $a$, the cochain complex
{ \renewcommand{\S}[3]{\displaystyle \bigoplus_{\dims{\delta}=#1}\bigwedge^{#2}\TT^\dual\delta\otimes\SF^{#3}(\conezero^\delta)}
\begin{align*}
\CCp{p}^{a,\bul}_0\colon \qquad 0 \to \S{a}0{p} &\to \S{a+1}1{p-1} \to \cdots \\ & \qquad\dots\to \S{a+p}p0 \to 0
\end{align*} }
is quasi-isomorphic to
{ \renewcommand{\S}[3]{\bigoplus_{\dims{\delta}=#1 \\ \dims{\delta_\infty}\leq#2} H^{#3}(\delta)}
\[ \STinf{p}^{a,\bul}_0\colon \qquad 0 \to \S{p+a}a0 \to \S{p+a-1}a2 \to \dots \to \S{a}a{2p} \to 0. \] }
Moreover, the induced isomorphisms in cohomology must commute with the differentials of degree $(1,0)$ on the respective first pages:
\[ \begin{tikzcd}
\CCp{p}_1^{a,\bul} \dar["\sim"{above,sloped}]\rar& \CCp{p}_1^{a+1,\bul} \dar["\sim"{above,sloped}]\\
\STinfI{p}_1^{a,\bul} \rar["\i^*+\pi^*"]& \STinfI{p}_1^{a+1,\bul}
\end{tikzcd} \]

In order to calculate the cohomology of the first cochain complex, we use the tropical Deligne resolution. Applying that exact sequence to the unimodular fan $\Sigma^\delta$ and to any integer $s$, we get the exact sequence
\[0\to \SF^s(\conezero^\delta) \to \bigoplus_{\zeta \in \Sigma^{\delta}\\ \dims \zeta =s} H^0(\zeta) \to \bigoplus_{\zeta \in \Sigma^{\delta}\\ \dims \zeta =s-1} H^2(\zeta) \to \dots \to \bigoplus_{\zeta \in \Sigma^{\delta}\\ \dims \zeta =1} H^{2s-2}(\zeta) \to H^{2s}(\conezero^\delta) \to 0.\]

Given the correspondence between cones $\zeta \in \Sigma^\delta$, and faces $\eta \in X$ which contain $\delta$, with the same sedentarity, we rewrite the exact sequence above in the form
{ \renewcommand{\S}[2]{\bigoplus_{\eta \supface \delta\\ \dims\eta=#1 \\ \makebox[0pt]{\scriptsize $\sed(\eta)=\sed(\delta)$}} H^{#2}(\eta)}
\[ 0\to \SF^s(\conezero^\delta) \to \S{s+\dims\delta}0 \to \S{s-1+\dims\delta}2 \to \dots \to \S{1+\dims\delta}{2s-2} \to H^{2s}(\delta) \to 0. \]
}

\medskip

Replacing now each $\SF^s(\conezero^\delta)$ in $\CCp{p}_0^{a,\bul}$ by the resolution given by the tropical Deligne complex, we get the double complex $\Da{a}^{b, b'}$ whose total complex has the same cohomology as $\CCp{p}^{a,\bul}$. More precisely, define the double complex $\Da{a}^{\bul, \bul}$ as follows.
\[ \Da{a}^{b, b'} := \Dab{a+b}b{p+a-b'}{2b'}. \]
The differential $\d'$ of bidegree $(0,1)$ comes from the Deligne sequence and is $\id\otimes\gys$. The differential of bidegree $(1,0)$ is defined as follows, thanks to the map $\i^*_\parr\colon \bigwedge^s\TT^\dual\gamma\to\bigwedge^{s+1}\TT^\dual\delta$ which has been defined right before Proposition~\ref{prop:grading}. The differential $\d$ of bidegree $(1,0)$ is chosen to be $\i^*_\parr\otimes\id$ (extended using our $\sign$ function as in Section \ref{sec:basic_maps}) on rows with even indices and to be $-\i^*_\parr\otimes\id$ on rows with odd indices. More concisely, we set
\[ \d:=(-1)^{b'}\i^*_\parr\otimes\id. \]

\begin{prop} \label{prop:D_commuting}
We have $\d\d'+\d'\d=0$. Moreover, the inclusion $(\CCp{p}_0^{a,\bul},\i^*_\parr\otimes\i^*_\perp) \hookrightarrow (\Da{a}^{\bul,0},\d')$ is a morphism of cochain complexes.
\end{prop}

\begin{proof}
This follows by a simple computation. The details are given in Appendix~\ref{sec:technicalities}.
\end{proof}

By the exactness of the tropical Deligne sequence that we proved in Proposition \ref{thm:deligne}, we know that the $b$-th column of $\Da{a}^{b,\bul}_0$ of $\Da{a}_0^{\bul,\bul}$ is a right resolution of $\CCp{p}^{a,b}_0$.

\medskip

We claim the following result.
\begin{prop} \label{prop:D_to_STinf}
The $b'$-th row $\Da{a}^{\bul,b'}_0$ of $\Da{a}_0^{\bul,\bul}$ is a right resolution of $\STinf{p}^{a,b'}_0$.
\end{prop}

The proof of this proposition is given in Section~\ref{sec:D_to_STinf}. Admitting this result for the moment, we explain how to finish the proof of the isomorphism between first pages. We will need the following Lemma, which seems to be folklore, for which we provide a proof.

\begin{lemma}[Zigzag isomorphism] \label{lem:zigzag}
Let $\AA^{\bul,\bul}$ be a double complex of differentials $\d$ and $\d'$ of respective degree $(1,0)$ and $(0,1)$. Assume that
\begin{itemize}
\item $\d\d'+\d'\d=0$,
\item $\AA^{b,b'}=0$ if $b<0$ or $b'<0$,
\item $\AA^{b,\bul}$ is exact if $b>0$,
\item $\AA^{\bul,b'}$ is exact if $b'>0$.
\end{itemize}
Then, there is a canonical isomorphism
\[ H^\bul(\AA^{0,\bul})\simeq H^\bul(\AA^{\bul,0}). \]

Moreover, if $\BB^{\bul,\bul}$ is another double complex with the same property, and if $\Phi\colon\AA^{\bul,\bul}\to\BB^{\bul,\bul}$ is a morphism of double complexes, then the following diagram commutes.
\[ \begin{tikzcd}
H^\bul(\AA^{0,\bul})  \dar["\sim"{above, sloped}]\rar["\Phi"]&  H^\bul(\BB^{0,\bul})  \dar["\sim"{above, sloped}] \\
H^\bul(\AA^{\bul,0})                             \rar["\Phi"]&  H^\bul(\BB^{\bul,0})
\end{tikzcd} \]
where $\Phi$ denotes the induced maps on the cohomologies. This also holds if $\Phi$ anticommutes with the differentials.
\end{lemma}
The isomorphism $ H^\bul(\AA^{0,\bul})\simeq H^\bul(\AA^{\bul,0})$ comes from the inclusions
\[(\AA^{0,\bul}, \d') \hookrightarrow (\Tot(\AA^{\bul,\bul}, \d+\d')) \hookleftarrow (\AA^{\bul,0}, \d)\]
which are both quasi-isomorphisms. In the following, we give a canonical description of this isomorphism.

\medskip

{
\renewcommand{\K}{\textnormal{\textsf{L}}}
\renewcommand{\L}{\textnormal{\textsf{R}}}

Let $b$ be a non-negative integer and $b'$ be any integer. Define
\[\K^{b,b'} := \frac{\ker(\d\d')\cap\AA^{b,b'}}{\bigl(\Im(\d)+\Im(\d')\bigr)\cap\AA^{b,b'}} \qquad \textrm{and} \qquad \L^{b,b'}:= \frac{\ker(\d)\cap\ker(\d')\cap\AA^{b,b'+1}}{\Im(\d\d')\cap\AA^{b,b'+1}}. \]

We claim the following.
\begin{claim} Notations as in Lemma~\ref{lem:zigzag}, the map $\d'$ induces an isomorphism
\[ \K^{b'b'} \underset{\d'}{\longsimto} \R^{b,b'}. \]
\end{claim}
\begin{proof}
Let $y$ be an element of $\ker(\d)\cap\ker(\d')\cap\AA^{b,b'+1}$. By the exactness of the $b$-th column, there exists a preimage $x$ of $y$ by $\d'$. Moreover, $x\in\ker(\d\d')\cap\AA^{b,b'}$. Hence, $\d'$ is surjective from $\ker(\d\d')\cap\AA^{b,b'}$ to $\ker(\d)\cap\ker(\d')\cap\AA^{b,b'+1}$. From the trivial identity $\Im(\d\d')=\d'\Im(\d)$, we get
\[ \d^{\prime-1}\bigl(\Im(\d\d')\cap\AA^{b,b'+1}\bigr)=\bigl(\Im(\d)+\ker(\d')\bigr)\cap\AA^{b,b'}=\bigl(\Im(\d)+\Im(\d')\bigr)\cap\AA^{b,b'}. \]
This gives the isomorphism $\d'\colon\K^{b,b'}\simto\L^{b,b'+1}$.
\end{proof}

\begin{proof}[Proof of Lemma~\ref{lem:zigzag}]
Let $b$ be a non-negative and $b'$ be any integer. By the previous claim, we have an isomorphism $\d'\colon\K^{b,b'}\simto\L^{b,b'+1}$.

By symmetry, if $b'$ is non-negative and $b$ is any integer, we get a second isomorphism $\d\colon\K^{b,b'}\simto\L^{b+1,b'}$. Thus, for any non-negative integer $k$, we obtain a zigzag of isomorphisms:
\[
\begin{tikzpicture}[baseline={([yshift=-.5ex]current bounding box.center)}, scale=.6, shorten <=-2pt, shorten >=-2pt]
\scriptsize
\foreach \i in {0,...,3} {
  \foreach \j in {1,...,4} {
    \node (\i\j) at (\i,\j) {$\bullet$};
  }
}
\node[fill=white] (01) at (0,1) {$\K^{0,k}$};
\node[fill=white] (34) at (3,4) {$\K^{k,0}$};
\begin{scope}[shift={(1,1)}]
\foreach \i/\j/\k in {0/1/2,1/2/3,2/3/4} {
  \draw[-latex] (\i\j)--(\j\j);
  \draw[latex-] (\j\j)--(\j\k);
}
\end{scope}
\end{tikzpicture}
\begin{tikzcd}
\K^{k,0} \rar["\sim", "\d'"']& \L^{k,1} &\lar["\sim"', "-\d"] \K^{k-1,1} \rar["\sim", "\d'"']& \L^{k-1,2} &\lar["\sim"', "-\d"] \cdots \rar["\sim", "\d'"']& \L^{1,k} &\lar["\sim"', "-\d"] \K^{0,k}.
\end{tikzcd} \]
Notice that we take $-\d$ and not $\d$. We refer to Remark \ref{rem:preserve_cohomology} for an explanation of this choice.

\medskip

Since $\d'$ is injective on $\AA^{\bul,0}$, and since $\Im(\d')\cap\AA^{k,0}=\{0\}$, we obtain
\[ \K^{k,0}=\frac{\ker(\d\d')\cap\AA^{k,0}}{\bigl(\Im(\d)+\Im(\d')\bigr)\cap\AA^{k,0}}=\frac{\ker(\d)\cap\AA^{k,0}}{\Im(\d)\cap\AA^{k,0}}=H^k(\AA^{\bul,0}). \]
By a symmetric argument, $\K^{0,k}=H^k(\AA^{0,\bul})$. This gives an isomorphism
\[ H^k(\AA^{\bul,0})\simto H^k(\AA^{0,\bul}). \]

The second part of the lemma is clear from the above claim and the preceding arguments. (In the case $\Phi$ anticommutes, the commutative diagram holds for $(-1)^{b+b'}\Phi$ which commutes.)
\end{proof}

\begin{remark} \label{rem:preserve_cohomology}
We now explain why we choose the isomorphisms $\d'$ and $-\d$ in the above proof. Let $a\in\L^{k-i,i+1}$. The zigzag map sends $a$ to $b\in\K^{k-i,i}$ such that $\d'b=a$, and also to $a'=-\d b\in\L^{k-i+1,i}$. Thus, $a-a'=(\d+\d')b\in\Im(\d+\d')$. Moreover, by the definition of $\L^{\bul,\bul}$, it is clear that $a$ and $a'$ belong to $\ker(\d+\d')$. Hence, $a$ and $a'$ are two representatives of the same element of
\[ H^{k+1}\bigl(\Tot(\AA^{\bul,\bul}, \d+\d')\bigr), \]
\ie, the zigzag isomorphism is just identity on the cohomology of $\AA^{\bul,\bul}$.
\end{remark}
}

\begin{proof}[Proof of Theorem \ref{thm:spectral_isomorphism} in page one]
Define the double complex $\AAa{a}^{b,b'}$ by shifting the double complex $\Da{a}^{b,b'}$ by the vector $(1,1)$ and by inserting the the complexes $\CCp{p}_0^{a,\bullet}[1]$ and $\STinf{p}_0^{a,\bullet}[1]$ as the zero-th row and the zero-th column, respectively. More precisely, we set
\[ \AAa{a}^{b,b'}=\begin{cases}
  0            & \text{if $b<0$ or $b'<0$,} \\
  \CCp{p}_0^{a,b-1} & \text{if $b'=0$} \\
  \STinf{p}_0^{a,b'-1} & \text{if $b=0$} \\
  \Da{a}^{b-1,b'-1} & \text{if $b>0$ and $b'>0$,}
\end{cases} \]
with the corresponding differentials $\d$ and $\d'$ of respective bidegree $(1,0)$ and $(0,1)$ such that
\begin{align*}
\AAa{a}^{\bul,0}  &= \CCp{p}_0^{a,\bul}[1], \\
\AAa{a}^{\bul,b'} &= \STinf{p}_0^{a,b'-1} \hookrightarrow \Da{a}^{\bul,b'-1}[1] \qquad \forall b'\geq 1, \\
\AAa{a}^{0,\bul}  &= \STinf{p}_0^{a,\bul}[1],\myand \\
\AAa{a}^{b,\bul}  &= \CCp{p}_0^{a,b-1}\hookrightarrow\Da{a}^{b-1,\bul}[1] \qquad \forall b\geq 1.
\end{align*}
One can easily extend Proposition \ref{prop:D_commuting} to get that $\d$ and $\d'$ commute (the details are given in Appendix~\ref{sec:technicalities}). Moreover, by the exactness of Deligne sequence and by Proposition \ref{prop:D_to_STinf}, all rows but the $0$-th one and all columns but the $0$-th one are exact. Thus we can apply Lemma \ref{lem:zigzag} to get an isomorphism
\[ H^\bul(\AAa{a}^{\bul,0})\simto H^\bul(\AAa{a}^{0,\bul}). \]
Looking at the definition of $\AAa{a}^{\bul,\bul}$, up to a shift by 1, this isomorphism is just what we want, \ie, the following isomorphism in page one.
\[ \CCp{p}_1^{a,\bul}\simto \STinf{p}_1^{a,\bul}. \]

\medskip

It remains to prove the commutativity of the following diagram.
\[ \begin{tikzcd}
\CCp{p}_1^{a,\bul} \dar["\sim"{above,sloped}]\rar["\i^*+\pi^*"]& \CCp{p}_1^{a+1,\bul} \dar["\sim"{above,sloped}]\\
\STinfI{p}_1^{a,\bul} \rar["\i^*+\pi^*"]& \STinfI{p}_1^{a+1,\bul}
\end{tikzcd} \]
To prove this, we construct two anticommutative morphisms of double complexes.

The first morphism corresponds to $\pi^*$ and is naturally defined as follows. On $\CCp{p}^{a,b}$ it is given by $\pi^*\otimes\id$, on $\STinf{p}^{a,b'}$ it equals $\pi^*$, and on $\Da{a}^{b,b'}$ it is defined by $\pi^*\otimes\id$. The anticommutativity properties are easy to check.

The second morphism, which we denote $\d^\i$, corresponds to $\i^*$ and is defined as follows.
To keep the reading flow, we postpone the proof that $\d^\i$ is indeed a morphism and more details about the construction to Appendix~\ref{sec:technicalities}.

\medskip

For each face $\delta\in X$, let $\p_\delta\colon N^{\sed(\delta)}_\R \to \TT_\delta$ be a projection as defined just after Proposition \ref{prop:grading1}. Assume moreover that the projections are compatible in the following sense:
\begin{itemize}
\item if $\gamma\ssubface\delta$ have not the same sedentarity, then $\pi_{\sed(\delta)\ssubface\sed(\gamma)} \p_\delta=\p_\gamma \pi_{\sed(\delta)\ssubface\sed(\gamma)}$, where $\pi_{\sed(\delta)\ssubface\sed(\gamma)}: N^{\sed(\delta)}_{\R}\to N^{\sed(\gamma)}_{\R}$ is the natural projection;
\item if $\gamma\subface\delta$, then $\p_\gamma\p_\delta=\p_\gamma$.
\end{itemize}

For instance, we can choose an inner product on $N_\R$ and extend it naturally on all the strata $N^\sigma_\R$ for $\sigma\in X_\infty$. Then, if the projections $\p_\delta$ are chosen to be the orthogonal projections with respect to this inner product, then they are compatible in the above sense.

\medskip

With these conventions, we define $\d^\i$ on $\STinf{p}^{a,b'}$ by $\d^\i:=\i^*$. On $\CCp{p}^{a,b}$, let $\alpha\otimes\beta \in \bigwedge^{b}\TT^\dual\gamma\otimes\SF^{p-b}(\conezero^\gamma)$, where $\delta\ssupface\gamma$ is a pair of faces of same sedentarity and $\dims\gamma=a+b$. Then the part of its image by $\d^\i$ in $\bigwedge^{b}\TT^\dual\delta\otimes\SF^{p-b}(\conezero^\delta)$ is defined by
\begin{equation} \label{eqn:d''}
\sign(\gamma,\delta) \p_\gamma\rest{\TT\delta}^*(\alpha) \otimes (\u \mapsto \pi_\delta^*\beta(\u')),
\end{equation}
where, for $\u\in\SF^{p-b}(\conezero^\delta)$, the multivector $\u'$ denotes the only element of $\bigwedge^{p-b}\ker(\p_\gamma)$ such that $\pi_\delta(\u')=\u$ (one can check that $\u\in \SF_{p-b}(\gamma)$).

\medskip

On $\Da{a}^{b,b'}$, let $\alpha\otimes x \in \bigwedge^b\TT^\dual\gamma \otimes H^{2b'}(\eta)$ where $\eta\supface\gamma$ is a face of dimension $\dims\gamma+p-b-b'$. Then the part of its image by $\d^\i$ in $\bigwedge^b\TT^\dual\delta \otimes H^{2b'}(\mu)$, with $\mu\ssupface\delta$, is given by
\[ \sign(\eta,\mu) \nvect^\dual_{\delta/\gamma}(\p_\delta(u))\p_\gamma\rest\delta^*(\alpha) \otimes \i^*_{\eta\ssubface\mu}(x), \]
where $u$ is a vector going from any point of $\TT\gamma$ to the vertex of $\mu$ which is not in $\eta$.

We will prove in Appendix~\ref{sec:technicalities} that this map is indeed a morphism, and that the induced map on the first page $\CCp{p}_1^{a,\bul}$ is equal to the differential corresponding to $\i^*$ given by the spectral sequence. This concludes the proof of the proposition.
\end{proof}

\subsection{Proof of Proposition~\ref{prop:D_to_STinf}} \label{sec:D_to_STinf} In order to conclude the proof of isomorphism between the first pages,
we are thus left to prove Proposition \ref{prop:D_to_STinf}. This proposition claims the exactness of the following sequence for any integer $b'$:
\begin{align*}
0 \to&
  \STinfab\eta{p+a-b'}{a}{2b'} \to
  \Dab{a}{0}{p+a-b'}{2b'} \to
  \Dab{a+1}{1}{p+a-b'}{2b'} \to \\[1em] & \hspace{3cm}
  \dots \to
  \Dab{a+p-b'}{p-b'}{p+a-b'}{2b'} \to
0,
\end{align*}
where we recall that the maps are $\i^*_\parr\otimes\id$ and that the relations $\eta \supface \delta$ only concern those faces which have the same sedentarity $\sed(\eta)=\sed(\delta)$.

\medskip

As the form of each term suggests, we can decompose this sequence as a direct sum of sequences over fixed $\eta$ of dimension $s:=p+a-b'$ as follows:
\[ \bigoplus_{\dims\eta=s}\ \Bigl(
  0 \to \R^{\epsilon_\eta} \to \Sxab{a}0 \to \Sxab{a+1}1 \to \dots \to \Sxab{s}{p-b'} \to 0
\Bigr)\otimes H^{2b'}(\eta), \]
where
\[ \epsilon_\eta=\begin{cases}
  1 & \text{if $\dims{\eta_\infty}\leq a$,} \\
  0 & \text{otherwise.}
\end{cases} \]

Therefore, the following result clearly implies Proposition~\ref{prop:D_to_STinf}.

\begin{prop}
For each pair of non-negative integers $r \leq n$, and for each (simplicial) unimodular polyhedron $\eta$ of dimension $n$, the cochain complex
\[ \Cx\eta{r}\colon 0 \longrightarrow \Sxab{r}0[r] \longrightarrow \Sxab{r+1}1[r+1] \longrightarrow \dots \longrightarrow \Sxab{n}{n-r}[n] \longrightarrow 0, \]
whose maps are given by $\i^*_\parr$, has cohomology
\[ \begin{cases}
  \R[r] & \text{if $\dims{\eta_\infty}\leq r$}, \\
  0 & \text{if $\dims{\eta_\infty}>r$.}
\end{cases} \]
\end{prop}
\begin{proof} We will reduce the statement to the computation of some appropriate cohomology groups of the tropical projective spaces.

We begin with the case $\dims{\eta_\infty}=0$, \ie, $\eta$ is a unimodular simplex. Let us study the dual complex instead, which using the duality between $\bigwedge^r\TT^\dual\delta $ and $\bigwedge^{\dims{\delta}-r} \TT^\dual\delta$ for each simplex $\delta$, has the following form
\[ 0 \longrightarrow \Sxab{n}r[n] \longrightarrow \Sxab{n-1}r[n-1] \longrightarrow \dots \longrightarrow \Sxab{r}r[r] \longrightarrow 0. \]
The maps in this complex are induced by restriction map $\bigwedge^r\TT^\dual\delta \to \bigwedge^r\TT^\dual\zeta$ for inclusion of faces $\zeta \ssubface \delta$, extended with $\sign$ according to our convention in Section \ref{sec:basic_maps}.

\medskip

For the ease of arguments, it will be convenient to choose an inner product $\pairing{\ccdot}{\rdot}$ on $\TT\eta$. This inner product restricts to $\TT\delta$ for each face $\delta$ and extends naturally to multivectors. More precisely, for any collection of vectors $u_1, \dots, u_k, v_1, \dots, v_k$, we set
\[ \pairing{u_1\wedge\dots\wedge u_k}{v_1\wedge\dots\wedge v_k} := \det\Bigl(\,\bigl(\pairing{u_i}{v_j}\bigr)_{1\leq i,j\leq k}\,\Bigr). \]

Let $\gamma \ssubface\delta$ be a pair of faces. We define $\p_{\delta\ssupface\gamma}\colon \TT\delta\to\TT\gamma$ to be the orthogonal projection and naturally extend it to multivectors in the exterior algebra.

\medskip

Let $\alpha\in\bigwedge^r\TT^\dual\delta$. We denote by $\alpha^\dual$ the dual multivector of $\alpha$, defined by the property
\[ \alpha( \u ) = \pairing{\alpha^\dual}{\u} \]
for any $\u \in \bigwedge^r\TT\delta$.

Moreover, we denote by $\alpha\rest\gamma$ the restriction of $\alpha$ on $\bigwedge^r\TT\gamma$. By adjunction property, we get the following equality:
\[ \p_{\delta\ssupface\gamma}(\alpha^\dual)=(\alpha\rest\gamma)^\dual. \]
Indeed, for any $\u\in\bigwedge^r\TT\gamma\subset\bigwedge^r\TT\delta$, we have
\[ \alpha\rest\gamma(\u) = \alpha(\u) = \pairing{\alpha^\dual}\u = \pairing{\p_{\delta\ssupface\gamma}(\alpha^\dual)}\u. \]

\medskip

This implies that the linear map which sends $\alpha\mapsto\alpha^\dual$ provides an isomorphism between the two following complexes: our original complex on one side,
\[ \begin{tikzcd}[column sep=small]
0 \rar& \Sxab{n}r \rar& \Sxab{n-1}r \rar& \cddots \rar& \Sxab{r}r \rar& 0,
\end{tikzcd} \]
and its dual
\[ \begin{tikzcd}[column sep=small]
0 \rar& \Syab{n}r \rar& \Syab{n-1}r \rar& \cddots \rar& \Syab{r}r \rar& 0.
\end{tikzcd} \]
Note that the maps of the second complex are given by the orthogonal projections $\p$.

\medskip

We claim that the second complex above computes $H^{r,\bul}_\trop(\TP^n)$, which leads the result. Indeed, the only non-trivial cohomology groups of $\TP^n$ are $H^{p,p}_\trop (\TP^n)\simeq \R$ for $0\leq p\leq n$.

\medskip

Denote by $\delta_0, \dots, \delta_n$ all the faces of codimension one of $\eta$. For each $\delta_i$, let $\ell_i$ be the affine map on $\eta$ which is identically equal to $1$ on $\delta_i$ and takes value $0$ on the opposite vertex to $\delta_i$. Let $\ell_i^0\in\TT^\dual\eta$ be the linear map corresponding to $\ell_i$. It is easy to see that $\sum_i\ell_i^0=0$. Set $u_i:=(\ell_i^0)^\dual\in\TT\eta$. Let $(e_0, \dots, e_n)$ be the standard basis of $\R^n$. We get a linear isomorphism $\TT\eta\simto N_\R=\rquot{\R^{n+1}}{(1,\dots,1)}$ by mapping the $u_i$ to $e_i$ seen in $N_\R$. The family $(e_i)_i$ induce a natural compactification of $N_\R$ into $\TP^n$, and each face $\delta$ of $\eta$ corresponds to a stratum $N_{\R, \delta}$ of dimension $\dims\delta$. Moreover, the linear isomorphism $\TT\eta\simto N_\R$ induces a linear isomorphism $\TT\delta\simto N_{\R,\delta}$. The projections $\pi$ on $\TP^n$ and $\p$ on $\eta$ commute with these isomorphisms. Thus we get an isomorphism between complexes:
\[ \begin{tikzcd}[column sep=small]
0 \rar& \Syab{n}r \rar& \Syab{n-1}r \rar& \cddots \rar& \Syab{r}r \rar& 0, \\
0 \rar& \SF^r(N_\R) \rar& \displaystyle\bigoplus_{\delta\subface\eta \\ \dims\delta = n-1}\SF^r(N_{\R,\delta}) \rar& \cddots \rar& \displaystyle\bigoplus_{\delta\subface\eta \\ \dims\delta = r}\SF^r(N_{\R,\delta}) \rar& 0.
\end{tikzcd} \]
The second complex is just the tropical simplicial complex of $\TP^n$ for the coarsest simplicial subdivision. Hence, its cohomology is $H^{r,\bul}_\trop(\TP^n)=\R[r]$. These concludes the case $\dims{\eta_\infty}=0$.

\medskip

It remains to generalize to any simplicial polyhedron $\eta$. If $\eta=\rho$ is a ray, by a direct computation, the proposition holds: $H^\bul(\Cx\rho{0})$ is trivial and $H^\bul(\Cx\rho{1})=\R[1]$. In general, $\eta$ is isomorphic to $\eta_\f\times\underbrace{\rho\times \dots \times \rho}_{m \textrm{ times }}$, where $\rho$ is a ray and $m:=\dims{\eta_\infty}$. Moreover,
\[ \Cx\eta{r}\simeq\bigoplus_{r_0,\dots,r_m\geq 0 \\ r_0+\dots+r_m=r} \Cx{\eta_\f}{r_0}\otimes\Cx\rho{r_1}\otimes\dots\otimes\Cx\rho{r_m}. \]
By the Künneth formula for cochain complexes of vector spaces, this decomposition also holds for cohomology. I.e., we have
\[ H^\bul(\Cx\eta{r})\simeq\bigoplus_{r_0,\dots,r_m\geq0 \\ r_0+\dots+r_m=r} H^\bul(\Cx{\eta_\f}{r_0})\otimes H^\bul(\Cx\rho{r_1})\otimes\dots\otimes H^\bul(\Cx\rho{r_m}). \]
Looking at the cohomology of the rays, all terms are trivial except the term of indices $r_1=\dots=r_m=1$ and $r_0=r-m$ if $m\leq r$. In this case, we get
\begin{align*}
H^\bul(\Cx\eta{r})
  &\simeq H^\bul(\Cx{\eta_\f}{r-m})\otimes H^\bul(\Cx\rho1)^{\otimes m} \\
  &\simeq \R[r-m]\otimes\R[1]^{\otimes m} \\
H^\bul(\Cx\eta{r}) &= \R[r].
\end{align*}
This concludes the proof of the proposition.
\end{proof}

We have so far proved that the first pages of $\CCp{p}^{\bul,\bul}$ and of $\STinf{p}^{\bul,\bul}$ coincide. Unfortunately, this is not sufficient in general to conclude that the two spectral sequences have the same abutment. Thus, we have to extend the isomorphism to further pages. A direct use of the zigzag lemma on further pages is tedious. Therefore, in the next section, we present a proof using a kind of generalization of the zigzag lemma.

\subsection{Isomorphisms on further pages}
\label{sec:every_pages}

{\renewcommand{\DI}{\prescript{\textup I}{}\Dnop}
\newcommand{\partialI}{\prescript{\textup I}{}\partial}
\newcommand{\CI}{\prescript{\textup I}{}C}
\newcommand{\CCpI}[1]{\prescript{\textup I}{#1}\CCnop}

In this section, we achieve the proof of Theorem \ref{thm:spectral_isomorphism}. To do so, start by stating a general lemma about spectral sequences which we call the \emph{spectral resolution lemma}. Roughly speaking, this lemma states that, under some natural assumptions, a \emph{resolution} of a spectral sequence results in another spectral sequence which coincides with the first one on any further pages.

To prove Theorem \ref{thm:spectral_isomorphism}, we will apply this lemma to the \emph{triply indexed differential complex} $(\Da{\bul}^{\bul,\bul}, \partial=\d^\i+\pi^*\d+\d')$ introduced above (actually a slight modification obtained by inserting the Steenbrink and tropical sequences). This triply indexed complex plays the role of a bridge between the two spectral sequences $\CCp{p}^{\bul,\bul}$ and $\STinfI{p}^{\bul,\bul}$. In fact, as we will show below, it provides a \emph{spectral resolution}, in a sense which will be made precise in a moment, of both the spectral sequences at the same time. Thus, applying the spectral resolution lemma allows to directly conclude the proof of Theorem \ref{thm:spectral_isomorphism}.

\medskip

Though the method presented below gives an alternate proof of the isomorphism in page one, independent of the one given in the previous section, the intermediate results of the previous section will be crucial in showing the stated resolution property, and so to verify the assumption of the spectral resolution lemma. In a sense, this lemma stands upon the basic case treated in the previous section. Moreover, as the proof of the spectral resolution lemma shows, this lemma can be regarded as a generalization of the zig-zag lemma used in the previous section.

\subsubsection{Triply indexed differential complexes} Before we state the lemma, we introduce some conventions. By a \emph{triply indexed differential complex $(\E^{\bul,\bul,\bul}, \partial)$}, we mean
\begin{itemize}
\item $\E^{a,b,c}=0$ unless $a, b, c\geq0$;
\item the map $\partial\colon \E^{\bul,\bul,\bul} \to \E^{\bul,\bul,\bul}$ is of degree one, \ie, it can be decomposed as follows:
\begin{gather*}
\partial=\bigoplus_{i,j,k\in\Z \\ i+j+k=1}\partial^{i,j,k}, \text{ with} \\
\forall i,j,k\in\Z, \qquad \partial^{i,j,k} = \bigoplus_{a,b,c\in\Z}\partial^{i,j,k}\rest{\E^{a,b,c}},
\end{gather*}
where $\partial^{i,j,k}\rest{\E^{a,b,c}}$ is a map from $\E^{a,b,c}$ to $\E^{a+i,b+j,c+k}$;
\item we have $\partial\partial=0$.
\end{itemize}

\medskip

With these assumptions, the total complex $\Tot^\bul(E^{\bul,\bul,\bul})$ becomes a differential complex.

\medskip

The \emph{filtration induced by the first index} on $E^{\bul, \bul,\bul}$ is by definition the decreasing filtration
\[ \E^{\bul,\bul,\bul} = \E^{\geq 0,\bul,\bul} \supseteq \E^{\geq 1,\bul,\bul} \supseteq \E^{\geq 2,\bul,\bul} \supseteq \dots, \]
where
\[ \E^{\geq a,\bul,\bul} = \bigoplus_{a'\geq a}\E^{a',\bul,\bul}. \]
Notice that the differential $\partial$ preserves the filtration induced by the first index if and only if $\partial^{i,j,k}=0$ for any $i<0$. In such a case, we denote by $\EI^{\bul,\bul}_0$ the 0-th page of the induced spectral sequence abutting to the cohomology of the total complex $\Tot^\bul(E^{\bul,\bul,\bul})$. It has the following shape:
\[ \EI^{a,b}_0 = \Tot^{a+b}(\E^{a,\bul,\bul})  = \bigoplus_{m, n\\ m+n=b} \E^{a, m,n}\]
endowed with the differential
\[ \partialI := \sum_{k,l\in\Z}\partial^{0,j,k}. \]
The $k$-th page of this spectral is denoted by $\EI^{\bul,\bul}_k$.

\medskip

We use analogous conventions as above to define bi-indexed differential complexes.

\subsubsection{Statement of the spectral resolution lemma} Let $(C^{\bul,\bul}, \d)$ be a bi-indexed differential complex. We assume that $\d$ preserves the filtration induced by the first index, \ie, we assume $\d^{i,j}=0$ if $i<0$.

\medskip

A \emph{spectral resolution} of $(C^{\bul,\bul}, \d)$ is a triply indexed differential complex $(\E^{\bul,\bul,\bul}, \partial)$ which verifies the following properties enumerated R1-R4:

\smallskip
\begin{enumerate}[label=(R\arabic*)]
\item \label{enum:R1} There is an inclusion $\i\colon C^{\bul,\bul}\to\E^{\bul,\bul,0} \subseteq \E^{\bul,\bul,\bul}$ which respects the bi-indices, \ie, such that $\i\d=\partial\i$ (here $\partial$ is \emph{not} restricted to $\E^{\bul,\bul,0}$).

\smallskip
\item \label{enum:R2} The differential $\partial$ of $E^{\bul,\bul,\bul}$ preserves the filtration induced by the first index, \ie, we have $\partial^{i,j,k}=0$ if $i<0$.

\smallskip
\item \label{enum:R3} We have $\partial^{i,j,k}=0$ if $k\geq 2$.
\end{enumerate}

\medskip

Denote by $\partial_3:=\partial^{0,0,1}$. From the fact that $\partial$ is a differential and the assumptions above, we deduce that $\partial_3\partial_3=0$. Finally, we assume the following resolution property.

\smallskip
\begin{enumerate}[resume*]
\item \label{enum:R4} For any $a,b$ in $\Z$, the differential complex $(\E^{a,b,\bul}, \partial_3\rest{\E^{a,b,\bul}})$ is a right resolution of $C^{a,b}$.
\end{enumerate}

\medskip

Denote by $\EI^{\bul,\bul}$ and $\CI^{\bul,\bul}$ the spectral sequences induced by the filtration by the first index on $E^{\bul,\bul,\bul})$. Then the inclusion $\i$ induces canonical compatible isomorphisms
\[ \CI_k^{\bul,\bul} \simeq \EI_k^{\bul,\bul} \]
between $k$-th pages for any $k\geq1$.

\begin{lemma}[Spectral resolution lemma] \label{lem:spectral_resolution}
Let $(C^{\bul,\bul}, \d)$ be a bi-indexed complex such that the differential $\d$ preserves the filtration induced by the first index. Let $(\E^{\bul,\bul,\bul}, \partial)$ be a spectral resolution of $C^{\bul, \bul}$.

Denote by $\EI^{\bul,\bul}$ and $\CI^{\bul,\bul}$ the spectral sequences induced by the filtration by the first indices on $E^{\bul,\bul,\bul}$ and $C^{\bul, \bul}$, respectively. Then, the inclusion $\i \colon C^{\bul,\bul} \hookrightarrow E^{\bul, \bul, \bul}$ induces canonical compatible isomorphisms
\[ \CI_k^{\bul,\bul} \simeq \EI_k^{\bul,\bul} \]
between $k$-th pages of the two spectral sequences for any values of $k\geq1$.
\end{lemma}

We assume for the moment this lemma and finish the proof of Theorem \ref{thm:spectral_isomorphism}.

\subsubsection{Proof of Theorem \ref{thm:spectral_isomorphism}}
We use the notations of Section~\ref{sec:CC_to_STinf}. We will apply the spectral resolution lemma twice in order to achieve the isomorphism of pages staged in the theorem.

\medskip

We start by gathering the double complexes $\Da{a}^{\bul,\bul}$ together for all $a$ in order to construct a triply indexed differential complex $(\D^{\bul,\bul,\bul}, \partial)$ as follows. First, set
\[ \D^{a,b,b'} := \Da{a}^{b,b'} \]
and let the two differentials of multidegrees $(0,1,0)$ and $(0,0,1)$ be equal to the differentials $\d$ and $\d'$ of the double complexes $\Da{a}^{\bul,\bul}$.

In the course of proving the isomorphism between the first pages of the spectral sequences in Theorem~\ref{thm:spectral_isomorphism}, we introduced in Section~\ref{sec:CC_to_STinf} two anticommutative morphisms $\pi^*$ and $\d^\i$ of multidegree $(1,0,0)$ from $\Da{a}^{\bul,\bul}\to \Da{a+1}^{\bul,\bul}$.

\medskip

We set $\partial =\d^\i+ \pi^* + \d+ \d'$.

\begin{prop}\label{prop:differential_triple_D}
We have $\partial \circ \partial =0$.
\end{prop}
\begin{proof}
The proof consists of calculating directly all the terms which appear in this decomposition. This will be explained in Appendix~\ref{sec:technicalities}.
\end{proof}

Denote by $\CCp{p}^{\bul,\bul}$ the bi-indexed complex $\CCp{p}_0^{\bul,\bul}$ which is endowed with the differential $\d+\d'$, where $\d$ and $\d'$ are the maps defined on $\AA^{\bul,0}=\Tot^\bul(\CCp{p}^{\bul,\bul})$ in the previous section that we restrict to its zero-th row.

The inclusion $\AA^{\bul,0}\to\AA^{\bul,1}$ gives an inclusion of complexes of $\CCp{p}^{\bul,\bul}$ into $\Da{p}^{\bul,\bul,\sqbullet}$. Moreover, for any integers $a$ and $b$, the complex $(\D^{a,b,\sqbullet}, \d')$ is the Deligne resolution of $\CCp{p}_0^{a,b}$. Thus we can apply the spectral resolution lemma. This gives canonical isomorphisms

\begin{equation} \label{eq:iso_pages_1}
\forall \,\,\, k\geq 1, \quad \CCpI{p}_k^{\bul,\bul} \simeq \DI_k^{\bul,\bul}.
\end{equation}

A priori, $\CCpI{p}_k^{\bul,\bul}$ could be different from $\CCp{p}_k^{\bul,\bul}$. So it might appear somehow surprising to see that they are actually equal thanks to the following proposition.
\begin{prop} \label{prop:isomorphism_filtrations}
The natural isomorphism $C^{p,\bul}_\trop \simto \Tot^\bul(\CCp{p}^{\bul,\bul})$ is an isomorphism of filtered differential complexes. Here $C^{p,\bul}_\trop$ comes with differential $\partial_\trop$ and weight filtration $W^\bul$, the differential on $\Tot^\bul(\CCp{p}^{\bul,\bul})$ is $\d+\d'$ and the filtration is induced by the first index.
\end{prop}
\begin{proof}
The proof of this proposition is given in Appendix~\ref{sec:technicalities}.
\end{proof}

In the same way, from Section~\ref{sec:CC_to_STinf}, we deduce an inclusion of complexes of $\STinf{a}^{\bul,\bul}$, with the usual differentials $\gys$ and $\i^*+\pi^*$, into $\D^{\bul,0,\bul}$. Moreover, by Proposition~\ref{prop:D_to_STinf}, proved in Section~\ref{sec:D_to_STinf}, for any pair of integers $a$ and $b'$, the complex $(\D^{a,\blacktriangle,b'}, \d)$ is a resolution of $\STinf{p}_0^{a,b'}$.

\medskip

Thus, by applying again the spectral resolution lemma, we get canonical isomorphisms

\begin{equation} \label{eq:iso_pages_2}
\forall \,\, k\geq 1, \quad \STinfI{p}_k^{\bul,\bul} \simeq \DI_k^{\bul,\bul}.
\end{equation}

\medskip

From \eqref{eq:iso_pages_1} and \eqref{eq:iso_pages_2}, we infer canonical isomorphisms

\[ \CCp{p}_k^{\bul,\bul} \simeq \STinfI{p}_k^{\bul,\bul} \]
between $k$-th pages for any $k\geq1$ of the two spectral sequences, which concludes the proof of Theorem \ref{thm:spectral_isomorphism}. \QED

\subsubsection{Proof of the spectral resolution lemma}  In this section we prove the resolution lemma. Assume that $\i\colon C^{\bul, \bul} \hookrightarrow \AA^{\bul,\bul, 0}\subseteq \AA^{\bul,\bul, \bul}$ is a spectral resolution of the bi-indexed differential complex $C^{\bul, \bul}$.

\begin{proof}[Proof of Lemma \ref{lem:spectral_resolution}]
Let $k\geq1$. By construction of the pages of the spectral sequences, we have
\[ \EI_k^{a,b}=\frac{\Bigl\{x\in\Tot^{a+b}(\E^{\bul,\bul,\bul})\cap\E^{\geq a,\bul,\bul} \st \partial x \in \E^{\geq a+k,\bul,\bul}\Bigr\}}{\E^{\geq a+1,\bul,\bul} + \partial\E^{\geq a-k+1,\bul,\bul}}. \]
(To be more rigorous, we should replace here the denominator by the intersection of the numerator with the denominator.)

\medskip

The idea is to prove that one can take a representative of the above quotient in the image of $\i$. To do so, we use the following induction which can be seen as a generalization of the zigzag lemma.

\begin{claim} Let $a$ be an integer. The following quotient
\[ \frac{\Bigl\{x\in\Tot^{a+b}(\E^{\bul,\bul,\bul})\cap\E^{\geq a,\bul,\bul} \st \partial x \in \E^{\geq a+k,\bul,\bul}\Bigr\}}  {\E^{\geq a+k,\bul,\bul} + \partial\E^{\geq a-k+1,\bul,\bul} + \i(C^{\geq a,\bul})}\]
is zero.
\end{claim}

Notice that, in the denominator, we did not use $\E^{\geq a+1,\bul,\bul}$ but $\E^{\geq a+k,\bul,\bul}$. Doing so we get a stronger result which will be needed later.

\medskip

Let $x$ be an element of the numerator of the following quotient. We first explain how to associate a vector in $\Z^2$ to any $x$ in the numerator. Endowing $\Z^2$ with the lexicographic order, and proceeding by induction on the lexicographical order of vectors associated to the elements $x$, we show that $x$ is zero in the quotient.

\medskip

We decompose $x$ as follows
\[ x = \hspace{-2ex} \sum_{\alpha,\beta,c\in\Z \\ \alpha\geq a \\ \alpha+\beta+c=a+b} \hspace{-2ex} x_{\alpha,\beta, c}, \qquad x_{\alpha,\beta,c}\in\E^{\alpha,\beta,c}. \]
Let $\alpha$ be smallest integer such that there exists $\beta,c\in\Z$ with $x_{\alpha,\beta,c}\neq 0$. Note that if no such $\alpha$ exists, then $x=0$ and we are done. Having chosen $\alpha$, let now $c$ be the largest integer such that $x_{\alpha,\beta,c}\neq0$ for $\beta=a+b-\alpha-c$. We associate to $x$ the vector $(-\alpha, c)$.

\medskip

In order to prove the claim for a fixed value of $a$, we proceed by induction on the lexicographic order of the associated vectors $(-\alpha, c)$. Given a vector $(-\alpha, c)$, in what follows, we set $\beta:=a+b-\alpha-c$.

\smallskip
\emph{Base of the induction: $\alpha\geq a+k$.} Writing $x =\sum x_{m,n,p}$ for $x_{m,n,p}\in E^{m,n,p}$ we see by the choice of $\alpha$ that $m \geq \alpha \geq a+k$ for all non-zero terms $x_{m,n,p}$ it follows that all these terms belong to $\E^{\geq a+k,\bul,\bul}$, and so $x$ obviously belongs to the denominator.

\medskip

\emph{first case: $\alpha<a+k$ and $c>0$.} Denote by $y=\sum y_{\alpha',\beta',c'}$ the boundary $\partial x$ of $x$. Then, since $\partial x$ belongs to $\E^{\geq a+k,\bul,\bul}$, the element $y_{\alpha,\beta,c+1}$ must be zero. Moreover, since $\partial^{a',b',c'}=0$ if $c'\geq2$ or $a'<0$, the part $y_{\alpha,\beta,c+1}$ of $y$ must be equal to $\partial_3( x_{\alpha,\beta,c})$. Thus $x_{\alpha,\beta,c}$ is in the kernel of $\partial_3$. Since $(\E^{\alpha,\beta,\bul}, \partial_3)$ is a resolution, we can take a preimage $z$ of $x_{\alpha,\beta,c}$ in $\E^{\alpha,\beta,c-1}$ for $\partial_3$. Then
\[ \partial z \in \partial\E^{\geq a,\bul,\bul} \subseteq \partial\E^{\geq a-k+1,\bul,\bul}. \]
Thus, $x$ is equivalent to $x'=x-\partial z$, the vector associated to $x'$ is strictly smaller than $(-\alpha, c)$, and we conclude again by the the hypothesis of our induction for $x'$.

\medskip

\emph{second case: $\alpha<a+k$ and $c=0$}. We can follow the arguments used in the second case above up to the moment $x_{\alpha,\beta,0}\in\ker(\partial_3)$. Then, $x_{\alpha,\beta,0}$ belongs to $\i(C^{\alpha,\bul})$. Since $\i(C^{\geq a,\bul})$ is in the denominator, we can subtract $x_{\alpha,\beta,0}$ from $x$ and conclude by induction.

\smallskip
In any case, we have proved $x$ is equal to zero and the claim follows.
As a consequence, we infer that
\begin{equation} \label{eqn:EI2}
\EI_k^{a,b} = \frac{\Bigl\{x\in\Tot^{a+b}(\E^{\bul,\bul,\bul})\cap\E^{\geq a,\bul,\bul} \st \partial x \in \E^{\geq a+k,\bul,\bul}\Bigr\} \cap \Im(\i)}  {\bigl(\E^{\geq a+1,\bul,\bul} + \partial\E^{\geq a-k+1,\bul,\bul}\bigr) \cap \Im(\i)}.
\end{equation}

\medskip

This has to be compared with
\[ \CI_k^{a,b} = \frac{\Bigl\{z\in\Tot^{a+b}(C^{\bul,\bul})\cap C^{\geq a,\bul} \st \d z \in C^{\geq a+k,\bul}\Bigr\}}  {C^{\geq a+1,\bul} + \d C^{\geq a-k+1,\bul}}. \]
From $\i\d=\partial\i$, we get that the induced map $\i\colon \EI_k^{a,b}\to \CI_k^{a,b}$ is well-defined. Moreover, it is now clear that this map is surjective. For the injectivity, we have to work a bit more on the denominator of $\EI_k^{a,b}$.

\medskip

Let $x$ be an element of the intersection of the denominator with the numerator of \eqref{eqn:EI2}. In particular, we have $x\in\E^{\geq a,\bul,\bul}$. Let $z\in\E^{\geq a-k+1,\bul,\bul}$ be such that $x-\partial z\in\E^{\geq a+1,\bul,\bul}$.

\medskip

Now we show that we can reduce to the case $z\in\Im(\i)+\E^{\geq a+1,\bul,\bul}$. Notice that $\partial z\in E^{\geq a,\bul,\bul}$. Thus, $z$ belongs to the numerator of $\EI_{k-1}^{a-k+1,b+k-2}$. If $k\geq 2$, we can apply the induction above to get that $z=z_1+\partial z_2+\i z_3$ with $z_1\in\E^{\geq a,\bul,\bul}$, $z_2\in\E^{\geq a-2k+3,\bul,\bul}$ and $z_3\in C^{\geq a-k+1,\bul}$. This decomposition also holds for $k=1$ setting $z_1=z$ and $z_2=z_3=0$. Since we are only concerned by $\partial z$, we can assume without loss of generality that $z_2=0$.

If $z_1\not\in E^{\geq a+1,\bul,\bul}$, let $c$ be the maximum integer such that $z_{a,\beta,c}\neq0$ for $\beta=b-1-c$. Denote by $y=\sum y_{\alpha',\beta',c'}$ the boundary $\partial z$ of $z$. Since $\partial\i z_3\in\E^{\bul,\bul,0}$, we get that $y_{a,\beta, c+1}=\partial_3 \, z_{a,\beta,c}$. Moreover $x-y=x-\partial z\in\E^{\geq a+1,\bul,\bul}$. Hence, $x_{a,\beta,c+1}$ must be equal to $y_{a,\beta,c+1}$. Since $x\in\Im(\i)\subseteq\E^{\bul,\bul,0}$, $x_{a,\beta,c+1}=0$. From this we get $y_{a,\beta,c+1}=\partial_3z_{a,\beta,c}=0$. Therefore, $z_{a,\beta,c}$ is in the kernel of $\partial_3$. Two cases can occur.

\begin{itemize}
\item If $c>0$, then we can find $w\in \E^{a,\beta,c-1}$ such that $z_{a,\beta,c}=\partial_1w$. In this case, $z':=z-\partial w$ has the same boundary as $z$, and the maximal $c'$ such that $z'_{a,b-1-c',c'}\neq 0$ verifies $c'<c$. Thus, we can apply again the previous argument until we get an element of $\E^{\geq a,\bul,\bul}$, or until we are in the following case.

\item If $c=0$, then $z_{a,\beta,0}\in\i(C^{a,\beta})$.
\end{itemize}

We deduce that, up to a cycle, we can assume that $z_1=z_1'+z_1''$ with $z_1'\in\i(C^{a,\bul})$ and $z_1''\in\E^{\geq a+1,\bul,\bul}$, \ie, $z\in\i(C^{\geq a-k+1,\bul})+\E^{\geq a+1,\bul,\bul}$.

\medskip

Therefore,
\begin{align*}
  &\hspace{-2cm}x \in \bigl(\E^{\geq a+1,\bul,\bul}+\partial\bigl(\i C^{\geq a-k+1,\bul}+\E^{\geq a+1,\bul,\bul}\bigr)\bigr)\cap\Im(\i) \\
  &= \bigl(\E^{\geq a+1,\bul,\bul} + \i\d C^{\geq a-k+1,\bul}\bigr)\cap\Im(\i) \\
  &= \E^{\geq a+1,\bul,\bul}\cap\Im(\i) + \i\d C^{\geq a-k+1,\bul} \\
  &= \i\bigl(C^{\geq a+1,\bul} + \d C^{\geq a-k+1,\bul}\bigr).
\end{align*}
We can replace Equation \eqref{eqn:EI2} by
\[ \EI_k^{a,b} = \frac{\Bigl\{x\in\Tot^{a+b}(\E^{\bul,\bul,\bul})\cap\E^{\geq a,\bul,\bul} \st \partial x \in \E^{\geq a+k,\bul,\bul}\Bigr\} \cap \Im(\i)}  {\i\bigl(C^{\geq a+1,\bul} + \d C^{\geq a-k+1,\bul}\bigr)}. \]
The injectivity of $\i\colon \EI_k^{a,b}\to\CI_k^{a,b}$ now follows from the injectivity of $\i\colon C^{\bul,\bul}\to\E^{\bul,\bul,0}$, which concludes the proof of the spectral resolution lemma.
\end{proof}
}

\vspace{.5cm}
The two remaining section contain the proofs of the two Propositions~\ref{prop:STinf_to_ST} and~\ref{prop:app1}.

\subsection{Vanishing of the cohomology of $\Tot^\bul{\STinf{p}^{\sigma,\bul,\bul}}$ for $\sigma\neq\conezero$} \label{sec:STinf_to_ST}

In this section, we prove Proposition \ref{prop:STinf_to_ST}, which, we recall, claims that the total cohomology of $\STinf{p}^{\sigma,\bul,\bul}_0$ is trivial for $\sigma\neq\conezero$, where
\[ \STinf{p}^{\sigma,a,b}:=
\begin{cases}
  \bigoplus_{\delta\in X \\ \maxsed(\delta)=\sigma \\ \dims\delta=p+a-b \\ \dims{\delta_\infty}\leq a} H^{2b}(\delta) & \text{if $b\leq p$}, \\
  0 & \text{otherwise.}
\end{cases} \]
The differential of bidegree $(1,0)$ is $\i^*+\pi^*$, and the differential of bidegree $(0,1)$ is $\gys$.

This double complex can be unfolded into a triple complex: $\Tot^\bul(\STinf{p}^{\sigma,\bul,\bul})=\Tot^\bul(\E^{\bul,\bul,\bul})$ where
\[ \E^{a,b,a'}:=
\begin{cases}
  \bigoplus_{\delta\in X \\ \maxsed(\delta)=\sigma \\ \dims\delta=p+a+a'-b \\ \dims{\delta_\infty}=a'} H^{2b}(\delta) & \text{if $b\leq p$ and $a\geq 0$}, \\
  0 & \text{otherwise.}
\end{cases} \]
where the differentials of multidegree $(1,0,0)$, $(0,1,0)$ and $(0,0,1)$ are respectively $\i^*$, $\gys$ and $\pi^*$.

To prove that the total cohomology of $\E^{\bul,\bul,\bul}$ is trivial, it suffices to show that $(\E^{a,b,\bul}, \pi^*)$ is exact for any integers $a$ and $b$.

\medskip

Let $\eta$ be a face of sedentarity $\conezero$ such that $\maxsed(\eta)=\sigma$. Notice that $\sed(\eta)=\conezero$ implies $\maxsed(\eta)=\eta_\infty$. For any cone $\tau\subface\sigma$, we denote by $\eta_\infty^\tau$ the face of $\eta$ of sedentarity $\tau$, which is given by
\[ \eta_\infty^\tau = \eta\cap N^\tau_\R. \]
Note that each face of $X$ of max-sedentarity $\sigma$ can be written in a unique way in this form. Notice also that
\begin{itemize}
\item $\pi_\sigma$ induces an isomorphism from $\eta_\f$ to $(\eta_\infty^\tau)_\f$.
\item $\Sigma^\eta\simeq\Sigma^{\eta_\infty^\tau}$. In particular $H^\bul(\eta_\infty^\tau)\simeq H^\bul(\eta)$.
\item $\pi^*_{\delta\ssupface\gamma}$ is nontrivial if and only if $\gamma=\eta_\infty^\tau$ and $\delta=\eta_\infty^\zeta$ for some $\eta$ of sedentarity $\conezero$, and some pair of cones $\zeta\ssubface\tau$ which are faces of $\maxsed(\eta)$.
\end{itemize}

\medskip

From all these observations, we deduce that the sequence $(\E^{a,b,\bul},\pi^*)$ can be decompose as a direct sum:
{ \renewcommand{\S}[1]{\bigoplus_{\tau\subface\sigma \\ \dims\tau=#1}H^{2b}(\eta_\infty^\tau)}
\[ \bigoplus_{\eta\in X^\conezero \\ \maxsed(\eta)=\sigma \\ \dims{\eta}=p+a+\dims\sigma-b} \Bigl( 0 \to H^{2b}(\eta_\infty^\sigma) \to \S{\dims\sigma-1} \to \dots \to \S{1} \to H^{2b}(\eta) \to 0 \Bigr). \] }
This sequence can be rewritten as
{ \renewcommand{\S}[1]{\bigoplus_{\tau\subface\sigma \\ \dims\tau=#1}\R}
\[ \bigoplus_{\eta\in X^\conezero \\ \maxsed(\eta)=\sigma \\ \dims{\eta}=p+a+\dims\sigma-b} \Bigl( 0 \to \R \to \S{\dims\sigma-1} \to \dots \to \S{1} \to \R \to 0 \Bigr) \otimes H^{2b}(\eta). \] }

Thus, it just remains to prove that these sequences are exact for every $\eta$. The complex is clearly isomorphic to the simplicial homology of $\sigma$ with real coefficients. This homology is trivial since it corresponds to the reduced homology of a simplex.

\medskip

This concludes the proof of Proposition \ref{prop:STinf_to_ST} and thus that of Theorem \ref{thm:steenbrink}. \QED

\subsection{Proof of Proposition~\ref{prop:app1}}
\label{sec:proofapp1}

In this final part of this section, we prove Proposition~\ref{prop:app1} which concerns the three basic relations concerning the composition of maps $\gys$ and $\i^*$.

\medskip

Let us start by proving that the composition $\gys\circ\gys$ is vanishing.

Let $a,b$ be a pair of integers, and let $s\geq\abs{a}$ be an integer with $s\equiv a\mod2$. Let $x\in\ST_1^{a,b,s}$. We have to prove that $\gys\circ\gys(x)=0$. By linearity of the Gysin maps, we can assume that $x\in H^{a+b-s}(\delta)$ for some $\delta\in X_\f$ of dimension $s$. If $s-1<\abs{a+1}$, or $s-2<\abs{a+2}$, then $\gys\circ\gys(x)=0$. Otherwise,
\[ \gys\circ\gys(x)=\sum_{\gamma\ssubface\delta}\sum_{\nu\ssubface\gamma}\sign(\nu, \gamma)\sign(\gamma, \delta)\gys_{\gamma\ssupface\nu}\circ\gys_{\delta\ssupface\gamma}(x). \]
This sum can be rewritten in the form
\[ \gys\circ\gys(x)=\sum_{\nu\subface\delta \\ \dims{\nu}=s-2}\sum_{\gamma\ssubface\delta \\ \gamma\ssupface\nu}\sign(\nu, \gamma)\sign(\gamma, \delta)\gys_{\gamma\ssupface\nu}\circ\gys_{\delta\ssupface\gamma}(x). \]
If $\nu$ is a subface of $\delta$ of codimension 2, then, by the diamond property, there exist exactly two faces $\gamma$ and $\gamma'$ such that $\nu\ssubface\gamma,\!\gamma'\ssubface\delta$. Moreover, the composition is just the natural Gysin map from $\delta$ to $\nu$. Thus, we have
\[ \gys_{\gamma'\ssupface\nu}\circ\gys_{\delta\ssupface\gamma'}(x)=\gys_{\gamma\ssupface\nu}\circ\gys_{\delta\ssupface\gamma}(x), \]
and
\[ \sign(\nu, \gamma)\sign(\gamma, \delta)=-\sign(\nu, \gamma')\sign(\gamma', \delta). \]
This shows that the previous sum vanishes, which proves $\gys\circ\gys=0$.

\medskip

One can prove in the same way that $\i^*\circ \i^*=0$, so we omit the details.

\medskip

We now study the map $\i^*\circ\gys+\gys\circ \i^*$. Let $a, b, s, \delta$ and $x\in H^{a+b-s}(\delta)\subseteq\ST_1^{a,b}$ be defined as above. If $s-1<\abs{a+1}$, then $a\geq0$ and $s<\abs{a+2}$ and the codomain of $\i^*\circ\gys+\gys\circ \i^*$ is trivial. Otherwise, as in the previous paragraph, we can decompose $\i^*\circ\gys(x)+\gys\circ \i^*(x)$ as a sum
\[ \sum_{\delta'\in X_\f \\ \dims{\delta'}=s} \Bigl(S_{\subface,\delta'}+S_{\supface,\delta'} \Bigr)\]
where
\begin{gather*}
S_{\subface,\delta'} =
  \sum_{\gamma\ssubface\delta \\ \gamma\ssubface\delta'}
    \sign(\gamma,\delta')\sign(\gamma,\delta)\i^*_{\gamma\ssubface\delta'}\circ\gys_{\delta\ssupface\gamma}(x)\myand \\
S_{\supface,\delta'} =
  \sum_{\eta\ssupface\delta \\ \eta\ssupface\delta'}
    \sign(\delta',\eta)\sign(\delta,\eta)\gys_{\eta\ssupface\delta'}\circ \i^*_{\delta\ssubface\eta}(x),
\end{gather*}
and the sums are over faces in $X_\f$. We claim that for each $\delta'$, the sum $S_{\subface,\delta'}+S_{\supface,\delta'}$ is zero. This will finish the proof of our proposition.

\medskip

\noindent Four cases can happen.

\begin{enumerate}[label=\defaultRoman, leftmargin=0pt]
\item \label{enum:gys_i*:1} Suppose $\delta'\neq\delta$ and $\dims{\delta\cap\delta'}<s-1$.
\end{enumerate}
In this case, both sums $S_{\subface,\delta'}$ and $S_{\supface,\delta'}$ have no terms, and they are both zero.

\medskip

\begin{enumerate}[resume*]
\item \label{enum:gys_i*:2} Suppose $\delta'\neq\delta$, the face $\gamma:=\delta\cap\delta'$ is of dimension $s-1$, and no face $\eta$ of dimension $s+1$ contains both $\delta$ and $\delta'$.
\end{enumerate}

In this case, $S_{\supface,\delta'}$ is a sum with no term. The other sum, $S_{\subface,\delta'}$, contains only one term which is $\i^*_{\gamma\ssubface\delta'}\circ\gys_{\delta\ssupface\gamma}(x)$. This term is zero because, the ray corresponding to $\delta$ and the ray corresponding to $\delta'$ in the fan $\Sigma^\gamma$, are not comparable.

\begin{enumerate}[resume*]
\item \label{enum:gys_i*:3} Suppose $\delta'\neq\delta$, the face $\gamma:=\delta\cap\delta'$ is of dimension $s-1$, and there exists a face $\eta\in X$ of dimension $s+1$ containing both $\delta$ and $\delta'$.
\end{enumerate}

In this case, notice that $\eta=\conv(\delta\cup\delta')$, thus $\eta\in X_\f$. Moreover the four faces form a diamond, \ie, we have $\gamma\ssubface\, \delta,\delta'\ssubface\eta$. We infer that both sums $S_{\ssubface,\delta'}$ and $S_{\ssupface,\delta'}$ contain a unique term, and the sum is given by
\[ \sign(\gamma,\delta')\sign(\gamma,\delta)\i^*_{\gamma\ssubface\delta'}\circ\gys_{\delta\ssupface\gamma}(x)
+ \sign(\delta',\eta)\sign(\delta,\eta)\gys_{\eta\ssupface\delta'}\circ \i^*_{\delta\ssubface\eta}(x). \]
This is zero because
\begin{gather*}
\i^*_{\gamma\ssubface\delta'}\circ\gys_{\delta\ssupface\gamma}(x)=\gys_{\eta\ssupface\delta'}\circ \i^*_{\delta\ssubface\eta}(x), \myand \\
\sign(\gamma,\delta')\sign(\gamma,\delta)=-\sign(\delta',\eta)\sign(\delta,\eta).
\end{gather*}

\medskip

It remains only the case $\delta =\delta'$.

\begin{enumerate}[resume*]
\item \label{enum:gys_i*:4} Suppose $\delta=\delta'$.
\end{enumerate}

This case is more technical. We need some notations. Let $v_0, \dots, v_s$ be the vertices of $\delta$. For $j\in\zint0s$ let $\gamma_j\ssubface\delta$ be the face opposite to $v_j$. Let $\eta_1, \dots, \eta_r$ be the compact faces of dimension $s+1$ containing $\delta$ and $\eta_{r+1}, \dots, \eta_{r+t}$ be the non-compact ones. For each $i\in\zint1r$ we denote by $w_i$ the only vertex of $\eta_i$ which is not in $\delta$. For each $i\in\zint{r+1}{r+t}$, we define $u_i$ to be the primitive vector of the ray $\eta_{i,\infty}$. We will work with the Chow rings. For each $i\in\zint1{r+t}$, the face $\eta_i$ corresponds to a ray $\rho_i$ in $\Sigma^\delta$ and to an element $x_i$ in $A^1(\delta)$. Note that every ray of $\Sigma^\delta$ is of this form. Moreover, for any $j\in\zint0s$, each ray $\rho_i$ corresponds to a ray in $\Sigma^{\gamma_j}$ which we will also denote by $\rho_i$. We also denote by $x_i$ the corresponding element in $A^1(\gamma_j)$. We denote by $\rho_{\delta,j}$ the ray of $\Sigma^{\gamma_j}$ corresponding to $\delta$, and by $x_{\delta,j}$ the associated element in $A^1(\gamma_j)$. If $\ell$ is a linear form on $N^\eta_\R$ for some face $\eta$, and if $\varrho$ is a ray in $\Sigma^\eta$, we will write $\ell(\varrho)$ for the value of $\ell$ on the primitive vector of $\varrho$.

\medskip

Since $\gys$ and $\i^*$ are ring homomorphisms, we can assume without loss of generality that the chosen element in $\ST_1^{a,b,s}(X)$ is the unity $1$ of $A^0(\Sigma)$. Then, for any $i\in\zint1r$,
\[ \gys_{\eta_i\ssupface\delta}\circ\i^*_{\delta\ssubface\eta_i}(1)=x_i. \]

\medskip

In this case, $S_{\supface,\delta}$ is given by the sum over \emph{compact} faces $\eta$ of dimension $s+1$ containing $\delta$ of
\[ \sign(\delta,\eta)^2\gys_{\eta\ssupface\delta}\circ\i^*_{\delta\ssubface\eta}(1). \]
It follows that
\[ S_{\supface,\delta}=\sum_{i=1}^rx_i. \]

\medskip

Let $j\in\zint0s$. We have $\gys_{\delta\ssupface\gamma_j}(1)=x_{\delta,j}\in A^1(\gamma_j)$. Choose an affine form $\ell_j$ on $N_\R$ that is zero on $\gamma_j$ and $-1$ on $v_j$. Let $\ell^0_j$ be the corresponding linear form. It induces a linear form $\lambda_j$ on $N^{\gamma_j}_\R$. Moreover, we have the following properties for $\lambda_j$:
\begin{gather*}
\lambda_j(\rho_{\delta,j})=\ell(v_j)=-1, \\
\forall i\in\zint1r, \,\, \lambda_j(\rho_i)=\ell_j(w_i), \, \myand \\
\forall i\in\zint{r+1}{r+t},\,\, \lambda_j(\rho_i)=\ell^0_j(u_i).
\end{gather*}
Thus, $x_{\delta,j}\in A^1(\gamma_j)$ is equal to
\[ \sum_{i=1}^r\ell_j(w_i)x_i+\sum_{i=r+1}^{r+t}\ell^0_j(u_i)x_i+x', \]
where $x'\in A^1(\gamma_j)$ is incomparable with $x_{\delta,j}$ in the sense that $x'x_{\delta,j}=0$. Hence, the image of $x_{\delta,j}$ by $\i^*_{\gamma_j\ssubface\delta}$ is
\[ \i^*_{\gamma_j\ssubface\delta}\circ\gys_{\delta\ssupface\gamma_j}(1)= \sum_{i=1}^r\ell_j(w_i)x_i+\sum_{i=r+1}^{r+t}\ell^0_j(u_i)x_i\in A^1(\delta). \]
Set $\ell=\ell_0+\dots+\ell_s$ and denote by $\ell^0$ the corresponding linear form. Summing over all $\gamma_j$, we get
\[ S_{\subface, \delta}=\sum_{i=1}^r\ell(w_i)x_i+\sum_{i=r+1}^{r+t}\ell^0(u_i)x_i. \]
Let $\cteun$ be the constant function on $N_\R$ which takes value $1$ everywhere. We can rewrite
\[ S_{\supface, \delta}=\sum_{i=1}^rx_i=\sum_{i=1}^r\cteun(w_i) x_i. \]

\medskip

Set $\ell'=\ell+\cteun$ and let $\ell^{\prime 0}=\ell^0$ be the corresponding linear form. We have
\[ S_{\subface, \delta}+S_{\supface, \delta}=\sum_{i=1}^r\ell'(w_i)x_i+\sum_{i=r+1}^{r+t}\ell^{\prime 0}(u_i)x_i. \]
Since $\ell$ equals $-1$ on each vertex of $\delta$, $\ell'$ is zero on $\delta$. Thus, it induces a linear form $\lambda'$ on $N^\delta_R$ which verifies
\begin{gather*}
\forall i\in\zint1r, \lambda'(\rho_i)=\ell'(w_i)\myand \\
\forall i\in\zint{r+1}{r+t}, \lambda'(\rho_i)=\ell^{\prime 0}(u_i).
\end{gather*}
Therefore, $S_{\subface, \delta}+S_{\supface, \delta}=0$ in $A^1(\delta)$, which is what we wanted to prove.

\medskip

In each of the four cases above we showed that $S_{\subface, \delta'}+S_{\supface, \delta'}=0$, thus $\gys\circ \i^*(x)+\i^*\circ\gys(x)=0$.

\medskip

The proof of Proposition~\ref{prop:app1} is complete.


\section{K\"ahler tropical varieties} \label{sec:kahler}
The aim of this section is to show that $\ST_1^{\bul, \bul}$ and the tropical cohomology $H^{\bul,\bul}_\trop$ admit a Hodge-Lefschetz structure. The precise meaning will only appear in Section~\ref{sec:differential_HL_structure} where we give a precise definition of Hodge-Lefschetz structures and treat the consequences of the materials introduced in this section in detail.

\medskip

We introduce the monodromy and Lefschetz operators on $\ST_1^{\bul, \bul}$ in this section. The definition of the latter is based on the definition of K\"ahler forms in tropical geometry which are introduced in this section as well.

\subsection{Monodromy operator $N$ on $\ST_1^{\bul, \bul}$} The monodromy operator is the data of maps $N^{a,b} \colon \ST_1^{a,b} \to \ST_1^{a+2, b-2}$ defined as follows. Writing
$\ST_1^{a,b} = \bigoplus \ST_1^{a,b,s}$ for $s \geq \abs{a}$ with $s \equiv a (\mod 2)$, then $N^{a,b}$ is defined by setting its restrictions $N^{a,b,s}$ on
$\ST_1^{a,b,s}$ to be

\begin{equation}
N^{a,b, s} = \begin{cases} \mathrm{id} \colon\ST_1^{a,b,s} \to \ST_1^{a+2,b-2,s}   & \textrm{if } s \geq \abs{a+2}, \\
0 & \textrm{otherwise}.
\end{cases}
\end{equation}

If there is no risk of confusion, we drop all the superscripts and denote simply by $N$ the operator $N^{a, b}$ and all its restrictions $N^{a,b,s}$.
\begin{remark} More precisely, we could view the collection of maps $N^{a,b}$ as a map $N$ of bidegree $(2, -2)$
\[N = \bigoplus_{a,b}N^{a,b} \colon \bigoplus_{a,b} \ST_1^{a,b} \longrightarrow \bigoplus_{a,b} \ST_1^{a,b}.\]
\end{remark}
We have the following proposition which summaries the main properties of the monodromy operator.

\begin{prop} \label{prop:N_commutes} We have
\begin{itemize}
\item $[N,\i^*] = [N, \gys] =0$
\item for $a\leq 0$, we have an isomorphism $N^{-a} \colon \ST_1^{a,b} \to \ST_1^{-a, b+2a} $
\item for $a \leq 0$, $\ker N^{-a+1} \cap \ST_1^{a,b} = \ST_1^{a, b, -a}$.
\end{itemize}
\end{prop}

\begin{proof}
For the first point, let $a, b, s$ be three integers with $a \equiv s \pmod 2$. We have the following commutative diagrams.
\[ \begin{tikzcd}
\ST^{a,b,s}_1 \rar["\id"]\dar["\i^*"] & \ST^{a+2,b-2,s}_1 \dar["\i^*"] \\
\ST^{a+1,b,s+1}_1 \rar["\id"] & \ST^{a+3,b-2,s+1}_1
\end{tikzcd} \qquad \begin{tikzcd}
\ST^{a,b,s}_1 \rar["\id"]\dar["\gys"] & \ST^{a+2,b-2,s}_1 \dar["\gys"] \\
\ST^{a,b+2,s}_1 \rar["\id"] & \ST^{a+2,b,s}_1
\end{tikzcd} \]
To get the commutativity in $\ST_1^{\bul,\bul}$, it remains to check that, in each diagram, if the top-left and bottom-right pieces appear in $\ST_1^{\bul,\bul}$, then the two other pieces appear as well. This comes from the following facts, which are easy to check.
\begin{gather*}
s \equiv a \equiv a+2 \pmod 2 \quad\myand\quad s+1 \equiv a+1 \equiv a+3 \pmod 2, \myand \\
s\geq\abs{a} \myand s+1\geq\abs{a+3} \Longrightarrow s\geq\abs{a+2} \myand s+1\geq\abs{a+1}.
\end{gather*}

\medskip

For the second point, for $a\leq0$, $s\geq\abs{a}$ implies $s\geq\abs{a+2k}$ for any $k\in\zint0{-a}$. Thus $N^{-a}$ induces the isomorphism $\ST_1^{a,b,s}\simeq\ST_1^{-a,b+2a,s}$, which implies the result.

\medskip

For the last point, for $a\leq0$, we have
\begin{align*}
\ker N^{-a+1}\cap\ST_1^{a,b}
  &= (N^{-a})^{-1}\ker\bigl(N\colon\ST_1^{-a,b+2a}\to\ST_1^{-a+2,b+2a-2}\bigr) \\
  &= (N^{-a})^{-1}\ST_1^{-a,b+2a,-a} \\
  &= \ST_1^{a,b,-a}. \qedhere
\end{align*}
\end{proof}

\begin{defi}[Primitive parts of the monodromy]
For $a\leq 0$, the kernel of $N^{-a+1}$ on $\ST_1^{a,b}$ is called the \emph{primitive part of bidegree $(a,b)$} for the monodromy operator $N$.
\end{defi}

\subsection{K\"ahler forms}

\subsubsection{Local ample classes} Let $\Sigma\subseteq N_\R$ be a (unimodular) Bergman fan. Recall that an element $\ell\in A^1(\Sigma)$ is called an ample class if it corresponds to a strictly convex cone-wise linear function on $\Sigma$.

\begin{prop} \label{prop:ample_stable}
Let $\Sigma$ be a Bergman fan. Let $\ell\in A^1(\Sigma)$ be an ample class of $\Sigma$. Then
\begin{itemize}
\item $\ell$ verifies the Hodge-Riemann relations $\HR(\Sigma,\ell)$,
\item for any cone $\sigma\in\Sigma$, $\i^*_\sigma(\ell)\in A^1(\Sigma^\sigma)$ is an ample class of $\Sigma^\sigma$.
\end{itemize}
\end{prop}

\begin{proof}
The first point is Theorem \ref{thm:mainlocal}. The second point follows from Proposition~\ref{prop:ell_f_independent} combined with \ref{thm:mainlocal}.
\end{proof}

\subsubsection{Kähler class} Let $X$ be a unimodular polyhedral complex structure on a smooth compact tropical variety.

\begin{defi}
A \emph{Kähler class} for $X$ is the data of ample classes $\ell^v \in H^2(v)$ for each vertex $v\in X_\f$ such that, for each edge $e$ of $X_\f$ of extremities $u$ and $v$, we have the compatibility
\[\i_{u \ssubface e}^*(\ell^u) = \i_{v \ssubface e}^*(\ell^v)\]
in $H^2(e)$.
\end{defi}

\begin{thm}\label{thm:KahlerClass}
Each Kähler class defines an element $\omega \in H^{1,1}(X)$.
\end{thm}

\begin{proof}
Let $\ell$ be a Kähler class. We have
\[ \ell \in \bigoplus_{v\in X_\f \\ \dims{v}=0}H^2(v) = \ST_1^{0,2,0}\subseteq \ST_1^{0,2}. \]
Thus, $\ell$ can be seen as an element of $\ST_1^{0,2}$. Moreover,
\begin{align*}
\d\ell
  &= \sum_{v\in X_\f \\ \dims{v}=0}\sum_{e\in X_\f \\ e\ssupface v}\sign(v,e)\i^*_{v\ssubface e}(\ell^v) \\
  &= \sum_{e\in X_\f \\ \dims{e}=1 \\ u,v\ssubface e}\sign(u,e)\i^*_{u\ssubface e}(\ell^u)+\sign(v,e)\i^*_{v\ssubface e}(\ell^v), \\
&= 0.
\end{align*}
Hence, $\ell$ is in the kernel of $\ST_1^{0,2}$ and thus induces an element $\omega\in \ST_2^{0,2}\simeq H^{1,1}(X)$.
\end{proof}

\begin{defi} \begin{itemize}[label=-]
\item An element $\omega$ in $H^{1,1}(X)$ is called a \emph{Kähler class} if there exist ample elements $\ell^v \in H^2(v)$ as above which define $\omega \in H^{1,1}(X)$.
\item A smooth compact tropical variety is called \emph{K\"ahler} if it admits a unimodular triangulation $X$ and a Kähler class $\omega\in H^{1,1}(X)$.
\end{itemize}

\end{defi}

\subsubsection{Projective tropical varieties are K\"ahler}
\label{sec:projective_Kahler}

In this section, we consider the case of a unimodular polyhedral complex $Y \subseteq \R^n$ and the compactification $X$ of $Y$ in the tropical toric variety $\TP_{Y_\infty}.$ In this situation $Y\subseteq\R^n$ is the part of $X$ of sedentarity $\conezero$.

We assume that the triangulation $Y$ is quasi-projective. Under this assumption, we will show that $X$ is Kähler, \ie, that there exists a Kähler class in $H^{1,1}(X)$.

\medskip

Let $f$ be a strictly convex piecewise linear function on $Y$. Let $\delta$ be a face of $X_\f$, and $\phi$ be an affine linear form with the same restriction to $\delta$ as $f$. Denote by $f^\delta$ the cone-wise linear function induced by $f-\phi$ on $\Sigma^\delta$.

\begin{propdefi} \label{propdefi:convex_subface}
The function $f^\delta$ is strictly convex. Denote by $\ell^\delta$ the corresponding element $\ell(f^\delta)$ of $H^2(\delta)$. This element $\ell^\delta$ does not depend on the chosen $\phi$. Moreover, if $\gamma$ is a face of $\delta$, then $\ell^\delta=\i^*_{\gamma\subface\delta}(\ell^\gamma)$.
\end{propdefi}

\begin{proof}
For both parts, one can adapt the proof given in Proposition~\ref{prop:ell_f_independent} in the local case.
\end{proof}

It follows that $\ell^\delta$ is an ample class in $A^1(\delta) = H^2(\delta)$. In particular, this holds for vertices of $X_\f$, giving an element $(\ell^v)_{v \text{ vertex of } X_\f}\in\ST^{0,2}_1$. We get the following theorem.

\begin{thm} \label{thm:projective_Kahler}
A strictly convex function $f$ on $Y$ defines a Kähler class of $X$.
\end{thm}
\begin{proof}
We just have to notice that, for any edge $e\in X_\f$ of extremities $u$ and $v$,
\[ \i^*_{v\ssubface e}(\ell^v)=\ell^e=\i^*_{u\ssubface e}(\ell^u). \qedhere \]
\end{proof}

\medskip

The two following results give a better understanding of the link between piecewise linear functions and elements of $H^{1,1}(\X)$. The corollary will be useful in Section \ref{sec:proofmaintheorem}.

Let $v$ be a vertex of $Y$. We define $\chi_v$ the \emph{characteristic function of $v$} to be the unique piecewise linear function on $Y$ which takes value $1$ at $v$, takes value $0$ at any other vertex of $Y$, and whose slope along any non-compact one dimensional face of $Y$ is zero.

\begin{prop}
Let $v$ be a vertex of $Y$ and let $\chi_v$ be the characteristic function of $v$. Then the class given by $\chi_v$ in $H^{1,1}(\X)$ is trivial.
\end{prop}

\begin{proof}
Let us also denote by $\chi_v$ the corresponding element of $\ST^{0,2}_1$. Then
\[ \chi_v = -\d\,\sum_{e \in Y_\f \\ e \ssupface v}\sign(v,e)1_e \]
where $1_e$ is the natural generator of $H^0(e)\subseteq\ST^{-1,2}$. Hence $\chi_v$ is a boundary, and its cohomological class is trivial.
\end{proof}

Recall that if $f$ is a piecewise linear function on $Y$, then we defined the asymptotic part $f_\infty\colon Y_\infty \to \R$ of $f$ in Section \ref{sec:X_infty_X}. We immediately get the following corollary.

\begin{cor} \label{cor:class_independent_finite_part}
Let $f$ be a piecewise linear function on $Y$ and let $\ell$ be the corresponding class in $H^{1,1}(X)$. Assume moreover that $f_\infty$ is well-defined. Then $\ell$ only depends on $f_\infty$.
\end{cor}

Another equivalent description of the class $\ell$ is given in Section \ref{sec:cl_D'}.

\subsection{The Lefschetz operator} Let $\ell$ be a Kähler class of $X$. For each face $\delta\in X_\f$, we define $\ell^\delta\in H^2(\delta)\simeq A^1(\delta)$ by
\[ \ell^\delta=\i^*_{v\subface\delta}(\ell^v), \]
where $v$ is any vertex of $\delta$. This definition does not depend on the choice of $v$ because, if $u$ is another vertex of $\delta$, we have
\[ \i^*_{v\subface\delta}(\ell^v)=\i^*_{e\subface\delta}\circ \i^*_{v\subface e}(\ell^v)=\i^*_{e\subface\delta}\circ \i^*_{u\subface e}(\ell^u)=\i^*_{u\subface\delta}(\ell^u), \]
where $e$ is the edge of extremities $u$ and $v$.

Let $Q_\delta$ be the bilinear form defined on $H^{2k}(\delta) = A^k(\delta)$, for $k \leq \frac{d-\dims\delta}2$, by
\[ \forall \, a,b\in H^{2k}(\delta), \qquad Q_\delta (a,b) = \deg(ab(\ell^\delta)^{d-\dims\delta-2k}).\]
Every $\ell^\delta$ is ample by Proposition \ref{prop:ample_stable}.

For each $\delta$, let
\[P^{2k}(\delta) := \ker\Bigl((\ell^\delta)^{d-\dims{\delta}-2k+1} \colon H^{2k}(\delta) \to H^{2d-2\dims{\delta}-2k+2}(\delta)\Bigr)\]
be the primitive part of $H^{2k}(\delta).$ Note that with the notations of Section~\ref{sec:local}, we have $P^{2k}(\delta) = A^k_{\prim, \ell^\delta}(\Sigma^\delta)$.

\medskip

The following properties hold.

\begin{itemize}
\item (Hard Lefschetz) For each $k \leq \frac{d-\dims{\delta}}2$, the map
\[(\ell^\delta)^{d-\dims\delta -2k} \colon H^{2k}(\delta) \to H^{2d-2\dims{\delta} -2k}(\delta)\]
is an isomorphism.
\item (Hodge-Riemann) The bilinear form $(-1)^k Q_\delta(\ccdot,\rdot)$ restricted to the primitive part $P^{2k}(\delta) = \ker\Bigl((\ell^\delta)^{d-\dims{\delta}-2k+1} \colon H^{2k}(\delta) \to H^{2d-2\dims{\delta}-2k+2}(\delta)\Bigr)$ is positive definite.
\item (Lefschetz decomposition) For $k\leq\frac{d-\dims\delta}2$, we have the orthogonal decomposition
\[ H^{2k}(\delta)=\bigoplus_{i=0}^k \ell^{k-i}P^{2i}(\delta), \]
where, for each $i\in\zint0k$, the map $\ell^{k-i}\colon P^{2i}(\delta)\simto \ell^{k-i}P^{2i}(\delta)$ is an isomorphism preserving $Q_\delta$.
\end{itemize}

\medskip

We now define the Lefschetz operator $\ell^{a,b} \colon \ST_1^{a,b} \to \ST_1^{a, b+2} $ by
\[\ell^{a,b} = \bigoplus \ell^{a,b,s} \colon \bigoplus_{s \geq \abs{a} \\
s \equiv a \pmod 2} \ST_1^{a,b, s} \to \bigoplus_{s \geq \abs{a} \\
s \equiv a \pmod 2} \ST_1^{a,b+2, s}\]
and the operator $\ell^{a,b,s}$ is defined by
\[\ell^{a,b,s} = \bigoplus \ell^\delta \colon \bigoplus_{\delta \in X_\f \\ \dims{\delta} =s} H^{a+b-s}(\delta) \longrightarrow    \bigoplus_{\delta \in X_\f \\ \dims{\delta} =s} H^{a+b+2-s}(\delta).\]

\medskip

We usually drop the indices $a,b$ and simply denote by $\ell$ the operator of bidegree $(0,2)$.

\begin{remark} More precisely, we could view the collection of maps $\ell^{a,b}$ as a map $\ell$
\[\ell = \bigoplus_{a,b}\ell^{a,b} \colon \bigoplus_{a,b} \ST_1^{a,b} \longrightarrow \bigoplus_{a,b} \ST_1^{a,b+2}.\]
\end{remark}

We now recast Lemma~\ref{lem:i_gys_basic_properties-local} in the following.

\begin{prop} \label{lem:i_gys_basic_properties}
Let $\gamma\ssubface\delta$ be a pair of faces, and let $x\in A^\bul(\gamma)$ and $y\in A^\bul(\delta)$. Denote by $\rho_{\delta/\gamma}$ the ray associated to $\delta$ in $\Sigma^\gamma$, and by $x_{\delta/\gamma}$ the associated element of $A^1(\gamma)$.
\begin{gather}
\text{$\i^*_{\gamma\ssubface\delta}$ is a surjective ring homeomorphism,} \label{eqn:i_surjective_homeo} \\
\gys_{\delta\ssupface\gamma}\circ\i^*_{\gamma\ssubface \delta}(x)=x_{\delta/\gamma}\cdot x, \label{eqn:gys_circ_i} \\
\gys_{\delta\ssupface\gamma}(\i^*_{\gamma\ssubface\delta}(x)\cdot y)=x\cdot\gys_{\delta\ssupface\gamma}(y). \label{eqn:gys_i_simplification}
\end{gather}
Moreover, if $\deg_\gamma\colon A^{d-\dims\gamma}\to\R$ and $\deg_\delta\colon A^{d-\dims\delta}\to\R$ are the natural degrees map, then
\begin{equation} \label{eqn:deg_circ_gys}
\deg_\delta=\deg_\gamma\circ\gys_{\delta\ssupface\gamma}.
\end{equation}
Finally, $\gys_{\delta\ssupface\gamma}$ and $\i^*_{\gamma\ssubface\delta}$ are dual in the sense that
\begin{equation} \label{eqn:i_gys_dual}
\deg_\gamma(x\cdot\gys_{\delta\ssupface\gamma}(y))=\deg_\delta(\i^*_{\gamma\ssubface\delta}(x)\cdot y).
\end{equation}
\end{prop}

\begin{proof} The proof is similar to the one given in the local case.
\end{proof}

\medskip

The following proposition summarizes basic properties of the Lefschetz operator $\ell$.
\begin{prop} We have
\begin{itemize}
\item $[\ell, N] = 0$, $[\ell, \i^*] = [\ell, \gys] =0$.
\item for a pair of integers $a,b$ with $a+b \geq d$, the map $\ell^{d-a-b} \colon \ST_1^{a,b} \longrightarrow \ST_1^{a, 2d - 2a-b}$ is an isomorphism.
\end{itemize}
\end{prop}
\begin{proof}
For the first point, let $\gamma\ssubface\delta$ be a pair of faces of $X_\f$. We recall that $\i^*_{\gamma\ssubface\delta}\colon A^\bul(\gamma)\to A^\bul(\delta)$ is a surjective ring homeomorphism. Thus, for any $a\in A^\bul(\gamma)$,
\[ \i^*_{\gamma\ssubface\delta}(\ell^\gamma a)=\i^*_{\gamma\ssubface\delta}(\ell^\gamma)\i^*_{\gamma\ssubface\delta}(a)=\ell^\delta\i^*_{\gamma\ssubface\delta}(a). \]
By Equation \eqref{eqn:gys_i_simplification}, for any $x\in A^\bul(\delta)$,
\[ \gys_{\delta\ssupface\gamma}(\ell^\delta x)=\gys_{\delta\ssupface\gamma}(\i^*_{\gamma\ssubface\delta}(\ell^\gamma) x)=\ell^\gamma\gys_{\delta\ssupface\gamma}(x). \]

Let $a,b,s$ be three integers with $a\equiv s\pmod2$. Multiplying the previous equations by $\sign(\gamma,\delta)$ and summing over couples $\gamma\ssubface\delta$, we get that the three following diagrams are commutative.
\[ \begin{tikzcd}[column sep=small]
\ST^{a,b,s}_1 \rar["\ell"]\dar["\id"] & \ST^{a,b+2,s}_1 \dar["\id"] \\
\ST^{a+2,b-2,s}_1 \rar["\ell"] & \ST^{a+2,b,s}_1
\end{tikzcd} \quad \begin{tikzcd}[column sep=small]
\ST^{a,b,s}_1 \rar["\ell"]\dar["\i^*"] & \ST^{a,b+2,s}_1 \dar["\i^*"] \\
\ST^{a+1,b,s+1}_1 \rar["\ell"] & \ST^{a+1,b+2,s+1}_1
\end{tikzcd} \quad \begin{tikzcd}[column sep=small]
\ST^{a,b,s}_1 \rar["\ell"]\dar["\gys"] & \ST^{a,b+2,s}_1 \dar["\gys"] \\
\ST^{a+1,b,s-1}_1 \rar["\ell"] & \ST^{a+1,b+2,s-1}_1
\end{tikzcd}
\]
To get the commutativity in $\ST_1^{\bul,\bul}$, it remains to check that, in each diagram, if the top-left and bottom-right pieces appear in $\ST_1^{\bul,\bul}$, then the two other pieces appear as well, which is immediate.

\medskip

For the second point, we have
\begin{gather*}
\ST_1^{a,b}
  = \bigoplus_{s\geq\abs a \\ s\equiv a\pmod2}\bigoplus_{\delta\in X_\f \\ \dims\delta=s} H^{a+b-s}(\delta) \\
\ST_1^{a,2d-2a-b}
  =\bigoplus_{s\geq\abs a \\ s\equiv a\pmod2}\bigoplus_{\delta\in X_\f \\ \dims\delta=s} H^{2d-a-b-s}(\delta)
\end{gather*}
and $\ell^{d-a-b}\colon \ST_1^{a,b}\to\ST_1^{a,{2d-2a-b}}$ can be rewritten as
\[ \bigoplus_{s\geq\abs a \\ s\equiv a\pmod2}\bigoplus_{\delta\in X_\f \\ \dims\delta=s} (\ell^\delta)^{d-a-b}\colon H^{a+b-s}(\delta)\longrightarrow H^{2d-a-b-s}(\delta). \]
For each $\delta$ of dimension $s$, since $\ell^\delta$ is ample, we get an isomorphism
\[ (\ell^\delta)^{d-a-b}\colon H^{a+b-s}(\delta)\longrightarrow H^{2d-a-b-s}(\delta)=H^{2(d-\dims\delta)-(a+b-s)}.\]
Thus, we conclude that $\ell^{d-a-b}\colon \ST_1^{a,b}\to\ST_1^{a,{2d-2a-b}}$ is an isomorphism.
\end{proof}

\subsection{The Hodge-Lefschetz-primitive part $P^{a,b}$}

For $a \leq 0$ and $b\leq d-a$, define the \emph{Hodge-Lefschetz} (or simply \emph{HL}) \emph{primitive part} $P^{a,b}$ of bidegree $(a,b)$ as
\[P^{a,b} \colon = \ST_1^{a, b} \cap \ker(N^{-a+1}) \cap \ker(\ell^{d-a-b+1}).\]

\begin{prop}
We have
\[P^{a,b} = \bigoplus_{\delta \in X_\f\\ \dims\delta = -a} P^{2a+b}(\delta).\]
\end{prop}
\begin{proof} For $b$ odd, both parts are zero. So suppose $b$ is even. In this case, we have
\[\ST_1^{a, b} \cap \ker(N^{-a+1}) = \bigoplus_{\delta \in X_\f\\ \dims\delta = -a} H^{2a+b}(\delta),\]
and thus
\[P^{a,b} = \ker \Bigl( \bigoplus_{\delta \in X_\f\\
\dims\delta=-a} (\ell^\delta)^{d-a-b+1} \colon   \bigoplus_{\delta \in X_\f\\ \dims\delta = -a} H^{2a+b}(\delta) \longrightarrow  \bigoplus_{\delta \in X_\f\\ \dims\delta = -a} H^{2d-b+2}(\delta)\Bigr).\]

To conclude, note that since $\dims\delta =-a$, we have $d-a-b+1 = d-\dims\delta -(2a+b)+1$, and thus
\[\ker\Bigl(\,(\ell^\delta)^{d-a-b+1} \colon H^{2a+b}(\delta) \rightarrow H^{2d-b+2}(\delta)\Bigr) = P^{2a+b} (\delta),\]
and the proposition follows.
\end{proof}

\subsection{Hodge diamond}

\begin{figure}
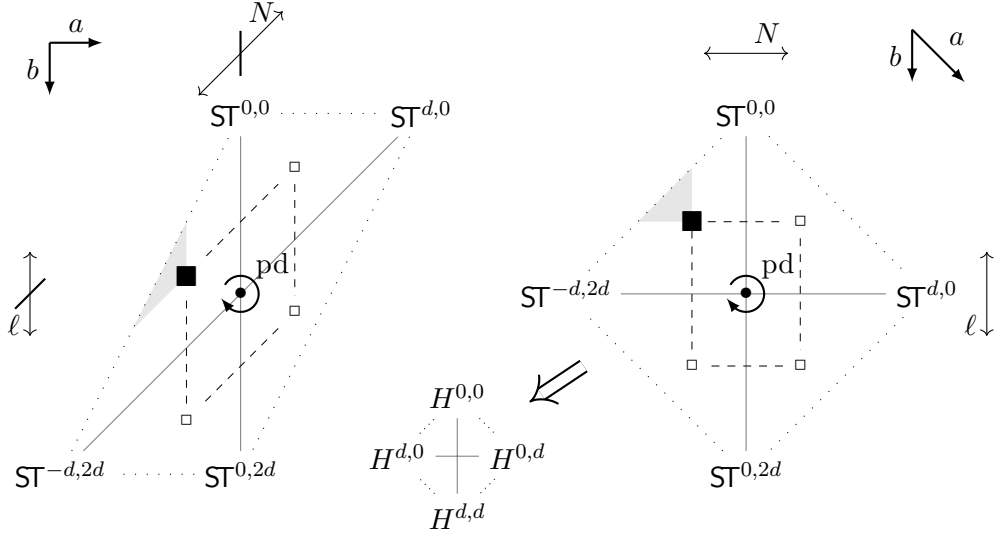

\caption{The tropical Hodge diamond} \label{fig:hodge_diamond}
\hodgediamond
\end{figure}

By what we proved in the previous sections, $\ST_1^{\bul,\bul}$ forms a \emph{Hodge diamond} as illustrated in Figure \ref{fig:hodge_diamond}. Let us explain what this means.

\medskip

We do a change of coordinates $(a,b) \to (a, d-a-b)$ to form
\[\HD^{a,d-a-b}:=\ST_1^{a,b}.\]
The result of this transformation is that the nonzero terms of $\HD^{\bul,\bul}$ now form a rhombus.

The horizontal symmetry maps a piece of $\HD^{\bul,\bul}$ onto an isomorphic piece by applying the monodromy operator $N$ enough number of times. More precisely, for an integer $a\geq0$, we get
\[ N^a\colon \HD^{-a,b} \simto \HD^{a,b} \]
Similarly, the vertical symmetry maps a piece of $\HD^{\bul,\bul}$ onto another isomorphic piece by applying the Lefschetz operator $\ell$ enough number of times: for $b\geq0$,
\[ \ell^b\colon \HD^{a,-b} \simto \HD^{a,b}. \]
In the figure, the four squares in both diagrams, in $\ST_1^{\bul,\bul}$ and $\HD^{\bul,\bul}$, represent four isomorphic pieces. Moreover, the central symmetry corresponds to the Poincaré duality $\HD^{-a,-b}\simeq H^{a,b}$.

To justify the name of tropical Hodge diamond, on the bottom of the figure is shown the tropical Hodge groups to which the Steenbrink diamond degenerates. We show in Theorem \ref{thm:main} that the symmetries described above also hold for the tropical Hodge decomposition, as in the usual Hodge diamond.

\medskip

Finally, notice that $\HD_{\bul,\bul}$ defined by $\HD_{a,b}=\HD^{-a,-b}$ is a HL-structure as defined in Section \ref{sec:differential_HL_structure}, \cf. Definition \ref{defi:HL-structure}.

\subsection{Polarization on $\ST_1^{\bul,\bul}$} We now describe a natural polarization on $\ST_1^{\bul, \bul}$. We prove in the next section that it induces a positive definite symmetric bilinear form on the HL-primitive parts $P^{a,b}$ defined above.

We first define a bilinear form $\psi$ on $\ST_1^{\bul, \bul}$ as follows. Let $x \in \ST_1^{a,b,s}$ and $y \in \ST_1^{a',b',s'}$. Write
\[ x=\sum_{\delta\in X_\f \\ \dims\delta=s} x_\delta \]
with $x_\delta\in H^{a+b-s}(\delta)$, and similarly for $y$. We define the \emph{degree map}
\[ \deg(x\cdot y) := \begin{cases}
\sum_{\delta\in X_\f \\ \dims\delta=s}\deg(x_\delta\cdot y_\delta) & \text{if $s=s'$}, \\
0 & \text{otherwise}.
\end{cases} \]
The bilinear form $\psi$ is then defined by
\begin{equation}
\psi(x,y) := \begin{cases}
\epsilon(a,b) \deg(x \cdot y)  & \textrm{if $a+a'=0$, \ $b+b'=2d$, and $s'=s$,} \\
0 & \textrm{otherwise},
\end{cases}
\end{equation}
where for a pair of integers $a,b$, with $b$ even,
\[\epsilon(a,b) := (-1)^{ a+ \frac b2},\]
and otherwise, if $b$ is odd, $\epsilon(a,b) =1$ (but in this case, we necessarily have $x=0$).

\medskip

The following proposition gathers basic properties of the the bilinear pairing $\psi$.

\begin{prop} \label{prop:properties_psi}
We have for all pairs of elements $x, y \in \bigoplus_{a,b} \ST_1^{a,b}$,
\begin{enumerate}[label=\defaultRoman]
\item \label{enum:psi:i} $\psi(x,y) = (-1)^d \psi(y,x)$.
\item \label{enum:psi:ii} $\psi(Nx, y) + \psi(x, Ny) =0$.
\item \label{enum:psi:iii} $\psi(\ell x, y) + \psi(x, \ell y) =0$.
\item \label{enum:psi:iv} $\psi(\i^*x, y) + \psi(x, \gys y) =0$.
\item \label{enum:psi:v} $\psi(\gys x, y) + \psi(x, \i^* y) =0$.
\item \label{enum:psi:vi} $\psi(\d x, y) + \psi(x, \d y) =0$.
\end{enumerate}
\end{prop}
\begin{proof}
Let $x\in\ST^{a,b,s}_1$ and $y\in\ST_1^{a',b',s'}$ be two nonzero elements for integers $a,b,s, a',b',s'$ with $s\geq\abs{a}$ and $s'\geq\abs{a'}$.
We proceed by verifying the statements point by point.

\medskip

\begin{enumerate}[label=\defaultRoman, leftmargin=0pt]
\item if $a+a'\neq 0$ or $b+b'\neq 2d$, then both $\psi(x,y)$ and $\psi(y,x)$ are zero. Otherwise,
\[ \epsilon(a',b')=\epsilon(-a,2d-b)=(-1)^{-a+\frac{2d-b}2}=(-1)^d\epsilon(a,b). \]
Thus, $\psi(x,y) = (-1)^d \psi(y,x)$.

\medskip

\item Recall that $N$ is a map of bidegree $(2,-2)$. If $a+a'+2\neq 0$ or $b+b'-2\neq 2d$ or $s\neq s'$, then both $\psi(Nx,y)$ and $\psi(x,Ny)$ are zero. Otherwise, $s=s'\geq\abs{a'}=\abs{a+2}$. Thus $Nx\neq 0$. By symmetry, $Ny\neq 0$. Since $N$ acts as identity on both $x$ and $y$, we get $\deg(Nx\cdot y)=\deg(x\cdot Ny)$. Moreover, $\epsilon(a+2,b-2)=-\epsilon(a,b)$. Hence, we obtain $\psi(Nx, y) + \psi(x, Ny) =0$.

\medskip

\item Recall that $\ell$ is a map of bidegree $(0,2)$. If $a+a'\neq 0$ or $b+b'+2\neq 2d$ or $s\neq s'$, then both $\psi(\ell x,y)$ and $\psi(x,\ell y)$ are zero. Otherwise,
\[ \deg(\ell x\cdot y)=\sum_{\delta\in X_\f \\ \dims\delta = s}\deg(\ell^\delta\cdot x_\delta\cdot y_\delta)=\deg(x\cdot \ell y). \]
Moreover, $\epsilon(a,b+2)=-\epsilon(a,b)$. Thus, we get $\psi(\ell x, y) + \psi(x,\ell y)=0$.

\medskip

\item If $a+a'+1\neq0$ or $b+b'\neq2d$ or $s+1\neq s'$, then $\psi(\i^*x,y)$ and $\psi(x,\gys y)$ are both zero. Otherwise, $s'-1=s\geq\abs a=\abs{a'+1}$. Thus, the Gysin map $\gys\colon\ST_1^{a',b',s'}\to\ST_1^{a'+1,b',s'-1}$ has a codomain belonging to $\ST_1^{\bul,\bul}$. This also holds for $\i^*$. Moreover, using Equation \eqref{eqn:i_gys_dual}, we get
\begin{align*}
\deg(\i^*(x)\cdot y)
  &= \sum_{\delta\in X_\f \\ \dims\delta=s+1}\sum_{\gamma\in X_\f \\ \gamma\ssubface\delta}\deg(\sign(\gamma,\delta)\i^*_{\gamma\ssubface\delta}(x_\gamma)\cdot y_\delta) \\
  &= \sum_{\gamma\in X_\f \\ \dims\gamma=s}\sum_{\delta\in X_\f \\ \delta\ssupface\gamma}\deg(x_\gamma\cdot\sign(\gamma,\delta)\gys_{\delta\ssupface\gamma}(y_\delta)) \\
  &=\deg(x\cdot\gys(y)).
\end{align*}
Since $\epsilon(a+1,b)=-\epsilon(a,b)$, we get $\psi(\i^*x,y)+\psi(x,\gys y)=0$.

\medskip

\item The equality $\psi(\gys x, y) + \psi(x, \i^* y) = 0$ follows from the points \ref{enum:psi:i} and \ref{enum:psi:iv}.

\medskip

\item Summing up the equalities in \ref{enum:psi:iv} and \ref{enum:psi:v} gives $\psi(\d x, y) + \psi(x, \d y) =0$ by the definition of the differential $\d$. \qedhere
\end{enumerate}
\end{proof}

\subsection{Induced polarization on the Hodge-Lefschetz primitive parts}

The bilinear form $\psi$ induces a positive definite symmetric bilinear form on the Hodge-Lefschetz primitive parts $P^{a,b}$ that we now describe.

\begin{prop} \label{prop:primitive_polarization}
For integer $a\leq 0$ and $b\leq d-a$, the bilinear form $\psi(\ccdot, \ell^{d-a-b}N^{-a}\ccdot)$ is symmetric positive definite on $P^{a,b}$. More precisely, it is the orthogonal sum of polarizations $(-1)^{a+b/2}Q_{\delta} (\ccdot,\rdot)$ for $\delta \in X_\f$ on the primitive part $P^{2a+b}(\delta)$, under the identification
\[P^{a,b} = \bigoplus_{\delta \in X_\f \\ \dims\delta =-a} P^{2a+b} (\delta).\]
\end{prop}

\begin{proof}
By the properties of the monodromy operator, the isomorphism
\[ N^{-a}\colon \ST_1^{a,2d-2a-b} \to \ST_1^{-a, 2d-b} \]
is the identity map on each isomorphic piece $H^{2a+b}(\delta)$ of $\ST_1^{a,2d-2a-b}$ and $\ST_1^{-a, 2d-b}$, for any $\delta \in X_\f$ with $\dims\delta =-a$.

By definition of $\psi$, for a pair of elements $x, y \in P^{2a+b}(\delta)$, we get
\begin{align*}
\psi(x,\ell^{d-a-b}N^{-a}y)
  &= \psi(x,N^{-a}\ell^{d-a-b}y) \\
  &= \epsilon(a,b) \deg(x\cdot\ell^{d-a-b}y) \\
  &= (-1)^{a+b/2}\deg(x\cdot\ell^{d-a-b}y) \\
  &= (-1)^{a+b/2} Q_\delta(x,y).
\end{align*}
This concludes the proof since $(-1)^{a+b/2}Q_\delta$ is symmetric positive definite by local Hodge-Riemann relations.
\end{proof}

\subsection{Hodge-Lefschetz decomposition}

From the preceding discussion, we get the following theorem. It can be illustrated as in Figure \ref{fig:hodge_diamond}: the black square can be decomposed as the direct sum of primitive parts lying in the grey area.

\begin{thm} \label{thm:ST_1_decomposition}
For each pair of integers $a\leq0$ and $b\leq d-a$, we have the decomposition
\[ \ST_1^{a,b} = \bigoplus_{s,r \geq 0} \ell^r N^s P^{a-2s, b+2s-2r}. \]
\end{thm}

\begin{proof}
We give here a direct proof based on the definition of $\ST_1^{\bul,\bul}$ and using local Lefschetz decompositions. We however note that this theorem can be proved thanks to the general theory of HL-structures, \cf. Proposition \ref{prop:HL_decomposition}. We have
\begin{align*}
\ST_1^{a,b}
  &= \bigoplus_{s \geq 0} \ST_1^{a,b,2s-a} \\
  &= \bigoplus_{s \geq 0} N^s \ST_1^{a-2s,b+2s,2s-a} \\
  &= \bigoplus_{s \geq 0} N^s \bigoplus_{\delta \in X_\f \\ \dims\delta = 2s-a} \underbrace{H^{2a+b-2s}(\delta)}_{\subseteq \ST_1^{a-2s, b+2s}} \\
  &= \bigoplus_{s \geq 0} N^s \bigoplus_{\delta \in X_\f \\ \dims\delta = 2s-a} \bigoplus_{r \geq 0} \ell^r \underbrace{P^{2a+b-2s-2r}(\delta)}_{\subseteq\ST_1^{a-2s, b+2s-2r}} \\
  &= \bigoplus_{r,s \geq 0} \ell^r N^s \bigoplus_{\delta \in X_\f \\ \dims\delta = 2s-a} P^{2a+b-2s-2r}(\delta), \\
\ST_1^{a,b}
  &= \bigoplus_{r,s \geq 0} \ell^r N^s P^{a-2s, b+2s-2r}. \qedhere
\end{align*}
\end{proof}

\subsection{Polarized Hodge-Lefschetz structure on the cohomology of $(\ST_1^{\bul, b}, \d)$}

In this section we show that when passing to the cohomology, the operators $N$, $\ell$ and the polarization $\psi$ induce a polarized Hodge-Lefschetz structure on the cohomology groups $L^{a,b} := H^{a}(\ST_1^{\bul, b}, \d)$. By this we mean the following.

\medskip

First, since $[\d, \ell] =0$ and $[\d,N]=0$, we get induced maps
\[\ell \colon L^{a,b} \to L^{a, b+2}\]
and
\[N\colon L^{a,b} \to L^{a+2, b-2}\]
such that $[N,\ell] = 0$.

Second, since $\psi(\d\ccdot,\rdot)=-\psi(\ccdot,\d\ccdot)$, we get an induced pairing
\[ \psi\colon L^{\bul,\bul}\times L^{\bul,\bul} \to \corps. \]

Moreover, we will show that
\begin{thm} \label{thm:tropical_differential_HL_structure}We have the following properties.
\begin{itemize}
\item For $a\leq 0$, the map $N^{-a} \colon L^{a,b} \to L^{-a,b}$ is an isomorphism.
\item For $a,b$ with $d-a-b \geq 0$, the map $\ell^{d-a-b} \colon L^{a,b} \to L^{a, 2d-2a-b}$ is an isomorphism.
\item Let $a\leq 0$ and $b \leq d-a$ and denote by $L_0^{a,b}$ the HL-primitive part of $L^{a,b}$ defined by
\[ L_0^{a,b} = L^{a,b} \cap \ker(N^{-a+1}) \cap \ker(\ell^{d-a-b+1}).\]
The polarization $\psi$ induces a symmetric bilinear form $\psi(\ccdot,\ell^{d-a-b}N^{-a} \ccdot)$ which is positive definite on $L_0^{a,b}$.
\item For integers $a\leq 0$ and $b\leq d-a$, we have the decomposition
\[ L^{a,b} = \bigoplus_{s,r \geq 0} \ell^r N^s L^{a-2s, b+2s-2r}_0.\]
Moreover, this decomposition is orthogonal for the pairing $\psi(\ccdot,\ell^{d-a-b}N^{-a}\ccdot)$.
\end{itemize}
\end{thm}

\begin{proof}
We can use the change of coordinates $\HD_{-a,-(d-a-b)}:=\ST_1^{a,b}$ for any integers $a,b$. By the previous discussion, we get that $(\HD_{-a,-(d-a-b)},N,\ell,\phi,\d)$ is a differential polarized HL-structure in the sense of Section~\ref{sec:differential_HL_structure}. We can therefore apply Theorem \ref{thm:differential_HL-structure} to conclude.
\end{proof}

\subsection{Singular monodromy operator}

We define the \emph{singular tropical monodromy operator} as follows. For each face $\delta\in X_\f$, let $o_\delta$ be a point of $\TT\delta\subseteq N_\R$. For instance, we can take $o_\delta$ to be the centroid of $\delta$. Let $\gamma$ be a face in $X_\f$ of dimension $q$. Then we set
\[ \begin{array}{rccc}
N\colon & C^{p,q}(X) & \to & C^{p-1,q+1}(X), \\
   & \alpha\in F^p(\gamma) & \mapsto & \underset{\substack{\delta\in X_\f\\
   \delta \ssupface \gamma}}{\bigoplus} \,\,\alpha\bigr(\ccdot \wedge (o_\delta-o_\gamma)\bigr).
\end{array}
\]

\begin{thm} \label{thm:tropical_monodromy}
The singular tropical monodromy operator corresponds to the monodromy operator on $\ST_1^{\bul,\bul}$ via the isomorphism of Theorem \ref{thm:steenbrink}.
\end{thm}

\begin{proof}
This is a direct consequence of Proposition \ref{prop:com:N_commutes}.
\end{proof}

\begin{cor} The operator $N\colon H^{p,q} \to H^{p-1, q+1}$ coincides with the eigenwave operator $\phi$ constructed in~\cite{MZ}.
\end{cor}
\begin{proof} One can show that the singular tropical monodromy operator coincides with the eigenwave operator. So the statement follows from Theorem~\ref{thm:tropical_monodromy}.
\end{proof}

\begin{remark} It follows that the monodromy operator $N$ coincides as well with the monodromy operator defined on the level of Dolbeault cohomology groups~\cites{Liu19}. This is a consequence of \cite{Jell19} which relates the eigenwave operator to the monodromy on Dolbeault cohomology.
\end{remark}

\subsection{HL, HR and Monodromy for tropical cohomology}

From the result in the previous section, we get the following theorem.

\begin{thm}\label{thm:main} We have
\begin{itemize}
\item \emph{(Weight-monodromy conjecture)} For $q>p$ two non-negative integers, we get an isomorphism
\[N^{q-p} \colon H_{\trop}^{q,p}(X) \to H_{\trop}^{p,q}(X).\]

\item \emph{(Hard Lefschetz)} For $p+q \leq d/2$, the Lefschetz operator $\ell$ induces an isomorphism
\[\ell^{d- p-q}\colon H_{\trop}^{p,q}(X) \to H^{d-q, d-p}_{\trop}(X).\]

\item \emph{(Hodge-Riemann)} The pairing $(-1)^p \bigl< \ccdot,\, \ell^{d-p-q} N^{q-p} \ccdot \bigr>$ induces a positive-definite pairing on the primitive part $P^{q,p}$, where $\bigl< \ccdot,\rdot \bigr>$ is the natural pairing

\[\bigl< \ccdot, \rdot\bigr> \colon H^{q,p}(X) \otimes H^{d-q,d-p}(X) \to H^{d,d}(X) \simeq \mathbb Q.\]
\end{itemize}
\end{thm}

\subsection{Hodge index theorem for tropical surfaces}
In this final part of this section, we explain how to deduce the Hodge index theorem for tropical surfaces from Theorem~\ref{thm:main}.
\begin{proof}[Proof of Theorem~\ref{thm:hodgeindex}] The primitive part decomposition theorem implies that we can decompose the cohomology group $H^{1,1}_\trop(X, \Q)$ into the direct sum
\[H^{1,1}_\trop(X, \Q) = \ell\, H^{0,0}_\trop(X, \Q) \oplus N\, H^{2,0}_\trop(X, \Q) \oplus H^{1,1}_\prim(X, \Q),\]
where $H^{1,1}_\prim(X, \Q) = \ker(\ell) \cap \ker(N)$. Moreover, by Hodge-Riemann, the pairing is positive definite on
$\ell\, H^{0,0}_\trop(X, \Q)$ and $N\, H^{2,0}_\trop(X, \Q)$, and it is negative definite on $H^{1,1}_\prim(X, \Q)$.

By Poinca\'re duality, we have $H^{2,0}_\trop(X, \Q) \simeq H^{0,2}_\trop(X, \Q)$. By definition of tropical cohomology, since $F^0$ is the constant sheaf $\Q$, we have $H^{0,2}_\trop(X, \Q) \simeq H^2(\X,\Q)$. Moreover, we have $H^{0,0}_\trop(X, \Q) \sim \Q$.
We conclude that the signature of the intersection pairing is given by $(1+ b_2, h^{1,1}- 1-b_2)$, as stated by the Hodge index theorem.
\end{proof}


\section{Hodge-Lefschetz structures}
\label{sec:differential_HL_structure}

In this section, we prove that the homology of a differential Hodge-Lefschetz structure is a polarized Hodge-Lefschetz structure, thus finishing the proof of Theorem~\ref{thm:tropical_differential_HL_structure} and the main theorem of the previous section.

\medskip

The main references for differential Hodge-Lefschetz structures are Saito's work on Hodge modules~\cite{Saito}, the paper by Guill\'en and Navarro Aznar on invariant cycle theorem~\cite{GNA90} and the upcoming book by Sabbah and Schnell~\cite{SabSch} on mixed Hodge modules, to which we refer for more information. We note however that our set-up is slightly different, in particular, our differential operator is skew-symmetric with respect to the polarization. The proof in~\cite{SabSch} is presented also for the case of mono-graded complexes and is based on the use of representation theory of $\mathrm{SL}_2(\R)$ as in~\cite{GNA90}. So we give a complete self-contained proof of the main theorem. Our proof is direct and does not make any recourse to representation theory, although it might be possible to recast in the language of representation theory the final combinatorial calculations we need to elaborate.

\medskip

In order to simplify the computations, it turns out that it will be more convenient to work with coordinates different from the ones in the previous section, and that is what we will be doing here. The new coordinates are compatible with the ones in~\cite{SabSch}.

\subsection{Hodge-Lefschetz structure}

\begin{defi} \label{defi:HL-structure}
A \emph{(bigraded) Hodge-Lefschetz structure}, or more simply an \emph{HL-structure}, $(H_{\bul,\bul}, N_1, N_2)$ is a bigraded finite dimensional vector space $H_{\bul,\bul}$ endowed with two endomorphisms $N_1$ and $N_2$ of respective bidegree $(-2,0)$ and $(0,-2)$ such that,
\begin{itemize}
\item $[N_1,N_2]=0$,
\item for any pair of integers $a,b$ with $a\geq 0$, we have $N_1^a\colon H_{a,b}\to H_{-a,b}$ is an isomorphism,
\item for any pair of integers $a, b$ with $b\geq0$, we have $N_2^b\colon H_{a,b}\to H_{a,-b}$ is an isomorphism.
\end{itemize}
\end{defi}

For the rest of this section, we fix a HL-structure $(H_{\bul,\bul}, N_1, N_2)$.

\medskip

We define the \emph{primitive parts} of $H_{\bul,\bul}$ by
\[\forall \, a,b\geq 0,\qquad P_{a,b}:=H_{a,b}\cap\ker(N_1^{a+1})\cap\ker(N_2^{b+1}). \]

We have the following decomposition into primitive parts.

\begin{prop} \label{prop:HL_decomposition}
For each pair of integers $a,b\geq0$, we have the decomposition
\[ H_{a,b} = \bigoplus_{s,r \geq 0} N_1^r N_2^s P_{a+2r, b+2s}. \]
\end{prop}

Denote by $H^+_{\bul,\bul}$ the bigraded subspace
\[ H^+_{\bul,\bul}=\bigoplus_{a,b\geq0}H_{a,b}. \]
Define $N_1^\dual\colon H^+_{\bul,\bul}\to H^+_{\bul,\bul}$ of bidegree $(2,0)$, by

\[ N_1^\dual\rest{H_{a,b}}=\Bigl(\bigl(N_1\rest{H_{a+2,b}}\bigr)^{a+2}\Bigr)^{-1}\Bigl(N_1\rest{H_{a,b}}\Bigr)^{a+1}\qquad\forall a,b\geq0. \]

\medskip

It verifies the following properties. Let $a,b\geq0$. Then,

\[ N_1^\dual N_1\rest{H_{a+2,b}}=\id\rest{H_{a+2,b}} \quad\myand\quad \ker(N_1^\dual)\cap H_{a,b}=\ker(N_1^{a+1})\cap H_{a,b}. \]

In particular, $N_1\colon H_{a+2,b}\to H_{a,b}$ is injective, and $N_1^\dual\colon H_{a,b}\to H_{a+2,b}$ is surjective.

From this, we can deduce that the map $N_1N_1^\dual\rest{H_{a,b}}$ is the projection along $\ker(N_1^{a+1})\cap H_{a,b}$ onto $N_1H_{a+2,b}$, \ie, for the decomposition
\[H_{a,b} =  N_1H_{a+2,b} \oplus \Bigl(\ker(N_1^{a+1})\cap H_{a,b}\Bigr).\]

\medskip

We define $N_2^\dual$ in a similar way, and note that it verifies similar properties. In particular, $N_2N_2^\dual\rest{H_{a,b}}$ is the projection along $\ker(N_2^{b+1})\cap H_{a,b}$ onto $N_2H_{a,b+2}$ for the decomposition
\[H_{a,b} =  N_2H_{a,b+2} \oplus \Bigl(\ker(N_2^{b+1})\cap H_{a,b}\Bigr).\]

\medskip

Moreover, from the commutativity of $N_1$ and $N_2$, more precisely, from
\[ [N_1^{a+2},N_2]\rest{H_{a+2,b+2}}=[N_1,N_2^{b+2}]\rest{H_{a+2,b+2}}=[N_1^{a+2},N_2^{b+2}]\rest{H_{a+2,b+2}}=0, \]
we infer that
\[ [N_1^\dual,N_2]\rest{H_{a,b+2}}=[N_1,N_2^\dual]\rest{H_{a+2,b}}=[N_1^\dual,N_2^\dual]\rest{H_{a,b}}=0. \]
In particular,
\[ N_1N_1^\dual N_2N_2^\dual=N_2N_2^\dual N_1N_1^\dual=N_1N_2N_1^\dual N_2^\dual. \]
Notice that $N_1N_2N_1^\dual N_2^\dual\rest{H_{a,b}}$ is a projection onto $N_1N_2H_{a+2,b+2}$.

\medskip

We come back to the proof of the proposition.
\begin{proof}[Proof of Proposition~\ref{prop:HL_decomposition}]
Clearly, $N_1\colon H_{a+2,b}\to H_{a,b}$ is injective for $a\geq0$. Similarly for $N_2$. Therefore, proceeding by induction, it suffices to prove that for any pair of integers $a,b\geq0$, we have
\begin{gather*}
H_{a,b}= \bigl(N_1H_{a+2,b} + N_2H_{a,b+2}\bigr) \oplus P_{a,b} \quad  \myand \quad N_1H_{a+2,b} \cap N_2H_{a,b+2} = N_1N_2H_{a+2,b+2}.
\end{gather*}
This can be deduced from the following lemma applied to $N_1N_1^\dual\rest{H_{a,b}}$ and $N_2N_2^\dual\rest{H_{a,b}}$ and the corresponding decompositions of $H_{a,b}$.
\end{proof}

\begin{lemma}
Let $V$ be a finite dimensional vector space. For $i\in\{1,2\}$, suppose we have decompositions $V = I_1 \oplus K_1 = I_2 \oplus K_2$ for subspaces $I_1, K_1, I_2, K_2 \subseteq V$. Denote by $\pi_i\colon V\to V$ the projection of $V$ onto $I_i$ along $K_i$. Then, if $\pi_1$ and $\pi_2$ commute, we get
\[ V=I_1\cap I_2 \woplus I_1\cap K_2 \woplus K_1\cap I_2 \woplus K_1\cap K_2\quad \myand \quad I_1\cap I_2=\Im(\pi_1\circ\pi_2).\]
\end{lemma}

\begin{proof}
Since $\pi_1$ and $\pi_2$ commute, $\pi_2(I_1)\subseteq I_1$. Thus, $\pi_2$ can be restricted to $I_1$. This restriction is a projection onto $I_2 \cap I_1$ along $K_2 \cap I_1$. Thus,
\[ I_1=I_2 \cap I_1 \woplus K_2 \cap I_1 \quad\myand\quad I_1\cap I_2=\Im(\pi_2\circ\pi_1). \]
Applying the same argument to the commuting projections $\id-\pi_1$ and $\id-\pi_2$, we get that
\[ K_1=K_2 \cap K_1 \woplus I_2 \cap K_1. \]
Hence,
\[ V = I_1 \oplus K_1 = I_1\cap I_2 \woplus I_1\cap K_2 \woplus K_1\cap I_2 \woplus K_1\cap K_2. \qedhere \]
\end{proof}

\subsection{Polarization $\psi$ and its corresponding scalar product $\phi$}

A \emph{polarization} on the HL-structure $H_{\bul,\bul}$ is a bilinear form $\psi(\ccdot,\rdot)$ such that:
\begin{itemize}
\item for any integers $a,b,a',b'$ and any elements $x\in H_{a,b}$ and $y\in H_{a',b'}$, the pairing $\psi(x,y)$ is zero unless $a+a'=b+b'=0$,
\item for any $a,b\geq0$, the pairing $\psi(\ccdot,N_1^aN_2^b\ccdot)$ is symmetric definite positive on $P_{a,b}$,
\item $N_1$ and $N_2$ are skew-symmetric for $\psi$.
\end{itemize}

For the rest of this section, we assume that $H_{\bul,\bul}$ is endowed with a polarization $\psi$.

\medskip

Denote as in the previous section $H^+_{\bul,\bul}$ the bigraded subspace
\[ H^{+}_{\bul,\bul}=\bigoplus_{a,b\geq0}H_{a,b}. \]

We define the operator $\w$ on $H_{\bul,\bul}$ as follows.

Let $a,b\geq0$ and $r,s\geq0$ be four integers. Let $x\in P_{a,b}$. Then, we set
\[ \w N_1^rN_2^sx:=\begin{cases}
(-1)^{r+s}\frac{r!}{(a-r)!}\frac{s!}{(b-s)!}N_1^{a-r}N_2^{b-s}x & \text{if $r\leq a$ and $s\leq b$} \\
0 & \text{otherwise.}
\end{cases} \]

Using the decomposition into sum of primitive parts from the previous section, one can verify that $\w$ is well-defined.

We will provide some explanations about the choice of the operator $\w$, in particular about the choice of the renormalization factors, at the beginning of Section \ref{subsec:Laplacian_commute}.

\medskip

Set $\phi(\ccdot,\rdot):=\psi(\ccdot,\w\ccdot)$.

\begin{prop}
We have the following properties.
\begin{itemize}
\item The operator $\w$ is invertible with its inverse verifying
\[ \w^{-1}N_1^rN_2^sx=(-1)^{a+b}\w N_1^rN_2^sx. \]
\item The pairing $\phi$ is symmetric positive definite.
\item The decomposition into primitive parts of $H_{\bul,\bul}$ is orthogonal for $\phi$. More precisely, let $a,b,r,s\geq0$ and $a',b',r',s'\geq0$ be eight integers and let $x\in P_{a,b}$ and $y\in P_{a',b'}$. Then,
\[ \phi(N_1^rN_2^sx, N_1^{r'}N_2^{s'}y)=0 \text{ unless $r=r', s=s', a=a'$ and $b=b'$}. \]
\end{itemize}
\end{prop}

\begin{proof}
The first point is an easy computation. For the two other points, take $x\in P_{a,b}$ and $y\in P_{a',b'}$. We can assume that $r\leq a$, $s\leq b$, $r'\leq a'$ and $s'\leq b'$. Set \[ C=\frac{r'!\,s'!}{(a'-r')!\,(b'-s')!}\,. \]
We have
\begin{align*}
\phi(N_1^rN_2^sx, N_1^{r'}N_2^{s'}y)
  &= \psi(N_1^rN_2^sx, \w N_1^{r'}N_2^{s'}y) \\
  &= (-1)^{r'+s'}C\psi(N_1^rN_2^sx, N_1^{a'-r'}N_2^{b'-s'}y) \\
  &= (-1)^{a'+b'}C\psi(N_1^{a'-r'+r}N_2^{b'-s'+s}x, y) \\
  &= (-1)^{r'+s'+r+s}C\psi(x, N_1^{a'-r'+r}N_2^{b'-s'+s}y).
\end{align*}
Suppose the value of the above pairing is nonzero. Then, since $x\in H_{a,b}$ and
\[ N_1^{a'+r-r'}N_2^{b'+s-s'}y\in H_{-a'+2r'-2r,-b'+2s'-2s}, \]
we should have $a-a'+2r'-2r=0$. Moreover, since $x\in\ker(N_1^{a+1})$ and $y\in\ker(N_1^{a'+1})$, we get
\[ a'-r'+r\leq a \quad\myand\quad a'-r'+r\leq a'. \]
Hence, we get
\[ r-r'\leq a-a'=2(r-r') \quad\myand\quad r-r'\leq 0,\]
which in turn implies that $r-r'=a-a'=0$, \ie, $r=r'$ and $a=a'$. By symmetry, we get $s=s'$ and $b=b'$. We infer that
\[ \phi(N_1^rN_2^sx, N_1^rN_2^sy) = C\psi(x, N_1^aN_2^by), \]
which concludes the proof since $\psi(\ccdot,N_1^aN_2^b\ccdot)$ is symmetric positive definite on $P_{a,b}$.
\end{proof}

\subsection{Differential polarized HL-structure}

A \emph{differential} on the polarized HL-structure $H_{\bul,\bul}$ is an endomorphism $\d$ of bidegree $(-1,-1)$ such that
\begin{itemize}
\item $\d^2=0$,
\item $[\d,N_1]=[\d,N_2]=0$,
\item $\d$ is skew-symmetric for $\psi$.
\end{itemize}

In the rest of this section, we fix a differential $\d$ on $H_{\bul,\bul}$.

\begin{prop}
The following holds for any pair of integers $a,b\geq0$:
\[ \d P_{a,b} \subseteq P_{a-1,b-1} \oplus N_1 P_{a+1,b-1} \oplus N_2 P_{a-1,b+1} \oplus N_1N_2 P_{a+1,b+1}, \]
where, by convention, $P_{a',b'}=0$ if $a'<0$ or $b'<0$.
\end{prop}

\begin{proof}
Since $\d$ and $N_1$ commute,
\[ \d P_{a,b} \subseteq \ker(N_1^{a+1})\cap H_{a-1,b-1} = \bigoplus_{r\geq 0} \ N_2^r P_{a-1,b-1+2r} \oplus N_1N_2^r P_{a+1,b-1+2r}. \]
Notice that this is true even if $a=0$ or $b=0$. We conclude the proof by combining this result with the symmetric result for $N_2$.
\end{proof}

Using this proposition, if $x\in P_{a,b}$, then for $i,j\in\{0,1\}$, we define $x_{i,j}\in P_{a-1+2i,b-1+2j}$ such that
\[ \d x = \sum_{i,j\in\{0,1\}}N_1^iN_2^jx_{i,j}. \]
Since $\d$ commutes with $N_1$ and $N_2$, more generally, for any pair of integers $r,s\geq0$, we get
\[ \d N_1^r N_2^s x = \sum_{i,j\in\{0,1\}}N_1^{r+i}N_2^{s+j}x_{i,j}. \]

With this notation, the fact that $\d^2=0$ can be reformulated as follows.

\begin{prop} \label{prop:dd=0}
With the above notations, for any pair $k,l\in\{0,1,2\}$,
\[ N_1^{\epsilon_1}N_2^{\epsilon_2}\sum_{i,j,i',j'\in\{0,1\} \\ i+i'=k \\ j+j'=l}(x_{i,j})_{i',j'}=0. \]
where,
\[ \epsilon_1 = \begin{cases}
  1 & \text{if $a=0$ and $k=1$} \\
  0 & \text{otherwise,}
\end{cases} \quad\myand\quad
\epsilon_2 = \begin{cases}
  1 & \text{if $b=0$ and $l=1$} \\
  0 & \text{otherwise.}
\end{cases} \]
\end{prop}

\begin{proof}
From the previous discussion, we get that
\[ \d^2 x \in \bigoplus_{k,l\in\{0,1,2\}} N_1^kN_2^lP_{a-2+2k,b-2+2l}. \]
Moreover, the piece corresponding to a fixed $k$ and $l$ is
\[ \sum_{i,j,i',j'\in\{0,1\} \\ i+i'=k \\ j+j=l}N_1^kN_2^l(x_{i,j})_{i',j'}. \]
Since, $\d^2=0$, this last sum is zero:
\[ N_1^kN_2^l\sum_{i,j,i',j'\in\{0,1\} \\ i+i'=k \\ j+j=l}(x_{i,j})_{i',j'}=0. \]

We know that, for any $b'$, $N_1^k\rest{P_{a-2+2k,b'}}$ is injective if $k\leq a-2+2k$. This is also true if $a-2+2k<0$ because, in this case, the domain is trivial. In both cases, we can remove $N_1^k$ in the above equation. The last case is
\[ 0\leq a-2+2k < k. \]
This implies $k=1$ and $a=0$, in which case we cannot a priori remove the $N_1^k=N_1$ in the previous equation. This concludes the proof of the proposition.
\end{proof}

\subsection{The Laplace operator $\Lap$} The \emph{Laplace operator} $\Lap\colon H_{\bul,\bul}\to H_{\bul,\bul}$ is defined by
\[ \Lap := \d \dd + \dd \d \]
where $\dd\colon H_{\bul,\bul}\to H_{\bul,\bul}$ denotes the \emph{codifferential} of bidegree $(1,1)$ defined by
\[ \dd:= -\w^{-1} \d \w. \]

Note that $\Lap$ is of bi-degree $(0,0)$. The following summarizes basic properties one expect from the Laplace operator and its corresponding Hodge decomposition.
\begin{prop} \label{prop:Lap}
We have the following properties.
\begin{itemize}
\item $\dd$ is the adjoint of $\d$ with respect to $\phi$.
\item $\Lap$ is symmetric for $\phi$.
\item On each graded piece $H_{a,b}$, we have the following orthogonal decomposition for $\phi$
\[ \ker(\d)\cap H_{a,b} = \ker(\Lap)\cap H_{a,b} \wooplus \Im(\d)\cap H_{a,b}. \]
\end{itemize}
\end{prop}

\begin{proof}To get the first point, note that for any $x,y\in H_{\bul, \bul}$, we have
\[ \phi(x,\dd y)=-\psi(x,\w \w^{-1}\d\w y)=\psi(\d x,\w y)=\phi(\d x, y). \]

It follows that the adjoint of $\Lap=\d\dd+\dd\d$ is $\Lap$ itself.

\medskip

Last point follows directly from the first two statements, as usual. Denote by $V:=(H_{a,b},\phi)$ the vector space $H_{a,b}$ endowed with the inner product $\phi$. In order to simplify the notations, we will write $\Im(\d)$ instead of $V\cap\Im(\d)$, and similarly for other intersections with $V$. We have to show that
\[ \ker(\d)=\ker(\Lap)\ooplus\Im(\d).\]
From the decomposition $V=\Im(\d)^\perp \oplus \Im(\d)$, and the inclusion $\Im(\d)\subseteq\ker(\d)$, we get
\[ \ker(\d)=\Im(\d)^\perp\cap\ker(\d) \wooplus \Im(\d). \]

By adjunction, we get $\Im(\d)^\perp=\ker(\dd)$. To conclude, note that $\ker(\Lap)=\ker(\d)\cap\ker(\dd)$. Indeed, if $x\in\ker(\Lap)$, then
\[ 0 = \phi(x,\Lap x) = \phi(x, \d\dd x)+\phi(x, \dd\d x) = \phi(\dd x, \dd x)+\phi(\d x, \d x), \]
which, by positivity of $\phi$, implies that $\dd x=\d x=0$, \ie, $x\in\ker(\d)\cap\ker(\dd)$. The other inclusion $\ker(\d)\cap\ker(\dd) \subseteq \ker(\Lap)$ trivially holds.  \qedhere
\end{proof}

\begin{prop}\label{prop:commutativity} The Laplace operator commutes with $N_1$ and $N_2$, \ie, we have $[\Lap,N_1]=[\Lap,N_2]=0$.
\end{prop}
The proof of this theorem is given in Section \ref{subsec:Laplacian_commute}.

\subsection{Polarized Hodge-Lefschetz structure on the cohomology of $H_{\bul,\bul}$}

In this section, we show that, when passing to the cohomology, the operators $N_1$, $N_2$ and the polarization $\psi$ induce a polarized Hodge-Lefschetz structure on the cohomology groups
\[ L_{a,b} := \frac{\ker(\d\colon H_{a,b}\to H_{a-1,b-1})}{\Im(\d\colon H_{a+1,b+1}\to H_{a,b})}. \]
By this we mean the following.

First, since $[\d, N_1] =0$ and $[\d,N_2]=0$, we get induced maps
\[ N_1 \colon L_{a,b} \to L_{a-1, b} \quad\myand\quad N_2\colon L_{a,b} \to L_{a, b-2}. \]

Second, since $\psi(\d\ccdot,\rdot)=-\psi(\ccdot,\d\ccdot)$, we get a pairing
\[ \psi\colon L_{\bul,\bul}\times L_{\bul,\bul} \to \R. \]

Then, we will show that
\begin{thm} \label{thm:differential_HL-structure}
For the corresponding induced map, $(L_{\bul,\bul}, N_1, N_2, \phi)$ is a polarized HL-structure. Moreover, the decomposition by primitive parts is induced by the corresponding decomposition on $H_{\bul,\bul}$.
\end{thm}

\begin{proof}
By Proposition \ref{prop:Lap}, the kernel of the Laplace operator gives a section of the cohomology $L_{\bul,\bul}$. Moreover, since the Laplacian commutes with $N_1$ and $N_2$, both $N_1$ and $N_2$ preserve the kernel of the Laplacian. Consequently, it suffices to prove that $\ker(\Lap)$ is a polarized HL-substructure of $H_{\bul,\bul}$. By abuse of notation, we denote by $L_{a,b}$ the kernel $\ker(\Lap\colon H_{a,b}\to H_{a,b})$.

Clearly, $N_1$ and $N_2$ commutes. Let $a\geq 0$ and $b\in\Z$. Let us prove that $N_1^a$ induces an isomorphism $L_{a,b}\simto L_{-a,b}$. Since $\Lap$ is symmetric for $\phi$,
\[ H_{a,b}=\ker(\Lap)\cap H_{a,b} \wooplus \Im(\Lap)\cap H_{a,b}. \]
Since $N_1$ preserves $\ker(\Lap)$ and $\Im(\Lap)$, the isomorphism $N_1^a\colon H_{a,b}\simto H_{-a,b}$ induces isomorphisms for both parts $\ker(\Lap)$ and $\Im(\Lap)$. In particular, we get an isomorphism $N_1^a\colon L_{a,b}\simto L_{-a,b}$. We use the symmetric argument for $N_2$.

The primitive part of $L_{a,b}$ is, by definition,
\[ L_{a,b}\cap\ker(N_1^{a+1})\cap\ker(N_2^{b+1})=L_{a,b}\cap P_{a,b}. \]
Thus, the decomposition by primitive parts is the one induced by the decomposition of $H_{a,b}$. Finally, the properties about $\psi$ on $L_{a,b}$ comes from the analogous properties on $H_{a,b}$.
\end{proof}

\subsection{Commutative property of the Laplacian}
\label{subsec:Laplacian_commute}

{
\renewcommand{\r}{{r}}
\renewcommand{\a}{{ a}}
\newcommand{\0}{{0}}

In this section, we prove Proposition~\ref{prop:commutativity}, that $[\Lap,N_1]=[\Lap,N_2]=0$.

\medskip

Let us start by justifying our definition of the operator $\w$.
\begin{remark} We feel necessary to provide an explanation of why we are not using the simpler more natural operator $\~\w$ which is defined by taking values for $a,b,r,s\geq0$ and $x\in P_{a,b}$ given by
\[ \~\w N_1^rN_2^sx=(-1)^{r+s}N_1^{a-r}N_2^{b-s}x, \]
where negative powers of $N_1$ and $N_2$ are considered to be zero.

With this operator, the bilinear form $\~\phi:=\psi(\ccdot,\~\w\rdot)$ will be still positive definite, and we can define the corresponding Laplace operator $\~\Lap$. However, this Laplacian will not commute with $N_1$ and $N_2$ in general. In order to explain where the problem comes from, let us try to proceed with a natural, but \emph{wrong}, proof that this Laplacian commutes $N_1$. More precisely, setting $\~\dd = -\~\w^{-1}\d\~\w$ the adjoint of $\d$, let us try to prove that $[\~\dd,N_1]=0$. Assume $0\leq r\leq a$ and $0\leq s\leq b$. Then we can write
\begin{align*}
\~\dd N_1^rN_2^sx
  &= -\~\w^{-1}\d\~\w N_1^rN_2^sx \\
  &= (-1)^{r+s+1}\~\w^{-1}\d N_1^{a-r}N_2^{b-s}x \\
  &= \sum_{i,j\in\{0,1\}}(-1)^{r+s+1}\~\w^{-1}N_1^{a-r+i}N_2^{b-s+j}x_{i,j}.
\end{align*}
We recall that $x_{i,j}\in P_{a-1+2i,b-1+2j}$. Notice that $a-1+2i-(a-r+i)=r-1+i$. Thus,
\begin{align*}
\~\dd N_1^rN_2^sx
  &= \sum_{i,j\in\{0,1\}}(-1)^{r+s+1}(-1)^{r-1+i+s-1+j}N_1^{r-1+i}N_2^{s-1+j}x_{i,j} \\
  &= \sum_{i,j\in\{0,1\}}(-1)^{i+j+1}N_1^{r-1+i}N_2^{s-1+j}x_{i,j}.
\end{align*}
Then, $\~\dd N_1(N_1^rN_2^sx)$ can be obtained replacing $r$ by $r+1$ and similarly for $N_1\~\dd(N_1^rN_2^sx)$. Thus, we would conclude that $N_1$ and $\~\dd$ commute!

Where is the mistake? First, notice that the formula for the derivation can only be applied if $a-r$ and $b-s$ are both non-negative. Thus, if $r=a$, we cannot just replace $r$ by $r+1$ to compute $\~\dd N_1(N_1^rN_2^sx)$. Second, notice that $0=N_1N_1^{-1}\neq N_1^{0}$. If $r=0$ and $i=0$, we get $r-1+i=-1$. Multiplying the corresponding term by $N_1$ is not equivalent to simply replacing $r$ by $r+1$. Notice moreover that we did not use the condition $\d^2=0$ at any moment. Thus, roughly speaking, everything goes well except on the boundaries $r=a, s=b, r=0$ and $s=0$. These problems still exist for the operator $\~\Lap$, which justify why $\~\Lap$ does not commute in general with $N_i$.

\medskip

The solution to remedy these boundary issues is to modify $\~\w$ in such a way that some zero factors on the boundary solve the above evoked problem: if the term $N_1^{r-1+i}N_2^{s-1+j}x_{i,j}$ is multiplied by some zero factor due to the new operator when $r=0$ and $i=0$, then multiplying by $N_1$ corresponds to replacing $r$ by $r+1$. To find the right form of the Laplacian, one can use tools of representation theory, where Lie brackets can be used to study commutation properties. This is what is done in \cites{GNA90, SabSch}. These proofs however take place in a slightly different context. For example, our differential is skew symmetric, and not symmetric, and we work in the real instead of the complex world. Nevertheless, one should be able to adapt their proofs to our setting. However, instead of doing this, we decided to provide a completely elementary proof. The arguments which follow are somehow technical but no knowledge in representation theory is needed. \end{remark}

\medskip

We now begin the proof of the commutativity of the Laplacian $\Lap$ with the operators $N_1$ and $N_2$.

\medskip

\paragraph{{\bf Convention.}} In order to simplify the calculations, we avoid to write symmetric factors by using the following notations. If $A,B,C$ are some sets, then we follow the following rules to extend the operations between these sets coordinatewise to operations involving their Cartesian products $A^2,B^2$ and $C^2$: if $\zeta\colon A\to C$ and $\theta\colon A\times B\to C$ are some operations, then we extend them as follows.
\[ \begin{array}{rrrll}
\zeta\colon  & A^2\to C^2, \qquad             & (a_1,a_2)                       & \mapsto (\zeta(a_1), \zeta(a_2)),                      \\
\theta\colon & A \times B^2 \to C^2, \qquad   & \bigl(a,(b_1,b_2)\bigr)         & \mapsto (\theta(a,b_1), \theta(a,b_2)), \myand \\
\theta\colon & A^2 \times B^2 \to C^2, \qquad & \bigl((a_1,a_2),(b_1,b_2)\bigr) & \mapsto (\theta(a_1, b_1), \theta(a_2, b_2)).
\end{array} \]
Moreover, if $A$ is endowed with a product structure, we define $\m\colon A^2\to A$ by $\m(a_1,a_2)=a_1a_2$.

\begin{example}
Let us give some examples. Let $a=(a_1, a_2)$ and $r=(r_1,r_2)$ for integers $a_1,a_2,r_1,r_2\geq 0$ with $r_1\leq a_1$ and $r_2\leq a_2$, and let $x\in P_{a,b}$. With our convention, we get $1+a = 1+(a_1, a_2) = (1+a_1, 1+a_2)$, $(-1)^{a} = (-1)^{(a_1, a_2)} = \bigl((-1)^{a_1}, (-1)^{a_2}\bigr)$, and $\m\bigl((-1)^{1+a}\bigr) =\m\bigl((-1)^{a}\bigr) $.

Let now $N=(N_1, N_2)$. The formula
\[ \w N_1^{r_1}N_2^{r_2}x=(-1)^{r_1+r_2}\frac{r_1!}{(a_1-r_1)!}\frac{r_2!}{(a_2-r_2)!}N_1^{a_1-r_1}N_2^{a_2-r_2}x, \]
can be rewritten with our notation as follows
\[ \w \m(N^\r)x=\m\Bigl((-1)^r\frac{r!}{(a-r)!}N^{a-r}\Bigr)x. \]
In the same way, we recall that
\[ \d N_1^{r_1} N_2^{r_2} x = \sum_{i=(i_1,i_2)\in\{0,1\}^2}N_1^{r_1+i_1}N_2^{r_2+i_2}x_{i_1,i_2}, \]
if $r_1$ and $r_2$ are \emph{non-negative}, where $x_{i_1,i_2}\in P_{a_1-1+2i_1,a_2-1+2i_2}$. With our notations, this formula reads
\[ \d \m(N^r) x = \sum_{i\in\{0,1\}^2}\m(N^{r+i})x_i. \]
\end{example}

We further use the following conventions. By $(a_1,a_2)\geq(b_1,b_2)$ we mean $a_1\geq b_1$ and $a_2\geq b_2$. Moreover, $0^0=1$ and $N_1^r=N_2^r=0$ if $r$ is negative. Notice that $N_1N_1^{-1}\neq N_1^0$, thus we have to be cautious with negative exponents. This is the main difficulty in the proof, besides the exception $k=1$ and $a=0$ in Proposition \ref{prop:dd=0}.

\medskip

Let $a=(a_1,a_2)\geq 0$. Let $r=(r_1,r_2)\in\Z^2$ with $0\leq r\leq a$. Let $x\in P_a=P_{a_1,a_2}$ and $y=N_1^{r_1}N_2^{r_2}x$. Set $N=(N_1,N_2)$. We recall that
\[ \Lap=-(\d\w^{-1}\d\w+\w^{-1}\d\w\d) \]

We now explicit $\d \w^{-1} \d y$. We omit ranges of indices in the sums which follow if they are clear from the context.
\begin{align*}
\d y
  &= \sum_{i\in\{0,1\}^2}\m\bigl(N^{r+i}\bigr) x_i. \\
\w^{-1} \d y
  &= \sum_{i \\ r+i\leq a-1+2i}\m\Bigl((-1)^{a-1+2i}(-1)^{r+i}\frac{(r+i)!}{((a-1+2i)-(r+i))!}N^{(a-1+2i)-(r+i)}\Bigr) x_i, \\
  &= \sum_{i \\ a-r-1+i\geq0}\m\Bigl((-1)^{a+r+i+1}\frac{(r+i)!}{(a-r-1+i)!}N^{a-r-1+i}\Bigr) x_i. \\
\d \w^{-1} \d y
  &= \sum_{i,i'\in\{0,1\}^2 \\ a-r-1+i\geq0} \m\Bigl((-1)^{a+r+i}\frac{(r+i)!}{(a-r-1+i)!}N^{a-r-1+i+i'}\Bigr) (x_i)_{i'}.
\end{align*}
where we used the fact that $\m\bigl((-1)^{a+r+i+1}\bigr)=\m\bigl((-1)^{a+r+i}\bigr)$ (since $\m\bigl((-1)^{(1,1)}\bigr) =1$).

\medskip

We can compute in the same way $\d \w \d y$, in which case, using
\[ \m\bigl((-1)^{a-1+2i}\bigr)=\m\bigl((-1)^a\bigr), \]
we get
\[ \d \w \d y=\m\bigl((-1)^a\bigr)\d \w^{-1} \d y. \]

\medskip

To get $\d \w^{-1} \d \w y$, it suffices to multiply by the good factor and to replace $r$ by $a-r$:
\begin{align*}
\d \w^{-1} \d \w y
  &= \m\Bigl((-1)^{r}\frac{r!}{(a-r)!}\Bigr) \sum_{i,i'\in\{0,1\}^2 \\ r+i-1\geq0}\m\bigl((-1)^{r+i}\frac{(a-r+i)!}{(r-1+i)!}N^{r-1+i+i'}\bigr) (x_i)_{i'} \\
  &= \sum_{i,i' \\ r+i-1\geq0}\m\bigl((-1)^{i}\frac{r!}{(r-1+i)!}\frac{(a-r+i)!}{(a-r)!}N^{r-1+i+i'}\bigr) (x_i)_{i'} \\
  &= \sum_{i,i' \\ r+i-1\geq0}\m\Bigl((-1)^{i}r^{1-i}(a-r+1)^i N^{r-1+i+i'}\Bigr) (x_i)_{i'}.
\end{align*}
Note that in the last sum we can remove the condition $r+i-1\geq0$. Indeed, if, for instance, $r_1+i_1-1<0$, then, since $r_1, i_1$ are non-negative, we should have $r_1=i_1=0$, and the factor $r_1^{1-i_1}$ is automatically zero.

\medskip

Now, since
\[ (a-2+2i+2i')-(a-r-1+i+i')=r-1+i+i', \]
we get
\begin{align*}
\w^{-1} &\d \w \d y \\
  &= \hspace{-5mm}\sum_{i,i' \\ a-r+i-1\geq0}\hspace{-5mm}\m\Bigl((-1)^{r+i}\frac{(r+i)!}{(a-r-1+i)!}(-1)^{r-1+i+i'}\frac{(a-r-1+i+i')!}{(r-1+i+i')!}N^{r-1+i+i'}\Bigr) (x_i)_{i'} \\
  &= \hspace{-5mm}\sum_{i,i' \\ a-r+i-1\geq0}\hspace{-5mm}\m\Bigl((-1)^{i'}(r+i)^{1-i'}(a-r+i)^{i'}N^{r-1+i+i'}\Bigr) (x_i)_{i'}.
\end{align*}
In the last sum, we can again remove the condition $a-r+i-1\geq0$. Indeed, if, for instance, $a_1-r_1+i_1-1<0$, then since $r_1\leq a_1$ and $i_1 \in\{0,1\}$, we get $a_1=r_1$ and $i_1=0$. Then, either $i'_1=1$ and $(a_1-r_1+i_1)^{i'_1}=0$, or $i'_1=0$ and
\[ N_1^{r_1-1+i_1+i'_1}(x_i)_{i'}=N_1^{a_1-1}(x_i)_{i'}=0, \]
because
\[ (x_i)_{i'}=(x_{0,i_2})_{0,i'_2}\in P_{a_1-2,a_2-2+i_2+i'_2}\subseteq\ker(N_1^{a_1-1}). \]

\smallskip
In any cases, we get
\begin{align*}
-\Lap y   &= \sum_{i,i'} \Bigl(\m\bigl((-1)^{i}r^{1-i}(a-r+1)^i\bigr) + \m\bigl((-1)^{i'}(r+i)^{1-i'}(a-r+i)^{i'}\bigr)\Bigr) \m(N^{r-1+i+i'})(x_i)_{i'} \\
  &= \sum_{i,i'} \Bigl(\m\bigl(r^{1-i}(r-a-1)^i\bigr) + \m\bigl((r+i)^{1-i'}(r-a-i)^{i'}\bigr)\Bigr) \m(N^{r-1+i+i'})(x_i)_{i'} \\
  &= \sum_{i,i'} \Bigl(\m\bigl(r+0^{1-i}(-a-1)^i\bigr) + \m\bigl(r+i^{1-i'}(-a-i)^{i'}\bigr)\Bigr) \m(N^{r-1+i+i'})(x_i)_{i'}.
\end{align*}

\medskip

For a fixed $a$, let $\kappa\colon \{0,1,2\}^2\to\Z^2$ be the function defined by taking for every $i,i'\in\{0,1\}^2$, the value
\[ \kappa(i+i'):=0^{1-i}(-a-1)^i+i^{1-i'}(-a-i)^{i'}. \]
One can verify that the function is well-defined, \ie, if $i+i' = j+j$ for $i,i',j,j' \in \{0,1\}^2$ then
\[0^{1-i}(-a-1)^i+i^{1-i'}(-a-i)^{i'} = 0^{1-j}(-a-1)^j+j^{1-j'}(-a-j)^{j'}.\]
The only nontrivial case is when $i_1+i'_1=1$ or when $i_2+i'_2=1$. For instance, if $i_1=0, i'_1=1$ and $j_1=1, j_1' =0$, then we get $\kappa(i+i')_1=0-a_1=-a_1$ which is equal to $\kappa(j+j')_1$.

\medskip

We also define $C\colon \{0,1\}^2\to\Z$ by
\[ C(i,i'):=\m\bigl(0^{1-i}(-a-1)^i\bigr) + \m\bigl(i^{1-i'}(-a-i)^{i'}\bigr).\]

\medskip

We will use the following specific values later:
\begin{gather*}
\kappa\bigl((1,k_2)\big)_1=-a_1, \\
C\bigl((0,i_2),(0,i'_2)\bigr)=0, \myand \\
C\bigl((1,i_2),(1,i'_2)\bigr)=(-a_1-1)\kappa\bigl((2, i_2+i'_2)\bigr)_2.
\end{gather*}

\medskip

Using these terminology, we can now write
\begin{align*}
\Lap y
  &= -\sum_{i,i'} (2r_1r_2+\kappa(i+i')_2r_1+\kappa(i+i')_1r_2+C(i,i')) \m(N^{r-1+i+i'})(x_i)_{i'}.
\end{align*}

Set
\begin{equation} \label{eqn:Lap'}
\Lap' y := -\sum_{i,i'}C(i,i')\m(N^{r-1+i+i'})(x_i)_{i'}.
\end{equation}

We claim that $\Lap y=\Lap' y$. Indeed, we can write
\begin{align*}
(\Lap'-\Lap) y
  &= \sum_{k\in\{0,1,2\}^2} \Bigl(2r_1r_2+\kappa(k)_2r_1+\kappa(k)_1r_2\Bigr)\m(N^{r-1+k}) \sum_{i,i'\in\{0,1\}^2 \\ i+i'=k} (x_i)_{i'}.
\end{align*}
Since $\d^2=0$, using Proposition \ref{prop:dd=0}, we get
\[\m(N^{r-1+k}) \sum_{i,i'\in\{0,1\}^2 \\ i+i'=k} (x_i)_{i'} = 0\]
in the above sum unless either, for $k_1=1$ and $a_1=0$, or for $k_2=1$ and $a_2=0$.

In the first case $k_1=1$ and $a_1=0$, we should have $r_1=0$ and the remaining factor is $\kappa((1,k_2))_1r_2$. We have seen that $\kappa((1,k_2))_1=-a_1$. Thus, here, the factor is zero and there is no contribution in the sum. The same reasoning applies by symmetry to the case $k_2=1$ and $a_2=0$. Therefore, in any case, we get $\Lap y=\Lap' y$ for any $y$, and the claim follows.

\medskip

Finally, we prove that $\Lap'$ commutes with $N_1$. Recall that $y=N_1^{r_1} N_2^{r_2} x$ and $x\in P_{a_1,a_2}$. We divide the proof into three cases, depending on whether $r_1=0, a_1$ or $0<r_1 <a_1$

\medskip

\begin{itemize}[leftmargin=0pt, itemindent=1.5em, labelsep=.5em]
\item If $0=r_1<a_1$, we can compute $\Lap'N_1y$ by replacing $r_1$ by $r_1+1$ in Equation \eqref{eqn:Lap'}. For $N_1\Lap'y$, this remains true except for the terms with $i_1=i'_1=0$. Indeed,
\[ N_1N_1^{r_1-1+0+0}=0\neq N_1^{r_1+1-1+0+0}=\id. \]
Luckily, for those terms, the factor is $C((0,i_2),(0,i'_2))$, which is zero. Thus, in fact, we can also do the replacement for those terms. Hence, $\Lap'N_1y=N_1\Lap'y$.

\medskip

\item If $0<r_1<a_1$, then we can compute $\Lap'N_1y$ by replacing $r_1$ by $r_1+1$, and $N_1\Lap'y$ in the same way. Thus, they are equal.

\medskip

\item Finally, if $r_1=a_1$, since we used the hypothesis $r_1\leq a_1$ in all the computations leading to the expression of $\Lap'$, we cannot apply those computations to $N_1y$. However, in this case, we have $N_1y=0$. Thus, we just need to prove that $N_1\Lap'y=0$. This can be done thanks to the defining Equation \eqref{eqn:Lap'}.

If $i_1+i'_1\leq 1$, then $1+r_1-1+i_1+i'_1\geq a_1-2+2(i_1+i'_1)+1$. Since $(x_i)_{i'}\in P_{a_1-2+2(i_1+i'_1), b'}$ for some $b'$, we get
\[ C(i,i')N_1\m(N^{r-1+i+i'})(x_i)_{i'}=0. \]
The only remaining terms are those where $i_1=i'_1=1$. For these terms, we have $C(i,i')=(-a_1-1)\kappa(i+i')_2$. Thus,
\[ N_1\Lap'y=(a_1+1)N_1\sum_{k_2\in\{0,1,2\}}\kappa((2,k_2))_2\m(N^{r-1+i+i'})\sum_{i_2,i'_2\in\{0,1\}^2 \\ i_2+i'_2=k_2} (x_{1,i_2})_{1,i'_2}. \]
As in the first case, we apply Proposition \ref{prop:dd=0} to infer that all the terms in this sum are zero unless for $k_2=1$ and $a_2=0$. For these values of $k_2$ and $a_2$, we use $\kappa((2,1))_2=-a_2=0$ to conclude.
\end{itemize}

\medskip

In all the cases, we see that $N_1$ commutes with $\Lap'$. Since we proved $\Lap=\Lap'$, we conclude with $[\Lap,N_1]=0$, and by symmetry, $[\Lap,N_2]=0$, as required.
}


\section{Projective bundle theorem}
\label{sec:projective_bundle_theorem}

Let $\X$ be a polyhedral space in $V = N_\R$. We define the \emph{asymptotic cone} of $\X$ which we denote by $\X_\infty$ to be the limit set of the rescaled spaces $\X/t$ when $t$ goes to $+\infty$. Exceptionally, the terminology cone here might refer to a non-necessarily convex cone. Note that, if $X$ is any polyhedral complex structure on $\X$, then $\X_\infty=\supp{X_\infty}$.

\medskip

Let $\Y$ be a smooth polyhedral space in $V = N_\R$. Let $\Delta$ be a unimodular fan subdividing $\Y_\infty$. Let $\X$ be the closure of $\Y$ in $\TP_\Delta$. For each cone $\eta \in \Delta$, we define $D^\eta$ as the closure of the intersection $\X \cap \TP_\Delta^\eta$.

\medskip

Using these notations, we see in particular that we have
\begin{itemize}
\item $D^{\conezero} = \X$.
\item $D^\eta$ are all smooth of dimension $d-\dims{\eta}$.
\item We have \[D^\eta = \bigcap_{\substack{\varrho \subface \eta \\ \dims{\varrho}=1}} D^\varrho.\]
\end{itemize}
In other words, the boundary $D = \X \setminus \Y$ is a \emph{simple normal crossing divisor} in $\X$.

\subsection{Cycle class associated to the cones in $\Delta$} In this section, we recall how to associate elements of the tropical cohomology groups to the strata $D^\eta$ defined in the preceding paragraph. So we fix the notations as above. Let $\Delta$ be a unimodular fan structure on the asymptotic cone $\Y_\infty$, and take $\eta \in \Delta$. The stratum $D^\eta \subset \X$ defines then an element $\class(D^\eta)$ of $H^{\dims\eta,\dims\eta}_\trop(\X)$ obtained as follows.

\medskip

Let $d$ be the dimension of $\Y$. From the inclusion $D^\eta \subseteq \X$, we get a map
\[
\begin{tikzcd}[row sep=tiny]
H_{d-\dims\eta, d-\dims\eta} (D^\eta) \rar & H_{d-\dims \eta, d-\dims\eta} (\X) \rar["\sim"] &
\Bigl(H^{d -\dims \eta, d-\dims \eta} (\X) \Bigr)^\dual \rar["\sim"]& H^{\dims\eta,\dims\eta}(\X)\\
{[D^{\eta}]}\arrow[mapsto]{rrr} &&&\class(D^\eta),
\end{tikzcd}
\]
\ie, $\class(D^\eta)$ is the image of $[D^\eta] \in H_{d-\dims\eta, d-\dims\eta} (D^\eta) $ in $H^{\dims\eta,\dims\eta}(\X)$ by the composition of above maps. Note that the last isomorphism in the above chain of maps follows from the Poincar\'e duality for $\X$.

\subsection{Explicit description of the divisor classes associated to the rays} \label{sec:cl_D'}

Consider now the situation where $\X$ comes with a unimodular polyhedral complex structure $X$ with $X_\infty =\Delta$. Moreover, we consider a ray $\rho$ of $\Delta$  and its associated stratum $D^\rho$ and its cycle class $\class(D^\rho) \in H^{1,1}_\trop(\X)$. In this case, there is a nice representative $c = (c_\zeta)_{\substack{\zeta\in X \\ \dims \zeta=1}}$ in $C^{1,1}_\trop(X)$ of $\class(D^\rho)$, with $c_\zeta \in \SF^1(\zeta)$, as we explain now. For details of this construction, we refer to~\cite{AP}.

\medskip

First we define the elements $\zeta \in X$ of dimension one which are \emph{parallel to $\rho$}.

So let $\zeta$ be a one dimensional face of $X$ and denote by $\tau \in \Delta$ the sedentarity of $\zeta$.
We say \emph{$\zeta$ is parallel to $\rho$} if the following two properties hold:
\begin{enumerate}
\item $\sigma:=\rho+\tau$ is a cone  of $\Delta$ which properly contains $\tau$; and
\item $\zeta^\tau=\sigma^\tau_\infty$.
\end{enumerate}
Here $\sigma^\tau_\infty$ is $\sigma/\tau$ seen in the stratum $N^\tau_{\infty, \R}$ of $\TP_\Delta$, and $\zeta^\tau$ is the intersection $\zeta \cap N^\tau_{\infty, \R}$.

\medskip

With this preliminary discussion, the element $c=(c_\zeta)_{\substack{\zeta\in X \\ \dims \zeta=1}}$ in $C^{1,1}_\trop(X)$ is defined as follows.

\medskip

First, consider a face $\zeta\in X$ which is not parallel to $\rho$ and which has sedentarity $\tau$ such that $\tau$ does not contain $\rho$ as a ray. For such $\zeta$,  we set $c_\zeta=0$.

\medskip

Second, consider a one dimensional cone $\zeta$ of $X$ which is parallel to $\rho$. Let $\tau$ be the sedentarity of $\zeta$ and choose a linear form $\phi_\zeta$ on $N^\tau_{\infty,\R}$ which takes value one on the primitive vector $\e_{\zeta_\infty}$ of $\zeta_\infty$. The elements $\phi_\zeta$ can be naturally seen as an element of $\SF^1(\zeta)$. We denote this element by $c_\zeta$.

\medskip

All the remaining $\zeta$ in $X$ are those whose sedentarity contains $\rho$ as a ray. In particular, they are not parallel to $\rho$. Any such $\zeta$ form a side of a unique quadrilateral in $X$ with two sides $\zeta_1, \zeta_2$ parallel to $\rho$, and the last side $\zeta_3$ which is neither parallel to $\rho$ nor with a sedentarity $\sed(\zeta_3)$ containing $\rho$. By the construction of the previous two cases, we have $c_{\zeta_3} =0$. Moreover,  the difference $c_{\zeta_1} - c_{\zeta_2}$ naturally defines a linear form on $N^\tau_{\infty, \R}$. We define $c_{\zeta} = \pm(c_{\zeta_1} - c_{\zeta_2})$, for a convenient choice of sign that we omit to precise here.

\medskip

The element  $c=(c_\zeta)_{\substack{\zeta\in X \\ \dims \zeta=1}}$ in $C^{1,1}_\trop(X)$ thereby defined is closed and its class in the cohomology $H^{1,1}_\trop(X)$ coincides with $\class(D^\rho)$.

\medskip

Moreover, to each face $\zeta$ with sedentarity $\conezero$ which is parallel to $\rho$, we can associate the corresponding element $x_\zeta$ of $A^1(\zeta_\f)=H^2(\zeta_\f)$. Note that $\zeta_\f$ is a vertex of $X_\f$. Since $\zeta_\f$ is a vertex in $X_\f$, we get an element of $\ST_1^{0,2,0}\subseteq\ST_1^{0,2}$. By summing all the $x_\zeta$, we get an element of $\ST_1^{0,2}$ which is closed in the Steenbrink sequence, and which corresponds to $c$ via the isomorphism of Theorem \ref{thm:steenbrink}.

\begin{remark} \label{rem:explicit_function}
We have seen in Section \ref{sec:projective_Kahler} how to associate an element of $\ST_1^{0,2,0}$ to a piecewise linear function $f$ on $Y:=X^\conezero$. Let us describe such a function giving the element $c$ considered in the preceding paragraph.

\medskip

Assume $X$ is simplicial. Denote by $\chi_\rho$ the characteristic function of $\rho$ on $\Delta$ which takes value one on $\e_\rho$, value zero on each $\e_\varrho$ for all ray $\varrho \neq \rho$ of $\Delta$, and which is linear on each cone of $\Delta$. Pick any linear form $\e^\dual_\rho$ on $N_\R$ which takes value $1$ on $\e_\rho$. Then $\chi_\rho$ can be naturally extended to $\comp\Delta$ as follows. For any cone $\sigma\in\Delta$ not containing $\rho$, we extend $\chi_\rho$ by continuity on $\Delta^\sigma$. For any cone $\sigma\in\Delta$ containing $\rho$, we define the extension $\chi_\rho$ on $\Delta^\sigma$ to be the continuous extension of $\chi_\rho-\e^\dual_\rho$ to this space.

Then $\chi_\rho$ naturally induces a piecewise linear function $\~\chi_\rho$ on $X$ defined by $x\mapsto\chi_\rho(x_\infty)$. Moreover, the element of $\ST_1^{0,2,0}$ associated to this function does not depend on the choice of $\e^\dual_\rho$ and is equal to the element $c$ constructed above.
\end{remark}

\subsection{Statement of the theorem}

Let now $\sigma$ be a cone in $\Delta$, and denote by $\Delta'$ the fan obtained by the barycentric star subdivision of $\sigma$ in $\Delta$. Denote by $\X'$ the closure of $\Y$ in $\TP_{\Delta'}$. By Theorem~\ref{thm:smoothness-compact} we have the following.
\begin{prop}
The variety $\X'$ is smooth.
\end{prop}

For each cone $\zeta$ in $\Delta'$, we define similarly as in the preceding section, the corresponding closed stratum $D'^\zeta$ of $\X'$. Let $\rho$ be the new ray in $\Sigma'$ obtained after the star subdivision of $\sigma \in \Delta$.

Note that $\rho$ is the ray generated by the sum of the primitive vectors of the rays of $\sigma$, \ie,
\[\rho = \R_+\bigl(\sum_{\varrho \subface \sigma \\ \dims{\varrho}=1} \e_\varrho\bigr).\]
Consider the divisor $D'^\rho \subseteq \X'$ and its corresponding cycle class $\class(D'^\rho)$ of $H^{1,1}_\trop(\X')$.

\medskip

Recall that for any smooth compact tropical variety $\W$, and for each integer $k$, we define the $k$-th cohomology of $\W$ by
\[H^k(\W) = \sum_{p+q=k} H_\trop^{p,q}(\W).\]

\medskip

For the inclusion $\i\colon\vZ\hookrightarrow \W$ of smooth compact tropical varieties, as in the local case, we get the corresponding restriction maps
\[\i^*\colon H^{k}(\W) \to H^{k}(\vZ),  \]
and the Gysin map
\[\gys\colon H^{k}(\vZ) \to H^{k+2\dim(\W) - 2\dim(\vZ)}(\W),\]
which is the Poincar\'e dual to the restriction map $H^{2\dim(\vZ) - k}(\W) \to H^{2\dim(\vZ)-k}(\vZ)$. The restriction maps respect the bi-degree $(p,q)$, while the Gysin map sends elements of bidegree $(p,q)$ to elements of bidegree $(p+\dim(\W)-\dim(\vZ), q+\dim(\W)-\dim(\vZ))$.

\medskip

\begin{prop}\label{prop:adjunction} Let $\vZ$ be a divisor in $\W$, and denote by $\class(\vZ)$ the class of $\vZ$ in $H^{1,1}(\W)$. For any $\beta \in H^k(\W)$, we have
\[\gys (\i^*(\beta)) = \class(\vZ) \cup \beta.\]
\end{prop}
\begin{proof} This follows for example by duality from the projection formula for tropical homology groups established in \cite{GS-sheaf}*{Proposition 4.18}.
\end{proof}

From the map of fans $\pi\colon  \Delta' \to \Delta $, we get a projection map $\pi \colon\TP_{\Delta'} \to \TP_{\Delta}$ which by restriction leads to a projection map
$\pi \colon \X' \to \X$. This map is \emph{birational} in the sense that it induces an isomorphism on the open set $\Y$ of $\X$ and $\X'$: in fact, the map $\pi$ restricts to an isomorphism from $\X' \setminus D'^\rho $ to $\X \setminus D^\sigma$. Moreover, it sends $D'^\rho$ to $D^\sigma$.

The compositions of $\gys$ and $\pi^*$ gives now a map
\[\gys \circ\,\pi^*\colon H^s(D^\sigma) \to H^{s+2}(\X').\]

With these preliminaries, we can state the main theorem of this section.
\begin{thm}[Projective bundle formula] \label{thm:keelglobal}
We have an isomorphism
\[H^k(\X' ) \simeq H^k(\X) \oplus T H^{k-2}(D^\sigma) \oplus T^2 H^{k-4}(D^\sigma) \oplus \dots \oplus T^{\dims{\sigma}-1}H^{k-2\dims{\sigma}+2}(D^\sigma).\]
\end{thm}
Here the map from the right hand side to the left hand side restricts to $\pi^*$ on $H^k(\X)$ and sends $T$ to $-\class(D'^\rho)$, so in view of Proposition~\ref{prop:adjunction}, it is given on each factor $T^s H^{k-2s}(D^\sigma)$ by
\[ \begin{array}{rcl}
T^s H^{k-2s}(D^\sigma) & \longrightarrow & H^k(\X') \\[1em]
T^s \alpha & \longmapsto & (-1)^s \class(D'^\rho)^{s-1} \cup \gys\circ\pi^*(\alpha).
\end{array} \]
The decomposition stated in the theorem gives a decomposition for each pair of non-negative integers $(p,q)$
\[H^{p,q}(\X') \simeq H^{p,q}(\X) \oplus T H^{p-1, q-1}(D^\sigma) \oplus T^2 H^{p-2, q-2}(D^\sigma) \oplus \dots \oplus T^{\dims{\sigma}-1}H^{p-\dims\sigma+1,q - \dims\sigma+1}(D^\sigma).\]

\begin{remark} In this case, we have
\[T^{\dims{\sigma}} + c_1(N_{D^\sigma})T^{\dims\sigma-1} +c_2(N_{D^\sigma})T^{\dims\sigma-2} +\dots + c_{\dims \sigma}(N_{D^\sigma}) =0 \]
where $c_i(N_{D^\sigma})$ are the Chern classes of the normal bundle $N_{D^\sigma}$ in $\X$.
\end{remark}
The proof of this theorem is given in the next section. We actually give two proofs, one based on the use of the materials developed in the previous sections, using triangulations and the Steenbrink spectral sequence, and the other one based on Mayer-Vietoris exact sequence for tropical cohomology groups.

It is worth noting  that the second proof is more classical in spirit, while the first one being more interesting, requires more refined results about the geometry of tropical varieties which are not developed here and will only appear in our future work.

\subsection{Proof of Theorem~\ref{thm:keelglobal} using triangulations}

In this section, we give a proof of Theorem~\ref{thm:keelglobal} using the tropical Steenbrink spectral sequence and the isomorphism theorem (Theorem~\ref{thm:steenbrink}). We need the following theorem.

\begin{thm} Suppose $\Y$ is a polyhedral space in $\R^n$, and let $\Delta$ be a unimodular fan structure on $\mathscr Y_\infty$. There is a unimodular triangulation $Y$ of $\Y$ with $Y_\infty=\Delta$.
\end{thm}
\begin{proof} In the case when $\Y$ comes with a polyhedral structure $Y_0$ with recession fan $\Delta$, this is proved in Section~\ref{sec:triangulation}, \cf. Theorem~\ref{thm:unimodular_preserving_recession}. The point is to show that any polyhedral space $\Y$ admits a polyhedral structure compatible with the choice of such a unimodular fan structure $\Delta$, \ie, with recession fan $\Delta$. This is more subtle and will be treated in our forthcoming work~\cite{AP-geom}.
\end{proof}

Let $Y$ be such a triangulation of $\Y$. Let $X$ be its extensions to $\X$.

\begin{prop} The star subdivision $\Delta'$ of $\Delta$ induces a polyhedral structure $Y'$ on $\Y$.
Moreover, this polyhedral structure is unimodular and induces a unimodular polyhedral structure $X'$ on $\X'$ with the same finite part as $X$, \ie, $X_\f = X'_\f$.
\end{prop}
\begin{proof} The polyhedral structure $Y'$ is the one obtained as follows. Consider a face $\eta$ in $Y$, and decompose $\eta =\eta_\f+ \eta_\infty$, with $\eta_\f$ in $X_\f$ and $\eta_\infty$ in the recession fan $Y_\infty$. We define the subdivision $\eta'$ of $\eta$ obtained as the result of the barycentric star subdivision of $\sigma$ in $\eta_\infty$: this means if $\sigma$ is not a face of $\eta_\infty$, then $\eta'=\eta$ remains unchanged. On the other hand, if $\sigma \subface \eta_\infty$, then we subdivide $\eta_\f+ \eta_\infty$ by first taking the star subdivision $\eta_\infty'$ of $\sigma$ in $\eta_\infty$, and define the subdivision $\eta'$ as the one consisting of the faces of the form $\gamma + \tau$ with $\gamma \subface \eta_\f$ and $\tau \subface \eta_\infty'$.
\end{proof}

The following proposition describes the star fans of $Y'$ in terms of those of $Y$. By an abuse  of the notation, we also denote by $D^\sigma$ the subdivision induced by $X$ on $D^\sigma$.

\begin{prop} Let $\eta \in X_\f =Y_\f$ be a compact polyhedron. The star fan $\Sigma'^{\eta}$ of $\eta$ in $X'$ is obtained from the star fan $\Sigma^\eta$ of $\eta$ as the result of star subdividing $\sigma$. In particular, if the star fan of $\eta$ does not contain $\sigma$, the two fans $\Sigma'^{\eta}$ and $\Sigma^\eta$ are isomorphic.
\end{prop}
In the above proposition, by star subdividing $\sigma$ in $\Sigma^\eta$ we mean the fan obtained by the star subdivision of the image of $\sigma$ in the space $N^\eta_\R = \rquot{N_\R} {N_{\eta, \R}}$, where we recall $N_{\eta, \R} =\TT \eta$.
\begin{proof} This follows directly from the description of the polyhedra in $Y'$.
\end{proof}

In the following, in order to distinguish the star fans and their cohomology in an ambient polyhedral complex, we use the following terminology. For a polyhedral complex $W$, and a face $\eta$ of $W$, we write $\Sigma^\eta(W)$ in order to emphasize that the star fan $\Sigma^\eta$ is taken in $W$. We furthermore use $H^\bul_W(\eta)$ for the cohomology groups $H^\bul\bigl(\,\comp{\Sigma^\eta(W)}\,\bigr)$. In particular, $H^\bul_{X}(\eta)$ and $H^{\bul}_{X'}(\eta)$ denote the cohomology groups $H^\bul\bigl(\,\comp{\Sigma^\eta}\,\bigr)$ and $H^\bul\bigl(\,\comp{\Sigma^{\prime\,\eta}}\,\bigr)$, respectively.

\medskip

By Keel's lemma (Theorem~\ref{thm:keel}) and the above proposition, for each face $\eta \in X_\f$ whose star fan $\Sigma^\eta$ contains the cone $\sigma$, we obtain a decomposition of the cohomology groups $H^k_{X'}(\eta)$.

\begin{prop} For each cone $\eta \in X_\f$ with $\sigma \in \Sigma^\eta$, we have
\[H^k_{X'}(\eta) = H^k_X(\eta) \oplus T \cdot H^{k-2}_{\Sigma^\eta}(\sigma) \oplus T^2 \cdot H^{k-4}_{\Sigma^\eta}(\sigma) \oplus \dots \oplus T^{\dims \sigma-1} \cdot H^{k- 2\dims \sigma+2}_{\Sigma^\eta}(\sigma),\]
where $T$ is the class of $-x_\rho \in A^1(\Sigma'^\eta) = H^2_{X'}(\eta)$.
\end{prop}

\medskip

For $W = X', X,$ or $D^\sigma$, and for any integer $l$, we denote by $\STp{l}_{W}^{\bullet, \bullet}$ the double complex we constructed in Section~\ref{sec:steenbrinkdoublecomplex} whose total complex is the shifted cochain complex $\ST_1^{\bul, 2l}[l]$, the $2l$-th line of the Steenbrink spectral sequence calculated in $W$.

Let $p$ be a non-negative integer. The previous proposition allows to get a decomposition of the double complex $\STp{p}_{X'}^{\bullet, \bullet}$ in terms of the Steenbrink double complexes associated to $X$ and $D^\sigma$.

\begin{prop} Notations as above, we have
{ \let\oldoplus\oplus \renewcommand{\oplus}{\,\,\oldoplus\,\,}
  \let\oldcdot\cdot \renewcommand{\cdot}{\!\oldcdot\!}
\[ \STp{p}_{X'}^{\bullet, \bullet} = \STp{p}_{X}^{\bullet, \bullet}  \oplus T \cdot \STp{p-1}_{D^\sigma}^{\bullet, \bullet-1} \oplus T^2\cdot\STp{p-2}_{D^\sigma}^{\bullet, \bullet-2}\oplus{} \cdots  \oplus T^{\dims\sigma-1}\cdot\STp{p-\dims{\sigma}+1}_{D^\sigma}^{\bullet, \bullet -\dims{\sigma}+1}. \]
}
\end{prop}

The map from the decomposition to $\STp{p}_{X'}^{\bullet, \bullet}$ is made clear in the proof.

\begin{proof}
By definition of the Steenbrink double complexes, we have
\[ \STp{p}^{a,b}_{X'}:=\begin{cases}
  \displaystyle\bigoplus_{\eta\in X'_\f \\ \dims\eta=p+a-b} H^{2b}_{X'}(\eta) & \text{if $a\geq0$ and $b\leq p$}, \\[2em]
  0 & \text{otherwise,}
\end{cases} \]
whose differential of bidegree $(1,0)$ is $\i^*$ and whose differential of bidegree $(0,1)$ is $\gys$. Note that we have $X_\f = X'_\f$. From the previous proposition, for any $\eta \in X_\f$ whose star fan $\Sigma^\eta$ in $X$ contains $\sigma$, we can decompose
\[H^{2b}_{X'}(\eta) = H^{2b}_X(\eta) \oplus \, T\cdot \,H^{2b-2}_{\Sigma^\eta}(\sigma) \oplus\, T^2\cdot H^{2b-4}_{\Sigma^\eta}(\sigma) \oplus \dots \oplus\, T^{\dims \sigma-1}\cdot H^{2b- 2\dims \sigma+2}_{\Sigma^\eta}(\sigma). \]
Here $T$ is the class of $-x_\rho$ in $H^2_{\Sigma^\eta}(\sigma)$.

Moreover, for any other polytope $\eta\in X_\f$, we have $H^{2b}_{X'}(\eta) = H^{2b}_{X}(\eta)$.

\medskip

Let $\eta \in X_\f$ whose star fan $\Sigma^\eta$ contains $\sigma$. Then the cone $\sigma\in\Sigma^\eta$ corresponds to the face $\~\eta:=\comp{\eta+\sigma}\in X$. We obtain a polytope $\~\eta\cap D^\sigma\in X$ in $(D^\sigma)_\f$ for the polyhedral structure induced by $X$ on $D^\sigma$. We denote this polytope by $\eta_\sigma$. The projection $\pi_\sigma\colon \Y\to D^\sigma$ along $\TT\sigma$ induces an isomorphism $\eta\simto\eta_\sigma$. Moreover, the star fan of $\sigma$ in $\Sigma^\eta$ is naturally isomorphic to the star fan of $\~\eta$ in $X$ and of $\eta_\sigma$ in $D^\sigma$.

Using this observation, we rewrite the above decomposition as
\[H^{2b}_{X'}(\eta) = H^{2b}_X(\eta) \oplus \, T\cdot \,H^{2b-2}_{D^\sigma}(\eta_\sigma) \oplus\, T^2\cdot H^{2b-4}_{D^\sigma}(\eta_\sigma) \oplus \dots \oplus\, T^{\dims \sigma-1}\cdot H^{2b- 2\dims \sigma+2}_{D^\sigma}(\eta_\sigma). \]

Moreover, the morphisms $\i^*$ and $\gys$ are compatible with this decomposition.
It follows that
{ \let\oldoplus\oplus \renewcommand{\oplus}{\,\,\oldoplus\,\,}
  \let\oldcdot\cdot \renewcommand{\cdot}{\!\oldcdot\!}
\[ \STp{p}_{X'}^{\bullet, \bullet} = \STp{p}_{X}^{\bullet, \bullet}  \oplus T \cdot \STp{p-1}_{D^\sigma}^{\bullet, \bullet-1} \oplus T^2\cdot\STp{p-2}_{D^\sigma}^{\bullet, \bullet-2}\oplus{}
\cdots  \oplus T^{\dims\sigma-1}\cdot\STp{p-\dims{\sigma}+1}_{D^\sigma}^{\bullet, \bullet -\dims{\sigma}+1}. \]
}
Here $T$ can be seen as the element of $\STp{1}^{0,1}_{X'}\subseteq\ST_1^{0,2}(X')$ defined by
\[ T=\sum_{v\in X_\f \\ \dims{v}=1 \\ \Sigma^v\ni\sigma}\underbrace{-x_\rho}_{\in H^2_{X'}(v)}\in\bigoplus_{v\in X_\f \\ \dims v=1}H^2_{X'}(v)\subseteq\STp{p}^{0,1}_{X'}. \]
The multiplication by $T$ is similar to the multiplication by a Kähler class. It multiplies an element of $H^\bul(\eta)$ by $-x_\rho\in H^2(\eta)$ if $\sigma\in\Sigma^\eta$, and by zero otherwise.
\end{proof}

\medskip

We are now ready to prove the projective bundle theorem.

\begin{proof}[Proof of Theorem~\ref{thm:keelglobal}]
Applying the previous proposition, and taking the total complexes on both sides, after shifting, gives a decomposition of cochain complexes
\[\ST_1^{\bul, 2p}(X') = \ST_1^{\bul, 2p}(X) \oplus T\cdot \ST_1^{\bul, 2p-2}(D^\sigma)  \oplus \dots \oplus T^{\dims\sigma-1}\cdot \ST_1^{\bul, 2p-2\dims\sigma+2}(D^\sigma).\]
Taking the cohomology of degree $q-p$ on both sides, and using our Isomorphism Theorem~\ref{thm:steenbrink}, we get the decomposition
\[H^{p,q}(\X') \simeq H^{p,q}(\X) \oplus T H^{p-1, q-1}(D^\sigma) \oplus T^2 H^{p-2, q-2}(D^\sigma) \oplus \dots \oplus T^{\dims{\sigma}-1}H^{p-\dims\sigma+1,q - \dims\sigma+1}(D^\sigma),\]
as stated in the theorem. Taking the sum over all $p,q$ with $p+q =k$, gives the decomposition for the $k$-th cohomology group of $X'$.

It remains to prove that $T$ corresponds to minus the cycle class associated to $D'^\rho$ in $H^{1,1}(X')$. To see that, it suffices to compare the definition of $T$ in the proof of the previous proposition with the description of $\class(D'^\rho)$ in Section \ref{sec:cl_D'}.
\end{proof}

\subsection{Proof of Theorem~\ref{thm:keelglobal} using Mayer-Vietoris}

Notations as in the beginning of the section, let $\X$ and $\X'$ be the closure of $\Y$ in $\TP_\Delta$ and $\TP_{\Delta'}$. We consider the open covering of $\X$ obtained by $\V:=\X\setminus D^\sigma=\X'\setminus D'^\rho$ and a tubular neighborhood $\U$ of $D^\sigma$. Let $\W$ be the intersection of $\V$ and $\U$. Let $\pi\colon \X'\to \X$ be the projection map, and let $\U' \colon \pi^{-1}(\U)$, which is a tubular neighborhood of $D'^\rho$. Note that $\W = \V \cap \U'$.

The Mayer-Vietoris exact sequences associated to these coverings of $\X'$ and $\X$ for tropical cohomology groups with coefficients in $\SF^p$ lead to the following pair of long exact sequences where the vertical arrows are induced by the pullbacks.

\[
\begin{tikzcd}
\cdots \rar& H^{p,q}(\X) \dar \rar & H^{p,q}(\V)\oplus H^{p,q}(\U) \dar \rar & H^{p,q}(\W) \dar\rar & H^{p,q+1}(\X) \dar\rar &  \cdots \\
\cdots \rar&  H^{p,q}(\X') \rar & H^{p,q}(\V)\oplus H^{p,q}(\U') \rar & H^{p,q}(\W) \rar & H^{p,q+1}(\X') \rar&  \cdots
\end{tikzcd} \]

The following is straightforward.
\begin{prop} The projection map from the tubular neighborhood to $D^\sigma$ induces an isomorphism
\[H^{p,q} (\U) \simeq H^{p,q}(D^\sigma)\]
for any $q$. In a similar way, we get
\[H^{p,q} (\U') \simeq H^{p,q}(D'^\rho).\]
\end{prop}

\begin{prop}\label{prop:proj_bundle_blowup} The projection map $D'^\rho \to D^\sigma$ is a fibration with tropical projective spaces of dimension $\dims\sigma -1$.
\end{prop}
\begin{proof} Define the normal vector bundle $N_{D^\sigma}$ associated to the embedding $D^\sigma \hookrightarrow \X$.
\[D^\rho \simeq \TP(N_{D^\sigma}), \]
the projectivization of $N_{D^\sigma}$, from which the result follows.
\end{proof}

\begin{thm}\label{thm:proj_bundle_thm} Let $E$ be a vector bundle of rank $r$ over a compact tropical variety $\vZ$, and let $\vZ' = \TP(E)$ be the projectivization of $E$. Then we get an isomorphism
\[H^{p,q}(\vZ') \simeq H^{p,q}(\vZ) \oplus T H^{p-1,q-1}(\vZ) \oplus T^2 H^{p-2,q-2}(\vZ) \oplus \dots \oplus T^{r-1}H^{p-r+1,q-r+1}(\vZ).\]
\end{thm}

\begin{proof} Applying the Mayer-Vietoris exact sequence, we reduce to the case where the vector bundle is trivial, \ie, $E = \R^r\times \vZ$. In this case, we have $\TP(E) = \TP^{r-1} \times \vZ$. The theorem now follows by K\"unneth decomposition theorem for tropical cohomology groups~\cite{GS-sheaf}.
\end{proof}

\begin{proof}[Proof of Theorem~\ref{thm:keelglobal}] The theorem follows by combining the two Mayer-Vietoris long exact sequences associated to covers $\{\V, \U\}$ and $\{\V, \U'\}$ of $\X$ and $\X'$, and Theorem~\ref{thm:proj_bundle_thm} applied to the fibration $D'^\rho \to D^\sigma$, by using Proposition~\ref{prop:proj_bundle_blowup}, which allows to decompose the cohomology of $D'^\rho$ in terms of the cohomology of $D^\sigma$.
\end{proof}


\section{Proof of Theorem~\ref{thm:main2}} \label{sec:proofmaintheorem}
In this section, we prove Theorem~\ref{thm:main2}. We start by recalling the statement of the theorem. Let $\Y$ be a smooth rational polyhedral space in $V = N_\R  \simeq \R^n$. Let $\Delta$ be a unimodular recession fan on $\Y_\infty$. Let $\X$ be the closure of $\Y$ in $\TP_\Delta$. Let $\ell$ be a strictly convex function on $\Delta$. Then, $\ell$ defines an element of $H^{1,1}(\X)$ given by
\[ \sum_{\varrho\in\Delta \\ \dims\varrho=1}\ell(\varrho)\class(D^\varrho). \]
By an abuse of the notation, we denote this element and the corresponding Lefschetz operator on $H^{\bul, \bul}(\X)$ by $\ell$ as well.

Working with $\Q$ coefficients, we have to prove the following statements.
\begin{itemize}
\item (Weight-monodromy conjecture) For $q>p$ two non-negative integers, we get an isomorphism
\[N^{q-p} \colon H_{\trop}^{q,p}(\X) \to H_{\trop}^{p,q}(\X).\]

\item (Hard Lefschetz) For $p+q \leq d/2$, the Lefschetz operator $\ell$ induces an isomorphism
\[\ell^{d- p-q}\colon H_{\trop}^{p,q}(\X) \to H^{d-q, d-p}_{\trop}(\X).\]

\item (Hodge-Riemann) The pairing $(-1)^p \bigl< \ccdot,\, \ell^{d-p-q} N^{q-p}\ccdot \bigr>$ induces a positive-definite pairing on the primitive part $P^{p,q}$ of $H_{\trop}^{p,q}$, where $\bigl< \ccdot,\rdot \bigr>$ is the natural pairing
\[\bigl< \ccdot, \rdot\bigr> \colon H^{q,p}_\trop(\X) \otimes H_\trop^{d-q,d-p}(\X) \to H_\trop^{d,d}(\X) \simeq \Q.\]
\end{itemize}

\subsection{Global ascent-descent} Situation as in the previous section, let $\Delta = \Y_\infty$ be unimodular and let $\sigma$ be a cone in $\Delta$. Let $\Delta' $ be the star subdivision of $\Delta$ obtained by star subdividing $\sigma$. Denote by $\rho$ the new ray in $\Delta'$. For an element $\omega$ in $H^{1,1}(\X)$, we denote by $\omega^\sigma$ the restriction of $\omega$ to $D^\sigma$, \ie, the image of $\omega$ under the restriction map $H^{1,1}(\X) \to H^{1,1}(D^\sigma)$, for the inclusion of tropical varieties $D^\sigma \hookrightarrow \X$.

Let $\ell$ be a strictly convex cone-wise linear function on $\Delta$. Denote by $\ell'$ the convex function induced by $\ell$ on $\Delta'$. As in the previous section, we denote by $\ell$ and $\ell'$ the corresponding elements of $H^{1,1}(\X)$ and $H^{1,1}(\X')$, respectively.

The strictly convex function $\ell$ induces a strictly convex function $\ell^\sigma$ on the star fan $\Sigma^\sigma$ of $\sigma$ in $\Delta$. By identification of $\Sigma^\sigma$ with the recession fan of $D^\sigma$, we get a strictly convex function $\ell^\sigma$ on the recession fan $D^\sigma_\infty$. This defines an element in $H^{1,1}(D^\sigma)$ that we denote with our convention above by $\ell^\sigma$.

\begin{prop} \label{prop:class_commutes}
The element $\ell^\sigma \in H^{1,1}(D^\sigma)$ coincides with the restriction of $\ell \in H^{1,1}(X)$ to $D^\sigma \hookrightarrow \X$.
\end{prop}
\begin{proof}
We can assume without loss of generality that $\ell$ is zero on $\sigma$ and that $\ell^\sigma$ is the pushforward of $\ell$ by the projection $\Sigma \to \Sigma^\sigma$. Denote by $D_\ell$ the Poincaré dual of $\ell$ in $H^\trop_{1,1}(\X)$, and by $D_{\ell^\sigma}$ the Poincaré dual of $\ell^\sigma$ in $H^\trop_{1,1}(D_\sigma)$. By Poincaré duality, it suffices to prove that the cap product $D_\ell \cap D^\sigma\in H^\trop_{1,1}(D^\sigma)$ equals $D_{\ell^\sigma}$. For a ray $\varrho$ of $\Sigma$, we write $\varrho\sim\sigma$ if $\varrho$ is in the link of $\sigma$, \ie, if $\sigma+\varrho$ is a face of $\Sigma$ and $\varrho\not\in\sigma$. From our assumption that $\ell$ is zero on $\sigma$, we get
\[ D_\ell \cap D^\sigma = \sum_{\varrho\in\Sigma \\ \dims\varrho=1}\ell(\varrho)(D^\varrho\cap D^\sigma) = \sum_{\varrho\sim\sigma \\ \dims\varrho=1}\ell(\varrho)(D^\varrho\cap D^\sigma) = D_{\ell^\sigma}. \qedhere \]
\end{proof}

We now state and prove the following global version of Theorem~\ref{thm:barycentric_subdivision}.

\begin{thm} \label{thm:global_ascent_descent} Notations as above, we have the following properties.
\begin{itemize}
\item \emph{(Ascent)} Assume the property $\HR(D^\sigma,\ell^\sigma)$ holds. Then $\HR(\X, \ell)$ implies $\HR\bigl(\X',\ell'-\epsilon\class(D^\rho)\bigr)$ for any small enough $\epsilon>0$.

\item \emph{(Descent)} We have the following partial inverse: if both the properties $\HR(D^\sigma,\ell^\sigma)$ and $\HL(\X,\ell)$ hold, and if we have the property $\HR\bigl(\X',\ell'-\epsilon\class(D^\rho)\bigr)$ for any small enough $\epsilon>0$, then we have $\HR(\X, \ell)$.
\end{itemize}
\end{thm}

\begin{proof}
The projective bundle formula (Theorem \ref{thm:proj_bundle_thm}) is the global analog of local Keel's lemma (Theorem \ref{thm:keel}). Moreover Poincaré duality holds for $\X$ and $\X'$. Thus, by Remark \ref{rem:keel}, the proof is identical to that of the corresponding local theorem, Theorem~\ref{thm:barycentric_subdivision}.
\end{proof}

\subsection{Proof of Theorem~\ref{thm:main2}} We proceed as in the proof of the properties $\HR$ and $\HL$ in the local case.

So let $\Y$ be a smooth tropical variety of dimension $d$ in $V \simeq \R^n$. Any unimodular recession fan $\Delta$ on $\Y_\infty$ provides a compactification of $\Y$ by taking its closure in $\TP_\Delta$. We denote by $\X_\Delta$ this compactification in order to emphasize the dependence in $\Delta$. Since $\Y$ is smooth, the variety $\X_\Delta$ for any unimodular fan $\Delta$ is a smooth tropical variety.

\medskip

From now on we restrict to fans $\Delta$ on the support of $\Y_\infty$ which are in addition quasi-projective.
We have to show that for any such fan $\Delta$ and any strictly convex function $\ell$ on $\Delta$, the pair $(H^\bul(X), \ell)$ verifies $\HR$, $\HL$, and the weight-monodromy conjecture.

\medskip

Proceeding by induction, we suppose that the theorem holds for all smooth projective tropical varieties of dimension strictly smaller than $d$.

\medskip

As in the local case, the following proposition is the base of our proof. It allows to start with a fan structure on the asymptotic cone of $\Y$ such that the corresponding compactification of $\Y$ verifies the Hodge-Riemann relations, and then use the global ascent and descent properties to propagate the property to any other smooth compactification of $\Y$, given by any other fan structure on $\Y_\infty$.

\begin{prop}\label{prop:baseHR-global}
Notations as above, there exists a quasi-projective unimodular fan $\Delta_0$ on $\Y_\infty$ such that, for any strictly convex element $\ell_0$ on $\Delta_0$, we have the properties $\HR(\X_{\Delta_0}, \ell_0)$ and $\HL(\X_{\Delta_0}, \ell_0)$ for the compactification $\X_{\Delta_0}$ of $\Y$ in $\TP_{\Delta_0}$.
\end{prop}
\begin{proof} By Theorem~\ref{thm:triangulation_unimodulaire_convexe}, replacing the lattice $N$ with $\frac 1k N$ for some integer $k$, we can find a projective unimodular triangulation $Y_0$ of $\Y$ such that the recession fan $\Delta_0$ of $Y_0$ is quasi-projective and unimodular.

Let $X_{\Delta_0}$ be the compactification of $Y_0$ in $\TP_{\Delta_0}$. Let $\~\ell_0$ be a strictly convex function on $Y_0$. By the second part of Theorem~\ref{thm:triangulation_unimodulaire_convexe}, we can assume that $\ell_0:=(\~\ell_0)_\infty$ is a well-defined and strictly convex function on $\Delta_0$. By theorem \ref{thm:projective_Kahler}, $\~\ell_0$ defines a Kähler class of $X_{\Delta_0}$. Moreover, by Corollary \ref{cor:class_independent_finite_part} and Remark \ref{rem:explicit_function}, this class equals $\class(\ell_0)\in H^{1,1}(\X_{\Delta_0})$. It thus follows from our Theorem~\ref{thm:main} that the pair $(X_{\Delta_0}, \ell_{0})$ verifies $\HL$ and $\HR$.
\end{proof}

\begin{prop} \label{prop:HR-oneall-global}
For any unimodular fan $\Delta$ on the support of $Y_\infty$, we have the equivalence of the following statements.
\begin{itemize}
\item $\HR(X_\Delta, \ell)$ is true for any strictly convex function $\ell$ on $\Delta$.
\item $\HR(X_\Delta, \ell)$ is true for one strictly convex element $\ell$ on $\Delta$.
\end{itemize}
\end{prop}
The proof uses the following proposition.
\begin{prop} \label{prop:global_HR}
If $\HR(D^\varrho, \ell^\varrho)$ holds for all rays $\varrho\in\Delta$, then we have $\HL(X, \ell)$.
\end{prop}
\begin{proof}
By Proposition \ref{prop:class_commutes}, we know that $\ell^\varrho=\i^*(\ell)$ where $\i$ is the inclusion $D^\varrho\hookrightarrow X$. One can now adapt the proof of \ref{prop:local_HR} to the global case to prove the proposition.
\end{proof}

\begin{proof}[Proof of Proposition~\ref{prop:HR-oneall-global}] Let $\ell$ be a strictly convex piecewise linear function on $\Delta$. For all $\varrho \in \Delta$, we get a strictly convex function $\ell^\varrho$ on the star fan $\Sigma^\varrho$ of $\varrho$ in $\Delta$. By the hypothesis of our induction, we know that $\HR(D^\varrho, \ell^\varrho)$ holds. By Proposition~\ref{prop:global_HR} we thus get $\HL(X_\Delta, \ell)$.

By Proposition \ref{prop:HRbis}, we get that $\HR(X_\Delta,\ell)$ is an open and closed condition on the set of all $\ell$ which satisfy $\HL(X_\Delta,\ell)$. In particular, if there exists $\ell_0$ in the cone of strictly convex elements which verifies $\HR(X_\Delta,\ell_0)$, any $\ell$ in this cone should verify $\HR(X_{\Delta},\ell)$.
\end{proof}

\medskip

Let $\Delta$ be a unimodular fan of dimension $d$ on $\Y_\infty$, let $\ell$ be a convex piecewise linear function on $\Delta$, and let $\Delta'$ be the fan obtained from $\Delta$ by star subdividing a cone $\sigma\in\Delta$. Denote by $\rho$ the new ray in $\Delta'$, and let $\ell'$ be the piecewise linear function induced by $\ell$ on $\Delta'$.

Denote by $\chi_\rho$ the characteristic function of $\rho$ on $\Delta'$ which takes value one on $\e_\rho$, value zero on each $\e_\varrho$ for all ray $\varrho \neq \rho$ of $\Delta'$, and which is linear on each cone of $\Delta'$.
The following is straightforward.

\begin{prop} For any small enough $\epsilon>0$, the function $\ell' - \epsilon \chi_\rho$ is strictly convex. \end{prop}

\begin{prop} \label{prop:HR-trans-global}
Notations as above, the following statements are equivalent.
\begin{enumerate}[label=\defaultRoman]
\item \label{enum:HR-trans-global:i} We have $\HR(\X_\Delta, \ell)$.
\item \label{enum:HR-trans-global:ii} The property $\HR(\X_{\Delta'},\ell'-\epsilon \chi_\rho)$ holds for any small enough $\epsilon>0$.
\end{enumerate}
\end{prop}

\begin{proof}
This follows from Theorem~\ref{thm:global_ascent_descent}. We first observe that the $(1,1)$-class associated to the convex function $\ell'-\epsilon \chi_\rho$ is precisely $\ell' - \epsilon \class(D^\rho)$, where in this second expression $\ell'$ represents the class of the convex function $\ell'$ in $H^{1,1}(X_{\Delta'})$. Moreover, by the hypothesis of our induction, we have $\HR(D^\sigma, l_\sigma)$. Now, a direct application of the ascent part of Theorem~\ref{thm:global_ascent_descent} leads to the implication $\ref{enum:HR-trans-global:i} \Rightarrow \ref{enum:HR-trans-global:ii}$.

We now prove $\ref{enum:HR-trans-global:ii} \Rightarrow \ref{enum:HR-trans-global:i}$. Proposition \ref{prop:global_HR} and our induction hypothesis imply that we have $\HL(X_{\Delta}, \ell)$. Applying now the descent part of Theorem~\ref{thm:global_ascent_descent} gives the result.
\end{proof}

\begin{proof}[Proof of Theorem~\ref{thm:main2}] By Proposition~\ref{prop:baseHR-global}, there exists a fan $\Delta_0$ with support $\Y_\infty$ and a strictly convex element $\ell_0$ on $\Delta_0$ such that the result holds for the pair $(\X_{\Delta_0}, \ell_0)$. By Propositions~\ref{prop:HR-trans-global} and~\ref{prop:HR-oneall-global}, for any quasi-projective unimodular fan $\Delta$ with support $\Y_\infty$ which can be obtained from $\Delta_0$ by a sequence of star subdivisions and star assemblies, and for any strictly convex function $\ell$ on $\Delta$, we get $\HR(\X_{\Delta}, \ell)$. The properties $\HR$ and $\HL$ both follow by applying Theorem~\ref{thm:equivalent_fan2}.

\medskip

The weight-monodromy conjecture is preserved by the projective bundle formula, Theorem~\ref{thm:keelglobal}. Thus the property propagates from $\Delta_0$ to all quasi-projective fans by Theorem~\ref{thm:equivalent_fan2}.
\end{proof}

\newpage


\section*{Appendix}

\section{Spectral resolution of the tropical complex}
\label{sec:technicalities}

{\renewcommand{\/}{\backslash}
\allowdisplaybreaks
\def\sign(#1,#2){\varepsilon_{#2/#1}}

In this section, we prove that the objects used to apply the spectral resolution lemma in Section~\ref{sec:steenbrink} verify the hypothesis of this lemma. We also prove Proposition~\ref{prop:isomorphism_filtrations}, which states that the isomorphism between the tropical complex $C^{p,\bul}$ and $\Tot^\bul(\CCp{p}^{\bul,\bul})$ is an isomorphism of filtered complexes for the corresponding differentials and filtrations. We recall all the objects below.

\medskip

We suppose $\X$ is a compact smooth tropical variety with a unimodular triangulation that we denote $X$. We suppose in addition that the underlying polyhedral space $\X$ is obtained by compactifying a polyhedral space $\Y$ in $\R^n$. This will allow later to fix a scalar product in $\R^n$ and define \emph{compatible} projections on subspaces of $\R^n$, in a sense which will be made clear.

\subsection{The triple complex $\AA^{\bul,\bul,\bul}$} We start by recalling the definition of the triple complex $\AA^{\bul, \bul, \bul}$, which is the main object of study in this section (and which is roughly obtained from $\Da{\bul}^{\bul,\bul}$ by inserting the Steenbrink and tropical chain complexes into it).

\medskip

Let $p$ be a non-negative integer. The triple complex $\AA^{\bul, \bul, \bul}$ is defined as follows. We set
\[\AA^{\bul, \sqbullet, -1} := \CCp{p}^{\bul,\sqbullet}, \qquad \AA^{\bul,-1,\blackdiamond} := \STinf{p}^{\bul,\blackdiamond}, \, \textrm{and} \qquad \AA^{\bul,\sqbullet,\blackdiamond} := \Da{\bul}^{\sqbullet,\blackdiamond}.\]
More explicitly, let $a,b,b'\geq 0$ be three non-negative integers. We set

\medskip

\begin{itemize}[label=-]
\item \( \displaystyle \AA^{a,b,-1} := \bigoplus_{\dims\delta = a+b} \AA^{a,b,-1}[\delta] \quad \text{ with } \quad   \AA^{a,b,-1}[\delta]  = \bigwedge^b\TT^\dual\delta \otimes \SF^{p-b}(\conezero^\delta). \)

\medskip

\item \( \displaystyle \AA^{a,-1,b'} = \begin{cases}
  \bigoplus_{\dims\eta = p+a-b' \\ \dims{\eta_\infty} \geq a} \AA^{a,-1, b'}[\eta] & \text{if $b'\geq$ p,} \\
  0 & \text{otherwise,}
\end{cases} \quad \text{ with } \quad \AA^{a,-1,b'}[\eta] := H^{2b'}(\eta).
\)

\medskip

\item \( \AA^{a,b,b'} = \D^{a,b,b'} = \bigoplus_{\delta\subface\eta \\ \sed{\delta}=\sed{\eta} \\ \dims{\delta}=a+b \\ \dims{\eta}=p+a-b' }\AA^{a,b,b'}[\delta,\eta] \quad \text{ with }\quad \AA^{a,b,b'}[\delta,\eta] := \bigwedge^b\TT^\dual\delta \otimes H^{2b'}(\eta).\)

\medskip

\item All the other pieces $\AA^{a,b,b'}$ are trivial.
\end{itemize}

\medskip

We now set $\AA_\trop^{\bul,\sqbullet}:=\AA^{\bul,\sqbullet,-1}$, $\AA_\ST^{\bul,\blackdiamond}:=\AA^{\bul,-1, \blackdiamond}$, and define $\AA_\D^{\bul,\bul,\bul}$ as the truncation of $\AA^{\bul,\bul, \bul}$ at non-negative indices, \ie, such that
\[ \AA_\D^{a,b,b'}=\begin{cases}
\AA^{a,b,b'} & \text{if $a,b,b'\geq0$}, \\
0 & \text{otherwise.}
\end{cases} \]
Note that $\AA_\D^{a,b,b'} = \Da{a}^{b,b'}$ for the double complex $\Da{a}^{b,b'}$ of Section~\ref{sec:steenbrink}.

\medskip

For $\sigma\in X_\infty$, we set
\[ \AAa{\sigma}^{\bul,\bul,\bul}:=\bigoplus_{\delta \text{ with }  \sed(\delta) = \sigma} \AA^{\bul,\bul,\bul}[\delta], \]
and define $\AAa{\sigma}^{\bul, \bul}_\trop, \AAa{\sigma}^{\bul, \bul}_\ST$, and $\AAa{\sigma}^{\bul, \bul, \bul}_\D$, similarly.

The differentials on $\AA^{a,b,b'}$ of respective degree $(1,0,0)$, $(0,1,0)$ and $(0,0,1)$ are denoted $\dfrak$, $\d$ and $\d'$.
We will recall all these maps in Section~\ref{sec:maps}. Their restrictions to $\AA^{\bul,\bul}_\trop$ are denoted by $\dfrak_\trop, \d_\trop$ and $\d'_\trop$ respectively. We only restrict the domain here, and not the codomain. Thus $\d'_\trop$ is a morphism from $\AA^{\bul,\bul}_\trop$ to $\AA^{\bul,\bul,\bul}$. We use the same terminology for the two other complexes $\AA^{\bul,\bul}_\ST$, and $\AA_D^{\bul,\bul,\bul}$.

\medskip

The notation $\d\bigl[\AA^{a,b,b'}[\eta] \to \AA^{a,b+1,b'}[\mu]\bigr]$ denotes restriction of the differential $\d$ to $\AA^{a,b,b'}[\eta]$ projected to the piece $\AA^{a,b+1,b'}[\mu]$, \ie, this is the composition $\pi\d\i$ where $\i\colon\AA^{a,b,b'}[\eta]\hookrightarrow\AA^{\bul,\bul,\bul}$ is the natural inclusion and $\pi\colon\AA^{\bul,\bul,\bul}\twoheadrightarrow\AA^{a,b+1,b'}[\mu]$ is the natural projection.

Moreover, the star $*$ in the expressions $\d[* \to \AA^{a,b+1,b'}]$ and $\d[\AA^{a,b,b'} \to *]$, denotes the natural domain and codomain of $\d$, \ie, $\AA^{\bul,\bul,\bul}$, respectively. Similarly, in the expressions $\d[\, {\prescript{}{\sigma}{*}} \to \AA^{a,b+1,b'}]$, $\d[\AA^{a,b,b'} \to {\prescript{}{\sigma}{*}}\,]$, and $\d[\,{\prescript{}{\sigma}{*}} \to {\prescript{}{\sigma}{*}}\,]$, the ${\prescript{}{\sigma}{*}}$ denotes $\AAa{\sigma}^{\bul, \bul,\bul}$.

We use similar notations for other differentials and other pieces of $\AA^{\bul,\bul,\bul}$.

\medskip

We decompose the differential $\dfrak$ as the sum $\d^\i+\d^\pi$ where
\begin{gather*}
\d^\i := \bigoplus_{\sigma\in X_\infty}\dfrak[\,{\prescript{}{\sigma}{*}} \to {\prescript{}{\sigma}{*}}\,], \myand \\
\d^\pi := \bigoplus_{\tau,\sigma\in X_\infty \\ \sigma\ssupface\tau}\dfrak[\,{\prescript{}{\sigma}{*}} \to {\prescript{}{\tau}{*}}\,].
\end{gather*}

The two other differential $\d$ and $\d'$ preserve $\AAa{\sigma}^{\bul,\bul,\bul}$ for any $\sigma\in X_\infty$.

\medskip

The first main result of this section is the following.
\begin{prop} \label{prop:tech_1}
The complex $\Tot^\bul(\AA_\trop^{\bul,\bul}, \dfrak_\trop+\d_\trop)$ is isomorphic as a complex to $(C_\trop^{p,\bul}, \partial_\trop)$ where $\partial_\trop$ is the usual tropical differential.
\end{prop}

The second result is Proposition~\ref{prop:differential_triple_D}, which, we recall, states the following.
\begin{prop} \label{prop:tech_2}
The map $\partial=\dfrak+\d+\d'$ is a differential, \ie,
\[ \partial\partial = 0. \]
\end{prop}

This section is supposed to be self-contained, in particular we will recall the maps below. We start by (re)introducing some notations and making some conventions. We note that a few of the notations below differ from the ones used in the previous sections.

\subsection{Notations and conventions}
Here is a set of conventions and notations we will use.

\subsubsection*{Naming convention on faces}
In what follows, we take faces $\gamma\ssubface\,\delta,\!\delta'\,\ssubface\chi$ forming a diamond and $\zeta\ssubface\,\eta,\!\eta'\,\ssubface\mu$ forming another diamond such that $\gamma \subface \zeta$, $\delta \subface \eta$, $\delta' \subface \eta',$ and $\chi \subface \mu$.

Moreover, we make the convention that if a relation involves both $\gamma$ and $\delta$, then it is assumed that $\gamma\ssubface\delta$. However, for this relation, we do \emph{not a priori assume} that there exists a face $\chi$ such that $\chi \ssupface\delta$: $\delta$ might be a face of maximum dimension. Unless stated otherwise, we assume that all these faces have the same sedentarity.

Moreover, $\tau\ssubface\sigma\subface\xi$ will be all the time three cones of $X_\infty$.

\medskip

\subsubsection*{Convention on maps $\pi, \iota, \i$ and $\p$}
If $\phi$ and $\psi$ are two faces of $X$ or cones of $\comp{X_\infty}$, we define, when this has a meaning, the map
$\pi^{\phi/\psi}\colon N^\psi_\R \twoheadrightarrow N^\phi_\R$ to be the natural projection. Recall that for $\phi$ a face of $X$ of sedentarity $\conezero$ or a cone of $X_\infty$, we define $N_{\phi,\R}=\TT\phi$ and $N^\phi_\R = \rquot{N_\R}{N_{\phi,\R}}$, and extend this to faces and cones of higher sedentarity by $N^\phi_\R = \rquot{N^{\sed(\phi)}_{\R}}{\TT\phi}$.

\medskip

For the map $\pi^{\phi/\psi}$, it might happen that we use the notation $\pi^{\psi\/\phi}$ in order to insist on the fact that the domain of the map concerns the face $\psi$ and the codomain $\phi$.

\medskip

To simplify the notations, if $\psi=\sed(\phi)$, we set $\pi^\phi:=\pi^{\psi\/\phi}$.
Sometimes, we will put the superscript as a subscript instead, for example $\pi_{\phi/\psi}^*$ is another notation for $\pi^{\phi/\psi *}$. Moreover, depending on the context and when this is clear, we use $\pi^{\psi\/\phi}$ to denote as well the natural projection $\bigwedge^k N^\phi_\R \twoheadrightarrow \bigwedge^k N^\psi_\R$, and also the restriction of this projection to some subspace.

To give an example of these conventions, for instance, the map $\SF^k(\delta) \twoheadrightarrow \SF^k(0^\delta)$ is represented by the notation $\pi^\delta$.

\medskip

We denote by $\eta^\delta$ the face corresponding to $\eta$ in $\Sigma^\delta$, \ie, $\eta^\delta=\R_+\pi^\delta(\eta)$. In the same way, if $\sigma \subface \delta_\infty$, we set
\[\delta_\infty^\sigma:=\pi^\sigma(\delta)=\rquot{\bigl(\delta + N_{\sigma,\R}\bigr)}{N_{\sigma,\R}}.
\]
If $\sed(\delta)=\conezero$, it is the intersection of the closure of $\delta$ in $\TP_{X_\sminfty}$ with the stratum of sedentarity $\sigma$.

\medskip

If $\gamma\ssubface\delta$ are two faces of different sedentarities, then $\pi^{\sed(\delta)\/\sed(\gamma)}$ induces a projection from $\TT\delta$ to $\TT\gamma$. We denote this projection by $\pi^{\sed}_{\delta/\gamma}$. Note that this is quite different from $\pi^{\gamma\/\delta}$; for such a pair, the projection $\pi^{\gamma\/\delta}$ is in fact an isomorphism $N^\delta_\R \simto N^\gamma_\R$ and induces an isomorphism $\Sigma^\delta \simto \Sigma^\gamma$. We denote these maps and the corresponding pushforwards by $\id^{\delta/\gamma}$, and the inverse maps and the corresponding pushforwards by $\id^{\gamma\/\delta}$.

\medskip

In the same way, we set $\iota_{\eta\/\delta}\colon \SF_\bul\eta \hookrightarrow \SF_\bul\delta$ and $\iota_{\delta}\colon \SF_\bul\delta \hookrightarrow \bigwedge^\bul N^{\sed(\delta)}_\R$ to be the natural inclusions. (By our convention on faces, we have $\delta\subface \eta$.) The inclusion concerning the canonical compactifications of Bergman fans is denoted by $\i^{\eta/\delta}\colon \comp{\Sigma^\eta} \hookrightarrow \comp{\Sigma^\delta}$. The same simplifications hold.

\medskip

We also introduce projections $\p_{\delta}\colon N^{\sed(\delta)}_\R \twoheadrightarrow \TT\delta$ as follows. Let $(\ccdot,\rdot)$ be an inner product on $N_\R$. If $\sed(\delta)=\conezero$, then we define $\p_{\delta}$ to be the orthogonal projection with respect to $(\ccdot,\rdot)$. We extend this definition to other strata as follows.

For any face $\phi$ of sedentarity $\sigma$, we observe that we can find a unique face $\delta$ of sedentarity $\conezero$ such that $\phi=\delta_\infty^\sigma$, and we define $\p_\phi$ as the quotient map
\[\p_\phi:=\pi^{\sigma}_*(\p_\delta),\]
making the following diagram commutative:
\[
\begin{tikzcd}
N_\R \dar["\pi^\sigma"] \rar["\p_\delta"] & \TT\delta \dar["\pi^\sigma"]  \\
N^\sigma_\R=\rquot{N_\R}{N_{\sigma, \R}} \rar["\p_\phi"] &  \TT\phi = \rquot{\TT \delta}{N_{\sigma, \R}}
\end{tikzcd} \]
In addition, if $\sed(\eta)=\sed(\delta)$, then we define the relative projection map $\p_{\eta/\delta}$ by setting
\[\p_{\eta/\delta}:=\p_\delta\rest{\TT\eta}.
\]
If $\phi$ is a subface of $\delta$ with the same sedentarity, we set $\p_\delta^\phi := \pi^\phi_*(\p_\delta)$ and $\p_{\eta/\delta}^\phi := \pi^\phi_*(\p_{\eta/\delta})$. We also set $\cp_\delta:=\id-\p_\delta$ the orthogonal projection onto $\ker(\p_\delta)$. We use the same notation for $\cp_\delta^\phi$, etc.

Notice that $\pi^\delta\rest{\ker(\p_\delta)}\colon \ker(\p_\delta) \to N^\delta_\R$ is an isomorphism. To simplify the notation, we set $\pi_\delta^{-1}:=(\pi^\delta\rest{\ker(\p_\delta)})^{-1}$. In the same way, $\pi_{\eta/\delta}^{-1}:=(\pi^{\eta/\delta}\rest{\ker(\p^\delta_{\eta})})^{-1}$.

\medskip

\subsubsection*{Operations on local cohomology groups}

Recall that $\i^{\eta/\delta}$ induces a surjective map $\i_{\delta\/\eta}^*$ from $H^\bul(\delta)$ to $H^\bul(\eta)$ of weight $0$, where $H^\bul(\delta):=H^\bul(\comp{\Sigma^\delta})$. The Poincaré dual map of $\i_{\delta\/\eta}^*$ is denoted $\gys^{\eta/\delta}$. It is an injective map of weight $2(\dims\eta-\dims\delta)$ from $H^\bul(\eta)$ to $H^\bul(\delta)$. We denote by $1^\delta\in H^0(\delta)$ the corresponding unit. We denote by $x_{\delta/\gamma} \in A^1(\delta)\simeq H^2(\delta)$ the element corresponding to the ray $\delta^\gamma$.

\medskip

\subsubsection*{Operations on faces}
If $\phi$ and $\psi$ are subfaces of a same face $\delta$, we denote by $\phi+\psi$ the smallest face containing $\phi$ and $\psi$. When this has a meaning, we denote by $\eta-\delta$ the smallest face $\phi$ such that $\phi+\delta=\eta$. By $\eta-\delta+\gamma$, for instance, we mean $(\eta-\delta)+\gamma$. Abusing of the notation, by $\p_{\eta-\delta+\gamma/\gamma}$ we mean $\p_{(\eta-\delta+\gamma)/\gamma}$. We might sometimes write $\p_{\eta-\delta/\gamma}$ instead of $\p_{\eta-\delta+\gamma/\gamma}$.

\medskip

\subsubsection*{Multivectors and multiforms $e$ and $\nu$} We have a sign function $\sign(\ccdot,\rdot)=\textrm{sign}(\ccdot,\rdot)$ which verifies the diamond property
\[\sign(\gamma,\delta)\sign(\delta,\eta)=\sign(\gamma,\delta')\sign(\delta',\eta),\]
here the faces are allowed to have different sedentarities, and $\sign(v,e)=-\sign(v',e)$ if $v$ and $v'$ are the two extremities of an edge $e$.

Let $\e_{\delta/\gamma} \in \TT{\delta^\gamma}$ be the primitive element of the ray $\delta^\gamma$. Let $e_{\delta/\gamma}:=\sign(\gamma,\delta)\e_{\delta/\gamma}$. Let $s:=\dims\eta-\dims\delta$ and let $\delta_0, \delta_1, \dots, \delta_s$ be a flag of faces such that
\[ \delta=\delta_0\ssubface\delta_1\ssubface\cdots\ssubface\delta_s=\eta. \]
We define the multivector $e_{\eta/\delta} \in \bigwedge^s\TT{\eta^\delta}$ by
\[
e_{\eta/\delta} :=\sign(\delta_0,\delta_1)\dots\sign(\delta_{s-1},\delta_s)\e_{\delta_1/\delta} \wedge \e_{(\delta_2-\delta_1)/\delta} \wedge \dots \wedge \e_{(\delta_s-\delta_{s-1})/\delta}.
\]
This definition is independent of the choice of the flag $\delta_0,\dots,\delta_s $.

\medskip

We set $\nvect_{\delta/\gamma}:=\pi_\gamma^{-1}(\e_{\delta/\gamma})$ and $\nu_{\eta/\delta}:=\pi_\delta^{-1}(e_{\eta/\delta})$. (Recall that by our convention, $\pi_\delta^{-1}$ denotes the inverse of the restriction of $\pi_\delta$ to the kernel of the projection map $\p_\delta$, which is an isomorphism.)

If $\phi$ is a subface of $\delta$ with the same sedentarity, we set $\nu^\phi_{\eta/\delta}:=\pi^\phi(\nu_{\eta/\delta})$. In particular, $\nu^\delta_{\eta/\delta}=e_{\eta/\delta}$. Similarly, for a subface $\phi$ of $\gamma$ with the same sedentarity, we set $\nvect^\phi_{\delta/\gamma}:=\pi^\phi(\nvect_{\delta/\gamma})$.

\smallskip
Let $e_{\eta\/\delta}^\dual\in\bigwedge^s\TT^\dual\eta^\delta$ be the dual of $e_{\eta/\delta}$. Set $\nu^\dual_{\delta\/\eta}:=\pi_\delta^*(e^\dual_{\eta\/\delta})\in\bigwedge^s\TT^\dual\eta$ and $\nu^{\phi\,*}_{\eta\/\delta}:=\pi_{\delta/\phi}^*(e^\dual_{\eta\/\delta})\in\bigwedge^s\TT^\dual\eta^\phi$.

\smallskip
We also denote by $1_\delta^\dual\in\bigwedge^0\TT^\dual\delta$ the corresponding unit.

\medskip

\subsubsection*{Operations on multivectors and on multiforms}
Let $V$ is a vector space. Let $\alpha\in\bigwedge^r V^\dual$ be a multiform on $V$ and $\u\in\bigwedge^s V$ be a multivector in $V$. Then we denote by $\alpha \vee \u \in \bigwedge^{r-s}V^\dual$ the \emph{right contraction of $\alpha$ by $\u$}, \ie,
\[ \alpha\vee\u := \begin{cases}
  \alpha(\ldot \wedge \u) & \text{if $s\geq r$}, \\
  0 & \text{otherwise}.
\end{cases} \]
We also define the \emph{left contraction of $\alpha$ by $\u$} by
\[ \u\vee\alpha := \begin{cases}
  \alpha(\u \wedge \ccdot) & \text{if $s\geq r$}, \\
  0 & \text{otherwise}.
\end{cases} \]

We denote by $\ker(\alpha)$ the subspace of $V$ defined by
\[ \ker(\alpha) := \{u \in V \st \alpha\vee u = 0 \}. \]
If $\alpha$ is nonzero, the kernel is of dimension $\dim(V)-r$. If $\beta \in \bigwedge^t \ker(\alpha)^\dual$ then we denote by $\beta\wedge\alpha \in \bigwedge^{t+r}V^\dual$ the multiform $\~\beta\wedge\alpha$ where $\~\beta$ is any extension of $\beta$ to $\bigwedge^t V$. The definition of $\beta\wedge\alpha$ does not depend on the chosen extension.

\subsection{Useful facts}
In this section, we list some useful facts which will be used in the proofs of Propositions~\ref{prop:tech_1} and~\ref{prop:tech_2}.

\begin{prop} Using the conventions introduced in the previous section, the following holds.
\newcommand{\sameline}{\ \vspace*{-1em}}
\let\olditem\item
\renewcommand{\item}{\stepcounter{equation}\olditem}
\begin{enumerate}[label=\textnormal{(\theequation)}, ref={\theequation}]
\item \label{com:eqn:maxsed_sed} If $\sed(\gamma)\ssupface\sed(\delta)$ then
  \[ \maxsed(\gamma)=\maxsed(\delta). \]
\item \label{com:eqn:maxsed/sed} \sameline
  \[ \rquot{\maxsed(\delta)}{\sed(\delta)}=\delta_\infty. \]
\item \label{com:eqn:zeta_pi} If $\sed(\eta)=\sed(\delta)\ssubface\sed(\gamma)$, then
  \[ \zeta=\eta_\infty^{\sed(\gamma)}. \]
\end{enumerate}

\bigskip

\noindent Let $\ell$ and $\ell'$ be two linear forms, $\alpha$ be a multiform and $\u$ and $\v$ be two multivectors on a vector space $V$. the following holds.
\begin{enumerate}[resume*]
\item \sameline
  \[ \alpha \vee (\u \wedge \v) = \alpha \vee \v \vee \u. \]
\item \label{com:eqn:wedge_vee} If $u\in\ker(\ell)$, then
  \[ \alpha \wedge \ell \vee u = -\,\alpha \vee u \wedge \ell. \]
\item \label{com:eqn:wedge_dual} If $u'\in\ker(\ell)$, then
  \[ (\ell\wedge\ell')(u\wedge u')=\ell(u)\cdot\ell'(u'). \]

\bigskip

\item \label{com:eqn:gys_x} \sameline
  \[ \i^*_{\zeta/\eta}\gys^{\eta/\zeta}(x)=\i^*_{\zeta/\eta}\gys^{\eta/\zeta}(1^\eta)x. \]
\item \label{com:eqn:i*_x} \sameline
  \[ \gys^{\eta/\zeta}\i^*_{\zeta/\eta}(x)=\gys^{\eta/\zeta}\i^*_{\zeta/\eta}(1^\eta)x. \]
\item \label{com:eqn:gys} \sameline
  \[ \gys^{\eta/\zeta}(1^\eta)=\i^*_{\gamma\/\zeta}(x_{\eta-\zeta/\gamma}). \]

\bigskip

\item \label{com:eqn:nu} \sameline
  \[ \nu_{\eta/\phi} = \nu_{\delta/\phi} \wedge \nu_{\eta/\delta}. \]
\item \label{com:eqn:nu*} \sameline
  \[ \nu^\dual_{\phi\/\eta} = \nu^\dual_{\phi\/\delta} \wedge \nu^\dual_{\delta\/\eta}. \]
\item \label{com:eqn:pi_nu*} If $\sed(\gamma)=\sed(\delta)\ssupface\sed(\delta')=\sed(\eta)$, then,
  \[ \pi^{\sed*}_{\chi/\delta}(\nu^\dual_{\gamma\/\delta})=\nu^\dual_{\delta'\/\chi}. \]
\item If $\beta$ is a multiform on $N^\delta_\R$,
  \[ \pi^{\delta\/\eta}_*(\beta \vee \nu^\delta_{\eta/\delta}) \text{ is well-defined.} \]
\item If $\phi$ is a subface of $\delta$ of same sedentarity, then
  \[ \cp_\delta(\nu_{\eta/\phi}) = \nu_{\eta/\delta}. \]
\item \label{com:eqn:nu*_nu} \sameline
  \[ \p_\delta(\nu_{\delta'/\gamma})=\nu^\dual_{\gamma\/\delta}(\p_\delta(\nu_{\delta'/\gamma}))\nu_{\delta/\gamma}. \]
\item If $\eta=\delta+\zeta$, even if $\sed(\delta)\neq\sed(\gamma)$ (in which case $\pi^{\gamma\/\delta}=\id^{\gamma\/\delta}$),
  \[ \pi^{\gamma\/\delta}(e_{\eta/\delta}) = (-1)^{\dims\eta-\dims\delta}\sign(\gamma,\delta)\sign(\zeta,\eta)e_{\zeta/\gamma}. \]

\bigskip

\item \label{com:eqn:p_nu*_nu} \sameline
  \[ \p^*_{\gamma\/\delta}=(\ldot\wedge\nu^\dual_{\gamma\/\delta})\vee\nu_{\delta/\gamma}. \]
\item \sameline
  \[ \p_\delta \p_\eta=\p_\delta \p_\eta = \p_\delta \quad\myand\quad \cp_\delta \p_\eta=\p_\eta \cp_\delta. \]
\item \label{com:eqn:p_pi} If $\sed(\gamma)=\sed(\delta)\ssupface\sed(\delta')=\sed(\eta)$, then
  \[ \pi^{\sed}_{\delta'/\gamma}\p_{\chi/\delta'} = \p_{\delta/\gamma}\pi^{\sed}_{\chi/\delta}. \]
\end{enumerate}
\end{prop}

\begin{proof} We omit the proof which can be obtained by straightforward computations.
\end{proof}

\subsection{Proof of Proposition~\ref{prop:tech_1}}
\label{sec:com:trop}

In this section, we define a (non canonical) isomorphism $\Phi\colon \Tot^\bul(\AA_\trop^{\bul,\bul}) \simto C_\trop^{p,\bul}$. We then prove that we can choose two differentials $\dfrak_\trop$ and $\d_\trop$ of bidegree $(1,0)$ and $(0,1)$ on $\AA_\trop^{\bul, \bul}$ such that we have the equality $\dfrak_\trop+\d_\trop = \Phi^{-1} \partial_\trop \Phi$ where $\partial_\trop$ is the usual tropical differential on $C_\trop^{p,\bul}$. The two differentials we obtain are those we have already seen in Section \ref{sec:steenbrink}. This will prove Proposition~\ref{prop:tech_1}.

\medskip

Recall that, if $\dims\delta = a+b$,
\[ \AA_\trop^{a,b}[\delta] = \bigwedge^b\TT^\dual\delta \otimes \SF^{p-b}(\conezero^\delta)  \quad\myand\quad C_\trop^{p,\dims\delta}[\delta] = \SF^p(\delta), \]
and other pieces are zero. The differential $\partial_\trop$ is equal to the sum $\iota^*_\trop+\pi^*_\trop$ where
\[ \iota^*_\trop \bigl[C_\trop^{p,\dims\gamma}[\gamma] \to C_\trop^{p,\dims\delta}[\delta]\bigr] = \sign(\gamma,\delta)\iota^*_{\gamma\/\delta} \]
and, if $\sed(\gamma)\ssupface\sed(\delta)$, then
\[ \pi^*_\trop \bigl[C_\trop^{p,\dims\gamma}[\gamma] \to C_\trop^{p,\dims\delta}[\delta]\bigr] = \sign(\gamma,\delta)\pi^*_{\gamma\/\delta}. \]

\medskip

Recall as well from Section~\ref{sec:steenbrink} that we have a filtration $W^\bul$ on $C_\trop^{p,\bul}$. With our notations introduced above, this is defined by
\[ W^sC_\trop^{p,\dims\delta}[\delta] := \Bigl\{\, \alpha\in\SF^p(\delta) \st \forall \u\in\bigwedge^{\dims\delta-s+1}\TT_\delta, \alpha\vee\u=0 \Bigr\}. \]
We proved in Section~\ref{sec:steenbrink} that this filtration is preserved by $\partial_\trop$.

\medskip

By an abuse of the notation, we also denote by $W^\bul$ the filtration on $\AA_\trop^{\bul,\bul}$ induced by the first index.

\medskip

We now define $\Phi\colon\AA_\trop^{\bul,\bul} \to C_\trop^{p,\bul}$ by
\[ \Phi\bigl[\AA_\trop^{a,b}[\delta] \to C_\trop^{p,a+b}[\delta]\bigr] = \p^*_\delta \wedge \pi^*_\delta\]
if $a+b=\dims\delta$, and all other pieces are zero. Set
\[ \Phi_a := \Phi[\AA_\trop^{a,\bul} \to *]. \]

We define $\Psi\colon C_\trop^{p,\bul} \to \AA_\trop^{\bul,\bul}$ by
\[ \Psi\bigl[C_\trop^{p,a+b}[\delta] \to \AA_\trop^{a,b}[\delta]\bigr](\alpha)(\u \otimes \v) = \alpha(\u \wedge \pi^{-1}_\delta(\v)) \]
if $a+b=\dims\delta$, and other pieces are zero. Set
\[ \Psi_a := \Psi[* \to \AA_\trop^{a,\bul}]. \]

\begin{prop}
The map $\Phi$ and $\Psi$ are inverse isomorphisms of weight zero, and both preserve the filtration $W^\bul$. Moreover, for any $a$ and $k$, $\Psi_{a+k}\partial_\trop\Phi_a = 0$ is zero except for $k\in\{0,1\}$.
\end{prop}

Set $\d_\trop := \bigoplus_{a}\Psi_a \partial_\trop \Phi_a$ and $\dfrak_\trop := \d^\i_\trop + \d^\pi_\trop$ where
\[ \d^\i_\trop := \Psi_{a+1} \iota^*_\trop \Phi_a \quad\myand\quad \d^\pi_\trop := \Psi_{a+1} \pi^*_\trop \Phi_a. \]
By the previous proposition it is now clear that $\Phi\colon (\Tot^\bul(\AA_\trop^{\bul,\bul}), \dfrak_\trop + \d_\trop) \simto (C_\trop^{p,\bul}, \partial_\trop)$ is an isomorphism of complexes which respects the filtrations. It remains to give an explicit formula of $\d_\trop$ and $\dfrak_\trop$. This can be shown by direct computations.

\begin{prop}
If $\dims\delta=a+b$, then
\[ \d_\trop\bigl[\AA_\trop^{a,b-1}[\gamma] \to \AA_\trop^{a,b}[\delta]\bigr](\alpha \otimes \beta) = \sign(\gamma,\delta)(\alpha \wedge \nu^\dual_{\gamma\/\delta}) \otimes (\pi^{\gamma\/\delta}_*(e_{\delta/\gamma} \vee \beta)), \]
as stated in Proposition \ref{prop:grading}, and
\[ \d^\i_\trop\bigl[\AA_\trop^{a-1,b}[\gamma] \to \AA_\trop^{a,b}[\delta]\bigr] = \begin{cases}
  \sign(\gamma, \delta)\p_{\gamma\/\delta}^* \otimes (\pi_{\delta/\gamma}^{-1})^*, & \text{if $\sed(\gamma)=\sed(\delta)$,} \\
  0 & \text{otherwise,}
\end{cases} \]
as stated by Equation \eqref{eqn:d''}. Moreover,
\[ \d^\pi_\trop\bigl[\AA_\trop^{a-1,b}[\gamma] \to \AA_\trop^{a,b}[\delta]\bigr] = \begin{cases}
  \sign(\gamma, \delta)\pi^{\sed*}_{\gamma\/\delta} \otimes \id & \text{if $\sed(\gamma) \ssupface \sed(\delta)$,} \\
  0 & \text{otherwise.}
\end{cases} \]
\end{prop}

It thus follows that $\d_\trop$ and $\dfrak_\trop$ defined above coincide with $\d_\trop$ and $\dfrak_\trop$ of Section~\ref{sec:steenbrink}.

\subsection{Differentials of $\AA^{\bul,\bul,\bul}$}
\label{sec:maps}

We now recall the definition of the differentials on $\AA^{\bul,\bul,\bul}$. Recall that, for $\dims\delta=a+b$ and $\dims\eta=p+a-b'$,
\begin{align*}
\AA_\trop^{a,b}[\delta] &= \bigwedge^b\TT^\dual\delta \otimes \SF^{p-b}(\conezero^\delta), \\
\AA_\ST^{a,b'}[\eta] &= H^{2b'}(\eta)\text{ if $\dims{\eta_\infty} \geq a$ and $b'\geq p$, and} \\
\AA^{a,b,b'}_\D[\delta,\eta] &= \bigwedge^b\TT^\dual\delta \otimes H^{2b'}(\eta).
\end{align*}

When the domain or the codomain of a map is clear, we only indicate the concerned faces. For $\d^\pi$, we assume that $\sed(\gamma)=\sed(\zeta)\ssupface\sed(\delta)=\sed(\eta)$. We have
\begin{align*}
\d_\trop[\gamma \to \delta]
  &= \sign(\gamma,\delta)(\ldot \wedge \nu^\dual_{\gamma\/\delta}) \otimes \pi^{\gamma\/\delta}_*(e_{\delta/\gamma} \vee \rdot), \\
\d'_\trop\bigl[\AA_\trop^{a,b}[\delta] \to \delta,\eta \bigr]
  &= (-1)^{(a+1)(b+p)+\dims{\sed(\delta)}}\id \otimes (\ldot \vee e_{\eta/\delta})1_\eta, \\
\d^\i_\trop[\gamma \to \delta]
  &= \sign(\gamma, \delta)\p_{\gamma\/\delta}^* \otimes (\pi_{\gamma\/\delta}^{-1})^*, \\
\d^\pi_\trop[\gamma \to \delta]
  &= \sign(\gamma, \delta)\pi^{\sed*}_{\gamma\/\delta} \otimes \id^{\gamma\/\delta},
\end{align*}
\begin{align*}
\d_\ST\bigl[\AA_\ST^{a,b'}[\eta] \to \delta,\eta\bigr]
  &= (-1)^{a+b'}1^\dual_\delta \otimes \id  \quad\text{if $\maxsed(\delta)=\maxsed(\eta)$}, \\
\d'_\ST[\mu \to \eta]
  &= \sign(\eta,\mu)\gys^{\mu/\eta}  \quad\text{if $\maxsed(\eta)=\maxsed(\mu)$}, \\
\d^\i_\ST[\eta \to \mu]
  &= \sign(\eta,\mu) \i^*_{\eta\/\mu}  \quad\text{if $\maxsed(\eta)=\maxsed(\mu)$}, \\
\d^\pi_\ST[\eta \to \mu]
  &= \sign(\eta,\mu)\id^{\eta\/\mu}  \quad\text{(by \eqref{com:eqn:maxsed_sed}, $\maxsed(\eta)=\maxsed(\mu)$)},
\end{align*}
\begin{align*}
\d_\D\bigl[\AA^{a,b,b'}[\gamma,\eta] \to \delta,\eta\bigr]
  &= (-1)^{a+b'}\sign(\gamma,\delta)(\ldot \wedge \nu^\dual_{\gamma\/\delta}) \otimes \id, \\
\d'_\D[\delta,\mu \to \delta,\eta]
  &= \sign(\eta,\mu)\id \otimes \gys^{\mu/\eta}, \\
\d^\i_\D[\gamma,\zeta \to \delta,\eta]
  &= \sign(\zeta,\eta) \nvect^\dual_{\gamma\/\delta}(\p_\delta(\nvect_{\eta-\zeta/\gamma}))\p_{\gamma\/\delta}^* \otimes \i^*_{\zeta\/\eta}, \\
\d^\pi_\D[\gamma,\zeta \to \delta,\eta]
  &= \sign(\zeta,\eta)\pi^{\sed*}_{\gamma\/\delta} \otimes \id^{\zeta\/\eta}.
\end{align*}

\subsection{Proof that $\partial=\dfrak+\d+\d'$ is a differential} Let $\partial = \dfrak+\d+\d'$ with $\dfrak=\d^\i+\d^\pi$. In this section, we prove that $\partial\circ\partial=0$. To do so, we separate thirty cases summarized in the following table. For example, the first case (denoted $a$ in the top-left corner) concerns the proof that $(\d^2)[\AA_\trop^{\bul,\bul} \to *]=0$. The case just after is the proof that $(\d\d'+\d'\d)[\AA_\trop^{\bul,\bul} \to *]=0$, etc. Some cases has been gathered (in particular the easy cases and those that have been already done). In each case, we fix a context by precising the faces involved. We let the reader check that the studied case is general: no other configuration of faces can be involved.

{
\newcommand\taa{\ref{com:trop}}
\newcommand\tab{\ref{com:tab}}
\newcommand\tac{\ref{com:trop}}
\newcommand\tad{\ref{com:trop}}
\newcommand\tbb{\ref{com:tbb}}
\newcommand\tbc{\ref{com:tbc}}
\newcommand\tbd{\ref{com:tbd}}
\newcommand\tcc{\ref{com:trop}}
\newcommand\tcd{\ref{com:trop}}
\newcommand\tdd{\ref{com:trop}}

\newcommand\saa{\ref{com:saa}}
\newcommand\sab{\ref{com:ST_pi}}
\newcommand\sac{\ref{com:sac}}
\newcommand\sad{\ref{com:ST_pi}}
\newcommand\sbb{\ref{com:ST}}
\newcommand\sbc{\ref{com:ST}}
\newcommand\sbd{\ref{com:ST_pi}}
\newcommand\scc{\ref{com:ST}}
\newcommand\scd{\ref{com:ST_pi}}
\newcommand\sdd{\ref{com:ST_pi}}

\newcommand\aaa{\ref{com:a_easy}}
\newcommand\aab{\ref{com:a_easy}}
\newcommand\aac{\ref{com:aac}}
\newcommand\aad{\ref{com:a_easy}}
\newcommand\abb{\ref{com:abb}}
\newcommand\abc{\ref{com:abc},\ref{com:abc'},\ref{com:abc''}}
\newcommand\abd{\ref{com:a_easy}}
\newcommand\acc{\ref{com:acc}}
\newcommand\acd{\ref{com:acd}}
\renewcommand\add{\ref{com:a_easy}}
\[ \bigarray \]
}

\begin{enumerate}[label={\bf(\alph*)}, ref=\alph*, leftmargin=0pt]
\newcommand{\Lar}{\Longleftarrow\quad}
\newcommand{\Rar}{\Longrightarrow\quad}
\newcommand{\LRar}{\Longleftrightarrow\quad}

\item \label{com:trop}$(\d_\trop+\dfrak_\trop)^2=0$. This was already done in Section \ref{sec:com:trop}.

\medskip

\item \label{com:tab} $\d\d'+\d'\d =0$ in $\AA^{\bul, \bul}_\trop$. We need to show that for integers $a,b$, a face $\gamma$ with $\dims \gamma =a+b$, and $\alpha \otimes \beta \in\bigwedge^b\TT^\dual\gamma \otimes \SF^{p-b}(\conezero^\gamma)$, we have $(\d\d'+\d'\d)\bigl[\AA_\trop^{a,b}[\gamma] \to \AA^{a,b,0}[\delta,\eta]\bigr](\alpha\otimes\beta) = 0$. This is obtained by the following chain of implications.
\begin{align*}
& \hspace{-2cm}\nu^\gamma_{\eta/\gamma}
  = \nu^\gamma_{\delta/\gamma} \wedge \nu^\gamma_{\eta/\delta}\\
\Rar& (-1)^a\beta\vee e_{\eta/\gamma}
  = -(-1)^{a+1}(e_{\delta/\gamma} \vee \beta) \vee \pi_{\delta/\gamma}^{-1}(e_{\eta/\delta})\\
\Rar& (-1)^{a+0}\sign(\gamma,\delta)(-1)^{(a+1)(b+p)+\dims{\sed(\gamma)}}\alpha\vee\nu^\dual_{\gamma\/\delta}\otimes (\beta\vee e_{\eta/\gamma})1_\eta
  \\&\qquad+ (-1)^{(a+1)(b+1+p)+\dims{\sed(\delta)}}\sign(\gamma,\delta)\alpha\vee\nu^\dual_{\gamma\/\delta} \otimes (\pi^{\gamma\/\delta}_*(e_{\delta/\gamma} \vee \beta) \vee e_{\eta/\delta}) 1_\eta = 0,
\end{align*}
which is the desired equality given the differentials.

\medskip

\item \label{com:tbb} We have $\d'\circ \d'=0$ in $\AA^{\bul, \bul}_\trop$. Let $a, b$ be two integers, let $\delta$ be a face of dimension $\dims\delta =a+b$, and choose $\alpha\otimes \beta \in\bigwedge^b\TT^\dual\delta \otimes \SF^{p-b}(\conezero^\delta)$.
Then we have to prove that $(\d'^2)\bigl[\AA_\trop^{a,b}[\delta] \to \AA^{a,b,1}[\delta,\eta]\bigr](\alpha\otimes\beta) = 0$. We have, using the linear relations in the cohomology ring $H^\bul(\eta)$ via the isomorphism with the Chow ring of $\Sigma^\eta$,
\begin{align*}
& \hspace{-2cm}\ssum_{\mu'} \bigl(\pi^{\delta\/\eta}_*(\beta \vee \nu^\delta_{\eta/\delta})\bigr)(\e_{\mu'/\eta}) x_{\mu'/\eta} = 0 \\
\Rar& (-1)^{\dims\eta-\dims\delta}\ssum_{\mu'} (\beta \vee \nu^\delta_{\eta/\delta} \vee \nvect^\delta_{\mu'/\eta}) x_{\mu'/\eta} = 0 \\
\Rar& \ssum_{\mu'} \sign(\eta,\mu') (\beta \vee \nu^\delta_{\mu'/\eta} \vee \nu^\delta_{\eta/\delta}) x_{\mu'/\eta} = 0 \\
\Rar& \sum_{\mu'\ssupface\eta \\ \sed(\mu')=\sed(\eta)} \sign(\eta,\mu')(-1)^{(b+p)(a+1)+\dims{\sed(\delta)}} \alpha \otimes (\beta \vee e_{\mu'/\delta}) \gys^{\mu'/\eta}(1_\mu') = 0.
\end{align*}

\medskip

\item \label{com:tbc} We have to show that $\d'\d^\i+\d^\i\d'=0$ in $\AA^{\bul, \bul}_\trop$. Let $a, b$ be two integers, let $\gamma$ be a face of dimension $\dims\gamma =a+b$, and choose $\alpha\otimes \beta \in\bigwedge^b\TT^\dual\gamma \otimes \SF^{p-b}(\conezero^\gamma)$. We need to prove
\[(\d'\d^\i+\d^\i\d')[\AA_\trop^{a,b}[\gamma] \to \AA^{a+1,b,0}[\delta,\eta]](\alpha \otimes \beta) = 0.\]
This is equivalent to
\begin{align*}
&\LRar \sign(\gamma,\delta) (-1)^{(a+2)(b+p)+\dims{\sed(\delta)}} \p^*_{\gamma\/\delta}(\alpha) \otimes ((\pi_{\gamma\/\delta}^{-1})^*(\beta) \vee e_{\eta/\delta})1_\eta
  \\&\quad + \sum_{\zeta' \ssubface \eta \\ \zeta' \supface \gamma} \sign(\zeta',\eta) (-1)^{(a+1)(b+p)+\dims{\sed(\gamma)}} \nvect^\dual_{\gamma\/\delta}(\p_\delta(\nvect_{\eta-\zeta'/\gamma})) \p^*_{\gamma\/\delta}(\alpha) \otimes (\beta \vee e_{\zeta'/\gamma}) \i^*_{\zeta'\/\eta}(1_{\zeta'}) =0 \\
&\LRar \sign(\gamma,\delta) (-1)^{p-(b+1)} \pi_{\delta/\gamma}^{-1}(e_{\eta/\delta})
  = {\ssum}_{\zeta'} \sign(\zeta',\eta) \nvect^\dual_{\gamma\/\delta}(\p_\delta(\nvect_{\eta-\zeta'/\gamma})) e_{\zeta'/\gamma} \\
&\LRar \nu_{\eta/\delta}
  = \sign(\gamma,\delta) (-1)^{\dims\eta-\dims\delta} {\ssum}_{\zeta'} \sign(\zeta',\eta) \nvect^\dual_{\gamma\/\delta}(\p_\delta(\nvect_{\eta-\zeta'/\gamma})) \nu_{\zeta'/\gamma}.
\end{align*}

Thus, we reduce to proving the last equation. Let $s=\dims\eta-\dims\delta$. Let $\delta_0,\dots,\delta_s$ such that
\[ \delta=\delta_0\ssubface\delta_1\ssubface\cdots\ssubface\delta_s=\eta. \]
Set $\epsilon=\delta-\gamma$, $\gamma_i=\delta_i-\epsilon$ and $\zeta_i=\eta-(\delta_{i+1}-\delta_{i})=\eta-(\gamma_{i+1}-\gamma_i)$. Then
\begin{align*}
\nu_{\eta/\delta}
  &= \oldbigwedge_{i=0}^{s-1} \sign(\delta_i,\delta_{i+1}) \nvect_{\delta_{i+1}-\delta_i/\delta} \\
  &= \bigwedge_{\,i} -\sign(\gamma_i,\delta_i)\sign(\gamma_i,\gamma_{i+1})\sign(\gamma_{i+1},\delta_{i+1}) \cp_\delta(\nvect_{\gamma_{i+1}-\gamma_i/\gamma}) \\
  &= (-1)^s\sign(\gamma_0,\delta_0)\sign(\gamma_s,\delta_s) \bigwedge_{\,i} \sign(\gamma_i,\gamma_{i+1}) \bigl(\nvect_{\gamma_{i+1}-\gamma_i/\gamma} - \p_\delta(\nvect_{\gamma_{i+1}-\gamma_i/\gamma})\bigr).
\end{align*}
By \eqref{com:eqn:nu*_nu},
\[ \p_\delta(\nvect_{\gamma_{i+1}-\gamma_i/\gamma})=\nvect^\dual_{\gamma\/\delta}(\p_\delta(\nvect_{\eta-\zeta_i}))\nvect_{\delta/\gamma}. \]
Hence,
\begin{align*}
&\nu_{\eta/\delta} = (-1)^s\sign(\gamma,\delta)\sign(\eta-\epsilon,\eta) \Bigl( \bigwedge_{\,i} \sign(\gamma_i,\gamma_{i+1}) \nvect_{\gamma_{i+1}-\gamma_i/\gamma} \\
&\qquad- \sum_i (-1)^{s-i-1}\sign(\gamma_i,\gamma_{i+1})\nvect^\dual_{\gamma\/\delta}(\p_\delta(\nvect_{\eta-\zeta_i/\gamma}))\oldbigwedge_{\substack{j=0 \\ j\neq i}}^{s-1} \sign(\gamma_j,\gamma_{j+1}) \nvect_{\gamma_{j+1}-\gamma_j/\gamma}\ \wedge \nvect_{\delta/\gamma}  \Bigr).
\end{align*}
In the big brackets, the first wedge product is simply $\nu_{\gamma_s/\gamma_0}=\nu_{\eta-\epsilon/\gamma}$. For the sum, set $\kappa_i:=\gamma_{i+1}-\gamma_i$. We get
\begin{align*}
&\sign(\gamma_j,\gamma_{j+1}) \nvect_{\gamma_{j+1}-\gamma_j/\gamma} =  \\
&\qquad-\sign(\gamma_j-\kappa_i,\gamma_j)\sign(\gamma_j-\kappa_i,\gamma_{j+1}-\kappa_i)\sign(\gamma_{j+1}-\kappa_i,\gamma_{j+1}) \nvect_{(\gamma_{j+1}-\kappa_i)-(\gamma_j-\kappa_i)/\gamma}.
\end{align*}
Hence,
\begin{align*}
&\oldbigwedge_{\substack{j=i+1}}^{s-1} \sign(\gamma_j,\gamma_{j+1}) \nvect_{\gamma_{j+1}-\gamma_j/\gamma} \\
  &\qquad= (-1)^{s-i-1}\sign(\gamma_{i+1}-\kappa_i,\gamma_{i+1})\sign(\gamma_s-\kappa_i,\gamma_s)\nu_{\gamma_s-\kappa_i/\gamma_{i+1}-\kappa_i} \\
  &\qquad= (-1)^{s-i-1}\sign(\gamma_i,\gamma_{i+1})\sign(\zeta_i-\epsilon,\eta-\epsilon)\nu_{\zeta_i-\epsilon/\gamma_i}.
\end{align*}
Thus, the $i$-th term of the sum equals
\begin{align*}
& \sign(\zeta_i-\epsilon,\eta-\epsilon) \nvect^\dual_{\gamma\/\delta}(\p_\delta(\nvect_{\eta-\zeta_i})) \nu_{\gamma_i/\gamma} \wedge \nu_{\zeta_i-\epsilon/\gamma_i} \wedge \nvect_{\delta/\gamma} \\
  &\qquad = \sign(\zeta_i-\epsilon,\eta-\epsilon) \nvect^\dual_{\gamma\/\delta}(\p_\delta(\nvect_{\eta-\zeta_i/\gamma})) \nu_{\zeta_i-\epsilon/\gamma} \wedge \sign(\zeta_i-\epsilon,\zeta_i) \nu_{\zeta_i/\zeta_i-\epsilon} \\
  &\qquad = -\sign(\zeta_i,\eta)\sign(\eta-\epsilon,\eta) \nvect^\dual_{\gamma\/\delta}(\p_\delta(\nvect_{\eta-\zeta_i/\gamma}))\nu_{\zeta_i/\gamma}.
\end{align*}
Therefore,
\begin{align*}
\nu_{\eta/\delta}
  &= (-1)^s\sign(\gamma,\delta) \Bigl(\sign(\eta-\epsilon,\eta)\nu_{\eta-\epsilon/\gamma} + \sum_i \sign(\zeta_i,\eta)\nvect^\dual_{\gamma\/\delta}(\p_\delta(\nvect_{\eta-\zeta_i}))\nu_{\zeta_i/\gamma} \Bigr).
\end{align*}
Notice that the $\zeta_i$ are exactly the faces $\zeta'$ such that $\delta\subface\zeta'\ssubface\eta$, and $\eta-\epsilon$ is the only face $\zeta'$ such that $\gamma\subface\zeta'\ssubface\eta$ and such that $\delta\not\subface\zeta'$. Moreover,
\[ \nu^\dual_{\gamma\/\delta}(\p_\delta(\nvect_{\eta-(\eta-\epsilon)/\gamma)})) = \nu^\dual_{\gamma\/\delta}(\p_\delta(\nvect_{\delta/\gamma)})) = 1. \]
Finally,
\[ \nu_{\eta/\delta}
  = (-1)^s\sign(\gamma,\delta) \sum_{\zeta' \ssubface \eta \\ \zeta' \supface \gamma} \sign(\zeta',\eta) \nvect^\dual_{\gamma\/\delta}(\p_\delta(\nvect_{\eta-\zeta'/\gamma})) \nu_{\zeta'/\gamma},\]
which is what we wanted to prove.

\medskip

\item \label{com:tbd} We have to show that $\d'\d^\pi+\d^\pi\d'=0$ in $\AA^{\bul,\bul}_\trop$. Let $a, b$ be two integers, let $\gamma$ be a face of dimension $\dims\gamma =a+b$, and choose $\alpha\otimes \beta \in\bigwedge^b\TT^\dual\gamma \otimes \SF^{p-b}(\conezero^\gamma)$. Suppose moreover that $\sed(\eta)=\sed(\delta)\ssubface\sed(\gamma)$. Let $\zeta=\eta_\infty^{\sed(\gamma)}$. Then, we have
\begin{align*}
&\hspace{-2cm}(\d'\d^\pi+\d^\pi\d')\bigl[\AA_\trop^{a,b}[\gamma] \to \AA^{a,b,0}[\delta, \eta]\bigr](\alpha \otimes \beta)=0\\
\LRar& (-1)^{(a+2)(b+p)+\dims{\sed(\delta)}}\sign(\gamma,\delta)\pi^{\sed*}_{\gamma\/\delta}(\alpha)\otimes(\id^{\gamma\/\delta}(\beta) \vee e_{\eta/\delta})1_\eta
  \\&\qquad+ \sign(\zeta,\eta)(-1)^{(a+1)(b+p)+\dims{\sed(\gamma)}}\pi^{\sed *}_{\gamma/\delta}(\alpha) \otimes (\beta \vee e_{\zeta/\gamma})\id^{\zeta\/\eta}(1_{\zeta}) = 0, \\
\LRar& (-1)^{b+p+1}\sign(\gamma,\delta)\id^{\delta/\gamma}(e_{\eta/\delta})
  = -\sign(\zeta,\eta)e_{\zeta/\gamma}, \\
\LRar& \id^{\delta/\gamma}(e_{\eta/\delta})
  = (-1)^{\dims\eta-\dims\delta}\sign(\gamma,\delta)\sign(\zeta,\eta)e_{\zeta/\gamma},
\end{align*}
which follows from the definition.

\medskip

\item \label{com:ST_pi} We need to show that $(\d^\pi_\ST)^2=\d^\i\d^\pi_\ST+\d^\pi\d^\i_\ST=\d'\d^\pi_\ST+\d^\pi\d'_\ST=\d\d^\pi_\ST+\d^\pi\d_\ST=\d^\i\d'_\ST+\d'\d^\i_\ST=0$ in $\AA^{\bul, \bul}_\ST$. This is clear by looking at the definitions and using \eqref{com:eqn:zeta_pi} and the defining property of the sign function.

\medskip

\item \label{com:ST} We have to show that $(\d'+\d^\i)^2\bigl[\AA^{a,b'}_\ST[\eta] \to *]=0$.

Assume $\sigma=\sed\eta$ and $\xi=\maxsed\eta$. If $\sigma=\xi$ the statement holds by applying Proposition \ref{prop:app1} to $(X^\sigma)_\f$. Otherwise, we can reduce to the previous case as follows. By \eqref{com:ST_pi}, the two differentials $\d'_\ST$ and $\d^\i_\ST$ anticommute with $\d^\pi_\ST$. Recall that $\d^\pi_\ST$ is, up to a sign, identity on each piece. Thus, setting
\[ \partial:=(\d^\pi)^{\dims\xi-\dims\sigma}[\AAa{\substack{\sed=\xi \\ \maxsed=\xi}}_\ST^{\bul,\bul}\to\AAa{\substack{\sed=\sigma \\ \maxsed=\xi}}_\ST^{\bul,\bul}], \]
we get
\[ \partial^{-1}(\d'+\d^\i)^2\partial=(-1)^{\dims\xi-\dims\sigma}(\d'+\d^\i)^2[\AAa{\substack{\sed=\xi \\ \maxsed=\xi}}_\ST^{\bul,\bul}\to *]. \]

\medskip

\item \label{com:saa} This says that $(\d^2)[\AA_\ST^{a,b'}[\eta] \to \AA^{a,1,b'}[\delta,\eta]](x) = 0.$ This is equivalent to
\begin{align*}
\LRar& \sum_{\gamma'\ssubface\delta \\ \maxsed(\gamma')=\maxsed(\eta) \\ \sed(\gamma')=\sed(\eta)}(-1)^{a+b'}(-1)^{a+1'}\sign(\gamma,\delta) (1^\dual_\delta \wedge \nu^\dual_{\gamma'\/\delta}) \otimes x = 0 \\
\LRar& \ssum_{\gamma'} \nvect^\dual_{\gamma'\/\delta} = 0.
\end{align*}
Notice that the sum has no term if $\maxsed(\delta)\neq\maxsed(\eta)$. Otherwise, the sum is over facets of $\delta$ of same sedentarity and such that $\gamma_\infty=\delta_\infty$. Thus, every linear forms of the sum is zero on $\delta_\infty$.
Let $u$ and $u'$ be two vertices of $\delta_\f$. For any $\gamma'$,
\[ \nvect^\dual_{\gamma'\/\delta}(u'-u) = \begin{cases}
  1 & \text{if $u'\in\gamma'$ and $u\not\in\gamma'$,} \\
  -1 & \text{if $u\in\gamma'$ and $u'\not\in\gamma'$} \\
  0 & \text{otherwise.}
\end{cases} \]
The first condition can be rewritten as \enquote{if $\gamma'=\delta-u$} and the second condition as \enquote{if $\gamma'=\delta-u'$}. Thus, the first condition and the second condition each occurs exactly once. Hence, the sum of form applied to $u-u'$ is zero. To conclude, notice that $\TT\delta$ is spanned by the union of the edges of $\delta_\f$ and of $\delta_\infty$.

\medskip

\item \label{com:sac} If $\maxsed(\delta)=\maxsed(\zeta)=\maxsed(\eta)$, we have to show that $(\d\d^\i+\d^\i\d)\bigl[\AA_\ST^{a,b'}[\zeta]\to\AA^{a+1,0,b'}[\delta,\eta]\bigr](x) = 0$. This is equivalent to
\begin{align*}
\LRar& (-1)^{a+1+b'}\sign(\zeta,\eta) 1^\dual_\delta \otimes \i^*_{\zeta\/\eta}(x)
  \\ &\qquad + \hspace{-2em}\sum_{\gamma'\ssubface\delta \\ \maxsed(\gamma')=\maxsed(\delta)}\hspace{-2em} \sign(\zeta,\eta)(-1)^{a+b'} \nvect^\dual_{\gamma'\/\delta}(\p_\delta(\nvect_{\eta-\zeta/\gamma'})) \p^*_{\gamma'\/\delta}(1^\dual_{\gamma'}) \otimes \i^*_{\zeta\/\eta}(x) = 0 \\
\LRar& 1 = \ssum_{\gamma'} \nvect^\dual_{\gamma'\/\delta}(\p_\delta(\nvect_{\eta-\zeta/\gamma'})).
\end{align*}
Let $v$ be any vertex of $\delta_\f$. We have already seen in \eqref{com:saa} that
\[ \ssum_{\gamma'} \nvect^\dual_{\gamma'\/\delta}=0. \]
Thus,
\[ \ssum_{\gamma'} \nvect^\dual_{\gamma'\/\delta}(\p_\delta(\nvect_{\eta-\zeta/\gamma'})) = \ssum_{\gamma'} \nvect^\dual_{\gamma'\/\delta}(\p_\delta(\nvect_{\eta-\zeta/\gamma'}-\nvect_{\eta-\zeta/v})). \]
For each $\gamma'$, set $u_{\gamma'}:=\nvect_{\eta-\zeta/\gamma'}-\nvect_{\eta-\zeta/v}$. Then $u_{\gamma'}$. is a vector going from some point of $\TT\gamma'$ to $v$. Thus, if $v\in\gamma'$, then $u_{\gamma'}\in\TT\gamma'$ and
\[ \nvect^\dual_{\gamma'\/\delta}(\p_\delta(u_{\gamma'})) = 0. \]
Otherwise, \ie, for $\gamma'=\delta-v$, up to an element of $\TT\gamma'$, $u_{\delta-v}$ will be equal to $\nu_{\delta/\gamma}$, and
\[ \nvect^\dual_{\gamma'\/\delta}(\p_\delta(u_{\delta-v})) = 1. \]
Thus, the sum equals $1$.

\medskip

\item \label{com:a_easy} We have to show that $\d_\D^2=\d_\D\d'_\D+\d'_\D\d_\D=\d_\D\d^\pi_\D+\d^\pi_\D\d_\D=\d_\D'\d_\D^\pi+\d_\D^\pi\d_\D'=(\d^\pi_\D)^2=0$. This follows directly from the definitions, the property of $\sign(\cdot,\cdot)$, \eqref{com:eqn:nu*}, \eqref{com:eqn:zeta_pi} and \eqref{com:eqn:pi_nu*}.

\medskip

\item \label{com:abb} We need to show that $\d'^2_\D=0$; see the proof that $\gys^2=0$ in Section \ref{sec:proofapp1}.

\medskip

\item \label{com:aac} This is the property $\d\d^\i+\d^\i\d\bigl[\AA_\D^{a,b,b'}[\gamma,\eta] \to \AA_\D^{a,b,b'}[\chi,\mu]\bigr](\alpha \otimes x) = 0$, which is equivalent to
\begin{align*}
\LRar& \sum_{\ddelta\in\{\delta,\delta'\}}(-1)^{a+1+b'}\sign(\ddelta,\chi)\sign(\eta,\mu)\nvect^\dual_{\gamma\/\ddelta}(\p_{\ddelta}(\nvect_{\mu-\eta/\gamma}))\p^*_{\gamma\/\ddelta}(\alpha)\wedge\nu^\dual_{\ddelta\/\chi} \otimes \i^*_{\eta\/\mu}(x)
  \\&\qquad + (-1)^{a+b'}\sign(\eta,\mu)\nvect^\dual_{\ddelta\/\chi}(\p_\chi(\nvect_{\mu-\eta/\ddelta}))\sign(\gamma,\ddelta) \p^*_{\ddelta\/\chi}(\alpha \wedge \nu^\dual_{\gamma\/\ddelta}) \otimes \i^*_{\eta\/\mu}(x) = 0 \\
\LRar&
  \ssum_\ddelta -\sign(\ddelta,\chi)\nvect^\dual_{\gamma\/\ddelta}(\p_{\ddelta}(\nvect_{\mu-\eta/\gamma})) (\alpha\wedge\nu^\dual_{\gamma\/\ddelta}\vee\nu_{\ddelta/\gamma})\wedge\nu^\dual_{\ddelta\/\chi}
    & \eqref{com:eqn:p_nu*_nu}
  \\&\qquad + \sign(\gamma,\ddelta)\nvect^\dual_{\ddelta\/\chi}(\p_\chi(\nvect_{\mu-\eta/\ddelta})) \alpha\wedge\nu^\dual_{\gamma\/\ddelta}\wedge\nu^\dual_{\ddelta\/\chi}\vee\nu_{\chi/\ddelta} = 0 \\
\LRar&
  \alpha\wedge\nu^\dual_{\gamma\/\chi}\Bigl( \ssum_\ddelta \sign(\ddelta,\chi)\nvect^\dual_{\gamma\/\ddelta}(\p_{\ddelta}(\nvect_{\mu-\eta/\gamma}))\nu_{\ddelta/\gamma}  +  \sign(\gamma,\ddelta)\nvect^\dual_{\ddelta\/\chi}(\p_\chi(\nvect_{\mu-\eta/\ddelta}))\nu_{\chi/\ddelta} \Bigr) = 0
    & \eqref{com:eqn:wedge_vee},\eqref{com:eqn:nu*} \\
\Lar&
  \ssum_\ddelta \sign(\ddelta,\chi)\sign(\gamma,\ddelta)\p_\ddelta(\nvect_{\mu-\eta/\gamma}) + \sign(\gamma,\ddelta)\sign(\ddelta,\chi)\p_\chi(\nvect_{\mu-\eta/\ddelta}) = 0
    & \eqref{com:eqn:nu*_nu} \\
\LRar&
  \p_\delta(\nvect_{\mu-\eta/\gamma}) + \p_\chi(\nvect_{\mu-\eta/\delta}) = \p_{\delta'}(\nvect_{\mu-\eta/\gamma}) + \p_\chi(\nvect_{\mu-\eta/\delta'}) \\
\LRar&
  \p_\chi \p_\delta(\nvect_{\mu-\eta/\gamma}) + \p_\chi\cp_\delta(\nvect_{\mu-\eta/\gamma}) = \p_\chi \p_{\delta'}(\nvect_{\mu-\eta/\gamma}) + \p_\chi\cp_{\delta'}(\nvect_{\mu-\eta/\gamma}) \\
\LRar& \p_\chi(\nvect_{\mu-\eta/\gamma}) = \p_\chi(\nvect_{\mu-\eta/\gamma}),
\end{align*}
which is obvious.

\vspace{.7cm}

For cases \eqref{com:abc}, \eqref{com:abc'} and \eqref{com:abc''}, recall our naming convention of faces: we have $\gamma \ssubface \eta, \eta'$. We have to prove that
\begin{gather*}
(\d'\d^\i+\d^\i\d')\bigl[\AA^{a,b,b'}_\D[\gamma,\eta'] \to \AA^{a+1,b,b'+1}_\D[\delta,\eta]\bigr], \textrm{ and } \\
(\d'\d^\i+\d^\i\d')\bigl[\AA^{a,b,b'}_\D[\gamma,\eta] \to \AA^{a+1,b,b'+1}_\D[\delta,\eta]\bigr](\alpha \otimes x) \text{ equals 0}.
\end{gather*}
For the first equation, we consider two cases depending on whether there exists $\mu$ such that $\eta,\eta' \ssubface \mu$ or not.

\medskip

\item \label{com:abc} If $\mu$ exists, then we have
\begin{align*}
&\hspace{-2cm}(\d'\d^\i+\d^\i\d')\bigl[\AA^{a,b,b'}_\D[\gamma,\eta'] \to \AA^{a+1,b,b'+1}_\D[\delta,\eta]\bigr](\alpha \otimes x)=0\\
\LRar&
  \sign(\eta,\mu)\sign(\eta',\mu)\nvect^\dual_{\gamma\/\delta}(\p_\delta(\nvect_{\mu-\eta'/\gamma}))\p^*_{\gamma\/\delta}(\alpha) \otimes \gys^{\mu/\eta}\i^*_{\eta'\/\mu}(x)
  \\&\qquad + \sign(\zeta,\eta)\sign(\zeta,\eta') \nvect^\dual_{\gamma\/\delta}(\p_\delta(\nvect_{\eta-\zeta/\gamma})) \p^*_{\gamma\/\delta}(\alpha) \otimes \i^*_{\zeta\/\eta}\gys^{\eta'/\zeta}(x) \\
\LRar&
  \nvect^\dual_{\gamma\/\delta}(\p_\delta(\nvect_{\mu-\eta'/\gamma}))
  \\&\qquad = \nvect^\dual_{\gamma\/\delta}(\p_\delta(\nvect_{\eta-\zeta/\gamma})) \\
\LRar& \mu-\eta'=\eta-\zeta,
\end{align*}
which is obviously the case.

\medskip

\item \label{com:abc'} If $\mu$ does not exist, then $(\d'\d^\i+\d^\i\d')\bigl[\AA^{a,b,b'}_\D[\gamma,\eta'] \to \AA^{a+1,b,b'+1}_\D[\delta,\eta]\bigr](\alpha \otimes x)=0$. We only have the second term from the case \eqref{com:abc}. However, here, $\i^*_{\zeta\/\eta}\gys^{\eta'/\zeta}(x)=0$ since $\eta_\infty^\zeta$ and $\eta'^\zeta_\infty$ are not comparable.

\medskip

\item \label{com:abc''} It remains to prove $(\d'\d^\i+\d^\i\d')\bigl[\AA^{a,b,b'}_\D[\gamma,\eta] \to \AA^{a+1,b,b'+1}_\D[\delta,\eta]\bigr](\alpha \otimes x)=0.$ This is equivalent to
\begin{align*}
\LRar&
  \sum_{\mu' \ssupface \eta \\ \sed(\mu')=\sed(\eta)}\sign(\eta,\mu')\sign(\eta,\mu')\nvect^\dual_{\gamma\/\delta}(\p_\delta(\nvect_{\mu'-\eta/\gamma}))\p^*_{\gamma\/\delta}(\alpha) \otimes \gys^{\mu'/\eta}\i^*_{\eta\/\mu'}(x)
  \\&\qquad + \sum_{\zeta' \ssubface \eta \\ \zeta' \supface \gamma \\ \sed(\zeta')=\sed(\eta)} \sign(\zeta',\eta)\sign(\zeta',\eta) \nvect^\dual_{\gamma\/\delta}(\p_\delta(\nvect_{\eta-\zeta'/\gamma})) \p^*_{\gamma\/\delta}(\alpha) \otimes \i^*_{\zeta'\/\eta}\gys^{\eta/\zeta'}(x) = 0 \\
\LRar&
  \ssum_{\mu'}\p_\delta^*(\nvect^\dual_{\gamma\/\delta})(\nvect_{\mu'-\eta/\gamma})\gys^{\mu'/\eta}\i^*_{\eta\/\mu'}(1^\eta)
  + \ssum_{\zeta'} \p_\delta^*(\nvect^\dual_{\gamma\/\delta})(\nvect_{\eta-\zeta'/\gamma})\i^*_{\zeta'\/\eta}\gys^{\eta/\zeta'}(1^\eta) = 0 \\
\LRar&
  \ssum_{\mu'}(\p^\gamma_\delta)^*(\e^\dual_{\gamma\/\delta})(\e_{\mu'-\eta/\gamma})\i^*_{\gamma\/\eta}(x_{\mu'-\eta/\gamma})
  + \ssum_{\zeta'} (\p^\gamma_\delta)^*(\e^\dual_{\gamma\/\delta})(\e_{\eta-\zeta'/\gamma})\i^*_{\gamma\/\eta}(x_{\eta-\zeta'/\gamma}) = 0.
\end{align*}

Notice now that rays of the form $\mu'^\gamma-\eta^\gamma$ are exactly the rays of $\Sigma^\gamma$ in the link of $\eta^\gamma$, and the rays of the form $\eta^\gamma-\zeta'^\gamma$ are exactly the rays of $\eta^\gamma$. Moreover, if $\varrho$ is a ray which is not in the link of $\eta^\gamma$ nor in $\eta^\gamma$, then $\i^*_{\gamma\/\eta}(x_\varrho)=0$. Thus, the above sum equals
\begin{align*}
\i^*_{\gamma\/\eta}\Bigl(\sum_{\varrho\in\Sigma^\gamma \\ \dims\varrho = 1} (\p^\gamma_\delta)^*(\e^\dual_{\gamma\/\delta})(\e_\varrho)x_\varrho\Bigr),
\end{align*}
which is clearly zero.

\medskip

\item \label{com:acc} We have to show that $(\d^\i)^2\bigl[\AA_\D^{a,b,b'}[\gamma,\zeta] \to \AA_\D^{a+2,b,b'}[\chi,\mu]\bigr](\alpha \otimes x)=0$. This is equivalent to
\begin{align*}
\LRar&
  \sum_{\eeta \in \{\eta,\eta'\}}\sum_{\ddelta \in \{\delta,\delta'\} \\ \ddelta \subface \eeta} \sign(\eeta,\mu)\nvect^\dual_{\ddelta\/\chi}(\p_\chi(\nvect_{\mu-\eeta/\ddelta})) \sign(\zeta,\eeta)\nvect^\dual_{\gamma\/\ddelta}(\p_\ddelta(\nvect_{\eeta-\zeta/\gamma})) \p^*_{\ddelta\/\chi}\p^*_{\gamma\/\ddelta}(\alpha) \otimes \i^*_{\eeta\/\mu}\i^*_{\zeta\/\eeta}(x)=0 \\
\LRar&
  \ssum_{\eeta,\ddelta} \sign(\eeta,\mu)\sign(\zeta,\eeta)\sign(\ddelta,\chi)\sign(\gamma,\ddelta) \nu^\dual_{\ddelta\/\chi}(\p_\chi\cp_\ddelta(\nvect_{\mu-\eeta/\gamma})) \nu^\dual_{\gamma\/\ddelta}(\p_\chi \p_\ddelta(\nvect_{\eeta-\zeta/\gamma})) =0 \\
\LRar&
  \ssum_{\eeta,\ddelta} \sign(\eeta,\mu)\sign(\zeta,\eeta)\sign(\ddelta,\chi)\sign(\gamma,\ddelta) \nu^\dual_{\gamma\/\chi}(\cp_\ddelta \p_\chi(\nvect_{\mu-\eeta/\gamma}) \wedge \p_\ddelta \p_\chi(\nvect_{\eeta-\zeta/\gamma})) =0.
\end{align*}

Notice that, if $\delta'\not\subface\eta$, then $\mu-\eta=\delta'-\gamma$ and $\cp_{\delta'}\p_\chi(\nvect_{\mu-\eta/\gamma})=0$. Using the same argument for $\delta$ and $\eta'$, this shows that we can now remove the condition $\ddelta\subface\eeta$ in the above sum.

\smallskip
\par Set $u=\p_\chi(\nvect_{\mu-\eta'/\gamma})=\p_\chi(\nvect_{\eta-\zeta/\gamma})$ and $u'=\p_\chi(\nvect_{\mu-\eta/\gamma})=\p_\chi(\nvect_{\eta'-\zeta/\gamma})$. The above sum equals
\begin{align*}
\sign(\eta,\mu)\sign(\zeta,\eta)\sign(\delta,\chi)\sign(\gamma,\delta)\nu^\dual_{\gamma\/\chi}\Bigl(\cp_\delta(u')\wedge \p_\delta(u) - \cp_\delta(u)\wedge \p_\delta(u') - \cp_{\delta'}(u')\wedge \p_{\delta'}(u) + \cp_{\delta'}(u)\wedge \p_{\delta'}(u')\Bigr).
\end{align*}
Then we have
\begin{align*}
&\cp_\delta(u')\wedge \p_\delta(u) - \cp_\delta(u)\wedge \p_\delta(u') \\
&\qquad= \bigl(\cp_\delta(u')+\p_\delta(u')\bigr)\wedge\bigl(\cp_\delta(u)+\p_\delta(u)\bigr) - \cp_\delta(u') \wedge \cp_\delta(u) - \p_\delta(u')\wedge \p_\delta(u) \\
&\qquad= u' \wedge u - 0 - \pi_\delta(u') \wedge \pi_\delta(u).
\end{align*}
Notice that $\nu^\dual_{\gamma\/\chi}(\v)=0$ if $\v\in\bigwedge^2\TT\delta$. Thus, up to a sign, the above sum equals
\[ \nu^\dual_{\gamma\/\chi}(u' \wedge u + u \wedge u') = 0. \]

\medskip

\item \label{com:acd} Finally, we have to show that $(\d^\i\d^\pi+\d^\pi\d^\i)\bigl[\AA_\D^{a,b,b'}[\gamma,\zeta] \to \AA_\D^{a+2,b,b'}[\chi,\mu]\bigr](\alpha \otimes x)=0$ with $\sed(\mu)=\sed(\chi)=\sed(\eta')=\sed(\delta')\ssubface\sed(\eta)=\sed(\delta)=\sed(\zeta)=\sed(\gamma)$. This is equivalent to
\begin{align*}
&
  \sign(\eta',\mu)\nvect^\dual_{\delta'\/\chi}(\p_\chi(\nvect_{\mu-\eta'/\delta'}))\sign(\zeta,\eta') \p^*_{\delta'\/\chi}\pi^{\sed*}_{\gamma/\delta'}(\alpha) \otimes \i^*_{\eta'\/\mu}\id^{\zeta\/\eta'}(x)
  \\*&\qquad + \sign(\eta,\mu)\sign(\zeta,\eta) \nvect^\dual_{\gamma\/\delta}(\p_\delta(\nvect_{\eta-\zeta/\gamma})) \pi^{\sed *}_{\delta\/\chi}\p^*_{\gamma\/\delta}(\alpha) \otimes \id^{\eta\/\mu}\i^*_{\zeta\/\eta}(x) = 0,
\end{align*}
which is implied by
\begin{align*}
&
  \e^\dual_{\delta'\/\chi}(\p^{\delta'}_\chi(\e_{\mu-\eta'/\delta'})) = \e^\dual_{\gamma\/\delta}(\p^{\gamma}_\delta(\e_{\eta-\zeta/\gamma}))\myand  \pi^{\sed}_{\gamma\/\delta'}\p_{\delta'\/\chi} = \p_{\gamma\/\delta}\pi^{\sed}_{\delta\/\chi}.
\end{align*}

\end{enumerate}

\medskip


\section{Tropical monodromy}

In this section, we prove Theorem \ref{thm:tropical_monodromy} that describes the simplicial tropical monodromy operator and claims that it corresponds to the monodromy operator on Steenbrink. We use the notations and conventions of Section \ref{sec:technicalities}. Recall that there is one triple complex $\AA^{\bul,\bul,\bul}$ for each even row $\ST^{\bul,2p}_1$ of the steenbrink spectral sequence. To distinguish these spectral sequences, we precise the corresponding $p$ as follows: $\AAp{p}^{\bul,\bul,\bul}$.

\medskip

For each face $\delta\in X_\f$, let $o_\delta$ be a point of $\Tan\delta\subseteq N_\R$, for instance the centroid of $\delta$. For other faces of $\delta$ of sedentarity $\conezero$, set $o_\delta:=o_{\delta_\f}$. For any faces $\delta$ of any sedentarity $\sigma$, set $o_\delta=\pi^\sigma(o_\eta)$, where $\eta$ is the only face of sedentarity $\conezero$ such that $\delta=\eta_\infty^\sigma$. We define $v_{\delta\/\eta}:=o_\eta-o_\delta$.

\medskip

Let $V$ and $W$ be any vector spaces such that $\TT\delta\subseteq V$ and $k$ be an integer. We define the morphism $\~N_{\gamma\/\delta}\colon \bigwedge^kV^\dual\otimes W\to\bigwedge^{k-1}V^\dual\otimes W$ by
\[ \~N_{\gamma\/\delta} = (\ldot \vee v_{\gamma\/\delta}) \otimes \id\!. \]

We define a map between triple complexes $N\colon\AAp{p}^{\bul,\bul,\bul} \to \AAp{p-1}^{\bul,\bul,\bul}$ as follows. Assume $\gamma$ is a face of dimension $a+b$ and $\eta$ of dimension $a+p-b'$, with $b'\leq p$, and such that $\dims{\eta_\infty}\leq a$. Then,
\begin{align*}
N\Bigl[\,\AAp{p}_\trop^{a,b}[\gamma] \to \AAp{p-1}^{\bul,\bul}_\trop[\delta] \,\Bigr]
  &= \~N_{\gamma\/\delta}(\d_\trop+\dfrak_\trop), \\
  &\\
N\Bigl[\, \AAp{p}_\D^{a,b,b'}[\gamma, \eta] \to \AAp{p-1}_\D^{\bul,\bul,\bul}[\delta, \eta]\, \Bigr]
  &= \~N_{\gamma\/\delta}(\d_\D+\dfrak_\D), \\
  &\\
N\Bigl[\,\AAp{p}_\ST^{a,b'}[\eta] \to \AAp{p-1}_\ST^{a+1,b'}[\eta]\,\Bigr]
  &= \begin{cases}
    \id & \text{if $b'\leq p-1$,} \\
    0 & \text{otherwise.}
  \end{cases}
\end{align*}

\begin{prop} \label{prop:com:N_commutes}
The map $N$ commutes with the differentials of $\AAp{p}^{\bul,\bul,\bul}$ and $\AAp{p-1}^{\bul,\bul,\bul}$.
\end{prop}

\begin{proof}
The fact that $N\d^\pi+\d^\pi N=0$ is easy to check.

\medskip

Let us prove that $(N\d-\d N)\bigl[\AAp{p}_{\trop}^{a,b}[\gamma] \to \AAp{p-1}_{\trop}^{a+1,b+1}[\chi]\bigr]=0$. The key point is that
\[ (\d\~N_{\gamma\/\delta} + \~N_{\gamma\/\delta}\d)\bigl [\AAp{p}_{\trop}[\eta] \to *\bigr] = 0, \]
which can easily been checked. We simplify the notations by writing $[ \gamma \to \chi]$ instead of $\bigl[\AAp{p}_{\trop}^{a,b}[\gamma] \to \AAp{p-1}_{\trop}^{a+1,b+1}[\chi]\bigr]$, or when the corresponding domain and codomain are clear.

\medskip

 Then we get
\begin{align*}
(N\d-\d N)[\gamma \to \chi]
  &= \sum_{\ddelta\in\{\delta,\delta'\}} N[\ddelta \to \chi]\d[\gamma \to \ddelta]-\d[\ddelta \to \chi]N[\gamma \to \ddelta] \\
  &= \sum_{\ddelta\in\{\delta,\delta'\}} (\~N_{\ddelta\/\chi}+\~N_{\gamma\/\ddelta})\d[\ddelta \to \chi]\d[\gamma \to \ddelta].
\end{align*}
From $\d^2=0$, we deduce that $\d[\delta' \to \chi]\d[\gamma \to \delta']=-\d[\delta \to \chi]\d[\gamma \to \delta]$. Thus,
\begin{align*}
(N\d-\d N)[\gamma \to \chi]
  &= (\~N_{\delta\/\chi}+\~N_{\gamma\/\delta}-\~N_{\delta'\/\chi}-\~N_{\gamma\/\delta'})\d[\delta\to\chi]\d[\gamma\to\delta].
\end{align*}
By definition of $\~N$, to prove this is zero, it suffices to check that
\[ v_{\delta\/\chi}+v_{\gamma\/\delta}-v_{\delta'\/\chi}-v_{\gamma\/\delta'}=0. \]
But this vector is by definition
\[ (o_\chi-o_\delta)+(o_\delta-o_\gamma)-(o_\chi-o_{\delta'})-(o_{\delta'}-o_\gamma), \]
which is clearly zero.

\medskip

One can prove with a similar argument that $N$ and $\d+\d'+\dfrak$ commute on the tropical part and also on the $\D$ part.

\medskip

For the Steenbrink part, Proposition \ref{prop:N_commutes} proves that $N$ commutes with $\d'_\ST$ and $\d^\i_\ST$ on $X_\f$. The proposition can easily been extended to all the faces, and the commutativity with $\d^\pi_\ST$ is also easy to check. Thus, it only remains to see the commutativity with $\d_\ST$, \ie, with the inclusion $\AA_\ST^{\bul,\bul} \hookrightarrow \AA_\D^{\bul,0,\bul}$.

\medskip

We have
\begin{align*}
&\hspace{-2em}(\d N-N\d)\bigl[\AAp{p}_\ST^{a,b'}[\eta] \to \AAp{p-1}_\D^{a,b',0}[\eta,\delta]\bigr](x) \\
  &= (-1)^{a+1+b'} 1^\dual_\delta \otimes x
    - \hspace{-2em}\sum_{\gamma' \ssubface \delta \\ \maxsed(\gamma')=\maxsed(\delta)}\hspace{-2em} (-1)^{a+b'}\sign(\gamma',\delta) v_{\gamma'\/\delta} \vee (1^\dual_{\gamma'}\wedge\nu^\dual_{\gamma'\/\delta}) \otimes x \\
  &= \bigl( -1 + \ssum_{\gamma'}\nvect^\dual_{\gamma'\/\delta}(v_{\gamma'\/\delta}) \bigr) 1^\dual_\delta \otimes x.
\end{align*}
The terms of the sum are the normalized barycentric coordinates of $o_\delta$ in $\delta_\f$. Hence, the terms sum up to one. This concludes the proof of the commutativity of $N$ with the differentials.
\end{proof}
}

\newpage

\bibliographystyle{alpha}
\bibliography{bibliographie}
\end{document}